# The Caporaso–Harris–Ran degeneration principle: proof and applications

**Abstract.** Severi varieties are the parameter spaces for curves with prescribed homology class and genus on a smooth surface. We describe their limits along degenerations of surfaces, with a view towards the enumeration of curves. This includes a complete proof of the Caporaso–Harris recursive formula, with all the necessary background on deformations of curves and singularities.

**Note.** This text is a coherent collection of interdependent contributions by various authors. It is one in a series[1] of texts issued from the seminar *Degenerations and enumeration of curves on surfaces*, held at the University of Roma "Tor Vergata" in the years 2015–2017.

---

[1] see `https://www.math.univ-toulouse.fr/~tdedieu/#enumTV`

# Contents

# Lecture I
# Limits of nodal curves: generalities and examples

by Ciro Ciliberto and Thomas Dedieu

## 1 – The basic question

The objective of this lecture is to introduce to the subject of interest in this series of seminars. This exposition is mainly based on [4]. In this text we will work over the complex field.

Let $f : S \to \mathbf{D}$ be a flat, projective family of surfaces of degree $d$ in $\mathbf{P}^r$ (with $r \geqslant 3$), with $S$ a smooth threefold, and $\mathbf{D}$ the complex unit disc, usually called a *degeneration* of the *general fibre* $S_t := f^{-1}(t)$, for $t \in \mathbf{D}^* := \mathbf{D} - \{0\}$, which we assume is a smooth irreducible surface, to the *central fibre* $S_0$, which may be singular and even reducible. This is also called a *smooth deformation* of $S_0$. One has $S \subseteq \mathbf{D} \times \mathbf{P}^r$, so that $S$ is endowed with the line bundle $\mathcal{L} = p^*(\mathcal{O}_{\mathbf{P}^r}(1))$, where $p : S \to \mathbf{P}^r$ is the projection to the second factor: this is called the *polarizing (hyperplane) bundle* on $S$.

**(1.1) Question.** *What are the limits of tangent, or in general of multitangent, hyperplanes to $S_t$, for $t \neq 0$, as $t$ tends to $0$?*

The same question makes sense, more generally, for varieties of any dimension. In the case of degeneration of (plane) curves we refer to [15, pp. 134–135] for a glimpse on this subject. In dimension higher than 2 the problem is quite complicated and few results are known (however, B. Segre's paper [23] is an interesting classical source of basic information on the subject).

One of the main interests in Question (1.1) arises form enumerative geometry, since, by degenerating a surface to a reducible one, whose components are much simpler than the original surface (e.g., a degeneration into a union of planes), we may hope that the configuration of limits of pluritangent hyperplanes is easier to understand and their number (if finite), or in general the degree of the variety parameterising them, may be computed. This is indeed the subject of the foundational work by L. Caporaso and J. Harris [1, 2], and independently by Z. Ran [18, 19, 20], which were both aimed to the study, with degeneration techniques, of the so–called *Severi varieties*, i.e., the families of irreducible nodal plane curves of a given degree.

Section 2 contains classical enumerative formulae for surfaces. Sections 3–5 present general facts about Question (1.1), following the presentation given in [4]. Finally, Sections 6–8 present several examples in detail in order to illustrate the general theory.

## 2 – Classical results

Enumerative results on the number of (pluri)tangent planes to a *general* surface $S$ of degree $d > 1$ in $\mathbf{P}^3$ are quite classical, and probably the first of them go back to G. Salmon, see [21, 22].

The set $\check{S} \subseteq \check{\mathbf{P}}^3$ of tangent planes to $S$ is the *dual* of $S$, which is birational to $S$. One has

$$\deg(\check{S}) = d(d-1)^2.$$

As soon as $d > 2$, $\check{S}$ is singular. How do the singularities of $\check{S}$ look like if $S$ is general enough? There is a *nodal curve* $D_b$, whose general points correspond to *bitangent planes* to $S$. Salmon computes

$$\deg(D_b) = \frac{1}{2}d(d-1)(d-2)(d^3 - d^2 + d - 12).$$

There is also a *cuspidal curve* $D_s$, whose points correspond to *stationary tangent planes* to $S$, i.e., tangent planes to $S$ in *parabolic points*, namely points $x \in S$ such that the tangent plane to $S$ at $x$ cuts $S$ in a curve with a cusp at $x$. Salmon computes

$$\deg(D_s) = 4d(d-1)(d-2),$$

whereas the curve $\Gamma_p$ of parabolic points on $S$ is the complete intersection of $S$ with its *Hessian surface*, hence

$$\deg(\Gamma_p) = 4d(d-2).$$

Finally, there are finitely many triple points of $\check{S}$, which are triple points of $D_b$ too. Salmon computes their number

(2.0.1) $$t = \frac{1}{6}d(d-2)(d^7 - 4d^6 + 7d^5 - 45d^4 + 114d^3 - 111d^2 + 548d - 960).$$

For all these classical formulae, the reader is invited to consult [**??**, **??**].

All these numbers have been computed in recent times by I. Vainsencher [25, 26] and by S. Kleiman and R. Piene [14, 15], via Chern class computations, see [**??**]. However several questions are still open, like:

**(2.1) Question.** *What are the singularities of $D_b$ and of $D_s$? Are these curves irreducible? How do they intersect? Is $\Gamma_p$ smooth? What is the monodromy of the covering $T \to U$, given by the triple points of $\check{S}$, when $S$ moves in the open subset $U \subseteq |\mathcal{O}_{\mathbf{P}^3}(d)|$ of smooth surfaces?*

Partial answers to some of these question have been given for $d = 4$ in [4].

## 3 – Semistable degenerations

A degeneration $f : S \to \mathbf{D}$ as in §1 is said to be *semistable* if the central fibre $S_0$ is reduced and with *local normal crossing* singularities. We will often assume that all irreducible components of $S_0$ are smooth, so the singularities of $S_0$ are: *double curves*, along which only two irreducible components of $S_0$ intersect transversally and *triple points*, where only three irreducible components of $S_0$ intersect transversally.

As in §1, we will assume there is a *polarising line bundle* $\mathcal{L}$ on $S$ (although, in principle, it could be neither ample nor relatively ample). Denote by $L_t$ the restriction of $\mathcal{L}$ to $S_t$ for all $t \in \mathbf{D}$. Then we will say that $(S_0, L_0)$ *is a limit* of $(S_t, L_t)$ (with $t \in \mathbf{D}^*$), when $t \to 0$ (or that $(S_t, L_t)$ with $t \in \mathbf{D}^*$ is a deformation of $(S_0, L_0)$).

The limit $(S_0, L_0)$ *is not unique* if $S_0$ is reducible. Indeed, if $W$ is an effective divisor supported on the central fibre $S_0$, consider the line bundle $\mathcal{L}(-W)$, which is said to be obtained from $\mathcal{L}$ by *twisting* by $W$. For $t \in \mathbf{D}^*$, its restriction to $S_t$ is the same as $L_t$, but in general this is not the case for $S_0$; any such line bundle $\mathcal{L}(-W)|_{S_0}$ is called a *limit line bundle* of $\mathcal{L}_t$ for $t \in \mathbf{D}^*$. If $Q$ is an irreducible component of $S_0$, the restriction $\mathcal{L}_{W,Q} := \mathcal{L}(-W)_{|Q}$ is called the *aspect* of $\mathcal{L}(-W)$ on $Q$.

**(3.1) Remark.** Since $\mathrm{Pic}(\mathbf{D})$ is trivial, the divisor $S_0 \subseteq S$ is linearly equivalent to 0. So if $W$ is a divisor supported on $S_0$, one has $\mathcal{L}(-W) \cong \mathcal{L}(mS_0 - W)$ for all integers $m$. In particular if $W + W' = S_0$ then $\mathcal{L}(-W) \cong \mathcal{L}(W')$.

Another family of surfaces $f' : S' \to \mathbf{D}$ as above is said to be a *model of* $f : S \to \mathbf{D}$ if there is a commutative diagram

$$\begin{array}{ccccccc}
S' & \longleftarrow & \bar{S}' & \stackrel{p}{\dashrightarrow} & \bar{S} & \longrightarrow & S \\
{\scriptstyle f'}\downarrow & \square & \downarrow & & \downarrow & \square & \downarrow {\scriptstyle f} \\
\mathbf{D} & \underset{t^{d'} \leftarrow t}{\longleftarrow} & \mathbf{D} & = & \mathbf{D} & \underset{t \mapsto t^d}{\longrightarrow} & \mathbf{D}
\end{array}$$

where the two squares marked with a $\square$ are Cartesian, and $p$ is a birational map, which is an isomorphism over $\mathbf{D}^*$. The family $f' : S' \to \mathbf{D}$, if semistable, is a *semistable model of* $f : S \to \mathbf{D}$ if in addition $d' = 1$ and $p$ is a morphism. The *semistable reduction theorem* of [12] asserts that any family $f : S \to \mathbf{D}$ as in §1 has a semistable model.

**(3.2) Example** (Families of surfaces in $\mathbf{P}^3$). Consider a *linear pencil* of degree $d$ surfaces in $\mathbf{P}^3$, generated by a *general* surface $S_\infty$ and a *special* one $S_0$. We will usually consider the case in which $S_0$ has local normal crossing singularities, and $S_\infty$ intersects transversally the double curve of $S_0$ at smooth points of it. This pencil gives rise to a flat, proper family $\varphi : \mathcal{S} \to \mathbf{P}^1$, with $\mathcal{S}$ a hypersurface of type $(d, 1)$ in $\mathbf{P}^3 \times \mathbf{P}^1$, isomorphic to the blow–up of $\mathbf{P}^3$ along the *base locus scheme* $S_0 \cap S_\infty$ of the pencil, and it has $S_0, S_\infty$ as fibres over $0, \infty \in \mathbf{P}^1$, respectively.

We will study the family obtained by restricting $\mathcal{S}$ to a disk $\mathbf{D} \subseteq \mathbf{P}^1$ centered at 0, that by abuse of notation we will still denote by $f : \mathcal{S} \to \mathbf{D}$, such that $S_t$ is smooth for all $t \in \mathbf{D}^*$, and we will consider a semistable model of $f : S \to \mathbf{D}$. To do so, we resolve the singularities of $S$ which occur only in the central fibre of $f$, at the points mapped by $\mathcal{S}_0 \to S_0 \subseteq \mathbf{P}^3$ to the intersection points of $S_\infty$ with the double curves of $S_0$ (they are the singular points of the curve $S_0 \cap S_\infty$). These are *ordinary double points* of $S$, i.e., singularities analytically equivalent to the one at the origin of the hypersurface $xy = zt$ in $\mathbf{A}^4$. Such a singularity is resolved by a single blow–up, which produces an exceptional divisor $E \cong \mathbf{P}^1 \times \mathbf{P}^1$, and then it is possible to contract $E$ in the direction of either one of its rulings without introducing any singularity: the result is called a *small resolution* of the ordinary double point.

Let $\tilde{f} : \tilde{S} \to \mathbf{D}$ be a semistable model thus obtained. One has $\tilde{S}_t \cong S_t$ for $t \in \mathbf{D}^*$. If $S_0$ has irreducible components $Q_1, \ldots, Q_r$, then $\tilde{S}_0$ consists of irreducible components $\tilde{Q}_1, \ldots, \tilde{Q}_r$ which are suitable blow–ups of $Q_1, \ldots, Q_r$, respectively. If $q$ is the number of ordinary double points of the original total space $S$, we will denote by $E_1, \ldots, E_q$ the exceptional curves on $\tilde{Q}_1, \ldots, \tilde{Q}_r$ arising from the small resolution process.

The simplest example of this situation is when $S_0$ is the union of a general surface $S'$ of degree $h$ with a general surface $S''$ of degree $d-h$ (and $1 \leqslant h \leqslant \lfloor \frac{d}{2} \rfloor$), intersecting along a smooth curve $R$ of degree $h(d-h)$. In this case $\tilde{S}_0$ consists of two components $Q_1, Q_2$, and we may assume that $Q_2 = S''$, whereas $Q_1 = \tilde{S}'$ is the blow-up of $S'$ at the $q = dh(d-h)$ intersection points of $R$ with $S_\infty$. Then $Q_1$ and $Q_2$ intersect along a curve which can be identified with $R$ and the exceptional curves $E_1, \ldots, E_q$ all lie on $Q_1$.

The case $h = 1$ in which $S' = P$ is a general plane, is particularly interesting since it shows the geometric significance of twisting the polarising line bundle $\mathcal{L} = p^*(\mathcal{O}_{\mathbf{P}^3}(1))$ on $\tilde{S}$. The polarising line bundle $L_0$ maps $\tilde{S}_0$ to the reducible surface $S'' \cup P$. If we twist by $Q_1 = \tilde{P}$ and consider $\mathcal{L}' = \mathcal{L}(-\tilde{P})$, then its aspects are:

- the trivial bundle on $Q_2 = S''$;
- $|\mathcal{O}_{\tilde{P}}(dH - \sum_{i=1}^{d(d-1)} E_i)|$ on $Q_1 = \tilde{P}$.

Hence $L'_0$ contracts $S''$ to a point $x$. Moreover $\dim(|\mathcal{O}_{\tilde{P}}(dH - \sum_{i=1}^{d(d-1)} E_i)|) = 3$ and this linear system maps $\tilde{P}$ to a *monoid* $\Sigma$, i.e., to a rational surface of degree $d$ with a point of multiplicity $d-1$ at $x$ (see Remark (7.3) for the case $d = 4$).

In conclusion, by twisting, we see that the degeneration of $(S_t, \mathcal{O}_{S_t}(1))$ to $(S'' \cup \tilde{P}, L_0)$ is also a semistable model of the degeneration of $(S_t, \mathcal{O}_{S_t}(1))$ to a general monoid $(\Sigma, \mathcal{O}_\Sigma(1))$.

**(3.3) Example** (Families of polarised K3 surfaces)**.** A special case of the previous example is the one of a general quartic surface in $\mathbf{P}^3$ degenerating to the general union of two quadrics. This is also a special case of degenerations of polarised K3 surfaces to a union of two scrolls (see [7]), which are type II degenerations according to the Kulikov-Persson-Pinkham classification of semistable degenerations of $K3$ surfaces (see [16, 17]).

Let $\mathcal{K}_p$ be the moduli space of primitively polarized $K3$ surfaces $(S, L)$ of genus $p \geqslant 3$, i.e., the dualising sheaf $\omega_S$ of $S$ is trivial, $h^1(S, \mathcal{O}_S) = 0$, and $L$ is big, nef and indivisible in $\mathrm{Pic}(S)$, with $L^2 = 2p - 2$. Write $p = 2l + \varepsilon$, with $\varepsilon = 0, 1$ and $l \in \mathbf{N}$. If $E' \subseteq \mathbf{P}^p$ is an elliptic normal curve of degree $p + 1$, set $L_{E'} = \mathcal{O}_{E'}(1)$. Consider two *general* line bundles $L_1, L_2 \in \mathrm{Pic}^2(E')$ with $L_1 \neq L_2$. Denote by $\Sigma'_i$ the rational normal scroll of degree $p - 1$ in $\mathbf{P}^p$ described by the secant lines to $E'$ spanned by the divisors in $|L_i|$, for $1 \leqslant i \leqslant 2$. One has

$$\Sigma'_i \cong \begin{cases} \mathbf{P}^1 \times \mathbf{P}^1 & \text{if } p = 2l+1 \text{ is odd,} \\ \mathbf{F}_1 & \text{if } p = 2l \text{ is even.} \end{cases}$$

The surfaces $\Sigma'_1$ and $\Sigma'_2$ are $\mathbf{P}^1$-bundles on $\mathbf{P}^1$. Denote by $\sigma_i$ and $f_i$ a minimal section and a fiber of the ruling of $\Sigma'_i$, respectively, so that $\sigma_i^2 = \varepsilon - 1$, $f_i^2 = 0$, and

(3.3.1) $$L_{\Sigma'_i} := \mathcal{O}_{\Sigma'_i}(1) \simeq \mathcal{O}_{\Sigma'_i}(\sigma_i + lf_i), \quad \text{for} \quad 1 \leqslant i \leqslant 2.$$

By [7, Thm. 1], $\Sigma'_1$ and $\Sigma'_2$ intersect transversally along $E'$, which is anticanonical on $\Sigma'_i$, i.e.

(3.3.2) $$E' \sim -K_{\Sigma'_i} \sim 2\sigma_i + (3 - \varepsilon)f_i \quad \text{for} \quad 1 \leqslant i \leqslant 2,$$

where $\sim$ is the linear equivalence. Hence $\Sigma' = \Sigma'_1 \cup \Sigma'_2$ has normal crossings and its dualising sheaf $\omega_{\Sigma'}$ is trivial. Set $L_{\Sigma'} := \mathcal{O}_{\Sigma'}(1)$. The first cotangent sheaf $T^1_{\Sigma'}$ (cf. [4, § 1]) is the degree 16 line bundle on $E'$

(3.3.3) $$T^1_{\Sigma'} \cong \mathbf{N}_{E'/\Sigma'_1} \otimes \mathbf{N}_{E'/\Sigma'_2} \cong L_{E'}^{\otimes 4} \otimes (L_1 \otimes L_2)^{\otimes(3-2l-\varepsilon)},$$

where $\mathbf{N}_{E'/\Sigma'_i}$ is the normal sheaf of $E'$ in $\Sigma'_i$, for $i = 1, 2$, and the last isomorphism comes from (3.3.1) and (3.3.2).

The surface $\Sigma'$ is a flat limit of smooth $K3$ surfaces in $\mathbf{P}^p$. Namely, if $\mathcal{H}_p$ is the component of the Hilbert scheme of surfaces in $\mathbf{P}^p$ containing $K3$ surfaces $S$ such that $(S, \mathcal{O}_S(1)) \in \mathcal{K}_p$, then $\Sigma'$ sits in $\mathcal{H}_p$ and, for general choices of $E', L_1, L_2$, the Hilbert scheme is smooth at $\Sigma'$ (see [7]). The fact that $T^1_{\Sigma'}$ is non-trivial implies that $\Sigma'$ is not a semi–stable limit of $K3$ surfaces: indeed, the total space $\mathfrak{S}'$ of every flat deformation of $\Sigma'$ to $K3$ surfaces in $\mathcal{H}_p$ is singular along a divisor $T \in |T^1_{\Sigma'}|$ (cf. [4, Prop. 1.11 and § 2]). More precisely (see again [7] for details), if

(3.3.4)
$$\begin{array}{ccccc} \Sigma' & \hookrightarrow & \mathfrak{S}' & \hookrightarrow & \mathbf{D} \times \mathbf{P}^p \\ \downarrow & & \downarrow \pi' & \searrow^{\mathrm{pr}_2} & \downarrow \\ 0 & \hookrightarrow & \mathbf{D} & & \mathbf{P}^p \end{array}$$

is a deformation of $\Sigma'$ in $\mathbf{P}^p$ whose general member is a smooth $K3$ surface, then $\mathfrak{S}'$ has double points at the support of a divisor $T \in |T^1_{\Sigma'}|$ associated to the tangent direction to $\mathcal{H}_p$ at $\Sigma'$ determined by the deformation (3.3.4), via the map

(3.3.5)
$$T_{\mathcal{H}_p, \Sigma'} \cong H^0(\Sigma', \mathbf{N}_{\Sigma'/\mathbf{P}^p}) \to H^0(T^1_{\Sigma'})$$

induced by the surjective sheaf map

$$\mathbf{N}_{\Sigma'/\mathbf{P}^p} \to T^1_{\Sigma'}.$$

If $T$ is reduced (this is the case if (3.3.4) is general enough), then the tangent cone to $\mathfrak{S}$ at each of the 16 points of $T$ has rank 4. If $\Pi : \mathfrak{S} \to \mathfrak{S}'$ is a small resolution of singularities, one gets a semistable degeneration $\pi := \pi' \circ \Pi : \mathfrak{S} \longrightarrow \mathbf{D}$ of $K3$ surfaces, with central fiber $\Sigma := \Sigma_1 \cup \Sigma_2$, where $\Sigma_i = \Pi^{-1}(\Sigma'_i)$, for $i = 1, 2$, and still $\omega_\Sigma \cong \mathcal{O}_\Sigma$. The polarising line bundle on $\mathfrak{S}$ is $\mathcal{L} := \Pi^*(\mathrm{pr}_2^*(\mathcal{O}_{\mathbf{P}^p}(1)))$.

Set $E = \Sigma_1 \cap \Sigma_2$; then $E \cong E'$ and

(3.3.6)
$$T^1_\Sigma \simeq \mathcal{O}_E.$$

Equation (3.3.6) is a particular case of the following general fact (see [1, 4] for a more general formulation).

**(3.4) Lemma** (Triple Point Formula). *Assume $f : S \to \mathbf{D}$ is semistable. Let $Q, Q'$ be smooth irreducible components of $S_0$ intersecting along the* double curve *$R$. Then*

$$\deg(\mathbf{N}_{R/Q}) + \deg(\mathbf{N}_{R/Q'}) + \mathrm{Card} \left\{ \begin{array}{c} \textit{triple points of } S_0 \\ \textit{along } R \end{array} \right\} = 0,$$

*where a* triple point *is the intersection $R \cap Q''$ with a component $Q''$ of $S_0$ different from $Q, Q'$.*

## 4 – Limit linear systems

Consider a semistable degeneration $f : \mathcal{S} \to \mathbf{D}$ as in §3. Suppose there is a polarising, ample, fixed components free line bundle $\mathcal{L}$ on the total space $\mathcal{S}$, such that $h^0(S_t, L_t)$ is constant for $t \in \mathbf{D}$.

Consider the subscheme $\mathrm{Hilb}(\mathcal{L})$ of the *relative Hilbert scheme* of curves of $\mathcal{S}$ over $\mathbf{D}$, which is the Zariski closure of the set of all curves $C \in |L_t|$, for $t \in \mathbf{D}^*$. Assume that $\mathrm{Hilb}(\mathcal{L})$ is a component of the relative Hilbert scheme, a condition satisfied if $\mathrm{Pic}(S_t)$ has no torsion for $t \in \mathbf{D}^*$. One has a natural projection morphism $\varphi : \mathrm{Hilb}(\mathcal{L}) \to \mathbf{D}$, which is a projective bundle over $\mathbf{D}^*$, because $\mathrm{Hilb}(\mathcal{L})$ is isomorphic to $\mathfrak{P} := \mathbf{P}(f_*(\mathcal{L}))$ over $\mathbf{D}^*$. We call the fibre of $\varphi$ over 0 the *limit linear system* of $|L_t|$ as $t \in \mathbf{D}^*$ tends to 0, and we will denote it by $\mathfrak{L}$.

**(4.1) Remark.** In general, *the limit linear system is not a linear system.* One would be tempted to say that $\mathfrak{L}$ coincides with $|L_0|$. This is the case if $S_0$ is irreducible, but it is in general no longer true when $S_0$ is reducible. Indeed in this case, there may be non–zero sections of $L_0$ whose zero–locus contains some irreducible component of $S_0$, and accordingly points of $|L_0|$ which do not correspond to points in the Hilbert scheme of curves. This is related to the non–uniqueness of the limit line bundle and to the twisting procedure mentioned in § 3 (see, e.g., Example (4.2) below). We will come back to this in a while.

The variety $\text{Hilb}(\mathcal{L})$ is birational to $\mathfrak{P}$, and $\mathfrak{L}$ is a suitable degeneration of the projective space $|L_t|$, $t \in \Delta^*$.

As the following example shows, in passing from $\mathfrak{P}$ to $\text{Hilb}(\mathcal{L})$, one has to perform a series of blow–ups along smooth centres contained in the central fibre, which correspond to spaces of non–trivial sections of some (twisted) line bundles which vanish on divisors contained in the central fibre. The exceptional divisors one gets in this way give rise to components of $\mathfrak{L}$, and may be identified with birational modifications of sublinear systems of twisted linear systems restricted to $S_0$.

**(4.2) Example** (See [9])**.** Consider a family of degree $d$ surfaces $f : \mathcal{S} \to \mathbf{D}$ arising, as in Example (3.2), from a linear pencil generated by a general surface $S_\infty$ and by $S_0 = S'' \cup P$, where $P$ is a plane and $S''$ a general surface of degree $d-1$. One has a semistable model $\tilde{f} : \tilde{\mathcal{S}} \to \mathbf{D}$ as described in Example (3.2), from which we keep the notation.

Let $\mathcal{L}$ be the polarising hyperplane bundle. The component $\text{Hilb}(\mathcal{L})$ of the Hilbert scheme is gotten from the projective bundle $\mathbf{P}(\mathcal{L})$, by blowing up the point of the central fibre $|\mathcal{O}_{S_0}(1)|$ corresponding to the 1–dimensional space of non–zero sections vanishing on the plane $P$. The limit linear system $\mathfrak{L}$ is the union of $\mathfrak{L}_1$, the blown–up $|\mathcal{O}_{S_0}(1)|$, and of the exceptional divisor $\mathfrak{L}_2 \cong \mathbf{P}^3$, identified as the twisted linear system $|\mathcal{O}_{\tilde{S}_0}(1) \otimes \mathcal{O}_{\tilde{S}}(-\tilde{P})|$, whose aspects have been described in Example (3.2).

The components $\mathfrak{L}_1$ and $\mathfrak{L}_2$ of $\mathfrak{L}$ meet along the exceptional divisor $\mathfrak{E} \cong \mathbf{P}^2$ of the morphism $\mathfrak{L}_1 \to |\mathcal{O}_{S_0}(1)|$. The elements of $\mathfrak{E} \subseteq \mathfrak{L}_1$ identify as the points of $|\mathcal{O}_R(1)| \cong |\mathcal{O}_P(1)|$, whereas the plane $\mathfrak{E} \subseteq \mathfrak{L}_2$ is the set of elements $C \in |\mathcal{O}_{\tilde{P}}(d) \otimes \mathcal{O}_{\tilde{P}}(-\sum_{i=1}^{d(d-1)} E_i)|$ containing the proper transform $\hat{R} \cong R$ of $R$ on $\tilde{P}$. The corresponding element of $|\mathcal{O}_R(1)|$ is cut out on $\hat{R}$ by the further component of $C$, which is the pull–back to $\tilde{P}$ of a line in $P$.

All this will be made explicit in the case $d = 4$ in §7.2.

## 5 – Severi varieties and their limits

Let $f : \mathcal{S} \to \mathbf{D}$ be a (not necessarily semistable) family as in §3, polarised by a line bundle $\mathcal{L}$ on $\mathcal{S}$. Fix a non–negative integer $\delta$, and consider the locally closed subset $\mathring{V}_\delta(\mathcal{S}, \mathcal{L})$ of $\text{Hilb}(\mathcal{L})$ formed by all curves $D \in |L_t|$, for $t \in \mathbf{D}^*$, such that $D$ is irreducible, nodal, and has exactly $\delta$ nodes. Define $V_\delta(\mathcal{S}, \mathcal{L})$ (resp. $V_\delta^{\text{cr}}(\mathcal{S}, \mathcal{L})$) as the Zariski closure of $\mathring{V}_\delta(\mathcal{S}, \mathcal{L})$ in $\text{Hilb}(\mathcal{L})$ (resp. in $\mathfrak{P} = \mathbf{P}(f_*(\mathcal{L}))$). This is the *relative Severi variety* (resp. the *crude relative Severi variety*). Sometimes we may write $\mathring{V}_\delta$, $V_\delta$, and $V_\delta^{\text{cr}}$, rather than $\mathring{V}_\delta(\mathcal{S}, \mathcal{L})$, $V_\delta(\mathcal{S}, \mathcal{L})$ and $V_\delta^{\text{cr}}(\mathcal{S}, \mathcal{L})$, respectively.

There is a natural map $f_\delta : V_\delta \to \mathbf{D}$. If $t \in \mathbf{D}^*$, the fibre $V_{\delta,t}$ of $f_\delta$ over $t$ is the *Severi variety* $V_\delta(S_t, L_t)$ of $\delta$–nodal curves in the linear system $|L_t|$ on $S_t$. The degree of $V_{\delta,t}$ as a subvariety of $|L_t|$ is independent on $t \in \mathbf{D}^*$, and will be denoted by $d_\delta(\mathcal{L})$ (or simply by $d_\delta$). Let $\mathfrak{V}_\delta(\mathcal{S}, \mathcal{L})$ (or simply $\mathfrak{V}_\delta$) be the central fibre of $f_\delta : V_\delta \to \mathbf{D}$; it is the *limit Severi variety* of $V_\delta(S_t, \mathcal{L}_t)$ as $t \in \mathbf{D}^*$ tends to 0. This is a subscheme of the limit linear system $\mathfrak{L}$, which has been studied by various authors (in the present setting by Z. Ran [18, 19, 20], then by L. Caporaso, J. Harris [1, 2], more recently by C. Galati and C. Galati and A. Knutsen [9, 11] following Z. Ran; this theme is however quite classical as B. Segre's paper [23] shows, and it appears in disguise in work of other authors, like S. Kleiman [13]).

In a similar way, one defines the *crude limit Severi variety* $\mathfrak{V}^{\mathrm{cr}}_\delta(\mathcal{S}, \mathcal{L})$ (or $\mathfrak{V}^{\mathrm{cr}}_\delta$), sitting in $|L_0|$.

**(5.1) Remark.** For $t \in \mathbf{D}^*$, the *expected dimension* of the Severi variety $V_\delta(S_t, \mathcal{L}_t)$ is $\dim(|L_t|) - \delta$ (see [24, Thm. 6.3]). We will always assume that the dimensions of (all components of) $V_\delta(S_t, L_t)$ equal the expected dimension for all $t \in \mathbf{D}^*$. This is a strong assumption, implying that $\mathring{V}_\delta$ is smooth; it will be satisfied in our applications.

Let us come now to the description of the limit Severi variety, under the assumption that the family $f : \mathcal{S} \to \mathbf{D}$ is semistable. We will freely use the notation introduced above and follow the approach in [9, 11]. We will suppose that the central fibre $S_0$ has smooth irreducible components $Q_1, \ldots, Q_r$, with double curves $R_1, \ldots, R_s$. We will assume, in addition, that there are $q$ exceptional curves $E_1, \ldots, E_q$ on $S_0$, arising from a small resolution of an original family with singular total space, as discussed in §3.

**(5.2) Notation.** Let $\underline{\mathbf{N}}$ be the set of sequences $\underline{\tau} = (\tau_m)_{m \geqslant 2}$ of non–negative integers with only finitely many non–vanishing terms. Define two maps $\nu$, $\mu : \underline{\mathbf{N}} \to \mathbf{N}$ as follows

$$\nu(\underline{\tau}) = \sum\nolimits_{m \geqslant 2} \tau_m \cdot (m-1), \quad \text{and} \quad \mu(\underline{\tau}) = \prod\nolimits_{m \geqslant 2} m^{\tau_m}.$$

Given a $p$-tuple $\underline{\boldsymbol{\tau}} = (\underline{\tau}_1, \ldots, \underline{\tau}_p) \in \underline{\mathbf{N}}^p$, set

$$\nu(\underline{\boldsymbol{\tau}}) = \nu(\underline{\tau}_1) + \cdots + \nu(\underline{\tau}_p), \quad \text{and} \quad \mu(\underline{\boldsymbol{\tau}}) = \mu(\underline{\tau}_1) \cdots \mu(\underline{\tau}_p),$$

defining two maps $\nu$, $\mu : \underline{\mathbf{N}}^p \to \mathbf{N}$. Given $\boldsymbol{\delta} = (\delta_1, \ldots, \delta_r) \in \mathbf{N}^r$, set

$$|\boldsymbol{\delta}| := \delta_1 + \cdots + \delta_r.$$

Given a subset $I \subseteq \{1, \ldots, q\}$, $|I|$ will denote its cardinality.

**(5.3) Definition.** *Consider a divisor $W$ on $\mathcal{S}$, supported on the central fibre $S_0$, such that the twist $\mathcal{L}(-W)$ is centrally effective, i.e., all the aspects of $\mathcal{L}(-W)$ on the components $Q_1, \ldots, Q_r$ of $S_0$ are effective. Fix*

$$\boldsymbol{\delta} \in \mathbf{N}^r, \quad \underline{\boldsymbol{\tau}} \in \underline{\mathbf{N}}^p \quad \text{and} \quad I \subseteq \{1, \ldots, r\}.$$

*Let $\mathring{V}(W, \boldsymbol{\delta}, I, \underline{\boldsymbol{\tau}})$ be the Zariski locally closed subset in $|\mathcal{L}(-W) \otimes \mathcal{O}_{S_0}|$ parametrizing curves $D$ such that:*
(i) *$D$ contains no double curve $R_l$, with $l \in \{1, \ldots, s\}$, and no triple point of $S_0$;*
(ii) *$D$ contains $E_i$, with multiplicity 1, if and only if $i \in I$, and has a node on it;*
(iii) *off the singular locus of $S_0$ the curve $D - \sum_{i \in I} E_i$ has only nodes as singularities, and their number is $\delta_s$ on $Q_s$, for $s \in \{1, \ldots, r\}$;*
(iv) *for every $l \in \{1, \ldots, s\}$ and $m \geqslant 2$, there are exactly $\tau_{l,m}$ points on $R_l$, off the intersections with $\sum_{i \in I} E_i$ and off the triple points of $S_0$ along $R_l$, at which $D$ has an $m$-tacnode (see below for the definition), with reduced tangent cone equal to the tangent line of $R_l$ there.*

*We let $V(W, \boldsymbol{\delta}, I, \underline{\boldsymbol{\tau}})$ be the Zariski closure of $\mathring{V}(W, \boldsymbol{\delta}, I, \underline{\boldsymbol{\tau}})$ in $|\mathcal{L}(-W) \otimes \mathcal{O}_{S_0}|$.*

Recall that an *m-tacnode* is an $A_{2m-1}$-double point, i.e., a plane curve singularity locally analytically isomorphic to the curve of $\mathbf{C}^2$ defined by the equation $y^2 = x^{2m}$ at the origin. Condition (iv) in Definition (5.3) requires that $D$ is a Cartier divisor on $S_0$, having $\tau_{l,m}$ $m$–th order tangency points with the curve $R_l$, at points of $R_l$ which are neither points on $\sum_{i \in I} E_i$ nor triple points of $S_0$.

**(5.4) Notation.** Rather than using the notation $V(W, \boldsymbol{\delta}, I, \boldsymbol{\tau})$, we may sometimes use a more expressive one like, e.g., $V(W, \delta_{Q_1} = 2, E_1, \tau_{R_1,2} = 1)$ for the variety parametrizing curves in $|\mathcal{L}(-W) \otimes \mathcal{O}_{S_0}|$, with two nodes on $Q_1$, one simple tacnode along $R_1$, and containing the exceptional curve $E_1$. Similarly $V(\delta_{Q_1} = 1)$ is the variety parametrizing curves in $|\mathcal{L} \otimes \mathcal{O}_{S_0}|$, with only one node on $Q_1$, etc.

**(5.5) Proposition** ([9, 11], see [VIII])**.** *Let $W, \boldsymbol{\delta}, I, \boldsymbol{\tau}$ be as above, and set*

$$\delta = |\boldsymbol{\delta}| + |I| + \nu(\boldsymbol{\tau}).$$

*Let $V$ be an irreducible component of $V(W, \boldsymbol{\delta}, I, \boldsymbol{\tau})$. If*
*(i) the linear system $|\mathcal{L}(-W) \otimes \mathcal{O}_{S_0}|$ has the same dimension as $|\mathcal{L}_t|$ for $t \in \mathbf{D}^*$, and*
*(ii) $V$ has (the expected) codimension $\delta$ in $|\mathcal{L}(-W) \otimes \mathcal{O}_{S_0}|$,*
*then $V$ is an irreducible component of multiplicity $\mu(V) := \mu(\boldsymbol{\tau})$ of the limit Severi variety $\mathfrak{V}_\delta(S, \mathcal{L})$.*

**(5.6) Remark.** Same assumptions as in Proposition (5.5). If there is at most one tacnode (i.e. all $\tau_{l,m}$ but possibly one vanish, and this is equal to 1), the relative Severi variety $V_\delta$ is smooth at the general point of $V$ (see [9, 11]), and thus $V$ belongs to only one irreducible component of $V_\delta$. There are other cases in which such a smoothness property holds (see [1] or, in this volume, [VI]).

If $V_\delta$ is smooth at the general point $D \in V$, the multiplicity of $V$ in the limit Severi variety $\mathfrak{V}_\delta$ is the minimal integer $m$ such that there are local analytic $m$–multisections of $V_\delta \to \mathbf{D}$, i.e. analytic smooth curves in $V_\delta$, passing through $D$ and intersecting the general fibre $V_{\delta,t}$, $t \in \mathbf{D}^*$, at $m$ distinct points.

Proposition (5.5) does not provide a complete picture of the limit Severi variety. For instance, curves passing through a triple point of $S_0$ could (and in fact do; see [4], and §8.1 below) play a role in this limit. It would be desirable to know that one can always obtain a semistable model of the original family, where *every* irreducible component of the limit Severi variety is realized as a family of curves of the kind stated in Definition (5.3).

**(5.7) Definition.** *Let $f : \mathcal{S} \to \mathbf{D}$ be a semistable family as in §3, with a polarising line bundle $\mathcal{L}$, and $\delta$ a positive integer. The* regular part of the limit Severi variety $\mathfrak{V}_\delta(\mathcal{S}, \mathcal{L})$ *is the cycle in the limit linear system $\mathfrak{L} \subseteq \mathrm{Hilb}(\mathcal{L})$ defined as*

$$(5.7.1) \qquad \mathfrak{V}_\delta^{\mathrm{reg}}(\mathcal{S}, \mathcal{L}) := \sum_W \sum_{|\boldsymbol{\delta}|+|I|+\nu(\boldsymbol{\tau})=\delta} \mu(\boldsymbol{\tau}) \cdot \left( \sum_{V \in \mathrm{Irr}^\delta(V(W, \boldsymbol{\delta}, I, \boldsymbol{\tau}))} V \right),$$

*(sometimes simply denoted by $\mathfrak{V}_\delta^{\mathrm{reg}}$) where:*
*(i) $W$ varies among all effective divisors on $S$ supported on the central fibre $S_0$, such that $\mathcal{L}(-W)$ is centrally effective and $h^0(\mathcal{L}_0(-W)) = h^0(L_t)$ for $t \in \mathbf{D}^*$;*
*(ii) $\mathrm{Irr}^\delta(Z)$ denotes the set of all codimension $\delta$ irreducible components of a scheme $Z$.*

Proposition (5.5) asserts that the cycle $Z(\mathfrak{V}_\delta) - \mathfrak{V}_\delta^{\mathrm{reg}}$ is effective where $Z(\mathfrak{V}_\delta)$ is the cycle associated to $\mathfrak{V}_\delta$. We call the irreducible components of the support of $\mathfrak{V}_\delta^{\mathrm{reg}}$ the *regular components* of the limit Severi variety.

Let now $\tilde{f} : \tilde{\mathcal{S}} \to \mathbf{D}$ be a semistable model of a (not necessarily semistable) degeneration $f : \mathcal{S} \to \mathbf{D}$, and $\tilde{\mathcal{L}}$ the pull–back on $\tilde{\mathcal{S}}$ of a polarising line bundle $\mathcal{L}$ on $\mathcal{S}$. There is a natural map $\mathrm{Hilb}(\tilde{\mathcal{L}}) \to \mathrm{Hilb}(\mathcal{L})$, which induces a morphism $\phi : \tilde{\mathfrak{L}} \to |L_0|$.

**(5.8) Definition.** *The semistable model $\tilde{f} : \tilde{\mathcal{S}} \to \mathbf{D}$ is a $\delta$–good model of $f : \mathcal{S} \to \mathbf{D}$ (or simply a good model, if it is clear which $\delta$ we are considering), if the following equality of cycles holds*

$$\phi_*\bigl(\mathfrak{V}_\delta^{\mathrm{reg}}(\tilde{\mathcal{S}}, \tilde{\mathcal{L}})\bigr) = \mathfrak{V}_\delta^{\mathrm{cr}}(\mathcal{S}, \mathcal{L}).$$

Note that the cycle $\mathfrak{V}_\delta^{\mathrm{cr}}(\mathcal{S}, \mathcal{L}) - \phi_*\bigl(\mathfrak{V}_\delta^{\mathrm{reg}}(\tilde{\mathcal{S}}, \tilde{\mathcal{L}})\bigr)$ is effective. The family $f : \mathcal{S} \to \mathbf{D}$ is said to be $\delta$–*well behaved* (or simply *well behaved* if $\delta$ is understood) if it has a $\delta$-good model.

**(5.9) Remark.** Suppose that $f : \mathcal{S} \to \mathbf{D}$ is $\delta$–well behaved, with $\delta$–good model $\tilde{f} : \tilde{\mathcal{S}} \to \mathbf{D}$. It is possible that some components in $\mathfrak{V}_\delta^{\mathrm{reg}}(\tilde{\mathcal{S}}, \tilde{\mathcal{L}})$ are contracted by $\phi$ to varieties of smaller dimension, and therefore their push–forwards to $\mathfrak{V}_\delta^{\mathrm{cr}}(\mathcal{S}, \mathcal{L})$ are zero. Hence these components of $\mathfrak{V}_\delta(\tilde{\mathcal{S}}, \tilde{\mathcal{L}})$ are *not visible* in $\mathfrak{V}_\delta^{\mathrm{cr}}(\mathcal{S}, \mathcal{L})$. They are however usually visible in the crude limit Severi variety of some other model $f' : \mathcal{S}' \to \mathbf{D}$, obtained from $\tilde{\mathcal{S}}$ via an appropriate twist of $\mathcal{L}$. The central fibre $S'_0$ is then a flat limit of $S_t$, as $t \in \mathbf{D}^*$ tends to $0$, different from $S_0$ (a situation met in Example (3.2)).

**(5.10) Conjecture** (See [4])**.** *Let $f : \mathcal{S} \to \mathbf{D}$ be a degeneration of surfaces, endowed with a line bundle $\mathcal{L}$ as above, and $\delta$ a positive integer. Then $f : \mathcal{S} \to \mathbf{D}$ is $\delta$–well behaved.*

The local computations in [9, 11] provide a criterion for being well behaved:

**(5.11) Proposition.** *Assume there is a semistable model $\tilde{f} : \tilde{\mathcal{S}} \to \mathbf{D}$ of $f : \mathcal{S} \to \mathbf{D}$, with a limit linear system $\tilde{\mathfrak{L}}$ which does not contain, in codimension $\delta + 1$, curves of the following types:*
*(i) curves containing double curves of $\tilde{S}_0$;*
*(ii) curves passing through a triple point of $\tilde{S}_0$;*
*(iii) non–reduced curves.*
*If in addition, for $W, \boldsymbol{\delta}, I, \boldsymbol{\tau}$ as in Definition (5.3), every irreducible component of $V(W, \boldsymbol{\delta}, I, \boldsymbol{\tau})$ has the expected codimension in $|L_0(-W)|$, then $\tilde{f} : \tilde{\mathcal{S}} \to \mathbf{D}$ is $\delta$–well behaved.*

**(5.12) Remark.** While there are good reasons to believe that Conjecture (5.10) has an affirmative answer, one cannot expect that, in the algebraic category, there is a good model of a given family, which *universally* works for all polarising line bundles $\mathcal{L}$. This in fact could require infinitely many birational modifications of the total space.

On the other hand, it is hopeless to ask for a semistable model on which all irreducible components of the limit Severi variety would parametrize nodal curves (i.e., such that there are no tacnodes as in Item (iv) of Definition (5.3)): we explain this in [B], together with C. Galati.

The three final sections are devoted to examples which will hopefully clarify the considerations and results in Section 5.

# 6 – Example: Limits of 1–nodal plane sections for surfaces in $\mathbf{P}^3$

Consider a degeneration of a general degree $d$ surface $S$ in $\mathbf{P}^3$ to a reducible surface $S' \cup S''$, as in Example (3.2), from which we keep the notation. Here $S_0 = S' \cup S''$, and $R = S' \cap S''$ is a smooth curve of degree $h(d - h)$. Then in the semistable model $\tilde{f} : \tilde{\mathcal{S}} \to \mathbf{D}$ constructed in Example (3.2), $\tilde{S}_0 = Q_1 \cup Q_2$, and $Q_1$ and $Q_2$ intersect transversally along a curve which may be identified with $R$.

**(6.1)** Assume $h \geqslant 2$. Then the only twist which is centrally effective is the trivial one and $\mathfrak{V}_1^{\text{reg}}$ consists of $3 + q$ components (remember that $q = hd(d-h)$), precisely:
(1) $V(\delta_{Q_i} = 1)$, with $i = 1, 2$, with multiplicity 1;
(2) $V(\tau_{R,2} = 1)$, with multiplicity 2;
(3) $V(E_i)$, for $i = 1, \ldots, q$, with multiplicity 1.

This tells us that in the above degeneration, the limit of the dual variety $\check{S}$, which sits in $\mathcal{L}_0 = \check{\mathbf{P}}^3$ and coincides with the crude limit of the Severi variety $\mathfrak{V}_1^{\text{cr}}$, contains:
(i) the dual varieties $\check{S}'$ and $\check{S}''$ with multiplicity 1 (corresponding to (1));
(ii) the dual variety $\check{R}$ with multiplicity 2 (corresponding to (2));
(iii) the $q$ planes $x_i^\perp \subseteq \check{\mathbf{P}}^3$, orthogonal to the points $x_i$, for $i = 1, \ldots, q$, with multiplicitiy 1 (corresponding to (3)).

One computes:
(a) $\deg(\check{S}') = h(h-1)^2$ and $\deg(\check{S}'') = (d-h)(d-h-1)^2$;
(b) $\deg(\check{R}) = h(d-h)(d-2)$.
So adding up, we see that the sum of the degrees of the components in (i)–(iii) above, counted with the appropriate multiplicities, is $d(d-1)^2$, hence the components in (1)–(3) exhaust the limit Severi variety $\mathfrak{V}_1$, i.e., the above degeneration is 1–well behaved.

**(6.2)** If $h = 1$, one needs to change the analysis, due to the presence of the non–trivial twist by $Q_1$ which is centrally effective. In this case $\mathfrak{V}_1^{\text{reg}}$ consists of the following components:
($\alpha$) $V(\delta_{Q_2} = 1)$, with multiplicity 1 (note that $V(\delta_{Q_1} = 1)$ is zero);
($\beta$) $V(Q_1, \delta_{Q_1} = 1)$, with multiplicity 1;
($\gamma$) $V(\tau_{R,2} = 1)$, with multiplicity 2;
($\delta$) $V(E_i)$, for $i = 1, \ldots, q$, with multiplicity 1.
However, only the components in $(\alpha), (\gamma)$ and $(\delta)$ contribute to the crude limit of the Severi variety $\mathfrak{V}_1^{\text{cr}}$, and exhaust it (same computation as above), because the component in $(\beta)$ is not visible in the crude limit, since it is contracted to a point (see Remark (5.9)). In conclusion also this degeneration is 1–well behaved and the limit of $\check{S}$ is the union of $\check{S}'$ plus $d(d-1)$ planes plus the dual of the plane curve $R$ counted with multiplicity 2.

# 7 – Example: Limits of 3–nodal plane sections for quartic surfaces in $\mathbf{P}^3$

The number of 3–nodal plane sections of a general quartic surface $S$ in $\mathbf{P}^3$ (or of triple points of the dual $\check{S}$) given by Salmon's formula (2.0.1) is 3200.

## 7.1 – Degeneration to two quadrics

Let us recover this first by considering a degeneration to two quadrics: thus we will look at the case $h = d - h = 2$ in Example (3.2) (from which we keep the notation; this is the same as the degeneration in Example (3.3) for $p = 3$).

In this case $R$ is an elliptic quartic curve, the semistable degeneration has central fibre $Q_1 \cup Q_2$, with $R = Q_1 \cap Q_2$, $Q_1$ is a quadric, and $Q_2$ a quadric blown up at the complete intersection points $x_1, \ldots, x_{16}$ of $R$ with a general quartic surface $S$ (see Figure 1).

For the computation, one has to consider the various cases for $V(W, \boldsymbol{\delta}, I, \boldsymbol{\tau})$ as in Definition (5.3). No non–trivial twist is centrally effective, so $W = 0$. If $|\boldsymbol{\delta}| = 3$, then $V(W, \boldsymbol{\delta}, I, \boldsymbol{\tau})$ is empty, because one cannot have two (or more) singular points on a plane section of a smooth quadric. Hence:

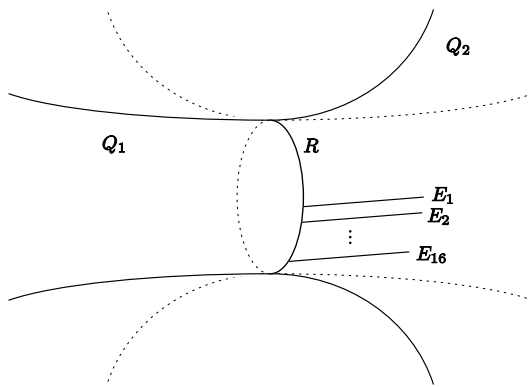

Figure 1: Degeneration into two quadrics

- $0 \leqslant |\boldsymbol{\delta}| \leqslant 2$, and if $|\boldsymbol{\delta}| = 2$ then $\delta_1 = \delta_2 = 1$;
- $\tau = (\tau_2, \tau_3, \tau_4)$ and $|I| + 2\tau_2 + 3\tau_3 + 4\tau_4 \leqslant 4$;
- $|\boldsymbol{\delta}| + |I| + \tau_2 + 2\tau_3 + 3\tau_4 = 3$.

Therefore the cases to be discussed are the following.

*Case 1:* $\delta_1 = \delta_2 = 1$, $|I| = 1$.

For each $i = 1, \ldots, 16$, a plane tangent to both $S'$ and $S''$, and passing through $x_i$ in $\mathbf{P}^3$ is entirely determined by the choice of a line passing through $x_i$ in both $Q_1$ and $Q_2$. There are exactly $2 \times 2 = 4$ such choices (two rulings on each quadric). We get $16 \times 4 = 64$ curves.

*Case 2:* $\delta_1 = \delta_2 = 1$, $\tau_2 = 1$.

The intersection $\check{S}' \cap \check{S}'' \cap \check{R}$ consists of $2 \times 2 \times 8 = 32$ points. From this, we want to remove the points corresponding to planes tangent to one quadric at some point lying on $R$. Each of them has multiplicity 2 in the intersection $\check{S}' \cap \check{S}'' \cap \check{R}$ by the:

**(7.1) Lemma** (Lemma 3.5, [4]). *Let $R$ be a smooth, irreducible curve contained in a smooth surface $S$ in $\mathbf{P}^3$. Let $\check{R}_S$ be the irreducible curve in $\check{\mathbf{P}}^3$ parametrizing planes tangent to $S$ along $R$. Then the dual varieties $\check{S}$ and $\check{R}$ both contain $\check{R}_S$, and do not intersect transversally at its general point.*

*Proof.* Clearly $\check{R}_S$ is contained in $\check{S} \cap \check{R}$. If either $\check{S}$ or $\check{R}$ are singular at the general point of $\check{R}_S$, there is nothing to prove. Assume that $\check{S}$ and $\check{R}$ are both smooth at the general point of $\check{R}_S$. One has to show that they are tangent there. Let $x \in R$ be general. Let $H$ be the tangent plane to $S$ at $x$. Then $H \in \check{R}_S$ is the general point. The biduality theorem (see , e.g., [10, Example 16.20]) says that the tangent planes to $\check{S}$ and $\check{R}$ at $H$ both coincide with the set of planes in $\mathbf{P}^3$ containing $x$, hence the assertion. □

The curve parametrizing planes tangent to $S'$ at some point of $R$ has degree 4 (to see this, intersect $R$ with a general polar plane of $S'$, cf. [**??**, Appendix **??**]). Hence this curve intersects $\check{S}''$ at 8 distinct points (without multiplicities). Eventually, we have to remove $2 \times 2 \times 8$ points from $\check{S}' \cap \check{S}'' \cap \check{R}$, and there does not remain anything.

*Case 3:* $\delta_i = 1$ for an $i = 1, 2$, $|I| = 2$.

For each $i = 1, 2$ and each set $\{r, s\} \subseteq I$, $r \neq s$, we count points of intersection of $\check{Q}_i$ with the line of $\check{\mathbf{P}}^3$ orthogonal to the line $\langle x_r, x_s \rangle$. There are $2 \times \binom{16}{2} \times 2 = 480$ such points.

*Case 4:* $\delta_i = 1$ for an $i = 1, 2$, $|I| = 1$, $\tau_2 = 1$.

For each $i = 1, 2$ and $r = 1, \ldots, 16$, there are two lines on $Q_i$ through the point $x_r$. They determine two pencils of planes tangent to $Q_i$ and passing through $x_r$.

The planes of each of these pencils cut out a $g_2^1$ on $R$: the intersection of such a plane with $R$ consists of 4 points, two of which are fixed (the two intersection points of the chosen line through $x_r$ with $R$).

Since a $g_2^1$ on an elliptic curve has 4 double points, we get $2 \times 4$ planes satisfying our conditions for each $i$ and $r$. These give a contribution of $2 \times 16 \times 2 \times 8 = 512$ to our computation.

*Case 5: $\delta_i = 1$ for an $i = 1, 2$, $\tau_2 = 2$.*

Planes satisfying these conditions are spanned by two lines belonging respectively to the two rulings of $Q_i$, both tangent to $R$. By genericity, two such lines do not meet on $R$.

Each ruling of $Q_i$ cuts out a $g_2^1$ on $R$, and thus contains exactly 4 tangent lines to $R$ (since a $g_2^1$ on $R$ has 4 ramification points). We thus get 16 planes satisfying our conditions for each $i$. This contributes $2 \times 4 \times 16 = 128$ to our computation.

*Case 6: $\delta_i = 1$ for an $i = 1, 2$, $\tau_3 = 1$.*

This requires planes spanned by two lines on $Q_i$, one of which is tangent to $R$, and meets the other line on $R$ at this tangency point. Such a plane is tangent to $Q_i$ at a point lying on $R$, and therefore cuts a curve with worse singularity than required. Hence they do not contribute to our computation.

*Case 7: $\tau_4 = 1$.*

Since $R$ is a smooth degree 4 elliptic curve in $\mathbf{P}^3$, the set of its hyperplane sections is

$$\{(a, b, c, d) \in \mathrm{Sym}^4(R) \mid a + b + c + d = 0 \in (R, +)\}.$$

Hence, the planes we are looking for are in 1-1 correspondence with 4-torsion points of $R$ and so there are 16 of them, each to be counted with multiplicity 4, so we get a contribution of 64.

*Case 8: $|I| = 1$, $\tau_3 = 1$.*

Fix $r = 1, \ldots, 16$, and consider the projection $\pi$ from $x_r$. It maps $R$ to a smooth cubic plane curve. Since this has 9 flexes, this gives us 9 planes through $x_r$ satisfying the required condition. Each of them has to be counted with multiplicity 3, and this gives a contribution of $16 \times 3 \times 9 = 432$ to our computation.

*Case 9: $|I| = 2$, $\tau_2 = 1$.*

For given $r, s$ distinct in $\{1, \ldots, 16\}$, the pencil of planes through $x_r$ and $x_s$ cuts out (off of $x_r + x_s$) a $g_2^1$ on $R$. It has 4 double points, hence the pencil contains 4 tangent planes to $R$. By genericity, the tangency points are distinct from $x_r$ and $x_s$. This gives a contribution of $\binom{16}{2} \times 2 \times 4 = 960$ to the computation.

*Case 10: $|I| = 3$.*

For each choice of three different indices $r$, $s$ and $l$ in $\{1, \ldots, 16\}$, there is one plane through $x_r$, $x_s$ and $x_l$. So we get a contribution of $\binom{16}{3} = 560$.

Adding up all the non–trivial contributions as in the following table, we get exactly 3200, which is the number predicted by Salmon's formula. This shows that this degeneration is 3–good.

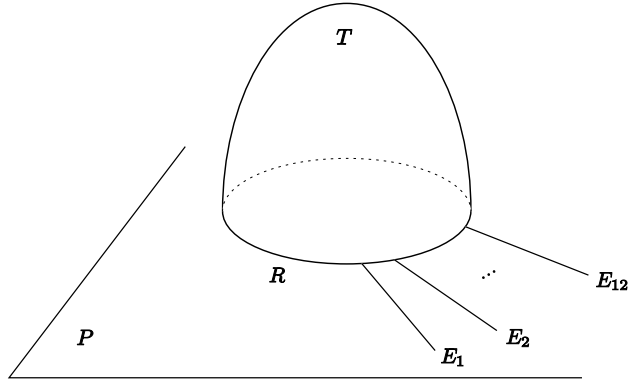

Figure 2: Degeneration into a cubic plus a plane

| Cases   | Contribution |
|---------|-------------:|
| Case 1  | 64           |
| Case 3  | 480          |
| Case 4  | 512          |
| Case 5  | 128          |
| Case 7  | 64           |
| Case 8  | 432          |
| Case 9  | 960          |
| Case 10 | 560          |
| **Total** | **3200**   |

## 7.2 – Degeneration to a cubic and a plane

We will now do the same enumeration as in Section 7.1, this time by considering a degeneration to a cubic and a plane. Thus, we want to count 3–nodal plane sections of a general quartic surface $S$ in $\mathbf{P}^3$, and we will look at the degeneration in Example (3.2), with $d=4$ and $h=1$. This will be instructive, because it will require to make one of the first non–trivial cases of computation of the degree of a Severi variety of plane curves with nodes.

We simplify a bit the notation of Example (3.2). We denote by $T+P_0$ the reducible surface $S_0$ in Example (3.2) (so $T$ is a general cubic, $P_0$ a general plane, with $P_0 \cap T = R$ a general plane cubic), and we denote by $T+P$ the central fibre of the semistable degeneration, where $P$ is the plane $P_0$ blown–up at 12 points $x_1, \ldots, x_{12}$ forming a general divisor in $|\mathcal{O}_R(4)|$ (see Figure 2).

Now we explain in some detail in this case what we stated in Example (3.2) about the limit linear system $\mathfrak{L}$ of the plane sections of $S_t$.

**(7.2) Proposition.** *Let $\tilde{\mathcal{H}}_{\mathbf{P}^3}$ be the blow up of $|\mathcal{O}_{\mathbf{P}^3}(1)| = \check{\mathbf{P}}^3$ at the point corresponding to the plane $P_0 \subseteq \mathbf{P}^3$. Let $D_{P_0} \cong \mathbf{P}^2$ be the exceptional divisor of this blow up.*

*Set $\mathcal{H}_P = |\mathcal{O}_P(4H - \sum_{i=1}^{12} E_i)| \cong \mathbf{P}^3$ and denote by $D_R \subseteq \mathcal{H}_P$ the plane parametrizing strict transforms of plane quartic curves that contain $R$ ($H$ denotes the pull back to $P$ of a general line of $P_0$).*

*Then the limit of the limit linear system $\mathfrak{L}$ of $|\mathcal{O}_{S_t}(1)|$ as $t$ tends to 0 is obtained by glueing $\tilde{\mathcal{H}}_{\mathbf{P}^3}$ and $\mathcal{H}_P$ along $D_{P_0}$ and $D_R$.*

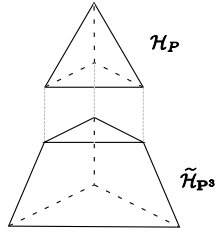

Figure 3: Limit hyperplane linear system for a degeneration into a cubic and plane

*Proof.* One has $\dim(|4H - \sum_{i=1}^{12} E_i|) = 3$ (see Example (3.2) and Remark (7.3) below).

There is only one non–trivial twist of the polarising hyperplane bundle $\mathcal{L}$, i.e., the twist by $P$, which is centrally effective. As we saw in Example (3.2), the twist $\mathcal{L}(-P)$ is trivial on $T$, whereas on $P$ it restricts to $\mathcal{O}_P(H + R) = \mathcal{O}_P(4H - \sum_{i=1}^{12} E_i)$.

Up to a multiplicative constant, there is only one non–zero section of $\mathcal{L}_{|T+P}$ vanishing on $P$, hence it gives no element in $\mathfrak{L}$. To see the corresponding curves in $\mathfrak{L}$ we have to blow up $|\mathcal{L}_{|T+P}| = |\mathcal{O}_{\mathbf{P}^3}(1)|$ at $P_0$. One has

$$T_{\check{\mathbf{P}}^3, [P_0]} = H^0(\mathbf{P}^3, \mathcal{O}_{\mathbf{P}^3}(1))/H^0(\mathbf{P}^3, \mathcal{O}_{\mathbf{P}^3}) \cong H^0(\mathbf{P}^3, \mathcal{O}_{P_0}(1)).$$

This says that $D_{P_0} \cong \mathbf{P}(T_{\check{\mathbf{P}}^3, [P_0]})$ identifies with the 2-dimensional linear system $|\mathcal{O}_{P_0}(1)|$. To give a geometric interpretation of this, let us look at the limits of the curves cut out on $T + P_0$ by a plane $\Pi$ of $\mathbf{P}^3$ tending to $P_0$. The limit curve $\Pi \cap T$ tends to $R$, while the limit curve of $\Pi \cap P_0$ is a line moving freely in $|\mathcal{O}_{P_0}(1)|$.

On the other hand, the plane $D_R$ is

$$|4H - R| \cong |\mathcal{O}_P(H)| \subseteq |4H - \sum_{i=1}^{12} E_i|.$$

There is therefore a natural identification between $D_{P_0}$ and $D_R$, with the glueing as in the statement. □

**(7.3) Remark.** As mentioned in Example (3.2), the surface $P$ is mapped to a quartic with a triple point in $\mathbf{P}^3$ by the linear system $|4H - \sum_{i=1}^{12} E_i|$. This mapping is an isomorphism outside $R$, and contracts $R$ to the triple point.

Indeed, one has an exact sequence

$$0 \to \mathcal{O}_P(H) \to \mathcal{O}_P(4H - \sum_{i=1}^{12} E_i) \to \mathcal{O}_R \to 0.$$

Since $h^1(P, \mathcal{O}_P(H)) = h^1(\mathbf{P}^2, \mathcal{O}_{\mathbf{P}^2}(1)) = 0$, the sequence is exact on global sections, which implies $h^0(P, \mathcal{O}_P(4H - \sum_{i=1}^{12} E_i)) = 4$.

The linear system $|4H - \sum_{i=1}^{12} E_i|$ has no base points and the morphism determined by it is birational onto its image (we leave this to the reader to check). The image of $P$ via this map has degree 4, since $(4H - \sum_{i=1}^{12} E_i)^2 = 4$. The map contracts $R \equiv 3H - \sum_{i=1}^{12} E_i$ because $(4H - \sum_{i=1}^{12} E_i) \cdot (3H - \sum_{i=1}^{12} E_i) = 0$, and its image is a triple point because $R^2 = (3H - \sum_{i=1}^{12} E_i)^2 = -3$.

### 7.2.1 Curves in the hyperplane bundle

According to the structure of the limit linear system $\mathfrak{L}$, one has to compute the limit of trinodal curves appearing in both $\tilde{\mathcal{H}}_{\mathbf{P}^3}$ and in $\mathcal{H}_P$. First we perform the computation for $\tilde{\mathcal{H}}_{\mathbf{P}^3}$.

Here one has to consider the various cases for $V(0, \boldsymbol{\delta}, I, \boldsymbol{\tau})$ as in Definition (5.3), with $I \subseteq \{1, \ldots, 12\}$. Since there cannot be any node on $P$ in this linear system, one has $|\boldsymbol{\delta}| = \delta$, the number of nodes on $T$. Then:
- $0 \leqslant \delta \leqslant 3$ and $\tau = (\tau_2, \tau_3)$;
- $|I| + 2\tau_2 + 3\tau_3 \leqslant 3$;
- $\delta + |I| + \tau_2 + 2\tau_3 = 3$.

Hence the cases to be discussed are the following.

*Case 1:* $\delta = 3$.

Tritangent planes to the general cubic surface $T$ correspond to the triangles contained in $T$. These are 45: the reader may apply Salmon's formula (2.0.1) or, better, directly compute this number as an exercise.

*Case 2:* $\delta = 2$, $|I| = 1$.

Binodal hyperplane sections of $T$ must contain a line. Since $T$ contains 27 lines, there are 27 pencils of such curves. In each of these pencils, and for every $r \in \{1, \ldots, 12\}$, there is exactly one curve passing in addition through the point $x_r$. This gives a contribution of $27 \times 12 = 324$ to the computation.

*Case 3:* $\delta = 2$, $\tau_2 = 1$.

Each of the 27 pencils of bitangent planes to $T$ cuts out a $g_2^1$ on $R$, off the fixed intersection of the line with $R$.

A $g_2^1$ on an elliptic curve possesses 4 ramification points. Each of this contributes to the computation, with multiplicity 2, hence the total contribution is $27 \times 4 \times 2 = 216$.

*Case 4:* $\delta = 1$, $|I| = 2$.

For any pair of distinct indices $r$ and $s$ in $\{1, \ldots, 12\}$, the pencil of planes passing through $x_r$ and $x_s$ contains $12 = \deg(\check{T})$ tangent planes to $T$. The contribution is then $\binom{12}{2} \times 12 = 792$.

*Case 5:* $\delta = 1$, $|I| = 1$, $\tau_2 = 1$.

For each $r = 1, \ldots, 12$, consider the projection $\pi_r$ from the point $x_r$. It gives a double cover $\tilde{T} \to \mathbf{P}^2$ branched along a plane quartic curve $B$, where $\tilde{T}$ is the blow up of $T$ at $x_r$ (the curve $B$ is the projection from $x_r$ of the curve $T \cap \mathrm{D}^{x_r}T$, which has degree 6 and multiplicity 2 at $x_r$, see [**??**, Section **??**]; by $\mathrm{D}^{x_r}T$ we denote the first polar of $T$ with respect to the point $x_r$). Being contained in the plane $P_0$, which passes through $x_r$, $R$ is mapped 2-1 to a line $\ell_R \subseteq \mathbf{P}^2$, which is nowhere tangent to $B$.

Tangent planes to $T$ containing $x_r$ in $\mathbf{P}^3$ map via $\pi_r$ to tangent lines to $B$. Tangent planes to $R$ containing $x_r$ map to lines passing through one of the 4 points of $\ell_R \cap B$.

Let $p$ be a point in $\ell_R \cap B$. The pencil of lines through $p$ in $\mathbf{P}^2$ contains $12 = \deg(\check{B})$ tangent lines to $B$. Among them there is the tangent line to $B$ at $q$, with multiplicity 2 (this can be seen as in Lemma (7.1)), which corresponds to a plane tangent to $T$ at a point lying on $R$ in $\mathbf{P}^3$. Such a plane gives a section of $T + P$ with a singularity worse than a tacnode, so it has to be discarded from computation.

Finally, the contribution in this case is $12 \times 4 \times 2 \times 10 = 960$.

*Case 6:* $\delta = 1$, $\tau_3 = 1$.

The plane cubic $R$ possesses 9 flexes. For each of them, there is the pencil of planes containing the tangent line to $R$ at that flex. Such a pencil is a line in $\check{\mathbf{P}}^3$, meeting $\check{T}$ at 12 points. One needs to remove from these 12 points the one corresponding to the tangent plane to $T$ at the flex of $R$. It is of multiplicity 3 (the proof of this fact is similar to the one of Lemma (7.1) and we leave it to the reader).

Planes satisfying the given conditions have to be counted with multiplicity 3. So we get a contribution of $3 \times 9 \times 9 = 243$.

The table below shows the result of the total summation.

| Cases  | Contribution |
|--------|--------------|
| Case 1 | 45           |
| Case 2 | 324          |
| Case 3 | 216          |
| Case 4 | 792          |
| Case 5 | 960          |
| Case 6 | 243          |
| **Total** | **2580**  |

### 7.2.2   Curves in the twisted hyperplane bundle: plane quartics

Since Salmon's formula (2.0.1) gives a total of 3200 tritangent planes to a general quartic, we are missing a total of $3200 - 2580 = 620$ limit trinodal curves. These have to be seen in $\mathcal{H}_P = |4H - \sum_{i=1}^{12} E_i|$. This system however contains the 2–dimensional system $D_R$ of trinodal curves, which is the plane along which $\mathcal{H}_P$ and $\tilde{\mathcal{H}}_{\mathbf{P}^3}$ glue. So these curves do not contribute to the computation. Hence, we have to compute the number of trinodal curves in $|4H - \sum_{i=1}^{12} E_i|$ off $D_R$. All such curves will be irreducible, since, by the genericity assumptions, the only reducible curves in $|4H - \sum_{i=1}^{12} E_i|$ are the ones in $D_R$ (the reader is invited to prove this).

The linear system $|4H - \sum_{i=1}^{12} E_i|$ is mapped to the plane $P_0$ to $|\mathcal{O}_{\mathbf{P}^2}(4) \otimes \mathcal{I}_D|$, where $D = \{x_1, \ldots, x_{12}\}$ and $x_1 + \cdots + x_{12} \in |\mathcal{O}_R(4)|$ is a general divisor. The varieties of $\delta$-nodal curves not containing $R$ in this linear system are logarithmic Severi varieties of the pair $(\mathbf{P}^2, R)$: on logarithmic Severi varieties, see [III], and in particular [III, Section 5.3] about the particular kind we are considering right now.

Since $\dim(|\mathcal{O}_{\mathbf{P}^2}(4)|) = 14$, $D$ imposes only 11 conditions to $|\mathcal{O}_{\mathbf{P}^2}(4)|$. Since passing through a point imposes at most one condition, there is a set $Y$ of 11 points in $D$ that imposes the 11 conditions (in fact, by genericity of $x_1 + \cdots + x_{12} \in |\mathcal{O}_R(4)|$, any subset of 11 points does). Hence
$$|\mathcal{O}_{\mathbf{P}^2}(4) \otimes \mathcal{I}_D| = |\mathcal{O}_{\mathbf{P}^2}(4) \otimes \mathcal{I}_Y|.$$
Finally, $Y$ can be seen as the limit of $|\mathcal{O}_{\mathbf{P}^2}(4) \otimes \mathcal{I}_Z|$, where $Z$ is a general set of 11 points in the plane.

We will compute the number of trinodal curves in $|\mathcal{O}_{\mathbf{P}^2}(4) \otimes \mathcal{I}_Z|$. Their limits, when $Z$ tends to $Y$, will be trinodal curves in $|\mathcal{O}_{\mathbf{P}^2}(4) \otimes \mathcal{I}_Y| = |\mathcal{O}_{\mathbf{P}^2}(4) \otimes \mathcal{I}_D|$, and this gives the required information on the number we want to compute.

Finally what we have to compute is the number of trinodal quartics passing through 11 general points in the plane, which is the degree of the codimension 3 family of trinodal curves in $|\mathcal{O}_{\mathbf{P}^2}(4)| \cong \mathbf{P}^{14}$, i.e., this is the degree of the appropriate Severi variety. This has been done, in general, by Caporaso and Harris in [1, 2]. Their result in this case reads as follows:

**(7.4) Proposition.** *If $Z$ is a general set of 11 points in $\mathbf{P}^2$, then the 3–dimensional linear system $|\mathcal{O}_{\mathbf{P}^2}(4) \otimes \mathcal{I}_Z|$ contains 675 trinodal curves:*
*(i) 620 of them are irreducible;*
*(ii) $55 = \binom{11}{2}$ are reducible in a line joining 2 of the points of $Z$ and the cubic through the remaining 9 points of $Z$.*

When $Z$ tends to $Y$ as above, the 55 reducible quartics as in (ii) of Proposition (7.4) tend to curves in $D_R$, which ought to be discarded. The limit of the 620 curve in (i) is the missing contribution we need.

It would be however unsatisfactory at this point to rely on (the quite difficult) general Caporaso–Harris' result to finish the computation. And in fact we can perform it *our way*, by using a new degeneration.

The degeneration we want to use is the following. Consider the blow up $\mathcal{S}$ of $\mathbf{D} \times \mathbf{P}^2$ along a general line $E$ of $\mathbf{P}^2$ in the *central fibre* over $0 \in \mathbf{D}$. Then there is a morphism $f : \mathcal{S} \to \mathbf{D}$ which is a semistable degeneration of $S_t = \mathbf{P}^2$ for $t \neq 0$ to $S_0 = P \cup F$, with $P \cong \mathbf{P}^2$ the proper transform of the original central fibre and $F$ the exceptional divisor of the blow up (see Figure 4).

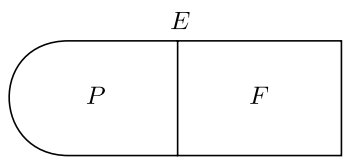

Figure 4: the degeneration of the plane

Note that $F \cong \mathbf{F}_1$, i.e., the plane blown up at a point $p$, and $P \cap F = E$, with $E^2 = -1$ on $F$ (i.e., the exceptional divisor corresponding to $p$), according to the Triple Point Formula in Lemma (3.4).

There is an obvious polarizing bundle $\mathcal{L}$ on $\mathcal{S}$, which restrict to $\mathcal{O}_{\mathbf{P}^2}(1)$ on the general fibre and to $L_0$ on the central fibre, with the aspects

$$\mathcal{L}_P = \mathcal{O}_{\mathbf{P}^2}(1), \quad \mathcal{L}_F = \mathcal{O}_F(f),$$

where $f$ is a fibre of the morphism $\pi : F \cong \mathbf{F}_1 \to \mathbf{P}^1$.

Now pick a set $A = \{\alpha_1, \ldots, \alpha_6\}$ of six general points on $P$ and a set $B = \{\beta_1, \ldots, \beta_5\}$ of five general points on $F$. Then $Z_0 := A \cup B$ is the limit of a set $Z_t$ of 11 general points on the general fibre $S_t$. Denote by $\mathcal{Z}$ the curve $\{Z_t\}_{t \in \mathbf{D}}$ in $\mathcal{S}$. Consider $\mathcal{L}^{\otimes 4} \otimes \mathcal{I}_{\mathcal{Z}|\mathcal{S}}$, which restricts to $\mathcal{O}_{\mathbf{P}^2}(4) \otimes \mathcal{I}_{Z_t}$ on the general fibre $S_t$. In order to have information about trinodal curves in $|\mathcal{O}_{\mathbf{P}^2}(4) \otimes \mathcal{I}_{Z_t}|$, we have to look at their limits in the limit linear system $\mathfrak{L}$ of $\mathcal{L}^{\otimes 4} \otimes \mathcal{I}_{\mathcal{Z}|\mathcal{S}}$.

**(7.5) Lemma.** *There is only one twist of $\mathcal{L}^{\otimes 4} \otimes \mathcal{I}_{\mathcal{Z}|\mathcal{S}}$ which is centrally effective, namely $\mathcal{L}^{\otimes 4} \otimes \mathcal{O}_\mathcal{S}(-F) \otimes \mathcal{I}_{\mathcal{Z}|\mathcal{S}}$.*

*Proof.* Any twist is of the form $\mathcal{L}^{\otimes 4} \otimes \mathcal{O}_\mathcal{S}(-aF) \otimes \mathcal{I}_{\mathcal{Z}|\mathcal{S}}$ and on $S_0$ this resticts to:
(i) $\mathcal{O}_{\mathbf{P}^2}(4-a) \otimes \mathcal{I}_A$ on $P \cong \mathbf{P}^2$;
(ii) $\mathcal{O}_{\mathbf{F}_1}(4f + aE) \otimes \mathcal{I}_B$ on $F \cong \mathbf{F}_1$. As $\mathbf{F}_1$ is the blow up of $\mathbf{P}^2$ at $p$, then (with a slight abuse of notation) one has $|\mathcal{O}_{\mathbf{F}_1}(4f + aE) \otimes \mathcal{I}_B| \cong |\mathcal{O}_{\mathbf{P}^2}(4) \otimes \mathcal{I}_B \otimes \mathcal{I}_p^{\otimes 4-a}|$.

Since $B$ is a general set of 5 points on $F$, for the effectivity on $F$ one needs $a \geqslant 1$, because $\dim(|\mathcal{O}_{\mathbf{F}_1}(4f)|) = 4$. Similarly, for the effectivity on $P$ we need $a \leqslant 1$, because $|\mathcal{O}_{\mathbf{P}^2}(h) \otimes \mathcal{I}_A|$ is empty for $h \leqslant 2$. This proves the assertion. $\square$

In conclusion, we have to look at limits of trinodal curves only in $|\mathcal{L}^{\otimes 4} \otimes \mathcal{O}_\mathcal{S}(-F) \otimes \mathcal{I}_{\mathcal{Z}|\mathcal{S}} \otimes \mathcal{O}_{S_0}|$, which is:
(i) $\mathcal{L}_P := |\mathcal{O}_{\mathbf{P}^2}(3) \otimes \mathcal{I}_A|$, the aspect on $P$, of dimension 3;
(ii) $\mathcal{L}_F := |\mathcal{O}_{\mathbf{F}_1}(4f + E) \otimes \mathcal{I}_B|$, the aspect on $F$, which is the same as $|\mathcal{O}_{\mathbf{P}^2}(4) \otimes \mathcal{I}_p^{\otimes 3} \otimes \mathcal{I}_B|$ in the identification of $F$ with the blow-up of $\mathbf{P}^2$ at a point $p$, also of dimension 3;
(iii) the two curves on $P$ and $F$ in the two aspects have to *match* along $E$, i.e., they have to cut out on $E$ the same divisor. This *matching condition* implies dimension 3 for $|\mathcal{L}^{\otimes 4} \otimes \mathcal{O}_\mathcal{S}(-F) \otimes \mathcal{I}_{\mathcal{Z}|\mathcal{S}} \otimes \mathcal{O}_{S_0}|$: indeed, given a curve in $|\mathcal{O}_{\mathbf{P}^2}(3) \otimes \mathcal{I}_A|$, which depends on 3 parameters, this fixes the degree 3 divisor it cuts on $E$ and there is only one curve in $|\mathcal{O}_{\mathbf{F}_1}(4f - E) \otimes \mathcal{I}_B|$ matching it.

Next we introduce the self-explanatory notation $\boldsymbol{\delta} = (\delta_P, \delta_F)$ and $\tau = (\tau_2, \tau_3)$ (in this case there is no $I$ to be considered) and

$$|\boldsymbol{\delta}| + \tau_2 + 2\tau_3 = 3.$$

It is useful to notice that a curve in $|\mathcal{O}_{\mathbf{F}_1}(4f + E)| \otimes \mathcal{I}_B$ with a node at a point $x \in F$, splits in the unique element $f_x$ of $|f|$ through $x$, and in a curve of $|\mathcal{O}_{\mathbf{F}_1}(3f + E) \otimes \mathcal{I}_B|$ (which is the same as $|\mathcal{O}_{\mathbf{P}^2}(3) \otimes \mathcal{I}_B \otimes \mathcal{I}_p^{\otimes 2}|$) passing through $x$.

Now the cases to be analyzed are the following. We will be sketchy, leaving the details to the reader. The corresponding analysis for the simpler enumeration of 2-nodal plane quartics is carried out in [VII, Section 4].

*Case 1:* $\delta_F = 3$.

A trinodal curve on $F$ contains three curves in $|f|$ with residual in $|\mathcal{O}_{\mathbf{F}_1}(f + E)|$, which is the same as $|\mathcal{O}_{\mathbf{P}^2}(1)|$ (with the curves in $|f|$ mapping to lines through $p$). This implies that the three splitting curves in $|f|$ have to contain three points of $B$, and the residual in $|\mathcal{O}_{\mathbf{F}_1}(f + E)|$ contains the remaining 2. This fixes the curve on $F$ and accordingly also the curve on $P$. This shows that the contribution in this case is $\binom{5}{2} = 10$.

*Case 2:* $\delta_F = 2, \delta_P = 1$.

Two curves in $|f|$ split from the curve on $F$ and either only one or exactly two points of $B$ lie on splitting elements of $|f|$. If two points $\beta_1, \beta_2$ of $B$ lie on splitting elements of $|f|$, these curves cut $E$ each in one point $b_1, b_2$. Then the matching curves on $P$ lie in the pencil of cubics $|\mathcal{O}_{\mathbf{P}^2}(3) \otimes \mathcal{I}_{A \cup \{b_1, b_2\}}|$ which, as well known, contains exactly 12 nodal curves (the reader may be asked to verify this as an exercise). Each such curve cuts a divisor on $E$ (which contains $b_1 + b_2$), and the matching uniquely determines the remaining component of the curve on $F$, which lies in the pencil $|\mathcal{O}_{\mathbf{F}_1}(2f + E) \otimes \mathcal{I}_{B - \{\beta_1, \beta_2\}}| = |\mathcal{O}_{\mathbf{P}^2}(2) \otimes \mathcal{I}_{B - \{\beta_1, \beta_2\} \cup \{p\}}|$. This gives a contribution of $\binom{5}{2} \times 12 = 120$ to our computation.

The reader may check, with similar arguments, that if only one point of $B$ lies on a splitting element of $|f|$ (and there are 5 possibilities for this), then one also has a well determined pencil of curves on $P$ with 12 nodal curves, hence we have $5 \times 12 = 60$ more limit trinodal curves of this type, for a total of 180.

*Case 3:* $\delta_F = 2, |\tau| \neq 0$.

Since, as above, two curves in $|f|$ split from the curve on $F$, it is not possible to have a tangency to $E$, so this case gives no contribution to the computation.

*Case 4:* $\delta_F = 1, \delta_P = 2$.

Only one curve in $|f|$ splits from the curve on $F$, and there are two possible cases:
(a) this curve does not contain any point of $B$;
(b) this curve contains one point of $B$ (this depends on 5 choices).

Also the curve on $P$ is reducible in a conic and a line, and there are two possible cases:
(a') 5 of the points in $A$ lie on the conic, and one lies on the line (this depends of 6 choices);
(b') 4 of the points in $A$ lie on the conic, and two lie on the line (this depends of 15 choices).

One has that:
Case (a,a') contributes with 12 curves,
Case (a,b') contributes with 30 curves,
Case (b,a') contributes with 30 curves,
Case (b,b') contributes with 75 curves, for a total contribution of 147.

Let us check case (b,b'), and leave the others to the reader. Each of the 15 choices of 2 points in $A$ determines a unique line, and then the residual conic through the 4 remaining points moves in a pencil. Each choice of one point in $B$ determines a unique curve in $|f|$, and then the residual curve in $|3f + E|$ through the 4 remaining points moves in a 2-dimensional linear system. Finally, the residual conic on $P$ is determined by its matching with the fixed curve in $|f|$ on $F$, and then the residual curve in $|3f + E|$ on $F$ is fixed with its matching with the conic on $P$ at its remaining point of intersection with $E$ and with the fixed line on $P$. Hence the contribution is of $5 \times 15 = 75$.

*Case 5:* $\delta_F = 1, \delta_P = 1, \tau_2 = 1$.

Again one has the two alternatives (a) and (b) as in Case 4. Case (a) contributes with 40 curves, case (b) with 160. Let us check case (b), and leave case (a) to the reader.

Once the splitting curve of $|f|$ has been fixed (which depends on 5 choices, i.e., the choice of the point in $B$), it intersects $E$ at a point $q$. The matching condition on $P$ gives the net of cubics through $A$ and $q$. This net maps $P$ to $\mathbf{P}^2$ as a double cover, branched along a quartic $D$ which, by genericity, is smooth. The curve $E$ is mapped to a conic $C$. We have to count the number of lines tangent to both, $D$ and $C$, which is the intersection number 24 of $\check{C}$ (a conic) and $\check{D}$ (a curve of degree 12), minus the number of points where $D$ and $C$ are tangent to each other (which is 4, the number of branch points of the double cover $E \to C$), counted with multiplicity 2 (this is similar to Lemma (7.1)). Since each of these curves contributes with multiplicity 2 to the total, we have a total of $2 \times 5 \times 16 = 160$.

*Case 6:* $\delta_P = 1, \tau_3 = 1$.

Each of the curves in question counts with multiplicity 3. The linear system $\mathcal{L}_P$ maps $P$ to a smooth cubic surface $\Phi$ in $\mathbf{P}^3$, and $E$ is mapped to a rational normal cubic $\Gamma \subseteq \Phi$. Consider the duals:
- $\check{\Phi}$, which is a surface of degree 12;
- $\check{\Gamma}$, which is a scroll of degree 4, with a cuspidal curve $\Gamma^*$ of degree 3, which is the rational normal cubic of $\check{\mathbf{P}}^3$ described by all osculating planes to $\Gamma$ in $\mathbf{P}^3$ (see [**??**]).

The surface $\check{\Gamma}$ is the projection of a smooth rational normal scroll $\Sigma$ of degree 4 in $\mathbf{P}^5$, and we denote by $H$ the hyperplane class in this embedding and $|\mathfrak{f}|$ the pencil of lines on $\Sigma$. One has $H^2 = 4, \mathfrak{f}^2 = 0, H \cdot \mathfrak{f} = 1$. [1]

The strict transform of $\Gamma^*$ on $\Sigma$ is clearly unisecant with the curves in $|\mathfrak{f}|$ and therefore $\Gamma^* \equiv H - \mathfrak{f}$.

The curve $\Gamma'$ described by the tangent planes to $\Phi$ at the points of $\Gamma$ sits in $\check{\Phi} \cap \check{\Gamma}$. It is a rational curve of degree 6 (to see this, intersect $\Gamma$ with a general polar quadric of $S'$), meeting the curves of $|\mathfrak{f}|$ in 1 point. Hence on $\Sigma$ one has $\Gamma' \equiv H + 2\mathfrak{f}$, and therefore $\Gamma' \cdot \Gamma^* = 5$.

We have to compute the intersection number of $\check{\Phi}$ and $\Gamma^*$, which is 36, and subtract the number of osculating planes to $\Gamma$ tangent to $\Phi$ at a point of $\Gamma$, which is $\Gamma' \cdot \Gamma^* = 5$, counted with multiplicity 3 (to see this, imitate the proof of Lemma (7.1)). In conclusion we get a contribution of $3 \times 21 = 63$.

*Case 7:* $\delta_P = 2, \tau_2 = 1$.

Each of the curves in question counts with multiplicity 2. The curve on $P$ splits, and we have the two alternatives (a') and (b') as in Case 4. However, case (a') gives no contribution to the computation, because in a pencil of lines there is no line tangent to $E$. Each of the situations in case (b') (which depend on 15 choices, i.e., the pairs of points of the splitting line $\ell$ on $P$) gives rise to 2 curves, corresponding to the conics of the residual pencil to $\ell$ tangent to $E$. In conclusion, one has a contribution of $2 \times 15 \times 2 = 60$.

*Case 8:* $\delta_P = 3$.

The curve on $P$ splits in 3 lines, which have to contain the points in $A$. The reader will readily check that their number is 15.

The following table shows the result of the total summation of non–trivial contributions.

---

[1] One has that $\Sigma \cong \mathbf{F}_0$ embedded in $\mathbf{P}^5$ via the morphism determined by the linear system $|H|$ of curves of type $(2, 1)$. Indeed it cannot be $\Sigma \cong \mathbf{F}_2$ because otherwise there would be a pencil of tangent planes to $\Gamma$ along moving points, which is not possible.

| Cases | Contribution |
|---|---:|
| Case 1 | 10 |
| Case 2 | 180 |
| Case 4 | 147 |
| Case 5 | 200 |
| Case 6 | 63 |
| Case 7 | 60 |
| Case 8 | 15 |
| **Total** | **675** |

From this total we have to subtract the limit of reducible quartic curves passing through 11 general points of the plane. Their number is $55 = \binom{11}{2}$, since they are the ones appearing in (ii) of Proposition (7.4). We leave it to the reader to check that the limits of these 55 curves are the ones in Case 1 and Case 4 (b,a') and (a,b').

**(7.6) Remark.** In conclusion, also the degeneration to a cubic plus a plane is 1–good. However the reader will have noticed that, both in this case and in the previous one, we do not have an a priori proof of this (as in the intentions of Conjecture (5.10)), but the proof relies on Salmon's formula (2.0.1), rather than proving it.

## 8 – A non-good case

### 8.1 – Degeneration of cubics to three planes

Consider a degeneration as in Example (3.2), where $d = 3$, with $S_0$ consisting of the union of three distinct planes $P_0, P_1, P_2$. Let us denote by $\ell_i$ the intersection line of $P_i$ and $P_{i+1}$ (for $i = 0, 1, 2$; the indices will now be considered to vary in $\mathbf{Z}/3\mathbf{Z}$), and by $p$ the intersection point of $P_0, P_1, P_2$. The central fibre of the corresponding semistable degeneration constructed in Example (3.2) is $\tilde{S}_0 = Q_0 \cup Q_1 \cup Q_2$ as in Figure 5, and we may assume that each of the surfaces $Q_i$ is the blow up of $P_i$ at three points along $\ell_i$.

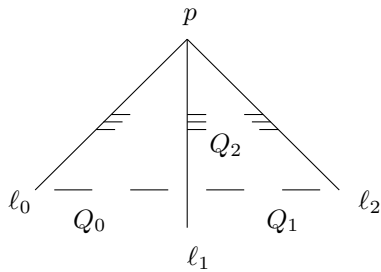

Figure 5: the degeneration of a cubic to the union of three planes

If we look at limits of nodal curves in the polarising hyperplane bundle, according to Proposition (5.5) we have to consider only planes passing by one of the aforementioned blown up points. This gives a contribution of 9 to the degree of the dual of a general cubic surface in $\mathbf{P}^3$, whereas we know that this degree is 12. So we are missing something of degree 3, which shows that this degeneration is not 1–good. In what follows we will construct a 1–good model, which will show that the degeneration under consideration is 1–well behaved.

**(8.1) Remark.** This degeneration is however 3–good. To see this, notice that the limit of the 27 lines on a general cubic can be seen here as the lines different from the $\ell_i$'s joining pairs of blown up points. Hence we see all triangles formed by these lines, which are 45, so we see also the limit of all 3–nodal planes sections of a general cubic surface.

For the same reason this degeneration is also 2–good.

The obvious guess is that the limit of the dual of the general cubic consists of the 9 aforementioned planes in $\check{\mathbf{P}}^3$ plus the plane $p^\perp$ counted with multiplicity 3, i.e., *there is something of degree 3 hidden at $p$*. This is indeed the case and to see it, one has to look at another model of the degeneration. This is explained in detail in [4] in the more complicated situation of a general quartic degenerating to a tetrahedron (i.e., the union of four linearly independent planes in $\mathbf{P}^3$). Let us now see how this 1–good model is gotten in the present case.

**(8.2) Construction of a good model.** Let $\bar{f}: \bar{S} \to \mathbf{D}$ be the family obtained from $\tilde{f}: \tilde{S} \to \mathbf{D}$ (notation as in Example (3.2)) by the base change $t \in \mathbf{D} \mapsto t^3 \in \mathbf{D}$. The central fibre $\bar{S}_0$ is isomorphic to $\tilde{S}_0$, so we will keep the above notation for it.

Analytically locally around $p$, the total space $\bar{S}$ is isomorphic to the hypersurface of $\mathbf{C}^4$ defined by the equation $xyz = t^3$ at the origin. Blow-up $\bar{S}$ at $p$. The blown–up total space locally sits in $\mathbf{C}^4 \times \mathbf{P}^3$. Let $[\xi : \eta : \zeta : \vartheta]$ be the homogeneous coordinates in $\mathbf{P}^3$. Then the new total space is locally defined in $\mathbf{C}^4 \times \mathbf{P}^3$ by the equations

$$(8.2.1) \qquad \mathrm{rk} \begin{pmatrix} x & y & z & t \\ \xi & \eta & \zeta & \vartheta \end{pmatrix} \leqslant 1 \quad \text{and} \quad \xi\eta\zeta = \vartheta^3.$$

Denote the exceptional divisor by $T$; it is isomorphic to the cubic surface with equation $\xi\eta\zeta = \vartheta^3$ in the $\mathbf{P}^3$ with coordinates $[\xi : \eta : \zeta : \vartheta]$. This cubic contains three lines, along which it intersects the proper transforms $\tilde{Q}_0, \tilde{Q}_1, \tilde{Q}_2$ of $Q_0, Q_1, Q_2$ respectively; thus, each line on $T$ identifies with the exceptional $(-1)$-curve on one of $\tilde{Q}_0, \tilde{Q}_1, \tilde{Q}_2$. The cubic $T$ has three $A_2$ double points, located at the intersections of $T$ with the proper transforms $\tilde{\ell}_0, \tilde{\ell}_1, \tilde{\ell}_2$ of $\ell_0, \ell_1, \ell_2$ respectively. See Remark (8.3) below for more about this cubic. The new central fibre is shown in Figure 6.

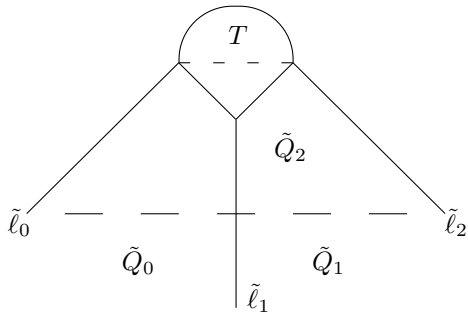

Figure 6: the central fibre after base change and blow–up of the vertex $p$

The model we have arrived at at this point is good enough to see all limits of 1-nodal curves, but it is not yet a good model because it is not semistable. Indeed, the total space is singular along the proper transforms $\tilde{\ell}_0, \tilde{\ell}_1, \tilde{\ell}_2$ of $\ell_0, \ell_1, \ell_2$: analytically locally around the general point of one of these curves, it is isomorphic to a neighbourhood of the origin in the hypersurface defined by $xy = t^3$ in $\mathbf{C}^4$; thus, the total space is locally the product of an $A_2$ surface singularity with a smooth curve.

To resolve these singularities, we blow up along the three (disjoint) curves $\tilde{\ell}_0, \tilde{\ell}_1, \tilde{\ell}_2$. This has the effect of replacing each of these curves by a chain of two ruled surfaces; correspondingly, all three $A_2$ double points of $T$ are resolved, each being replaced by a chain of two $(-2)$-curves. Thus, $T$ is replaced by its minimal resolution, which we will denote by $\tilde{T}$. The three surfaces $\tilde{Q}_0, \tilde{Q}_1, \tilde{Q}_2$, on the other hand, remain unchanged, and we will denote their proper transforms by the same symbols.

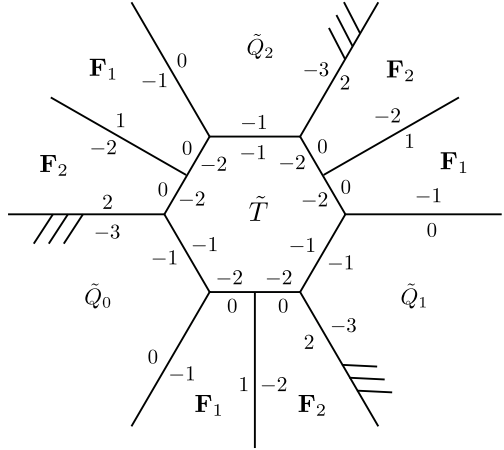

Figure 7: the central fibre of the good model (seen from above)

The new central fibre is shown in Figure 7. We leave to the reader the verification of the self-intersections of the double curves of the central fibre indicated on the figure, as an exercise in the triple point formula. With our choice of distribution of the $(-1)$-curves in the small resolution at the beginning, the chains of rational ruled surfaces are made of an $\mathbf{F}_1$ and an $\mathbf{F}_2$ surface, as indicated on the figure.

**(8.3) Remark.** The cubic surface $T$ is the image of the plane via the linear system $\mathcal{C}$ of cubic curves with two base points $x_1, x_2$, and infinitely near base points such that the cubics $C \in \mathcal{C}$ have two common flexes along lines $L_1$ and $L_2$. One has to blow up the plane at $x_1$ and $x_2$ three consecutive times to make the proper transform $\tilde{\mathcal{C}}$ of $\mathcal{C}$ base point free. This is shown in Figure 8, where we denote by $\tilde{C}$ the proper transform of the general curve $C \in \mathcal{C}$, by $\tilde{L}_i$ the proper transforms of the lines $L_i$, for $i = 1, 2$, and by $E_{ij}$, for $i = 1, 2, j = 1, 2, 3$, the exceptional curves of the blow–ups, and one has $E_{i3}^2 = -1$ for $i = 1, 2$ and $E_{ij}^2 = -2$, for $i, j = 1, 2$. [2]

The cubic $T$ has three $A_2$–double points, at the vertices of the triangle $\vartheta = 0$, $\xi\eta\zeta = 0$. They are the images of $\tilde{L}_1 + \tilde{L}_2$ and of $E_{i1} + E_{i2}$, for $i = 1, 2$. The three lines pairwise joining the $A_2$–double points, lie on $T$ and are easily seen to be the only lines on $T$. They correspond to the two exceptional divisors $E_{13}$ and $E_{23}$ and to the line $L_3$ joining $x_1$ and $x_2$.

The cubic $T$ is *self–dual*, i.e., $\check{T}$ is projectively equivalent to $T$, hence it has degree 3. This (heuristically) explains the contribution by 3 to the degree of the dual hidden at $p$.

**(8.4) Goodness of the new model.** The degeneration we have just constructed is 1–good. Indeed now there is a new centrally effective twist, which corresponds to a degeneration of the cubics $S_t$, $t \in \mathbf{D} \setminus \{0\}$, to the cubic $T \subseteq \mathbf{P}^3$. It provides the missing contribution by 3, as we will now briefly explain.

Let $\mathcal{L}$ be the pull-back to the new degeneration of the polarising hyperplane bundle on the initial degeneration. On the central fibre, it restricts to the pull-back $H$ of the line class on the surfaces $\tilde{Q}_0, \tilde{Q}_1, \tilde{Q}_2$, to the trivial class on $\tilde{T}$, and to the class of the ruling on the chains of ruled surfaces. The new relevant twist is $\mathcal{L}(-\tilde{T})$: it restricts to the class $H - E_i$ on $\tilde{Q}_i$ for all $i = 0, 1, 2$ (recall the notation from (8.2)), to the pull-back of the hyperplane class of $T \subseteq \mathbf{P}^3$ on $\tilde{T}$, and it is trivial on the chains of ruled surface; thus, it maps $\tilde{T}$ to $T$, contracts $\tilde{Q}_i$ to the line $E_i$ on $T$ for all $i = 0, 1, 2$, and contracts the three chains of ruled surfaces to the three $A_2$ double points of $T$.

---

[2] In [4] we gave a different plane linear system of cubics representing the same surface $T$. It is the linear system of cubics that pass through three independent points $y_1, y_2, y_3$ in the plane and are tangent there to the lines $\langle y_i, y_{i+1} \rangle$, with $i = \in \mathbf{Z}/3\mathbf{Z}$. To pass from $\mathcal{C}$ to this linear system, just make a quadratic transformation based at $x_1$ and $x_2$ and at the point infinitely near to $x_1$ along $L_1$.

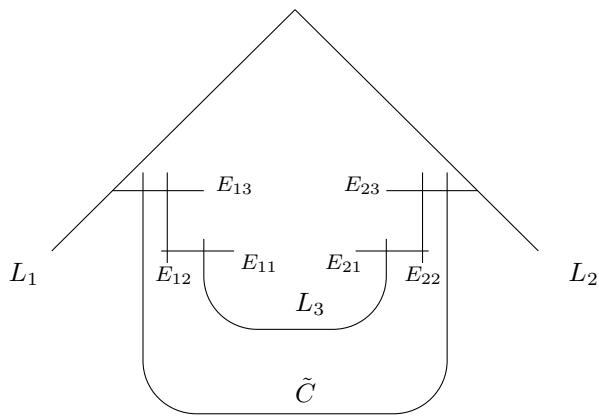

Figure 8: the linear system of curves $\tilde{C}$

Now, let us fix a general pencil of planes in the $\mathbf{P}^3$ spanned by $S_0$, and let us count the "tangent planes to $S_0$" in this pencil. Nine of them are the planes in the pencil passing through each of the special points on the lines $\ell_i$, $0 \leqslant i \leqslant 2$. As for the remaining three, consider the plane in the pencil passing through the triple point; it cuts out three coplanar lines meeting at the point $p$, the triple point of $S_0$. Thus its proper transform determines three points on $T$, one on each line of $T$, and these three points are collinear. The planes in the $\mathbf{P}^3$ spanned by $T$ which contain these three points form a pencil, and in this pencil there are three planes which are tangent to $T$ in simple points. This provides us with the remaining three "tangent planes to $S_0$". For further details, see [4].

## 8.2 – Degeneration of quartics to four planes

We conclude this text by a brief discussion of degenerations, as in Example (3.2) and Section 8.1 above, of smooth quartics $S_t$ to the union $S_0$ of four planes $P_0, P_1, P_2, P_3$ in linear general position, which we call a *tetrahedron*. This is one of the main topics in [4], and many examples in this volume are inspired by this situation. An interesting feature is that this is a degeneration of $K3$ surfaces.

There are four different kinds of alternative degenerations of the $S_t$'s to be considered in order to reveal all limits of $\delta$-nodal hyperplane sections, $1 \leqslant \delta \leqslant 3$. Each corresponds to a kind of twist in some suitable good model of the initial degeneration.

— The contributions of the *faces* of the tetrahedron (i.e., the planes $P_i$ themselves) are visible in degenerations to monoid quartic surfaces, similarly to what we have seen for the degeneration to a cubic and a plane in Section 7.2.
— The contributions of the sections by a hyperplane containing a *vertex* (i.e., a triple point $P_i \cap P_j \cap P_k$) are visible in degenerations to the union of a cubic $T \subseteq \mathbf{P}^3$ as in the previous Section 8.1 and the plane spanned by the three lines of $T$; the plane corresponds to the face of the tetrahedron opposite to the vertex in consideration.
— The contributions of the sections by a hyperplane containing an *edge* (i.e., a double line $P_i \cap P_j$) are visible in degenerations to a rational quartic surface with a double line and two triple points on the latter; the double line is the image of the opposite edge, and the two triple points are the respective images of the two faces adjacent to the edge.

More can be found on the rational surfaces appearing in these degenerations in [**??**, Section **??**].

The construction of a good model follows the same lines as in Section 8.1; it is rather complicated in practice, but can be organized so as to be manageable, see [4]. Starting from the degeneration to $S_0 = P_0 + P_1 + P_2 + P_3$ as in Example (3.2), one performs the base change $t \mapsto t^6$ and then resolve the singularities. Morally, a degree 3 base change is needed to unravel the contributions involving the vertices, as in Section 8.1, and a degree 2 base change is needed for those involving the edges. The central fibre of this good model is represented in Figure 9; the reader may also see [4, p. 143] for another picture, with a different design. The dotted curve represents (in a tropical fashion, see [**??**, Section **??**]) the pull-back of a hyperplane section of the initial tetrahedron.

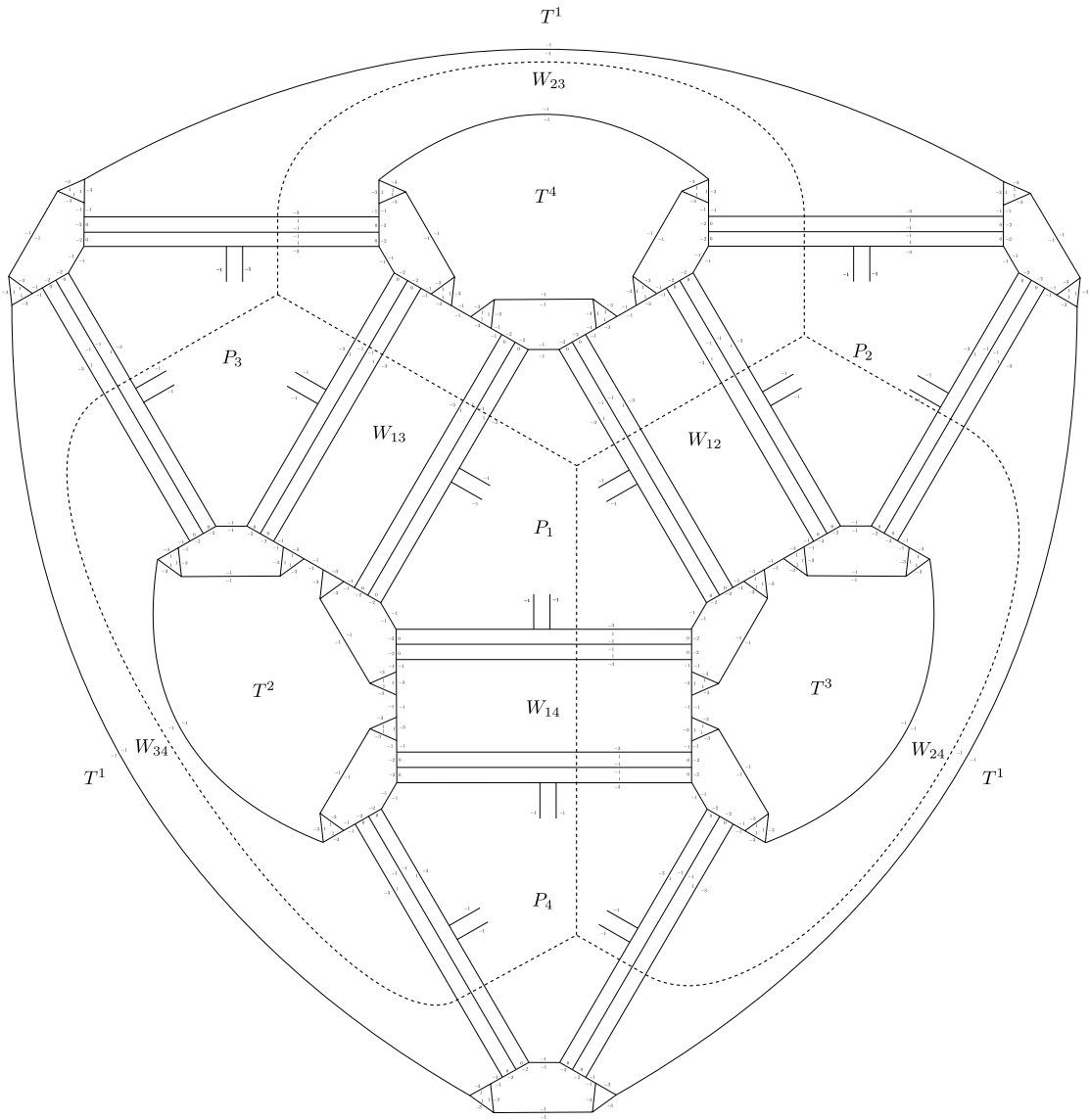

Figure 9: Central fibre of a good model of the degeneration of quartic surfaces to a tetrahedron

# Lecture II
# Deformations of curves on surfaces

## by Thomas Dedieu and Edoardo Sernesi

This text surveys some foundational material from the deformation theory of singular curves on surfaces, of ubiquitous use and necessity in this whole collection. The included topics are best described by the above table of contents.

## 1 – Algebroid plane curves

A *curve singularity* is a scheme of the form $\mathrm{Spec}(B)$, where $(B, \mathfrak{m})$ is the local ring at a singular point of a reduced algebraic curve. If the curve is contained in a nonsingular algebraic surface then we speak of a *planar curve singularity*. We will be mostly interested in this second case and in some of our considerations it will be convenient to replace such a $B$ by its $\mathfrak{m}$-adic completion. In this section we will focus on this case and therefore we will consider local rings of the form

$$B = \mathbb{C}[[X,Y]]/(f)$$

for some power series $f$ without constant term and without multiple factors. The scheme $C := \mathrm{Spec}(B)$ for such a ring is called an *algebroid plane curve*. If $f$ is irreducible then $C$ is called an *analytic branch*.

Assume that $f$ is irreducible. Then the normalization $\widetilde{B}$ of $B$ is a complete DVR and therefore it is isomorphic to a formal power series ring in one variable. Thus there exists $t \in K(B)$, the field of fractions of $B$, such that $\widetilde{B} = \mathbb{C}[[t]]$ and $K(B) = \mathbb{C}((t))$, the fraction field of $\mathbb{C}[[t]]$, which coincides with the field of formal Laurent power series in $t$.

We can write:
$$f(X,Y) = f_n(X,Y) + f_{n+1}(X,Y) + \cdots$$

for some $n \geq 1$, where $f_i$ is a homogeneous polynomial of degree $i$. The integer $n$ is the *multiplicity* of $f$ (or of $C$) at $0 = [\mathfrak{m}]$. The polynomial $f_n(X,Y)$ is called the *initial form* of $f$.

We have two ideals associated to the algebroid plane curve $C$:

$$\left(\frac{\partial f}{\partial X}, \frac{\partial f}{\partial Y}\right) = (f_X, f_Y)$$

is called the *Milnor ideal* or *gradient ideal* of $f$. The corresponding algebra:

$$M_f := \mathbb{C}[[X,Y]]/(f_X, f_Y)$$

is the *Milnor algebra* of $f$.

The ideal

$$(f, f_X, f_Y)$$

is called the *Tjurina ideal* of $f$ and

$$T_f := \mathbb{C}[[X,Y]]/(f, f_X, f_Y)$$

is the *Tjurina algebra* of $f$. It will be convenient to describe the Tjurina algebra as follows:

$$T_f = B/J$$

where $J = (f, f_X, f_Y)/(f) \subseteq B$ is the so-called *jacobian ideal*. Viewed as a $B$-module $T_f$ is also called the *first cotangent module* of $B$, and denoted by $T_B^1$. Since $C$ is either nonsingular or has an isolated singularity at $[\mathfrak{m}]$, the Tjurina ideal is either $(1)$ or supported at $[\mathfrak{m}]$, and therefore $T_f$ is a finite dimensional $\mathbb{C}$-vector space. Then

$$\tau(f) := \dim(T_f) < \infty$$

is called the *Tjurina number* of $f$ (or of $B$). On the other hand it is also true that $M_f$ is finite dimensional because of the following:

**(1.1) Lemma** (Risler [12]). *If $f \in (X,Y)^2 \subseteq \mathbb{C}[[X,Y]]$ has no multiple factors then $f_X, f_Y$ form a regular sequence.*

*Proof.* It suffices to prove that $f_X$ and $f_Y$ have no common factors because they are contained in $(X,Y)$. Assume that they have a non-constant irreducible factor $p$ and consider the algebra $A = \mathbb{C}[[X,Y]]/(p)$. Let $\overline{f} \in A$ be the image of $f$. Then $\overline{f}$ is annihilated by every $\mathbb{C}$-derivation $D : K(A) \to K(A)$. Let $t \in K(A)$ be such that $K(A) = \mathbb{C}((t))$ and consider the derivation $D = \frac{d}{dt}$. Then $\ker(D) = \mathbb{C}$. Therefore $\overline{f} \in \mathbb{C}$ and therefore $\overline{f} = 0$, which means that $f = pq$ for some $q \in (X,Y)$. Therefore, since $f_X = p_X q + p q_X$ is divisible by $p$ but $p_X$ is not, it follows that $q$ is divisible by $p$, and therefore $f$ is divisible by $p^2$, a contradiction. $\square$

For a proof of Lemma (1.1) in the case of convergent power series see [8], Lemma 2.3 p. 113. The dimension

$$\mu(f) = \dim(M_f)$$

is called the *Milnor number* of $f$ (or of $B$). By construction we have:

$$\tau(f) \leq \mu(f)$$

and equality holds if and only if $f \in (f_X, f_Y)$. Notice also that

$$\mu(f) = i(f_X, f_Y)$$

where $i(-,-)$ denotes intersection multiplicity.

**(1.2) Definition.** *A polynomial $P(X,Y) = \sum_\alpha a_\alpha X^{\alpha_1} Y^{\alpha_2}$ is called* quasi-homogeneous *(qh) of weights $w = (w_1, w_2) \in \mathbb{N}^2$ and degree $d$ if*

$$\alpha_1 w_1 + \alpha_2 w_2 = d$$

*for all $\alpha$ such that $a_\alpha \neq 0$. An algebroid plane curve $\mathrm{Spec}\,(\mathbb{C}[[X,Y]]/(f)))$ is quasi homogeneous if there is a q.h. polynomial $P$ such that $\mathbb{C}[[X,Y]]/(f) \cong \mathbb{C}[[X,Y]]/(P)$.*

If $P(X,Y)$ is qh then we have the generalized *Euler relation* in $\mathbb{C}[X,Y]$:

$$dP = w_1 X P_X + w_2 Y P_Y$$

and therefore $P \in (P_X, P_Y)$. Therefore an algebroid qh plane curve satisfies $\mu(f) = \tau(f)$. The converse is true, by the following:

**(1.3) Theorem** (Saito [13])**.** *The algebroid plane curve $\mathrm{Spec}\,(\mathbb{C}[[X,Y]]/(f))$ is qh if and only if $\mu(f) = \tau(f)$.*

**(1.4) Example** ([8, p. 111])**.** The polynomial $P(X,Y) = X^5 + Y^5 + X^2 Y^2$ is not qh. One can compute that a $\mathbb{C}$-basis of $T_f$ is defined by the monomials $1, X, \ldots, X^4, XY, Y, \ldots, Y^4$ and a $\mathbb{C}$-basis of $M_f$ has the additional monomial $Y^5$. Thus $\tau(P) = 10 < \mu(P) = 11$.

**(1.5) Example.** The polynomial $f = X^a \pm Y^b$ is qh of weights $(b, a)$ and degree $ab$. Therefore:

$$\tau(f) = \mu(f) = (a-1)(b-1)$$

**(1.6) Example.** The *simple singularities* (also called *ADE singularities*) are qh. They have the following equations:
- $A_k$: $x^{k+1} + y^2 = 0$, $\quad k \geq 1$.
- $D_k$: $x(y^2 + x^{k-2}) = 0$, $\quad k \geq 4$.
- $E_6$: $x^3 + y^4 = 0$.
- $E_7$: $x(x^2 + y^3) = 0$.
- $E_8$: $x^3 + y^5 = 0$.

We have:

$$\mu(A_k) = \mu(D_k) = \mu(E_k) = k.$$

## 2 – The conductor

Let $B$ be a ring and let $\overline{B} \subseteq \mathrm{Frac}(B)$ be a subring of the total ring of fractions containing $B$. The homomorphism

$$\mathrm{Hom}_B(\overline{B}, B) \to B, \quad \varphi \mapsto \varphi(1)$$

induces an isomorphism

$$\mathrm{Hom}_B(\overline{B}, B) \cong \mathrm{Ann}_B(\overline{B}/B).$$

Then $\mathrm{Ann}_B(\overline{B}/B)$ is an ideal both in $B$ and in $\overline{B}$. It is called the *conductor* of $\overline{B}$ in $B$. We define the *$\delta$-invariant* of $B \subseteq \overline{B}$ as

$$\delta(\overline{B}/B) := \dim_{\mathbb{C}}(\overline{B}/B).$$

If $\overline{B}$ is the integral closure of $B$ in $\mathrm{Frac}(B)$ then $\mathrm{Ann}_B(\overline{B}/B)$ is simply called the *conductor of $B$* and denoted by $\mathcal{A}(B)$; $\delta(\overline{B}/B)$ is denoted by $\delta(B)$ and called *$\delta$-invariant of $B$*.

In case $C = \mathrm{Spec}(B)$ is either a curve singularity or an algebroid plane curve, $\delta(B)$ is denoted by $\delta(C)$ or $\delta(f)$ in case $B = \mathbb{C}[[X,Y]]/(f)$. It is called the *$\delta$-invariant* of $B$ (of $C$).

Notice that if $C = \mathrm{Spec}(B)$ is any curve singularity then $\delta(B) = \delta(\widehat{B})$, by the flatness of $\widehat{B}$ over $B$.

**(2.1) Definition.** *Let $\varphi : \widetilde{C} \longrightarrow C$ be a birational morphism of reduced curves. The* conductor *of $\varphi$ is the sheaf of $\mathcal{O}_C$-ideals:*

$$\mathcal{A}(\varphi) := \mathcal{A}nn(\varphi_*\mathcal{O}_{\widetilde{C}}/\mathcal{O}_C)$$

*or, equivalently:*

$$\mathcal{A}(\varphi) = \mathcal{H}om(\varphi_*\mathcal{O}_{\widetilde{C}}, \mathcal{O}_C).$$

*For each $x \in C$, the finite nonnegative integer*

$$\delta(\varphi, x) := \delta\big((\varphi_*\mathcal{O}_{\widetilde{C}})_x/\mathcal{O}_{C,x}\big) = \dim_\mathbb{C}\big((\varphi_*\mathcal{O}_{\widetilde{C}})_x/\mathcal{O}_{C,x}\big)$$

*is called the (local) $\delta$-invariant of $\varphi$ at $x$, and*

$$\delta(\varphi) := h^0(\varphi_*\mathcal{O}_{\widetilde{C}}/\mathcal{O}_C) = \sum_{x \in C} \delta(\varphi, x)$$

*is the (global) $\delta$-invariant of $\varphi$. In case $\varphi$ is the normalization of $C$ (i.e., $\widetilde{C}$ is nonsingular) we write $\mathcal{A}(\varphi) = \mathcal{A}(C)$ and $\delta(\varphi) = \delta(C)$ and call them the* conductor *of $C$ and the $\delta$-invariant of $C$ respectively.*

The exact sequence

$$0 \longrightarrow \mathcal{O}_C \longrightarrow \nu_*\mathcal{O}_{\widetilde{C}} \longrightarrow \nu_*\mathcal{O}_{\widetilde{C}}/\mathcal{O}_C \longrightarrow 0$$

gives

$$\chi(\mathcal{O}_{\widetilde{C}}) = \chi(\nu_*\mathcal{O}_{\widetilde{C}}) = \chi(\mathcal{O}_C) + \delta(\varphi),$$

which can be rephrased as:

(2.1.1) $$\delta(\varphi) = p_a(C) - p_a(\widetilde{C})$$

where as usual $p_a(-)$ denotes *arithmetic genus* of $-$. In case $\varphi$ is the normalization of $C$, the above identity takes the form:

(2.1.2) $$\delta(C) = p_a(C) - p_g(C)$$

which shows that the $\delta$-invariant of $C$ measures the difference between the arithmetic genus and the geometric genus of the curve.

For the next lemma, observe that $\mathcal{A}(\varphi)$ is naturally a sheaf of $\varphi_*\mathcal{O}_{\widetilde{C}}$-modules. Therefore, since $\varphi$ is affine, we have $\widetilde{C} = \mathrm{Spec}(\varphi_*\mathcal{O}_{\widetilde{C}})$ and $\mathcal{A}(\varphi)$ corresponds to a sheaf of $\mathcal{O}_{\widetilde{C}}$-modules which we denote by $\widetilde{\mathcal{A}}(\varphi)$.

**(2.2) Lemma.** *Let $\varphi : \widetilde{C} \to C$ be a birational morphism of reduced projective curves. Then there is a canonical isomorphism:*

(2.2.1) $$\varphi_*\omega_{\widetilde{C}} \cong \mathcal{H}om_{\mathcal{O}_C}(\varphi_*\mathcal{O}_{\widetilde{C}}, \omega_C)$$

*If moreover $C$ is Gorenstein, then*

(2.2.2) $$\varphi_*\omega_{\widetilde{C}} \cong \mathcal{A}(\varphi) \otimes_{\mathcal{O}_C} \omega_C$$

*and*

(2.2.3) $$\omega_{\widetilde{C}} = \widetilde{\mathcal{A}}(\varphi) \otimes_{\mathcal{O}_{\widetilde{C}}} \varphi^*\omega_C.$$

*Proof.* The proof of (2.2.1) follows from standard facts and is valid for finite morphisms of projective schemes. Recall the two following facts from [9, III, Ex. 6.10, p. 239]. For each quasi-coherent sheaf $\mathcal{G}$ on $C$, the module $\mathcal{H}om_{\mathcal{O}_C}(\varphi_*\mathcal{O}_{\widetilde{C}}, \mathcal{G})$ has a natural structure of quasi-coherent sheaf on $\widetilde{C}$, and as such it is denoted by $\varphi^!\mathcal{G}$. For each quasi-coherent sheaf $\mathcal{F}$ on $\widetilde{C}$, there is a natural isomorphism:
$$\varphi_*\mathcal{H}om_{\mathcal{O}_{\widetilde{C}}}(\mathcal{F}, \varphi^!\mathcal{G}) \cong \mathcal{H}om_{\mathcal{O}_C}(\varphi_*\mathcal{F}, \mathcal{G}).$$
Taking $\mathcal{G} = \omega_C$ and $\mathcal{F} = \mathcal{O}_{\widetilde{C}}$ gives
$$\varphi_*(\varphi^!\omega_C) \cong \mathcal{H}om_{\mathcal{O}_C}(\varphi_*\mathcal{O}_{\widetilde{C}}, \omega_C).$$
Therefore it suffices to prove that $\varphi^!\omega_C = \omega_{\widetilde{C}}$ (which is [9, III, Ex. 7.2(a)]. For all coherent $\mathcal{F}$ on $\widetilde{C}$ we have:
$$H^1(\widetilde{C}, \mathcal{F})^\vee = H^1(C, \varphi_*\mathcal{F})^\vee = \mathrm{Hom}_{\mathcal{O}_C}(\varphi_*\mathcal{F}, \omega_C) \cong \mathrm{Hom}_{\mathcal{O}_{\widetilde{C}}}(\mathcal{F}, \varphi^!\omega_C),$$
and this isomorphism is functorial in $\mathcal{F}$. This implies that $\varphi^!\omega_C$ is a dualizing sheaf for $\widetilde{C}$ in the sense of [9, Definition p. 241] (see [18, Exercise 29.1.A]) , and thus (2.2.1) is proved.

If $\omega_C$ is invertible, then we have:
$$\varphi_*\omega_{\widetilde{C}} = \mathcal{H}om_{\mathcal{O}_C}(\varphi_*\mathcal{O}_{\widetilde{C}}, \omega_C) \cong \mathcal{H}om_{\mathcal{O}_C}(\varphi_*\mathcal{O}_{\widetilde{C}}, \mathcal{O}_C) \otimes \omega_C = \mathcal{A}(\varphi) \otimes \omega_C,$$
which is (2.2.2). Then, by the Projection Formula we have
$$\varphi_*\bigl(\widetilde{\mathcal{A}}(\varphi) \otimes \varphi^*\omega_C\bigr) = \varphi_*\widetilde{\mathcal{A}}(\varphi) \otimes \omega_C = \mathcal{A}(\varphi) \otimes \omega_C,$$
where the equality $\varphi_*\widetilde{\mathcal{A}}(\varphi) = \mathcal{A}(\varphi)$ follows from [9, II, Ex. 5.17e]. Thus
$$\varphi_*\omega_{\widetilde{C}} = \varphi_*\bigl(\widetilde{\mathcal{A}}(\varphi) \otimes \varphi^*\omega_C\bigr),$$
hence $\omega_{\widetilde{C}} = \widetilde{\mathcal{A}}(\varphi) \otimes \varphi^*\omega_C$, again by [9, II, Ex. 5.17e]. $\square$

We refer the reader to [18, §29.3–4] for versions of the formula $\varphi^!\omega_C = \omega_{\widetilde{C}}$ valid in a more general framework.

**(2.3) Corollary.** *(a) Let $C = \mathrm{Spec}(B)$ be a Gorenstein reduced curve singularity. Then:*
$$\delta(B) = \dim_{\mathbb{C}}(B/\mathcal{A}(B)).$$

*(b) Let $C$ be a Gorenstein reduced projective curve and $\varphi : \widetilde{C} \to C$ a birational morphism. Then*
$$\delta(\varphi) = h^0(\mathcal{O}_C/\mathcal{A}(\varphi)).$$

*Proof.* It suffices to prove (b). Let $\mathcal{A} = \mathcal{A}(\varphi)$. By definition we have an exact sequence:
$$0 \longrightarrow \mathcal{H}om_{\mathcal{O}_C}(\varphi_*\mathcal{O}_{\widetilde{C}}, \mathcal{O}_C) \longrightarrow \mathcal{O}_C \longrightarrow \mathcal{O}_C/\mathcal{A} \longrightarrow 0 .$$
Tensoring by $\omega_C$, using the fact that $C$ is Gorenstein and recalling (2.2.1), we obtain:
$$0 \longrightarrow \varphi_*\omega_{\widetilde{C}} \longrightarrow \omega_C \longrightarrow [\mathcal{O}_C/\mathcal{A}] \otimes \omega_C \longrightarrow 0 .$$

From this exact sequence we deduce:

$$\begin{aligned} h^0(\mathcal{O}_C/\mathcal{A}) &= \chi([\mathcal{O}_C/\mathcal{A}] \otimes \omega_C) \\ &= \chi(\omega_C) - \chi(\varphi_*\omega_{\widetilde{C}}) \\ &= \chi(\omega_C) - \chi(\omega_{\widetilde{C}}) \\ &= p_a(C) - p_a(\widetilde{C}) \\ &= \delta(\varphi) \end{aligned}$$

and comparing with (2.1.1) we conclude. □

**(2.4) Corollary.** *Let $\varphi : \widetilde{C} \to C$ be a birational morphism of projective reduced curves and assume that $C \subseteq Y$, a projective nonsingular surface. Then*

$$\mathcal{A}(\varphi) \otimes \mathcal{O}_C(C) = \varphi_*\omega_{\widetilde{C}} \otimes \omega_Y^{-1}$$

*and*

$$\widetilde{\mathcal{A}}(\varphi) \otimes \varphi^*\mathcal{O}_C(C) = \omega_{\widetilde{C}} \otimes \varphi^*\omega_Y^{-1}.$$

*Proof.* The assumptions imply that $C$ is Gorenstein. Then the corollary is an easy consequence of Lemma (2.2). □

**(2.5) Lemma.** *Let $C$ be a reduced curve contained in a nonsingular algebraic surface $Y$. Let $p \in C$ be a point of multiplicity $m$ and $b : Z \to Y$ the blow-up of $Y$ at $p$. Then, letting $\overline{C} \subseteq Z$ be the proper transform of $C$ and $E \subseteq Z$ the exceptional curve we have:*

$$b_*\omega_{\overline{C}} = \omega_C \otimes b_*\mathcal{O}_Z\big(-(m-1)E\big).$$

*Proof.* Left to the reader. □

The Milnor number and the $\delta$-invariant are related as follows:

**(2.6) Proposition** (Milnor formula)**.** *Let $f = f_1 \cdots f_r \in \mathbb{C}[[X,Y]]$, where $f_1, \ldots, f_r$ are pairwise distinct analytic branches. Then:*

$$\delta(f) = \frac{1}{2}\big[\mu(f) + r - 1\big].$$

*Proof.* See [8, Prop. 3.35, p. 208]. □

**(2.7) Remark.** Let $\varphi : \widetilde{C} \longrightarrow C$ be a birational morphism of reduced curves, and denote by $\Delta \subseteq C$ and $\widetilde{\Delta} \subseteq \widetilde{C}$ the two subschemes defined by $\mathcal{A}(\varphi)$ and $\widetilde{\mathcal{A}}(\varphi)$ respectively. If both $C$ and $\widetilde{C}$ are Gorenstein, then it follows from (2.2.3) that $\widetilde{\Delta}$ is Cartier. However it is not true in general that $\Delta$ is a Cartier divisor (for instance, at a node $xy = 0$ the conductor is the ideal generated by $x$ and $y$).

**(2.8) Remark.** In general the two sheaves of $\mathcal{O}_{\widetilde{C}}$-modules $\widetilde{\mathcal{A}}(\varphi)$ and $\varphi^*\mathcal{A}(\varphi)$ are different, even though, in the words of [11, (1.6), p. 350], "many writers who should know better write $\varphi^*\mathcal{A}(\varphi)$ for the ideal $\mathcal{A}(\varphi) \cdot \mathcal{O}_{\widetilde{C}}$" (the latter being $\widetilde{\mathcal{A}}(\varphi)$ in our situation).

For the remainder of this remark, we write $\mathcal{A}$ and $\widetilde{\mathcal{A}}$ for $\mathcal{A}(\varphi)$ and $\widetilde{\mathcal{A}}(\varphi)$ respectively. By definition, $\varphi^*\mathcal{A} = \varphi^{-1}\mathcal{A} \otimes_{\varphi^{-1}\mathcal{O}_C} \mathcal{O}_{\widetilde{C}}$, and then the multiplication map induces an exact sequence of $\mathcal{O}_{\widetilde{C}}$-modules

$$0 \longrightarrow \mathcal{K} \longrightarrow \varphi^*\mathcal{A} \longrightarrow \widetilde{\mathcal{A}} \longrightarrow 0.$$

We shall see in the two following examples that in general $\mathcal{K}$ is non-trivial and $\varphi^*\mathcal{A}$ has a non-trivial torsion part, which is a manifestation of the fact that $\varphi$ is not flat. The conclusion is that the object we want to consider is indeed $\widetilde{\mathcal{A}}$ and not $\varphi^*\mathcal{A}$.

Our first example is the normalization of the node $xy = 0$. The adjoint ideal is $(x, y)$. Locally on the branch $y = 0$ in the normalization, the multiplication map $\varphi^{-1}\mathcal{A} \otimes_{\varphi^{-1}\mathcal{O}_C} \mathcal{O}_{\widetilde{C}} \to \mathcal{O}_{\widetilde{C}}$ corresponds to the multiplication map

$$(x, y) \otimes_{\mathbb{C}[x,y]/(xy)} \mathbb{C}[x] \longrightarrow \mathbb{C}[x],$$

which has non-trivial kernel, generated by $y \otimes 1$ and thus killed by multiplication by $x$. The upshot is that the torsion part of $\varphi^*\mathcal{A}$ is $\mathcal{O}_{\widetilde{\Delta}}$, where $\widetilde{\Delta}$ consists of the two points in the preimage of the node.

Our second example is the normalization of the cusp $y^2 - x^3 = 0$, which we see as $t \mapsto (t^2, t^3)$. The adjoint ideal is again $(x, y)$, equivalently $(t^2, t^3)$. The multiplication map $\varphi^{-1}\mathcal{A} \otimes_{\varphi^{-1}\mathcal{O}_C} \mathcal{O}_{\widetilde{C}} \to \mathcal{O}_{\widetilde{C}}$ corresponds to the multiplication map

$$(t^2, t^3) \otimes_{\mathbb{C}[t^2,t^3]} \mathbb{C}[t] \longrightarrow \mathbb{C}[t],$$

which has non-trivial kernel generated by $t^2 \otimes t - t^3 \otimes 1$ and thus killed by multiplication by $t^2$. Again the torsion part of $\varphi^*\mathcal{A}$ is $\mathcal{O}_{\widetilde{\Delta}}$, where this time $\widetilde{\Delta}$ is defined by the ideal $(t^2)$ of $\mathbb{C}[t]$.

## 3 − Equivalence of singularities

We say that two algebroid plane curves:

$$C = \mathrm{Spec}(\mathbb{C}[[X, Y]]/(f)), \quad D = \mathrm{Spec}(\mathbb{C}[[X, Y]]/(g))$$

are *isomorphic* or *analytically equivalent* if there is an isomorphism of $\mathbb{C}$-algebras

$$\mathbb{C}[[X, Y]]/(f) \cong \mathbb{C}[[X, Y]]/(g)$$

Because of the following result it is not restrictive to assume that $f$ is a polynomial when studying isomorphism classes of algebroid plane curves.

**(3.1) Theorem** (Samuel [14]). *Given a reduced non-constant $f \in \mathbb{C}[[X, Y]]$ there is a polynomial $g$ such that $\mathbb{C}[[X, Y]]/(f) \cong \mathbb{C}[[X, Y]]/(g)$.*

A weaker equivalence relation has been introduced by Zariski as follows. Consider and analytic branch $C = \mathrm{Spec}(\mathbb{C}[[X, Y]]/(f))$ and let $n = e(C)$ be its multiplicity (at 0). The blowing up of $\mathrm{Spec}(\mathbb{C}[[X, Y]])$ at 0 is defined by the substitution

$$X = X', \quad Y = X'Y'$$

Then

$$f(X, Y) = f(X', X'Y') = X'^n f'(X', Y')$$

and $f'$ is the *proper transform* of $f$. One easily shows (see [20], p. 3) that $f'$ is irreducible, i.e., defines an analytic branch $C'$. Let $e(C')$ its multiplicity. Then $e(C') \leq e(C)$. The theorem of resolution of singularities guarantees that after a finite number of iterations of this operation we obtain a nonsingular branch $\overline{C}$, i.e., such that $e(\overline{C}) = 1$. Let

$$C, C', \ldots, C^{(i)}, \ldots$$

be the successive proper transforms of $C$ and $e_i = e(C^{(i)})$. Let $N \geq 1$ be such that $e_N = 1$ but $e_{N-1} \neq 1$ and let

$$e^*(C) = \{e_0, e_1, \ldots, e_{N-1}\}$$

The sequence $e^*(C)$ is uniquely determined by $C$.

**(3.2) Definition.** *Two analytic branches $C, D$ are said to be* equivalent *or* equisingular, *if*

$$e^*(C) = e^*(D)$$

It is possible to show that two simple singularities are equisingular if and only if they are analytically equivalent. The simplest examples for which the two notions differ are ordinary $n$-fold points, with $n \geq 4$. For example, any two ordinary 4-fold points are equisingular, but they are isomorphic if and only if their respective 4-tuples of principal tangent lines have equivalent cross ratios. For more details about equisingularity we refer to [20].

## 4 – Generalities on deformations

Deformation theory studies local deformations of a given algebro-geometric object $X$ (a projective scheme, a singularity, etc.). Under good hypotheses it produces a family containing all sufficiently small deformations of $X$, up to isomorphism, and having certain functorial properties. Here we are interested in deformations of planar curve singularities, a case where the theory simplifies consistently. Before adventuring into the details we need a general preamble about what we mean by (local) deformations.

Consider the category $\mathfrak{Loc}$ of local $\mathbb{C}$-algebras with residue field $\mathbb{C}$. A *deformation* of $X$ parametrized by $\mathrm{Spec}(R)$, where $R \in \mathrm{Ob}(\mathfrak{Loc})$ is a pullback diagram:

$$\zeta : \quad \begin{array}{ccc} X & \hookrightarrow & \mathcal{X} \\ \downarrow & & \downarrow \pi \\ \mathrm{Spec}(\mathbb{C}) & \longrightarrow & \mathrm{Spec}(R) \end{array}$$

where $\pi$ is flat. Given another deformation of $X$ parametrized by $\mathrm{Spec}(R)$:

$$\eta : \quad \begin{array}{ccc} X & \hookrightarrow & \mathcal{Y} \\ \downarrow & & \downarrow \rho \\ \mathrm{Spec}(\mathbb{C}) & \longrightarrow & \mathrm{Spec}(R) \end{array}$$

we say that $\zeta$ and $\eta$ are *isomorphic* if there is an isomorphism of $R$-schemes

$$\begin{array}{ccc} \mathcal{X} & \xrightarrow{\cong} & \mathcal{Y} \\ & \pi \searrow \quad \swarrow \rho & \\ & \mathrm{Spec}(R) & \end{array}$$

compatible with the identifications of $X$ with the fibres over $\mathrm{Spec}(\mathbb{C})$.

Deformation theory associates to $X$ a covariant functor $M_X : \mathfrak{Loc} \to \mathfrak{Sets}$ taking values in the category of sets. It is defined by:

$$M_X(R) = \{\text{isom. classes of def.s of } X \text{ parametrized by } \mathrm{Spec}(R)\}$$

for all $R \in \mathrm{Ob}(\mathfrak{Loc})$, and $M_X(f : R \to S)$ is defined by pullback. This functor can be extended to the category of pointed $\mathbb{C}$-schemes in an obvious way. Since deformation theory is unable to give a satisfactory information on either of these functors what can be done is to restrict $M_X$ to more amenable subcategories of $\mathfrak{Loc}$. Two full subcategories are relevant in deformation theory. The first one is $\widehat{\mathfrak{Loc}}$, the category of complete local $\mathbb{C}$-algebras with residue field $\mathbb{C}$. The second one is $\mathfrak{Art}$, whose objects are the local artinian $\mathbb{C}$-algebras with residue field $\mathbb{C}$.

Covariant functors $F : \mathfrak{Art} \to \mathfrak{Sets}$ are called *functors of Artin rings*. The restriction of $M_X$ to $\mathfrak{Art}$ is such a functor. If $A \in \mathrm{Ob}(\mathfrak{Art})$ then the elements of $M_X(A)$ are called *infinitesimal deformations* of $X$ parametrized by $A$.

A *formal deformation* of $X$ is a pair $(R, \widehat{\xi})$ where $(R, \mathfrak{m})$ is a complete local $\mathbb{C}$-algebra, i.e., $R \in \mathrm{Ob}(\widehat{\mathfrak{Loc}})$, and $\widehat{\xi} = \{\xi_n \in M_X(R/\mathfrak{m}^{n+1})\}_{n \geq 1}$ is a sequence of infinitesimal deformations which are compatible, i.e., such that $\xi_n \mapsto \xi_{n-1}$ under the map

$$M_X(R/\mathfrak{m}^{n+1}) \to M_X(R/\mathfrak{m}^n)$$

for each $n \geq 1$. Associating to every $R \in \mathrm{Ob}(\widehat{\mathfrak{Loc}})$ the set of formal deformations $(R, \widehat{\xi})$ we define a new functor

$$\widehat{M_X} : \widehat{\mathfrak{Loc}} \to \mathfrak{Sets}.$$

A formal deformation $(R, \widehat{\xi})$ is called *effective* if there is a deformation $\xi \in M_X(R)$ such that $\xi \mapsto \xi_n$ under the map:

$$M_X(\pi_n) : M_X(R) \to M_X(R/\mathfrak{m}^{n+1})$$

induced by the projection $\pi_n : R \to R/\mathfrak{m}^{n+1}$, for each $n \geq 1$.

*Important point:* Not every formal deformation is effective, i.e., the functor $\widehat{M_X}$ is not isomorphic to the restriction of $M_X$ to $\widehat{\mathfrak{Loc}}$. Namely in general $\widehat{M_X}(R) \neq M_X(R)$.

A formal deformation $(R, \widehat{\xi})$ is called *algebraizable* if there is a deformation $\zeta$ of $X$ parametrized by a pointed algebraic scheme $(S, s_0)$ and an isomorphism $u : R \cong \widehat{\mathcal{O}_{S, s_0}}$ and such that the sequence of deformations

$$\widehat{\zeta} = \{\zeta_n \in M_X(\mathcal{O}_{S, s_0}/\mathfrak{m}^{n+1})\}$$

obtained by pulling back $\zeta$ is mapped to $\widehat{\xi}$ under $u$. Every algebraizable formal deformation is effective because $\zeta$ can be pulled back to $\widehat{\mathcal{O}_{S, s_0}}$, but the converse is not true.

To every $R \in \mathrm{Ob}(\widehat{\mathfrak{Loc}})$ one can associate the functor of Artin rings:

$$h_R : \mathfrak{Art} \to \mathfrak{Sets}, \quad A \mapsto \mathrm{Hom}(R, A).$$

Functors of this form are called *prorepresentable*. Consider a formal deformation $\widehat{\xi} \in \widehat{M_X}(R)$, for some $R \in \mathrm{Ob}(\widehat{\mathfrak{Loc}})$. Let $A \in \mathrm{Ob}(\mathfrak{Art})$ and $f \in h_R(A) = \mathrm{Hom}(R, A)$. Since $A$ is artinian there is an $n \gg 0$ such that $f$ factors as

$$\begin{array}{c} R \longrightarrow R/\mathfrak{m}^{n+1} \\ {}_{f} \searrow \downarrow {}^{f_n} \\ A \end{array}$$

One can then associate to $f$ the element $M_X(f_n)(\xi_n) \in M_X(A)$. Since this is independent of $n$ we have defined a morphism of functors of Artin rings that we denote with the same symbol:

$$\widehat{\xi} : h_R \to M_X.$$

**(4.1) Definition.** *If a formal deformation $\widehat{\xi} \in \widehat{M}_X(R)$ is such that $\widehat{\xi}$ is an isomorphism of functors then we say that the pair $(R, \widehat{\xi})$* prorepresents *$M_X$. In this case $M_X$ is* prorepresentable.

In practice if a pair $(R, \widehat{\xi})$ prorepresents $M_X$ it means that for every $A \in \mathrm{Ob}(\mathfrak{Art})$ every $\xi \in M_X(A)$ is induced as above by a unique $f : R \to A$.

The prorepresentability of $M_X$ is a strong property which is often not satisfied even by very simple objects $X$. For example, if $X$ is the singularity $XY = 0$ then $M_X$ is not prorepresentable (see [15, Example 2.6.8, p. 95]). The condition can be weakened in two ways as follows.

**(4.2) Definition.** *A formal deformation $\widehat{\xi} \in \widehat{M}_X(R)$ is called* versal *if for every $A \in \mathrm{Ob}(\mathfrak{Art})$ the induced map*
$$\widehat{\xi} : h_R(A) \to M_X(A)$$
*is surjective.*

*A formal deformation $\widehat{\xi} \in \widehat{M}_X(R)$ is called* miniversal *(or* semiuniversal*) if it is versal and in addition*
$$\widehat{\xi}(\mathbb{C}[\varepsilon]) : h_R(\mathbb{C}[\varepsilon]) \to M_X(\mathbb{C}[\varepsilon])$$
*is bijective.*

The map $\widehat{\xi}(\mathbb{C}[\varepsilon])$ is called the *differential* of $\widehat{\xi}$. Versality and miniversality are satisfied in most natural geometrical situations. Finally we give the following:

**(4.3) Definition.** *An algebraic deformation $\zeta$ of $X$ parametrized by a pointed algebraic scheme $(S, s_0)$ is called* formally versal *(resp.* formally semiuniversal*) if the associated algebraizable formal deformation $(\widehat{\mathcal{O}}_{S,s_0}, \widehat{\zeta})$ is versal (resp. semiuniversal).*

For more details on the contents of this section we refer to [15] and [10].

## 5 – The semiuniversal deformation of an isolated singularity

The case we will consider is $X = \mathrm{Spec}(B)$ where $(B, \mathfrak{m})$ is a noetherian local $\mathbb{C}$-algebra with residue field $B/\mathfrak{m} = \mathbb{C}$. Assume that $B$ is either essentially of finite type (eft), i.e., a localization of a $\mathbb{C}$-algebra of finite type, or complete, and that $\mathrm{Spec}(B)$ has an isolated singularity at $[\mathfrak{m}]$. Under these conditions deformation theory ensures the existence of a formal semiuniversal family of deformations of $B$, or of $\mathrm{Spec}(B)$ (see [15]). In the cases of interest for us such a deformation can be easily constructed explicitly. Note firstly the following:

**(5.1) Proposition.** *If $B$ and $B'$ are local $\mathbb{C}$-algebras as above such that $\widehat{B} = \widehat{B}'$, then their deformation functors are isomorphic.*

*Proof.* [10, ex. 18.6, p. 127]. □

This means that it is not restrictive to study deformations of a complete local $\mathbb{C}$-algebra. Since we are interested in planar curve singularities, i.e., in the local ring $B$ of a reduced curve $C$ contained in a nonsingular surface $Y$ at a point $p \in C$, we can replace $B$ by its completion and therefore we may assume that $= \mathbb{C}[[X,Y]]/(f)$, i.e., that $C$ is an algebroid plane curve.

Consider indeterminates $t_1, \ldots, t_N$ where $N = \tau(f)$ and let
$$p_1(X,Y), \ldots, p_N(X,Y) \in \mathbb{C}[[X,Y]]$$

be polynomials inducing a basis of the Tjurina algebra $T_f$. Consider the following family:

(5.1.1)
$$\begin{array}{ccc} C & \hookrightarrow & \mathrm{Spec}(\mathcal{B}) =: \mathcal{C} \\ \downarrow & & \downarrow \pi \\ \mathrm{Spec}(\mathbb{C}) & \longrightarrow & \mathrm{Spec}(\mathbb{C}[[t_1,\ldots,t_N]]) =: M(B) \end{array}$$

where
$$\mathcal{B} := \mathbb{C}[[t_1,\ldots,t_N,X,Y]]/(f(Y,Y) + \sum_1^N t_i p_i(X,Y)).$$

**(5.2) Proposition.** *The family* (5.1.1) *is a formally semiuniversal deformation of $C$.*

*Replacing $\mathbb{C}[[t_1,\ldots,t_N]]$ by $\mathbb{C}[t_1,\ldots,t_N]$ we obtain an algebraic formally semiuniversal family of deformations of $C$:*

(5.2.1)
$$\begin{array}{ccc} C & \hookrightarrow & \mathrm{Spec}\left(\mathbb{C}[t_1,\ldots,t_N][[X,Y]])/(f(Y,Y) + \sum_1^N t_i p_i(X,Y))\right) \\ \downarrow & & \downarrow \pi_a \\ \mathrm{Spec}(\mathbb{C}) & \longrightarrow & \mathrm{Spec}(\mathbb{C}[t_1,\ldots,t_N]) \end{array}$$

*Therefore, letting*
$$S_n := \mathrm{Spec}(\mathbb{C}[[t_1,\ldots,t_N]]/(t_1,\ldots,t_N)^{n+1}), \quad n \geq 1,$$

*the restrictions of* (5.1.1) *(or, what is the same, of* (5.2.1)*) to $S_n$ define a sequence of (infinitesimal) deformations of $C$:*

$$\begin{array}{ccc} C & \hookrightarrow & \mathcal{C}_n \\ \downarrow & & \downarrow \pi_n \\ \mathrm{Spec}(\mathbb{C}) & \longrightarrow & S_n \end{array}$$

*which is an algebraizable formal semiuniversal deformation $\widehat{\pi} = \{\pi_n\}$ of $C$.*

*Proof.* See [10], §14, for a detailed explicit proof. □

The support of $\mathcal{C}$ is different from the support of $\pi_n$, which is $C$, and a priori there is no information on how to recover $\mathcal{C}$ starting from $\widehat{\pi}$. Nevertheless from general properties of formal schemes it follows that $\pi$ is uniquely determined by $\widehat{\pi} := \{\pi_n\}$ (see e.g. [15], Theorem 2.5.11, p. 81). This fact implies that Proposition (5.2) can be strengthened as follows.

**(5.3) Corollary.** *a) For every deformation*

$$\begin{array}{ccc} C & \hookrightarrow & \mathcal{D} \\ \downarrow & & \downarrow \\ \mathrm{Spec}(\mathbb{C}) & \longrightarrow & \mathrm{Spec}(R) \end{array}$$

*of $C$ where $R \in \mathrm{Ob}(\widehat{\mathfrak{Loc}})$, there is a morphism $\varphi : \mathrm{Spec}(R) \to M(B)$ (not unique, but whose differential is unique) inducing an isomorphism*

$$\mathcal{D} \cong \mathrm{Spec}(R) \times_{M(B)} \mathcal{C}$$

*compatible with the isomorphisms of the closed fibres with $C$.*

   *b) Given a family of deformations of $C$:*

$$\begin{array}{ccc} C & \hookrightarrow & \mathcal{X} \\ \downarrow & & \downarrow \\ \mathrm{Spec}(\mathbb{C}) & \xrightarrow{s_0} & S \end{array}$$

*parametrized by an* algebraic *scheme $S$ (or by the spectrum of an eft local $\mathbb{C}$-algebra) there is a morphism (not unique, but whose differential is unique)*

$$\mathrm{Spec}(\widehat{\mathcal{O}}_{S,s_0}) \to \mathrm{Spec}(\mathbb{C}[[t_1,\ldots,t_N]]) = M(B)$$

*which induces an isomorphism of deformations over $\mathrm{Spec}(\widehat{\mathcal{O}}_{S,s_0})$:*

$$\mathrm{Spec}(\widehat{\mathcal{O}}_{S,s_0}) \times_S \mathcal{X} \cong \mathrm{Spec}(\widehat{\mathcal{O}}_{S,s_0}) \times_{M(B)} \mathcal{C}.$$

   For another discussion of this property see [10], Prop. 15.2, p. 108.

**(5.4) Remarks.** (1) If $C$ is nonsingular then $T_B^1 = 0$ and the semiuniversal deformation is:

$$\begin{array}{ccc} C & = & C \\ \downarrow & & \downarrow \\ \mathrm{Spec}(\mathbb{C}) & = & \mathrm{Spec}(\mathbb{C}) \end{array}$$

(see §1 for the definition of $T_B^1$).

   (2) If $C : f = 0$ is singular then it is immediate to verify that $f + \varepsilon = 0$ is nonsingular for small values of $\varepsilon$. This implies that $C$ has nonsingular deformations, i.e., that it is *smoothable*. The same proof shows that any germ of hypersurface in $\mathbb{C}^n$ with an isolated singularity at the origin is smoothable. An obvious modification of the argument shows that a germ of complete intersection isolated singularity is smoothable.

   (3) Versality is an open property, but not semiuniversality. In other words, the semiuniversal family (5.1.1) is versal in a neighborhood of 0. Obviously it is not semiuniversal, e.g. where the fibres are smooth, because there the fibres are rigid and their semiuniversal deformation just consists of the fibre itself (see Remark (1) above). The standard reference for openness of versality is [7].

**(5.5) Example.** (i) The semiuniversal deformation of $A_1 : X^2 - Y^2 = 0$ is $X^2 - Y^2 + t = 0$. The only singular fibre is for $t = 0$.

   (ii) Consider $A_2 : X^2 - Y^3 = 0$. Then the semiuniversal deformation is

$$F(t, u, X, Y) = X^2 - Y^3 + tY + u = 0$$

The singular fibres are determined by the conditions $F = F_X = F_Y = 0$ which are equivalent to $t = 3Y^2$, $u = -2Y^3$. We obtain the locus $\Delta : 27u^2 - 4t^3 = 0$ inside $M(B)$. The fibres are the cusp at $(0,0)$ and one node over the other points of $\Delta$. The fibres over $M(B) \setminus \Delta$ are nonsingular.

## 6 – Equigenericity and equisingularity

It is possible to define two interesting loci inside $M(B)$, the equigeneric locus and the equisingular locus.

The *equigeneric locus* $EG \subseteq M(B)$ is supported on the points of $M(B)$ whose fibre is a singularity imposing the same number of conditions to adjoints as $C$ does, i.e., on the points having the same $\delta$-invariant as $C$. The locus $EG \subseteq M(B)$ is only defined set-theoretically, it does not have a functorial definition. As a locally closed subset of $M(B)$ it is analytically irreducible near $0 \in M(B)$ and its tangent cone is a vector subspace of $T_f$. Identifying $T_f = B/J$ the tangent cone to $EG$ at $0$ is identified with $A/J$. In particular $A$ contains $J$ and $EG$ has codimension $\delta(C)$ in $M(B)$. In Example (5.5)(ii) we have $EG = \Delta$.

The *equisingular locus* $ES \subseteq M(B)$ is roughly the locus where the topological type of the singularity is the same as that of $C$. In particular $ES \subseteq EG$. In example (5.5)(ii) we have $ES = \{0\}$.

The locus $ES$ has a scheme structure because it has a (complicated) functorial definition that we do not need to recall here. $ES$ is nonsingular and its tangent space at $0 \in M(B)$ is the vector space $I/J \subseteq B/J = T_f$ where $I \subseteq A \subseteq B$ is an ideal called the *equisingular ideal*.

**(6.1) Lemma.** *(i) $I = A$ if and only if $B$ is an $A_1$-singularity (a node).*

*(ii) All the simple singularities satisfy $I = J$, i.e., they have no nontrivial equisingular deformations. Equivalently $ES = \{0\}$.*

*(iii) The simplest example such that $\dim(ES) > 0$, i.e., $J \neq I$, is an ordinary 4-ple point.*

For detailed information on equisingularity and many examples we refer to [8]. The following results due to Teissier are important in this context.

**(6.2) Theorem** (Teissier [17]). *Let $\mathcal{C} \to S$ be an equigeneric family of reduced curves, $S$ reduced. Then there is a Zariski dense open set $U \subseteq S$ such that $\mathcal{C} \times_S U \to U$ is equisingular.*

**(6.3) Theorem** (Teissier [17]). *Let $p : \mathcal{C} \to S$ be a flat family of reduced curves with $\mathcal{C}, S$ reduced of finite type. If $S$ is normal then the following are equivalent:*

*(a) $p : \mathcal{C} \to S$ is equigeneric.*

*(b) The composition*

$$\mathcal{C}' \xrightarrow{\nu} \mathcal{C} \xrightarrow{p} S$$

*where $\nu$ is the normalization, is smooth and for all $s \in S$ the induced morphism on the fibres $\nu_s : \mathcal{C}'_s \to \mathcal{C}_s$ is the normalization.*

**Note:** if $S$ is not assumed to be normal the theorem is false: see [5], Ex. (2.6).

## 7 – Families of curves on surfaces

We now fix a projective nonsingular algebraic surface $Y$, a class $\xi \in \mathrm{NS}(Y)$ and the scheme of curves of class $\xi$ in $Y$, denoted $\mathrm{Curves}_Y^\xi$. We have a universal family:

$$\begin{array}{ccc} \mathcal{C} & \hookrightarrow & \mathrm{Curves}_Y^\xi \times Y \\ \downarrow & & \\ \mathrm{Curves}_Y^\xi & & \end{array}$$

whose fibres have arithmetic genus $p_a(\xi)$ depending only of $\xi$. If $[C] \in \text{Curves}_Y^\xi$ then

$$p_a(\xi) = \frac{(C + K_Y) \cdot C}{2} + 1.$$

For a given $0 \leq g \leq p_a(\xi)$ one can consider the (possibly empty) locally closed locus $V_g^\xi \subseteq \text{Curves}_Y^\xi$ parametrizing reduced curves having geometric genus $g$. Inside $V_g^\xi$ one has the locus $V^{\xi,\delta}$ parametrizing reduced curves having $\delta = p_a(\xi) - g$ nodes and no other singularities. We want to compare these two loci, exploring the conditions under which we have $\overline{V^{\xi,\delta}} = \overline{V_g^\xi}$.

Let $C \subseteq Y$ be a reduced curve parametrized by a point of $V_g^\xi$. Dualizing the conormal sequence we obtain an exact sequence of sheaves on $C$:

$$0 \longrightarrow T_C \longrightarrow T_{Y|C} \longrightarrow \mathcal{O}_C(C) \longrightarrow T_C^1 \longrightarrow 0$$

where $T_C^1$ is the first cotangent sheaf of $C$, which has been described locally in §1. Note that in §1 we considered the case of algebroid plane curves. The module $T_B^1$ has an analogous definition for eft local rings and, being a torsion module in our case (isolated singularities), it is invariant under completion.

Locally at the point $v \in V_g^\xi$ parametrizing $C$ the restricted universal family

$$p : \mathcal{C}_g \to V_g^\xi$$

induces a family of deformations of the singular points of $C$. The semiuniversal families induce a morphism from $\text{Spec}(\widehat{\mathcal{O}}_{V_g^\xi,v})$ to the product of the bases of the semiuniversal families of deformations of the singularities of $C$. The differential of this morphism turns out to be precisely the map:

$$H^0(C, \mathcal{O}_C(C)) \to H^0(T_C^1)$$

appearing in the above exact sequence. Let's denote, as before, by $\mathcal{J} \subseteq \mathcal{I} \subseteq \mathcal{A} \subseteq \mathcal{O}_C$ respectively the jacobian, equisingular, and conductor ideal sheaves. Then this observation implies the following:

**(7.1) Proposition.** *If we identify $T_v(\text{Curves}_Y^\xi) = H^0(C, \mathcal{O}_C(C))$ then:*

(i) $H^0(C, \mathcal{J} \otimes \mathcal{O}_C(C))$ *is the tangent space to the subscheme of formally locally trivial deformations of $C$ in $Y$.*

(ii) $H^0(C, \mathcal{I} \otimes \mathcal{O}_C(C))$ *is the tangent space to the subscheme of equisingular deformations of $C$ in $Y$.*

(iii) $H^0(C, \mathcal{A} \otimes \mathcal{O}_C(C))$ *contains the reduced tangent cone to the subscheme of equigeneric deformations of $C$ in $Y$.*

## 8 – Deformations of morphisms

Given $C \subseteq Y$ as in the previous section we may consider a partial normalization morphism:

$$\varphi : \widetilde{C} \to C \subseteq Y$$

where we assume that $\widetilde{C}$ has at worst local complete intersection singularities. There is a well defined deformation theory for such a morphism, which produces a semiuniversal deformation

of $\varphi$:

$$\begin{CD} \widetilde{C} @>\Phi>> M(\varphi) \times Y \\ @VVV \\ M(\varphi) \end{CD}$$

where $M(\varphi) = \mathrm{Spec}(R)$ for some complete local $\mathbb{C}$-algebra. Semiuniversality here means a property analogous to the one considered in §5 for deformations of local rings; in particular $\varphi$ is the fibre of $\Phi$ over the the closed point $0 := [\mathfrak{m}_R] \in M(\varphi)$. The properties of $M(\varphi)$ are determined by the deformation functor $\mathrm{Def}_\varphi$ of $\varphi$. For the infinitesimal study of $M(\varphi)$ one considers the complex:

$$\Omega_\varphi^\bullet := [\varphi^* \Omega_Y^1 \to \Omega_{\widetilde{C}}^1].$$

The hyperext spaces $\mathrm{Ext}^1_{\mathcal{O}_{\widetilde{C}}}(\Omega_\varphi^\bullet, \mathcal{O}_{\widetilde{C}})$ and $\mathrm{Ext}^2_{\mathcal{O}_{\widetilde{C}}}(\Omega_\varphi^\bullet, \mathcal{O}_{\widetilde{C}})$ are respectively the tangent space and the obstruction space of $\mathrm{Def}_\varphi$. Then deformation theory implies the following:

**(8.1) Proposition.** *There is a canonical identification*

$$T_0 M(\varphi) = \mathrm{Ext}^1_{\mathcal{O}_{\widetilde{C}}}(\Omega_\varphi^\bullet, \mathcal{O}_{\widetilde{C}})$$

*and inequalities:*

$$\mathrm{ext}^1_{\mathcal{O}_{\widetilde{C}}}(\Omega_\varphi^\bullet, \mathcal{O}_{\widetilde{C}}) - \mathrm{ext}^2_{\mathcal{O}_{\widetilde{C}}}(\Omega_\varphi^\bullet, \mathcal{O}_{\widetilde{C}}) \leq \dim(M(\varphi)) \leq \mathrm{ext}^1_{\mathcal{O}_{\widetilde{C}}}(\Omega_\varphi^\bullet, \mathcal{O}_{\widetilde{C}}).$$

For details we refer to [6, Sec. 3.2]. Our assumptions imply that $\varphi$ is unramified at the general point of every component of $\widetilde{C}$. Therefore a simple computation (see [3, Lemma 11]) implies that there is a sheaf $N_\varphi$ which fits into an exact sequence:

(8.1.1) $\quad 0 \to \mathcal{H}om(\Omega_{\widetilde{C}}^1, \mathcal{O}_{\widetilde{C}}) \to \mathcal{H}om(f^*\Omega_Y^1, \mathcal{O}_{\widetilde{C}}) \to N_\varphi \to \mathcal{E}xt^1(\Omega_{\widetilde{C}}^1, \mathcal{O}_{\widetilde{C}}) \to 0$

and such that

$$H^0(C, N_\varphi) = \mathrm{Ext}^1_{\mathcal{O}_{\widetilde{C}}}(\Omega_\varphi^\bullet, \mathcal{O}_{\widetilde{C}}), \quad H^1(C, N_\varphi) = \mathrm{Ext}^2_{\mathcal{O}_{\widetilde{C}}}(\Omega_\varphi^\bullet, \mathcal{O}_{\widetilde{C}}).$$

The sheaf $N_\varphi$ is called *normal sheaf* of $\varphi$.

We will be interested in two special cases:

(a) $\widetilde{C}$ is nonsingular, i.e., $\varphi$ is the normalization map.

(b) $\varphi$ is unramified.

In case (a) the sheaf $\Omega_{\widetilde{C}}^1$ is locally free and therefore the sequence (8.1.1) becomes:

(8.1.2) $\quad 0 \longrightarrow T_{\widetilde{C}} \longrightarrow \varphi^* T_Y \longrightarrow N_\varphi \longrightarrow 0.$

The sheaf $N_\varphi$ has torsion in general. More precisely we have an exact sequence:

(8.1.3) $\quad 0 \longrightarrow \mathcal{H}_\varphi \longrightarrow N_\varphi \longrightarrow \overline{N}_\varphi \longrightarrow 0$

where $\mathcal{H}_\varphi = \mathcal{O}_Z$ for some effective divisor $Z$ supported on the locus where the differential $d\varphi$ degenerates. Therefore $\varphi(Z) \subseteq C$ is the set of points where $C$ has a non-linear branch. The exact sequences (8.1.2) and (8.1.3) imply that:

$$c_1(N_\varphi) = \omega_{\widetilde{C}} \otimes \varphi^* \omega_Y^{-1}$$

and
$$\overline{N}_\varphi \cong \omega_{\widetilde{C}} \otimes \varphi^* \omega_Y^{-1}(-Z).$$

We have the following useful:

**(8.2) Lemma** (Arbarello-Cornalba [4]). *If $\theta \in H^0(\mathcal{H}_\varphi)$ then the corresponding infinitesimal (first order) deformation of $\varphi$ leaves the image fixed.*

In case (b) we have an exact sequence

(8.2.1) $$0 \to \mathcal{K} \longrightarrow \varphi^*\Omega_Y^1 \longrightarrow \Omega_{\widetilde{C}}^1 \to 0$$

where $\mathcal{K}$ is invertible because $\Omega_{\widetilde{C}}^1$ has homological dimension 1 since $\widetilde{C}$ has l.c.i. singularities (which implies that $\Omega_{\widetilde{C}}^1$ is resolved by the conormal exact sequence of a local embedding). Then in the sequence (8.1.1) we have an identification $N_\varphi = \mathcal{K}^\vee$.

**(8.3) Lemma.** *In case (b) we have*
$$N_\varphi \cong \omega_{\widetilde{C}} \otimes \varphi^* \omega_Y^{-1}.$$

*Proof.* Since $\widetilde{C}$ is a local complete intersection, we have
$$\omega_{\widetilde{C}} = \det(\Omega_{\widetilde{C}}^1)$$
by [9, Thm. III.7.11]. On the other hand the exact sequence (8.2.1) yields
$$\det(\Omega_{\widetilde{C}}^1) = \det(\varphi^*\Omega_Y^1) \otimes \mathcal{K}^{-1} = \varphi^*\omega_Y \otimes N_\varphi.$$
The result follows since $\varphi^*\omega_Y$ is invertible. □

**(8.4) Corollary.** *In case (b) we have, in the notation of Section 2 for the conductor:*
$$N_\varphi \cong \varphi^* N_{C/Y} \otimes \widetilde{\mathcal{A}}(\varphi)$$
*and*
$$\varphi_* N_\varphi \cong N_{C/Y} \otimes \mathcal{A}(\varphi).$$

*Proof.* It is a direct consequence of the previous lemma and Corollary (2.4). □

## 9 – A criterion for the density of nodal curves

We keep the same notations of the previous section.

**(9.1) Theorem.** *Let $V \subseteq V_g^\xi$ be an irreducible component and let $[C] \in V$ be a general point with normalization $\varphi : \widetilde{C} \to Y$. Consider the following conditions:*

(a) $\omega_{\widetilde{C}} \otimes \varphi^* \omega_Y^{-1}$ *is globally generated.*

(b) $\dim(V) \geq h^0(\widetilde{C}, \omega_{\widetilde{C}} \otimes \varphi^* \omega_Y^{-1})$.

(c) $h^0(\omega_{\widetilde{C}} \otimes \varphi^* \omega_Y^{-1}(-\mathfrak{a})) = h^0(\omega_{\widetilde{C}} \otimes \varphi^* \omega_Y^{-1}) - 3$ *for every effective divisor $\mathfrak{a}$ of degree 3 on $\widetilde{C}$.*

*If (a) and (b) are satisfied then $C$ is immersed, i.e., all its local branches are linear. If also (c) is satisfied then $C$ is nodal.*

For the proof we refer to [10]. Let's see an application.

**(9.2) Corollary** (Arbarello-Cornalba [3], Zariski [19])**.** *The general element of any component of the Severi variety $V_{d,g}$ of integral plane curve of degree $d$ and genus $g$ is a nodal curve.*

*Proof (outline).* Let $[C] \in V_{d,g}$ be a general element of a component. We have

$$\deg(\omega_{\widetilde{C}} \otimes \varphi^* \omega_{\mathbb{P}^2}^{-1}) = 2g - 2 + 3d$$

Since $d \geq 3$ (otherwise there is nothing to prove) conditions (a) and (c) of the Theorem are clearly satisfied. One can construct a family of curves in $V_{d,g}$ as the image of the semiuniversal deformation of $\varphi$. Using (8.1.3) one can easily deduce that this family satisfies (b) as well. $\square$

# Lecture III
# Geometry of logarithmic Severi varieties at a general point

## by Thomas Dedieu

This lecture is devoted to the study of logarithmic Severi varieties of a pair $(S, R)$, where $S$ is a surface and $R$ is a curve on $S$: these are the families of curves on $S$ with prescribed homology class and geometric genus, and prescribed contact pattern with $R$, meaning that contact orders with $R$ are prescribed both at fixed assigned points and at unassigned, a priori mobile, points. Logarithmic Severi varieties are often referred to as relative Severi varieties.

Logarithmic Severi varieties naturally appear as irreducible components of limit Severi varieties in degenerations of smooth surfaces, see [I]. The emblematic example is that of logarithmic Severi varieties of the pair $(\mathbf{P}^2, \text{line})$ appearing in degenerations of Severi varieties of plane curves of a given degree and genus, see [VI] and [VII] in this collection. On the other hand, in degenerations of $K3$ surfaces we will rather encounter logarithmic Severi varieties of pairs $(S, R)$ such that $K_S + R = 0$, see Section 5.3.

The central result of this text is Theorem (1.7). On the one hand it gives the expected dimension of logarithmic Severi varieties, together with a condition ensuring that this is indeed the actual dimension. On the other hand it gives various regularity properties of the general members of logarithmic Severi varieties, under suitable assumptions on the positivity of $-K_S$ with respect to the prescribed homology class and contact conditions with $R$.

Another important result is Proposition (4.3), which characterizes logarithmic Severi varieties as essentially the families of curves of maximal dimension with respect to the homology class, the genus, and the number of moving contact points with $R$, under a positivity assumption on the logarithmic canonical divisor $-(K_S + R)$. It will be used in subsequent lectures to prove (in specific situations) that the irreducible components of a degeneration of logarithmic Severi varieties are again logarithmic Severi varieties (or rather, irreducible components thereof).

The main source for this lecture is Caporaso and Harris' [1, §2]; slightly earlier important works are [18] and [4]. I have adapted the presentation given in [1], which is restricted to the pair consisting of the projective plane and a line, in order to make it fit to the more general situation studied here, and needed further in this volume. My treatment is very much inspired by [10], which I worked out together with Edoardo Sernesi, see also [II] earlier in this volume. I have included many examples throughout, with the two purposes of clarifying the subtleties of the various definitions and statements, and of putting these results in perspective.

The organization of the text is as follows. Section 1 is devoted to the definition of logarithmic Severi varieties and the statement of the main theorem (1.7), and Section 2 contains the material in deformation theory necessary for the proof of the main theorem, which is given in Section 3. Section 4 is devoted to the aforementioned Proposition (4.3) characterizing logarithmic Severi varieties by the maximality of their dimensions, and its proof. Finally, Section 5 is devoted to expanded examples, featuring many examples of superabundant logarithmic Severi varieties, and highlighting some specificities of the case when the logarithmic canonical divisor $K_S + R$ is trivial; we also point out some interesting applications of logarithmic Severi varieties (e.g., to the tropical vertex group, and to $\mathbf{A}^1$-curves on open surfaces obtained as the complement of an anticanonical curve).

## 1 – Definitions and main results

**(1.1)** Let $S$ be a nonsingular, projective, connected, algebraic surface over the field $\mathbf{C}$ of complex numbers, and let $R \subseteq S$ be a reduced curve.

By a curve on $S$, we always mean a closed subscheme of $S$ of pure dimension 1. The *geometric genus* of a reduced curve $C$ is the (arithmetic) genus of its normalization $\bar{C}$, namely

$$1 - \chi(\mathcal{O}_{\bar{C}}) = \sum_{i=1}^{n} g_i - n + 1$$

where $g_1, \ldots, g_n$ are the respective genera of the connected components $\bar{C}_1, \ldots, \bar{C}_n$ of $\bar{C}$. Note that with this definition, denoting by $g$ the geometric genus of $C$, the canonical bundle of the normalization $\bar{C}$ has degree $2g - 2$, and one has the usual Riemann–Roch Formula on $\bar{C}$, i.e.,

$$\chi(L) = 1 - g + \deg(L)$$

for all line bundles on $\bar{C}$.

For every $\xi \in \mathrm{NS}(S)$ (thus, $\xi$ is a homology class on $S$ that can be represented by a divisor) and integer $g$, we consider the space $M_g^{\xi,\mathrm{bir}}(S)$ parametrizing morphisms

$$\phi : C \to S$$

from a smooth genus $g$ curve $C$ (projective, but possibly disconnected) that are birational on their image, and such that $\phi_*[C] = \xi$. The space $M_g^{\xi,\mathrm{bir}}(S)$ is defined more precisely and studied in Subsection 2.1.

We also consider the locally closed subscheme $V_g^\xi(S)$ of $\mathrm{Curve}(S)$, consisting of those points $[C]$ such that $C$ is a reduced curve having geometric genus $g$ and homology class $\xi$; by $\mathrm{Curve}(S)$, we denote the Hilbert scheme of curves on $S$.

**(1.2)** We denote by $\underline{\mathbf{N}}$ the set of sequences $\alpha = [\alpha_1, \alpha_2, \ldots]$ of non-negative integers with all but finitely many $\alpha_i$ non-zero. In practice, we may omit the infinite sequence of zeros at the end. We may also drop the brackets if only $\alpha_1$ is non-zero.

For all $\alpha \in \underline{\mathbf{N}}$, we let

$$|\alpha| = \alpha_1 + \alpha_2 + \cdots;$$
$$I\alpha = \alpha_1 + 2\alpha_2 + \cdots + n\alpha_n + \cdots.$$

For all $\alpha, \alpha' \in \underline{\mathbf{N}}$, we say that $\alpha \geqslant \alpha'$ if $\alpha_i \geqslant \alpha'_i$ for all $i \geqslant 1$.

By a set $\Omega$ of cardinality $\alpha \in \underline{\mathbf{N}}$, we mean a sequence of sets $\Omega = (\Omega_1, \Omega_2, \ldots)$ such that each $\Omega_i$ has cardinality $\alpha_i$.

**(1.3) Definition.** *Let $g \in \mathbf{Z}$, $\xi \in \mathrm{NS}(S)$, and $\alpha, \beta \in \underline{\mathbf{N}}$ be such that*

$$I\alpha + I\beta = \xi \cdot R,$$

*and consider a set $\Omega = \bigl(\{p_{i,j}\}_{1 \leqslant j \leqslant \alpha_i}\bigr)_{i \geqslant 1}$ of $\alpha$ points on $R$.*[1]

*We define $M_g^\xi(\alpha, \beta)(\Omega)$ as the locally closed subset of $M_g^{\xi,bir}(S)$ consisting of those $[\phi : C \to S]$ such that the intersection $\phi(C) \cap R$ is proper and contained in the smooth locus of $R$, and there exist $\alpha$ points $q_{i,j} \in C$, $1 \leqslant j \leqslant \alpha_i$, and $\beta$ points $r_{i,j} \in C$, $1 \leqslant j \leqslant \beta_i$, such that*

(1.3.1) $$\forall 1 \leqslant j \leqslant \alpha_i : \quad \phi(q_{i,j}) = p_{i,j} \quad \text{and}$$

(1.3.2) $$\phi^* R = \sum_{1 \leqslant j \leqslant \alpha_i} i\, q_{i,j} + \sum_{1 \leqslant j \leqslant \beta_i} i\, r_{i,j}.$$

Note that this definition is functorial. In other words, $M_g^\xi(\alpha, \beta)(\Omega)$ represents a certain functor, see (2.10).

**(1.4) Remark.** We emphasize that we do not assume the points in $\Omega$ to be general, and neither do we in the Main Theorem (1.7) below. It will be important for some applications to have this flexibility.

**(1.5) Definition.** *Continuing with the situation of Definition (1.3), we let $V_g^\xi(\alpha, \beta)(\Omega)$ be the locally closed subscheme of $V_g^\xi(S)$ consisting of those points $[C]$ such that the normalization $\nu : \bar{C} \to C \subseteq S$ belongs to $M_g^\xi(\alpha, \beta)(\Omega)$. We call it a* logarithmic[2] *Severi variety of the pair $(S, R)$.*

**(1.6) Notation, examples, and comments.** In practice we will try to find a balance between rigorous and decipherable notation. For instance we will frequently drop the $\Omega$, and replace $\xi$ with an adequate shorthand.

Let us consider the emblematic case when $S = \mathbf{P}^2$ and $R$ is a line. Let $\xi = d[H]$ with $[H]$ the hyperplane class; we shall simply denote $\xi$ by $d$. For example, for all $g$, we may write $V_g^d\bigl([0,1], d-2\bigr)$ for the family of plane $d$-ics of genus $g$, tangent to the line $R$ at some prescribed general point $p \in R$ (in this case, $\Omega = \{p\}$).

Let $\mathrm{p}(d)$ be the arithmetic genus of plane curves of degree $d$. The logarithmic Severi variety $V_{\mathrm{p}(d)}^d\bigl([0,1], d-2\bigr)$ parametrizes smooth $d$-ics tangent to $R$ at $p$; it has codimension 2 in the linear system $|dH|$. The logarithmic Severi variety $V^d_{\mathrm{p}(d)-1}([0,1], d-2)$ parametrizes plane $d$-ics of cogenus 1, tangent to $R$ at $p$; it has codimension 3 in the linear system $|dH|$ (hence $V^d_{\mathrm{p}(d)-1}([0,1], d-2)$ is a divisor in $V^d_{\mathrm{p}(d)}([0,1], d-2)$), and its general member is a curve with one node at a general point of $\mathbf{P}^2$, as we will see.

---

[1]If $R$ is reducible one should consider more precise data, specifying the distribution of the contact points with $R$ on its various components: let $R_1, \ldots, R_n$ be the irreducible components of $R$; one should consider $\alpha_1, \ldots, \alpha_n \in \underline{\mathbf{N}}$ and $\beta_1, \ldots, \beta_n \in \underline{\mathbf{N}}$ such that $I\alpha_i + I\beta_i = \xi \cdot R_i$ for all $i = 1, \ldots, n$, and $\Omega_1, \ldots, \Omega_n$ sets of points on $R_1, \ldots, R_n$ respectively, with cardinalities $\alpha_1, \ldots, \alpha_n$. For the purposes of this chapter we can safely ignore this, and we will do so for obvious reasons of simplicity. (There will be, however, some examples with reducible $R$, in Section 5 most notably, and then of course we will consider data with the necessary precision).

[2]in [1] it is a called a *generalized* Severi variety; the terminology *relative* Severi variety is nowadays quite popular; I prefer the term *logarithmic* which refers to the fact that we consider pairs instead of "absolute" surfaces, whereas *relative* often refers to the consideration of a map (or a family) instead of a surface.

We emphasize that Definition (1.3) requires that curves in $V^d_{\mathrm{p}(d)-1}([0,1], d-2)$ have one local branch with the prescribed order two of contact with $R$ at $p$; we illustrate this below.

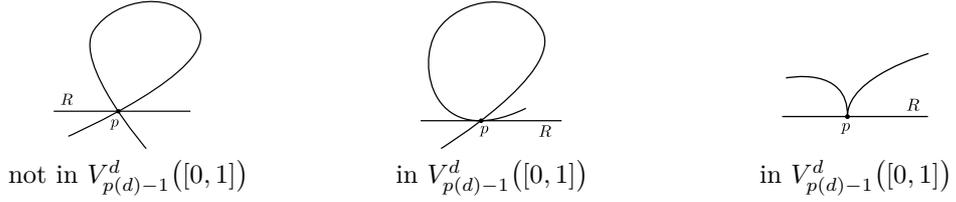

not in $V^d_{\mathrm{p}(d)-1}([0,1])$      in $V^d_{\mathrm{p}(d)-1}([0,1])$     in $V^d_{\mathrm{p}(d)-1}([0,1])$

Another subtlety of Definition (1.3) worth noting is the following. Consider the logarithmic Severi variety $V = V^d_{\mathrm{p}(d)}(1, [d-3, 1])$: loosely speaking, it parametrizes curves passing through an assigned point on $R$, call it $p$, and tangent to $R$. However, looking closely at the definition, one sees that a curve that is simply tangent to $R$ at $p$ and otherwise intersects it transversely does not belong $V$; a curve having an order 3 contact at $p$ and otherwise transverse does belong to $V$, on the other hand.

The following is the main result of this chapter. It is a vast generalization of the result proved independently by Arbarello and Cornalba [3] and Zariski [18] for plane curves of fixed degree and genus.

**(1.7) Theorem.** *Let $g \in \mathbf{Z}$, $\xi \in \mathrm{NS}(S)$, $\alpha, \beta \in \underline{\mathbf{N}}$, and $\Omega = \{p_{i,j}\}_{1 \leqslant j \leqslant \alpha_i} \subseteq R$ be as in Definition (1.3), and consider an irreducible component $V$ of $V_g^\xi(\alpha, \beta)(\Omega)$. Let $[C]$ be a general member of $V$, $\phi: \bar{C} \to C \subseteq S$ its normalization, $q_{i,j}$ ($1 \leqslant j \leqslant \alpha_i$), $r_{i,j}$ ($1 \leqslant j \leqslant \beta_i$) points on $\bar{C}$ such that* (1.3.1) *and* (1.3.2) *hold. Set*

$$D = \sum_{1 \leqslant j \leqslant \alpha_i} i\, q_{i,j} + \sum_{1 \leqslant j \leqslant \beta_i} (i-1)\, r_{i,j}.$$

*(1.7.00)* One has

$$\dim V \geqslant -(K_S + R) \cdot \xi + g - 1 + |\beta|.$$

*(1.7.0)* *If $-K_S \cdot C_i - \deg \phi_* D|_{C_i} \geqslant 1$ for every irreducible component $C_i$ of $C$, then the inequality in (1.7.00) above is an equality.*

*(1.7.1)* *If $-K_S \cdot C_i - \deg \phi_* D|_{C_i} \geqslant 2$ for every irreducible component $C_i$ of $C$, then:*

(a$^\flat$) *the normalization map $\phi$ is an immersion, except possibly at the points $r_{i,j}$;*

(b) *the points $q_{i,j}$ and $r_{i,j}$ of $\bar{C}$ are pairwise distinct;*

(c$^\flat$) *none of the points $s_{i,j} := \phi(r_{i,j})$ belongs to $\Omega$;*

(d) *for every curve $G \subseteq S$ and finite set $\Gamma \subseteq S$ such that $(G \cup \Gamma) \cap \Omega = \varnothing$, if $[C]$ is general with respect to $G$ and $\Gamma$, then $C$ intersects $G$ transversely and does not intersect $\Gamma$.*

*(1.7.2)* *If $-K_S \cdot C_i - \deg \phi_* D|_{C_i} \geqslant 3$ for every irreducible component $C_i$ of $C$, then:*

(a$^\flat$) *the normalization map $\phi$ is an immersion;*

(c) *the points $p_{i,j}$ and $s_{i,j} = \phi(r_{i,j})$ on $C$ are pairwise distinct;*

(e) *the curve $C$ is smooth at its intersection points with $R$.*

*(1.7.3)* If $-K_S \cdot C_i - \deg \phi_* D|_{C_i} \geqslant 4$ *for every irreducible component $C_i$ of $C$, then*
  (a) *the curve $C$ is nodal.*

The tangent space of $V_g^\xi(\alpha, \beta)(\Omega)$ at the general point $[C]$ of $V$ is described in Remark (3.9) under the assumptions of (1.7.2), and in Paragraph (3.6) under the weaker assumptions of (1.7.0).

We emphasize once again that it is not required in the above statement that the set $\Omega$ be general. The only requirement is that the points in $\Omega$ be smooth points of $R$, which in turn is imposed by the requirement, in Definitions (1.3) and (1.5), that the intersection with $R$ be proper and contained in the smooth locus of $R$.

We will see examples in Section 5 of the effect that imposing the passing through a set $\Omega$ of non-general points may have. We will also see examples of components of Severi varieties parametrizing reducible curves, for which the assumption of (1.7.0) fails for one irreducible component $C_i$, and having dimension strictly larger than expected (i.e., the inequality in (1.7.00) is strict).

## 2 – Background from deformation theory

Throughout this section, we consider a nonsingular, connected, projective surface $S$, defined over $\mathbf{C}$.

### 2.1 – Deformations of maps with fixed target

**(2.1)** Let $\phi : C \to S$ be a non-constant morphism from a smooth projective curve $C$. A *deformation of $\phi$ with fixed target* over a pointed base $(B, 0)$ is the data of a deformation $\mathcal{C} \xrightarrow{\pi} B$ of $C$ over $(B, 0)$ together with a morphism $\Phi : \mathcal{C} \to S \times B$ of $B$-schemes, such that the restriction of $\Phi$ over $0$ equals $\phi$.

This defines the *deformation functor* $\mathrm{Def}_{\phi/S}$ of $\phi$ with fixed target $S$. It is prorepresented by a complete local $\mathbf{C}$-algebra $R_\phi$, see [15, Thm. 3.4.8].

**(2.2)** The deformations of $\phi$ with fixed target $S$ are controlled by the *normal sheaf of $\phi$*, i.e., the sheaf $N_\phi$ of $\mathcal{O}_C$-modules defined by the exact sequence on $C$

(2.2.1) $$0 \to T_C \xrightarrow{d\phi} \phi^* T_S \to N_\phi \to 0 \quad :$$

the spaces $H^0(C, N_\phi)$ and $H^1(C, N_\phi)$ are respectively the Zariski tangent space and an obstruction space for the deformations of $\phi$ with fixed target $S$. In particular, we have

(2.2.2) $$\chi(N_\phi) \leqslant \dim R_\phi \leqslant h^0(N_\phi).$$

**(2.3)** The rank 1 sheaf $N_\phi$ may have torsion. We denote by $\mathcal{H}_\phi$ its torsion part and by $\bar{N}_\phi$ its maximal torsion-free quotient; they fit in an exact sequence

(2.3.1) $$0 \to \mathcal{H}_\phi \to N_\phi \to \bar{N}_\phi \to 0.$$

The torsion sheaf $\mathcal{H}_\phi$ is supported on the divisor $Z$ of zeroes of the differential $d\phi$, and it is zero if and only if $Z = 0$, i.e., if and only if the differential $d\phi$ is everywhere non-vanishing; in this case, we say that $\phi$ is an immersion. Moreover, there is an exact sequence of locally free sheaves on $C$

(2.3.2) $$0 \to T_C(Z) \to \phi^* T_S \to \bar{N}_\phi \to 0,$$

which readily implies the identification of line bundles on $C$

$$(2.3.3) \qquad \bar{N}_\phi \cong \omega_C \otimes \phi^* \omega_S^{-1}(-Z).$$

**(2.4) The scheme of morphisms.** We will pretend in this text that the global deformation functor of maps from smooth genus $g$ curves to the surface $S$ is represented by a scheme, which we will denote by $M_g(S)$. This comports the existence of a family $\mathcal{U} \to M_g(S)$ of smooth genus $g$ curves defined over $M_g(S)$, and of a universal morphism $\Phi : \mathcal{U} \to S \times M_g(S)$.

It is not strictly true that the global deformation functor of genus $g$ maps to $S$ is representable in the category of schemes, but it will be harmless for our purposes to pretend it is: one possible way to make things rigorously work out is to use Hartshorne's modular families [12, Def. p.171], as in [10, Subsec. 2.1]. A modular family of curves of genus $g$ is a surrogate of a universal family of curves of genus $g$. Let us consider such a modular family $\mathcal{C}_g/\mathcal{M}_g$, with the abuse of notation of denoting its base $\mathcal{M}_g$ although it is not the moduli space of genus $g$ curves. Then we may work with the scheme $\mathrm{Hom}_{\mathcal{M}_g}(\mathcal{C}_g, S \times \mathcal{M}_g)$ as a surrogate for $M_g(S)$, with Hom defined as in [14, Thm. I.1.10]. We will not enter in these technicalities here, and refer the reader to [10, Subsec. 2.1] for the complete details.

In this text, we will be mostly interested in the following subspace of $M_g(S)$. For all $\xi \in \mathrm{NS}(S)$, we consider $M_g^{\xi,\mathrm{bir}}(S)$ the subspace of $M_g(S)$ parametrizing morphisms $\phi$ that are birational on their image, and such that $\phi_*[C] = \xi$.

## 2.2 – Comparison of the spaces of maps and curves

Let $\xi \in \mathrm{NS}(S)$, and $g$ be a non-negative integer.

**(2.5) From maps to curves.** Consider the universal morphism $\Phi : \mathcal{U}_M \to S \times M_g^{\xi,\mathrm{bir}}$ defined over $M_g^{\xi,\mathrm{bir}}$. Let $\bar{M}$ be the semi-normalization of the reduced scheme underlying $M_g^{\xi,\mathrm{bir}}$, $\bar{\mathcal{U}}_M := \mathcal{U}_M \times_{M_g^{\xi,\mathrm{bir}}} \bar{M}$ the base change of the universal family $\mathcal{U}_M$, and $\bar{\Phi} : \bar{\mathcal{U}}_M \to S \times \bar{M}$ the induced morphism of $\bar{M}$-schemes. I claim that the scheme-theoretic image $\bar{\Phi}(\bar{\mathcal{U}}_M)$ is flat over $\bar{M}$.

Indeed, the morphism $\varpi := \mathrm{pr}_2 : \bar{\Phi}(\bar{\mathcal{U}}_M) \to \bar{M}$ is a well-defined family of codimension 1 algebraic cycles of $S$ in the sense of [14, I.3.11]. Since $\bar{M}$ is normal, the claim follows from [14, I.3.23.2].

It follows that there is a morphism from $\bar{M}$ to the Hilbert scheme of curves on $S$. It factorizes through $V_g^\xi(S)$, and actually through the normalization $\bar{V}_g^\xi \to V_g^\xi$ by the universal property of the normalization. Since by definition the semi-normalization morphism $\bar{M} \to M_g^{\xi,\mathrm{bir}}$ is $1:1$, two points $[\phi : C \to S], [\phi' : C' \to S] \in \bar{M}$ are mapped to the same point $[\Gamma] \in V_g^\xi$ if and only if there exists an isomorphism $\iota : C \cong C'$ such that $\phi = \phi' \circ \iota$.

**(2.6) From curves to maps.** On the other hand, consider the universal family $\mathcal{U}_V \to V_g^\xi$ of curves gotten from the universal family over the Hilbert scheme of curves on $S$. Let $\bar{V}$ be the normalization of $V_g^\xi$, and $\bar{\mathcal{U}}_V$ the normalization of $\mathcal{U}_V \times_{V_g^\xi} \bar{V}$. Teissier's Résolution Simultanée Theorem [17] asserts that $\bar{\mathcal{U}}_V \to \bar{V}$ is a family of smooth genus $g$ curves; it comes with a morphism of $\bar{V}$-schemes

$$\bar{\mathcal{U}}_V \to \mathcal{U}_V \times_{V_g^\xi} \bar{V} \subseteq S \times \bar{V}.$$

It follows that there is a morphism from $\bar{V}$ to the space $M_g^{\xi,\mathrm{bir}}$. It is generically injective, because the universal family of curves over $V_g^\xi$ is nowhere isotrivial.

From the considerations in (2.5) and (2.6) above, one deduces the following.

**(2.7) Proposition.** *Let $[\phi : C \to S] \in M_g^{\xi,bir}(S)$ be a general point (i.e., a general point of any irreducible component of $M_g^{\xi,bir}(S)$). Let $\Gamma := \phi(C)$, $\xi \in \mathrm{NS}(S)$ the homology class of $\Gamma$, and $g$ its geometric genus. Then $[\Gamma]$ belongs to a unique irreducible component of $V_g^\xi$ and*

$$\dim_{[\Gamma]} V_g^\xi = \dim R_\phi.$$

(Recall that $R_{\phi_0}$ is the complete local **C**-algebra that prorepresents $\mathrm{Def}_{\phi/S}$).

The next result provides a sharper upper bound on the dimension of the Severi varieties than that given by the inequality $\dim R_\phi \leqslant h^0(N_\phi)$ in (2.2.2).

Let $\phi : C \to S$ be a morphism from a smooth projective curve $C$, birational onto its image $\Gamma$. Let $\xi \in \mathrm{NS}(S)$ be the homology class of $\Gamma$, and $g$ its geometric genus. We consider $\Phi : \mathcal{C} \to S \times B$ a deformation of $\phi$ over a pointed normal connected scheme $(B, 0)$. Then $\Phi(\mathcal{C}) \subseteq S \times B$ is a deformation of $\Gamma$ over $(B, 0)$, see (2.5). There are thus two classifying morphisms $\kappa$ and $\gamma$ from $(B, 0)$ to $M_g^{\xi,\mathrm{bir}}(S)$ (or $R_\phi$) and $\mathrm{Curve}(S)$ respectively, with differentials

$$d\kappa : T_{B,0} \to H^0(C, N_\phi) \quad \text{and} \quad d\gamma : T_{B,0} \to H^0(\Gamma, N_{\Gamma/X}).$$

**(2.8) Lemma.** *The inverse image by $d\kappa$ of the torsion $H^0(C, \mathcal{H}_\phi) \subseteq H^0(C, N_\phi)$ is contained in the kernel of $d\gamma$.*

*Proof.* Given a non-zero section $\sigma \in H^0(C, N_\phi)$, the first order deformation of $\phi$ defined by $\sigma$ can be described in the following way: consider an affine open cover $\{U_i\}_{i \in I}$ of $C$, and for each $i \in I$ consider a lifting $\theta_i \in C(U_i, \phi^*T_X)$ of the restriction $\sigma_{|U_i}$. Each $\theta_i$ defines a morphism

$$\tilde\phi_i : U_i \times \mathrm{Spec}(\mathbf{C}[\varepsilon]) \to S$$

extending $\phi_{|U_i} : U_i \to X$. The morphisms $\tilde\phi_i$ are then made compatible after gluing the trivial deformations $U_i \times \mathrm{Spec}(\mathbf{C}[\varepsilon])$ into the first order deformation of $C$ defined by the coboundary $\partial(\sigma) \in H^1(C, T_C)$ of the exact sequence (2.2.1). In case $\sigma \in \mathrm{H}^0(C, \mathcal{H}_\phi)$, everyone of the maps $\tilde\phi_i$ is the trivial deformation of $\sigma_{|U_i}$ over an open subset. This implies that the corresponding first order deformation of $\phi$ leaves the image fixed, hence the vanishing of $dq_0(\sigma)$. □

**(2.9) Corollary.** *Let $g \in \mathbf{Z}$, $\xi \in \mathrm{NS}(S)$. Let $[C]$ be a general point of $V_g^\xi$, and $\phi : \bar C \to C \subseteq S$ its normalization. Then*

$$\dim V_g^\xi \leqslant h^0(\bar C, \bar N_\phi).$$

*Proof.* By generality we may assume that $[C]$ is a smooth point of $V_g^\xi$. Then $\dim V_g^\xi = \dim T_{[C]} V_g^\xi$, and by (2.6) there is a map

$$d\kappa_{[\phi]} : T_{[C]} V_g^\xi \to H^0(\bar C, N_\phi).$$

It is injective because to every tangent vector $\theta \in T_{[C]} V$ corresponds a non-trivial deformation of $C$. On the other hand, it follows from Lemma (2.8) that $\mathrm{Im}\, d\kappa_{[\phi]} \subseteq H^0(\bar C, \bar N_\phi)$. □

Lemma (2.8) (see also [II, Lemma (8.2)]) is a crucial observation (and indeed one of the cornerstones of the proof of Theorem (1.7)) that was made by Arbarello and Cornalba [4, p. 26], who deemed it a *fenomeno assai curioso*. They write: « nel caso in cui $\phi$ sia una birazionalità tra $C$ e $\Gamma$, la presenza di "cuspidi" su $\Gamma$, comporta l'esistenza, dal punto di vista infinitesimo, di più di un modello liscio della curva $\Gamma$, se così ci possiamo esprimere. »[3]

---

[3] a very curious phenomenon; in the case when $\phi$ is birational between $C$ and $\Gamma$, the presence of "cusps" on $\Gamma$ comports the existence, at the infinitesimal level, of more than one smooth model of the curve $\Gamma$, if we may say so.

Next, paraphrasing them, in order to use this phenomenon constructively they establish [4, Cor. 6.11]: in the above notation, if $B$ is the complex unit disk and if the family of curves is *equisingular*, then $d\kappa(\partial/\partial t)$ belongs to $H^0(C, \mathcal{H}_\phi)$ if and only if it is zero. Later, Caporaso and Harris (together with J. de Jong, they write) state and prove [1, Lem. 2.3]. They add the remark that this is linked to the notion of equisingularity, even though they make absolutely no use of that notion, neither in their statement nor in its proof.

The treatment I give above is that of Sernesi and myself in [10]. Although essentially equivalent to that of [1], it slightly differs in its formulation. This formulation, I hope, sheds some light on what is actually going on, and in particular shows that it is not necessary to invoke equisingularity in order to prove Corollary (2.9).

The standard application of Corollary (2.9) is to the proof that curves in a given family are immersed, i.e., they have no cusps: we will use it in Paragraph (3.7) below to prove assertion ($a^\flat$) of Theorem (1.7).

## 2.3 – Tangency conditions with respect to a fixed curve

We consider $R$ a fixed reduced curve on $S$. In this subsection we study deformations of curves on $S$ satisfying some tangency conditions with $R$; it follows from our Definition (1.3) that it suffices to treat the case when $R$ is smooth.

**(2.10)** Let $m$ be a non-negative integer. Let $\phi : C \to S$ be a non-constant morphism from a smooth projective curve $C$. A deformation of $\phi$ with fixed target *preserving a tangency of order $m$ with $R$* over a pointed base $(B, 0)$ is a deformation $\Phi : \mathcal{C} \to S \times B$ of $\phi$ with fixed target as in (2.1), such that there exists a section $Q$ of $\mathcal{C} \to B$ such that the pulled-back divisor $\Phi^* R$ contains $Q$ with multiplicity $m$ (i.e., $\Phi^* R - mQ \geqslant 0$).

$$\begin{array}{c} \mathcal{C} \xrightarrow{\Phi} S \times B \\ Q \left( \pi \downarrow \swarrow_{\mathrm{pr}_2} \right. \\ B \end{array}$$

The tangency is said to be respectively at a *fixed point* $p \in R$ if $\Phi(Q) = \{p\} \times B$, and at a *variable point* if $\mathrm{pr}_1 \circ \Phi(Q)$ is a curve.

We say that a family of maps $\Phi : \mathcal{C} \to S \times B$ preserves a tangency of order $m$ with $R$ if for all $b \in B$ it is locally around $b$ a deformation of maps preserving a tangency of order $m$.

The following result displays the additional conditions the class of a deformation of maps has to meet for this deformation to preserve a tangency with $R$. It is [1, Lem. 2.6]. Let $B$ be a reduced scheme, and $\Phi : \mathcal{C} \to S \times B$ be a family of maps preserving a tangency with $R$ of order exactly $m$. Let $b \in B$ be a general point, and $\phi : C \to S$ be the corresponding map. This comes with a classifying map $\kappa$, with differential $d\kappa : T_0 B \to H^0(C, N_\phi)$; we call $\bar{d\kappa}$ its composition with the projection $H^0(C, N_\phi) \to H^0(C, \bar{N}_\phi)$.

Let $q := Q \cap C$, and $l - 1$ be the order of vanishing of the differential $d\phi$ at $q$ (i.e., $l$ is the multiplicity of the point $p := \phi(q)$ in the local branch of $\phi(C)$ corresponding to $q$). Note that necessarily $l \leqslant m$.

**(2.11) Lemma.** *Let $\sigma \in \mathrm{Im}(\bar{d}\kappa) \subseteq H^0(C, \bar{N}_\phi)$ be a non-zero section, and denote by $v_q(\sigma)$ its order of vanishing at $q = Q \cap C$.*
(a) *One has $v_q(\sigma) \in \{m - l\} \cup [\![m, +\infty[\![$.[4]*
(b) *If the tangency is at a fixed point of $R$, then actually $v_q(\sigma) \in [\![m, +\infty[\![$.*

---
[4] I use the notation $[\![a, b]\!] = [a, b] \cap \mathbf{Z}$ for $a, b \in \mathbf{R} \cup \{\pm\infty\}$; albeit being standard probably only in France, it is convenient.

*Proof.* This is a local computation. Let $(x,y)$ be (analytic) local coordinates on $S$ at $p = \phi(q)$, such that $R$ is defined by the equation $y = 0$. Then the vector fields $\partial/\partial x$ and $\partial/\partial y$ generate $T_S$ near $p$, and their pull-backs generate $\phi^* T_S$ near $q$; by abuse of notation we shall denote them by $\partial/\partial x$ and $\partial/\partial y$ as well.

We may assume that $B$ is a curve. Let $\varepsilon$ be a local coordinate on $B$ centered at $b$, and $t$ be a local equation of the section $Q$ near $q$. Thus $(t, \varepsilon)$ are local coordinates on $\mathcal{C}$ at $q$. We may assume that $t$ gives a local coordinate on $C$ at $q$, in such a way that $\phi$ is given locally by

$$\phi(t) = \begin{cases} (t^l + a_{l+1}t^{l+1} + \cdots, t^m), & \text{if } l < m \\ (t^n + a_{n+1}t^{n+1} + \cdots, t^m) \text{ for some } n \geqslant m, & \text{if } l = m. \end{cases}$$

From now on we assume $l < m$ and leave the other, similar, case to the reader. Then the differential of $\phi$ at $q$ is

$$\tfrac{\partial}{\partial t} \longmapsto (lt^{l-1} + (l+1)a_{l+1}t^l + \cdots) \cdot \tfrac{\partial}{\partial x} + mt^{m-1} \cdot \tfrac{\partial}{\partial y}$$
$$= t^{l-1}\big((l + (l+1)a_{l+1}t + \cdots) \cdot \tfrac{\partial}{\partial x} + mt^{m-l} \cdot \tfrac{\partial}{\partial y}\big).$$

We see that around $q$ on $C$, the torsion part $\mathcal{H}_\phi$ of $N_\phi$ is a skyscraper sheaf of length $l-1$ concentrated at $q$, generated by the image in $N_\phi$ of the local section

$$\tau : t \mapsto \big(l + (l+1)a_{l+1}t + \cdots\big) \cdot \tfrac{\partial}{\partial x} + mt^{m-l} \cdot \tfrac{\partial}{\partial y}$$

of $\phi^* T_S$; thus the torsion-free quotient $\bar{N}_\phi$ is generated by the class of $\partial/\partial y$ modulo $\tau$. Observe that modulo $\tau$, $\partial/\partial x$ equals $t^{m-l} \cdot \partial/\partial y$ times an invertible, hence the image of $\partial/\partial x$ in $\bar{N}_\phi$ vanishes to the order exactly $m - l$ at $q$.

In turn the family $\Phi$ may be written locally as

$$\Phi(t, \varepsilon) = \big((t^l + a_{l+1}t^{l+1} + \cdots) + \varepsilon(u_0 + u_1 t + \cdots) + O(\varepsilon^2), t^m, \varepsilon\big),$$

since $\Phi^* R$ contains $Q$ with multiplicity $m$. By definition, the corresponding section $\bar{d}\kappa(\partial/\partial \varepsilon)$ of $\bar{N}_\phi$ is

$$\bar{d}\kappa(\tfrac{\partial}{\partial \varepsilon}) = (u_0 + u_1 t + \cdots) \cdot \tfrac{\partial}{\partial x} \bmod \tau.$$

Since $\partial/\partial x$ itself vanishes modulo $\tau$ to the order $m - l$ as we have seen, one has $v_q\big(\bar{d}\kappa(\partial/\partial \varepsilon)\big) \geqslant m - l$ in any event.

Moreover, by generality of $b \in B$ we may assume that all $\Phi(\cdot, \varepsilon)$ have their differential vanishing to the order $l$ at $Q \cap \mathcal{C}_\varepsilon$, which translates into the fact that $u_1 = \cdots = u_{l-1} = 0$. Then $\bar{d}\kappa(\partial/\partial \varepsilon)$ vanishes either to the order $m - l$, if $u_0 \neq 0$, or to some order larger than $m$, if $u_0 = 0$. When the tangency is maintained at a fixed point we are necessarily in the latter case. The lemma is proved. □

Let $\mathcal{C} \subseteq S \times B$ be a family of a reduced curves over a reduced base $B$. It is said to preserve a tangency of order $m$ with $R$ if the corresponding family of normalization maps over the normalization $\bar{B}$ of $B$ (see (2.6)) does.

**(2.12) Corollary.** *Let $V \subseteq \mathrm{Curve}(S)$ be a family of curves of genus $g$ having a tangency of order $m$ with the divisor $R$. Let $[C]$ be a general member of $V$, $\phi : \bar{C} \to C \subseteq S$ be the normalization of $C$, $q \in \bar{C}$ the tangency point, and $l - 1$ the order of vanishing of the differential $d\phi$ at $q$. Then*

$$\dim V \leqslant h^0\big(\bar{C}, \bar{N}_\phi(-(m-l)_+ \cdot q)\big),$$

*where $(m-l)_+$ denotes $\max(m-l, 0)$. If the tangency is at a fixed prescribed point on $R$, then actually*

$$\dim V \leqslant h^0\big(\bar{C}, \bar{N}_\phi(-m \cdot q)\big).$$

*Proof.* As in the proof of Corollary (2.9), we have $\dim V = \dim T_{[C]}V$ by generality of $[C]$, and there is an injective map $T_{[C]}V \to H^0(\bar{C}, \bar{N}_\phi)$. By Lemma (2.11) the image of this map is contained in $H^0(\bar{C}, \bar{N}_\phi(-a \cdot q))$ with $a = (m-l)_+$ if the tangency is mobile, and $a = m$ if the tangency is fixed. This ends the proof. □

## 3 – Proof of the Main Theorem

In this section we prove Theorem (1.7). The proof itself is in Subsection 3.2, after we give some lemmas in Subsection 3.1.

### 3.1 – Applications of the Riemann–Roch formula

Lemma (3.1) below is standard, but we will also use the more clever Lemma (3.2), which amounts to [1, Observation 2.5].

**(3.1) Lemma.** *Let $X$ be a smooth (possibly disconnected) projective curve, and $L$ a line bundle on $X$. Let $k \in \mathbf{N}^*$. If $\deg(L \otimes \omega_X^{-1}|_{X_i}) > k$ for all irreducible component $X_i$ of $X$, then the linear system $|L|$ separates any $k$ points on $X$.*

*Proof.* Let $Z$ be a subscheme of $X$ of length $k$, and let $Z'$ be another subscheme of $X$ such that $Z' \subsetneq Z$. The assumption on the degree of $L$ ensures that both $L(-Z)$ and $L(-Z')$ are non-special, hence $h^0(L(-Z)) < h^0(L(-Z'))$ by the Riemann–Roch formula. □

**(3.2) Lemma.** *Let $X$ be a smooth (possibly disconnected) projective curve of genus $g = 1 - \chi(\mathcal{O}_X)$, and $L, M$ two line bundles on $X$ such that*

$$\forall\, X_i \text{ irreducible component of } X: \quad \deg(M|_{X_i}) \leqslant \deg(L|_{X_i}).$$

(a). *If $\deg L \otimes \omega_X^{-1}|_{X_i} > 0$ for every component $X_i$ of $X$, then*

$$(3.2.1) \qquad h^0(X, M) \leqslant h^0(X, L) = \deg(L) - g + 1.$$

(b). *If $\deg L \otimes \omega_X^{-1}|_{X_i} > 1$ for every component $X_i$ of $X$, then equality holds in (3.2.1) if and only if $\deg M = \deg L$.*

*Proof.* Assumption (a) ensures that $L$ is non-special, hence the right-hand-side equality in (3.2.1), by the Riemann–Roch formula. If $M$ is non-special as well, then $h^0(L) - h^0(M) = \deg(L) - \deg(M)$ again by the Riemann–Roch formula, which gives the result in both cases (a) and (b).

We thus assume from now on that $M$ is special. Then,

$$h^0(X, M) \leqslant h^0(X, \omega_X) = g + n - 1.$$

Under Assumption (a),

$$\deg(L) = \sum_{i=1}^n \deg(L|_{X_i}) \geqslant \sum_{i=1}^n (2g_i - 1) = 2g - 2 + n,$$

where $n$ is the number of components of $X$ and $g_i$ is the genus of $X_i$ for all $i = 1, \ldots, n$; recall that $g = \sum g_i - n + 1$. By the Riemann–Roch formula, one has

$$h^0(X, L) = \deg(L) - g + 1 \geqslant g - 1 + n$$

hence $h^0(X, L) \geqslant h^0(X, M)$.

Under Assumption (b) one has in the same fashion $\deg(L) \geqslant 2g - 2 + 2n$, hence

$$h^0(X, L) \geqslant g + 2n - 1 > h^0(X, M).$$

This ends the proof, as the specialty of $M$ implies $\deg(M) < \deg(L)$ under the assumption of (a), and a fortiori under the assumption of (b). □

## 3.2 – Proof of the main theorem

This section is dedicated to the proof of Theorem (1.7). We consider an irreducible component $V$ of $V_g^\xi(\alpha, \beta)(\Omega)$ and a general point $[C]$ of $V$, as in the theorem. We use freely the notation introduced in the statement.

**(3.3)** *We start by proving that $V$ has expected dimension*

(3.3.1) $$\operatorname{expdim} V_g^\xi(\alpha, \beta)(\Omega) = -(K_S + R) \cdot \xi + g - 1 + |\beta|.$$

By (2.2.2) the expected dimension of $V_g^\xi$ is $\chi(N_\phi)$, which by the Riemann–Roch formula and the exact sequence (2.2.1) equals

(3.3.2) $$\operatorname{expdim} V_g^\xi = \deg \omega_{\bar{C}} - \deg \phi^* \omega_S + 1 - g = -K_S \cdot C + g - 1.$$

Now requiring that a curve $C$ have tangency of order $m$ with $R$ at a specified point $p$ is $m$ linear conditions on the coefficients of the equation of $C$, and if we let the point $p$ vary along $R$ the expected codimension of the corresponding locus of curves $C$ is one less, i.e. $m - 1$. We thus end up with

$$\operatorname{expdim}\bigl(V_g^\xi(\alpha, \beta)(\Omega)\bigr) = \operatorname{expdim}\bigl(V_g^\xi\bigr) - \sum_i \sum_{1 \leqslant j \leqslant \alpha_i} i - \sum_i \sum_{1 \leqslant j \leqslant \beta_i} (i - 1)$$
$$= \operatorname{expdim}\bigl(V_g^\xi\bigr) - I\alpha - (I\beta - |\beta|),$$

which together with (3.3.2) gives (3.3.1) after one remarks that $I\alpha + I\beta = R \cdot \xi$.

Note that this proves that, in any event,

(3.3.3) $$\dim(V) \geqslant -(K_S + R) \cdot \xi + g - 1 + |\beta|,$$

which is assertion (1.7.00). □

**(3.4)** *We now turn to the proof that the dimension of $V$ equals its expected dimension under assumption (1.7.0).*

Note that the points $q_{i,j}$ are necessarily pairwise distinct because they have distinct images $p_{i,j} \in R$. Let us first assume in addition that the points $q_{i,j}$ and $r_{i,j}$ are all together pairwise distinct; the case when this does not hold will be dealt with in (3.5).

We set

(3.4.1) $$D := \sum_{1 \leqslant j \leqslant \alpha_i} i \, q_{i,j} + \sum_{1 \leqslant j \leqslant \beta_i} (i - 1) \, r_{i,j},$$

the divisor on $\bar{C}$ of "infinitesimal tangency conditions with $R$" (compare (3.3)), and

(3.4.2) $$D_0 := \sum_{1 \leqslant j \leqslant \beta_i} (l_{i,j} - 1) \, r_{i,j},$$

where $l_{i,j} := v_{r_{i,j}}(d\phi)$ is the order of vanishing of the differential $d\phi$ at the point $r_{i,j}$, i.e., $D_0$ is the "ramification divisor of $\phi$ in the points $r_{i,j}$". We then decompose the difference of these two divisors as

$$(3.4.3) \qquad D - D_0 = D_1 - D_1'$$

where $D_1$ and $D_1'$ are non-negative divisors on $\bar{C}$ with disjoint supports; note that $D_1'$ may be nonzero only at the points $r_{i,j}$, and that it is so if and only if $l_{i,j} > i$.

It follows from Corollary (2.12) that

$$(3.4.4) \qquad \dim(V) \leqslant h^0\bigl(\bar{C}, \bar{N}_\phi(-D_1)\bigr).$$

Let $Z_0$ be the non-negative divisor on $\bar{C}$ such that the ramification divisor of $\phi$ is $D_0 + Z_0$. Then by (2.3.3) we have $\bar{N}_\phi \cong \omega_{\bar{C}} \otimes \phi^* \omega_S^{-1}(-D_0 - Z_0)$, and therefore (3.4.4) above reads

$$(3.4.5) \qquad \dim(V) \leqslant h^0\bigl(\bar{C}, \omega_{\bar{C}} \otimes \phi^* \omega_S^{-1}(-D_0 - Z_0 - D_1)\bigr)$$
$$(3.4.6) \qquad = h^0\bigl(\bar{C}, \omega_{\bar{C}} \otimes \phi^* \omega_S^{-1}(-D - D_1' - Z_0)\bigr)$$
$$(3.4.7) \qquad \leqslant h^0\bigl(\bar{C}, \omega_{\bar{C}} \otimes \phi^* \omega_S^{-1}(-D)\bigr).$$

Now by assumption (1.7.0) the line bundle $\omega_{\bar{C}} \otimes \phi^* \omega_S^{-1}(-D)$ is non-special, hence

$$(3.4.8) \qquad h^0\bigl(\bar{C}, \omega_{\bar{C}} \otimes \phi^* \omega_S^{-1}(-D)\bigr) = h^0\bigl(\bar{C}, \omega_{\bar{C}} \otimes \phi^* \omega_S(R)^{-1}(\phi^* R - D)\bigr)$$
$$(3.4.9) \qquad = 2g - 2 - (K_S + R) \cdot \xi + |\beta| + 1 - g$$
$$(3.4.10) \qquad = \operatorname{expdim} V_g^\xi(\alpha, \beta)(\Omega).$$

We thus have $\dim V \leqslant \operatorname{expdim} V_g^\xi(\alpha, \beta)(\Omega)$ which, together with (3.3.3) implies that $V$ has the expected dimension, if indeed the points $q_{i,j}$ and $r_{i,j}$ are all together pairwise distinct.

**(3.5)** Now if it is not true that the points $q_{i,j}$ and $r_{i,j}$ are all together pairwise distinct, then $V$ is actually a component of some Severi variety $V_g^\xi(\alpha', \beta')(\Omega')$ with $|\beta'| < |\beta|$ for which the corresponding points $q_{i,j}'$ and $r_{i,j}'$ are indeed pairwise disjoint (as sets $\Omega = \Omega'$, i.e., $\bigcup_i \Omega_i = \bigcup_i \Omega_i'$, and $\Omega_i \subseteq \bigcup_{k \geqslant i} \Omega_k'$).

Then, setting correspondingly $D'$ as in (3.4.1), one gets $\dim V \leqslant h^0(\bar{C}, \omega_{\bar{C}} \otimes \phi^* \omega_S^{-1}(-D'))$ exactly as in (3.4). Now $\deg D' > \deg D$ because $|\beta'| < |\beta|$, and it therefore follows from Lemma (3.2), part (a) that

$$(3.5.1) \qquad h^0\bigl(\bar{C}, \omega_{\bar{C}} \otimes \phi^* \omega_S^{-1}(-D')\bigr) \leqslant h^0\bigl(\bar{C}, \omega_{\bar{C}} \otimes \phi^* \omega_S^{-1}(-D)\bigr),$$

so that it still holds that $\dim V \leqslant \operatorname{expdim} V_g^\xi(\alpha, \beta)(\Omega)$, hence $V$ has the expected dimension. □

**(3.6)** *Note that the above proof gives the additional fact that, under the assumption of (1.7.0), the tangent space of $V_g^\xi(\alpha, \beta)(\Omega)$ at the general point $[C]$ of $V$ identifies with*

$$H^0\bigl(\bar{C}, \omega_{\bar{C}} \otimes \phi^* \omega_S^{-1}(-D)\bigr) \cong H^0\bigl(\bar{C}, \bar{N}_\phi(-D_1)\bigr) \cong H^0\bigl(\bar{C}, \bar{N}_\phi(D_0 - D)\bigr) \subseteq H^0\bigl(\bar{C}, \bar{N}_\phi\bigr).$$

**(3.7)** *We now prove that under Assumption (1.7.1) the assertions* (a$^\flat$), (b), (c$^\flat$), *and* (d) *hold.*

Suppose by contradiction that (b) doesn't hold. Then we argue as in (3.5). In this case, part (b) of Lemma (3.2) applies thanks to Assumption (1.7.1), and we get that the inequality (3.5.1) is strict, which is in contradiction with (3.3.3).

The same argument shows that none of the points $\phi(r_{i,j})$ can be fixed on $R$. This implies in particular that (c$^\flat$) holds.

The proof of (d) is similar: if $C$ were tangent to $G$, then it would belong to an irreducible component $W$ of some Severi variety of the pair $(S, R+G)$. Assumption (1.7.1) implies that $W$ is liable for part (1.7.0) of Theorem (1.7), hence $\dim(W) < \dim(V)$, in contradiction with the fact that $[C]$ is a general member of $V$. The same argument shows that $C$ avoids $\Gamma$ (pick some random curve on $S$ containing $\Gamma$).

Finally, we note that equality holds in (3.4.7) if and only if $D'_1 = Z_0 = 0$ since the line bundle $\omega_{\bar{C}} \otimes \phi^* \omega_S^{-1}(-D)$ is globally generated by assumption (1.7.1). Now it is indeed the case that equality holds in (3.4.7) since we have proved that $\dim(V) = \mathrm{expdim}(V_g^\xi(\alpha, \beta))$. We conclude that $D'_1 = Z_0 = 0$, which means that (i) $\phi$ is an immersion outside of the points $r_{i,j}$ (this is assertion (a$^\natural$)) and (ii) $l_{i,j} \leqslant i$ for $1 \leqslant j \leqslant \beta_i$. □

*Remark.* It is not enough to assume that $\omega_{\bar{C}} \otimes \phi^* \omega_S^{-1}(-D)$ is non-special and globally generated because of the argument we made to assume that the points $q_{i,j}$ and $r_{i,j}$ are pairwise distinct. We actually need to know something about every possible $\omega_{\bar{C}} \otimes \phi^* \omega_S^{-1}(-D')$, where $D' = \sum i q'_{i,j} + \sum (i-1) r'_{i,j}$ in the notation used for this argument.

**(3.8)** *We now prove that, under the assumption (1.7.2), $\phi$ is an immersion also at the points $r_{i,j}$, i.e. that $l_{i,j} = 1$ for $1 \leqslant j \leqslant \beta_i$, thus completing the proof of assertion* (a$^\flat$).

Let $i \geqslant 1$ and $1 \leqslant j \leqslant \beta_i$. It follows from the assumption that the linear series $\big|\omega_{\bar{C}} \otimes \phi^* \omega_S^{-1}(-D)\big|$ separates any two points, so there exists a section $\sigma \in H^0(\bar{C}, \omega_{\bar{C}} \otimes \phi^* \omega_S^{-1}(-D))$ with vanishing order $v_{r_{i,j}}(\sigma) = 1$ at the point $r_{i,j}$. Seen as a section $\tilde{\sigma} \in H^0(\bar{C}, \bar{N}_\phi)$, it vanishes at $r_{i,j}$ with order $v_{r_{i,j}}(\tilde{\sigma}) = 1 + (i - l_{i,j})$ (see (2.3.3)). By Lemma (2.11) this implies

$$1 + i - l_{i,j} \in \{i - l_{i,j}\} \cup [\![i, +\infty[\![$$

and therefore $l_{i,j} = 1$ as required. □

**(3.9) Remark.** Under the conditions of (1.7.2) the map $\phi$ is an immersion, hence $N_\phi = \bar{N}_\phi \cong \omega_C \otimes \phi^* \omega_S^{-1}$, see Paragraph (2.2). Therefore, the tangent space of $V_g^\xi(\alpha, \beta)(\Omega)$ at the general point $[C]$ of $V$ identifies with

$$H^0\big(\bar{C}, N_\phi(-D)\big) \cong H^0\big(\bar{C}, \omega_{\bar{C}} \otimes \phi^* \omega_S^{-1}(-D)\big)$$

(the definition of $D$ is in (3.4.1) above).

**(3.10)** *Let us prove that Assumption (1.7.2) implies Assertion* (c), *i.e., the points $p_{i,j}$ and $s_{i,j} = \phi(r_{i,j})$ are pairwise distinct.*

By (3.7), we already know that (c$^\flat$) holds, i.e., none of the $s_{i,j} = \phi(r_{i,j})$ belongs to $\Omega = \{p_{i,j}\}$. We thus only need to show that no two of the points $s_{i,j}$ coincide.

Suppose there exist $(i, j)$ and $(i', j')$ distinct such that $\phi(r_{i,j}) = \phi(r_{i',j'})$. Assumption (1.7.2) implies that the series $|\omega_{\bar{C}} \otimes \phi^* \omega_S^{-1}(-D)| = |T_{[C]}V|$ separates any two points. Therefore there exists a section $\sigma \in T_{[C]}V \cong H^0(\bar{N}_\phi(-D))$ such that

$$v_{r_{i,j}}(\sigma) = 1 \quad \text{and} \quad v_{r_{i',j'}}(\sigma) = 0.$$

This implies the existence of a deformation of $C$ in which the points $\phi(r_{i,j})$ and $\phi(r_{i',j'})$ no longer coincide, a contradiction to the generality of $[C]$ in $V$. □

**(3.11)** *Let us prove that Assumption (1.7.2) implies Assertion* (e), *i.e., $C$ is smooth at its intersection points with $R$.*

At this point we know that ($a^\flat$) and (c) hold under Assumption (1.7.2), i.e., the curve $C$ is immersed and the points $p_{i,j}$ and $s_{i,j}$ are pairwise distinct. Because the intersection $C \cap R$ is set-theoretically the union of all the points $p_{i,j}$ and $s_{i,j}$, this implies that $C$ is smooth at its intersection points with $R$. □

**(3.12)** *We eventually prove that under the Assumption (1.7.3) the curve $C$ is nodal, which is Assertion* (a) *of Theorem (1.7).*

Since we already know that the curve $C$ is immersed, it is enough to show that for all point $p \in C$, $C$ has neither three or more local branches, nor two or more tangent local branches.

To exclude the former possibility, suppose by contradiction that there exist $a, b, c \in \bar{C}$ pairwise distinct such that $\phi(a) = \phi(b) = \phi(c)$. The assumption (1.7.3) implies that the linear series $\left|\omega_{\bar{C}} \otimes \phi^*\omega_S^{-1}(-D)\right|$ separates any three points, so there exists a section $\sigma \in H^0\left(\bar{C}, \omega_{\bar{C}} \otimes \phi^*\omega_S^{-1}(-D)\right)$ such that

$$\sigma(a) = \sigma(b) = 0 \quad \text{and} \quad \sigma(c) \neq 0;$$

it corresponds to a first-order deformation of $\phi$ leaving both $\phi(a)$ and $\phi(b)$ fixed while moving $\phi(c)$. By generality of $[C]$, there is correspondingly an actual deformation of the curve $C$ for which the three local branches under consideration are no longer intersecting in a common point, a contradiction to the generality of $[C]$ in $V$ (see, e.g., [10, Prop. 1.4]).

We exclude the second possibility in a similar fashion. Suppose by contradiction that there exist $a, b \in \bar{C}$ distinct such that $\phi(a) = \phi(b)$ and $\text{Im}\, d\phi_a = \text{Im}\, d\phi_b$. Again since $\left|\omega_{\bar{C}} \otimes \phi^*\omega_S^{-1}(-D)\right|$ separates any three points, there exists $\sigma \in H^0\left(\bar{C}, \omega_{\bar{C}} \otimes \phi^*\omega_S^{-1}(-D)\right)$ such that $\sigma(a) = 0$ and $\sigma(b) \notin \text{Im}\, d\phi_a$. It corresponds to a first-order deformation of $\phi$ leaving $\phi(a)$ fixed while moving $\phi(b)$ in a direction transverse to the common tangent to the two local branches of $C$ under consideration. This contradicts the generality of $[C]$ as before. □

This ends the proof of Theorem (1.7).

## 4 – A dimensional characterization of logarithmic Severi varieties

We consider a pair $(S, R)$ consisting of a smooth surface $S$ and a reduced boundary divisor $R$, as described in 1. In this Section we give an upper bound for the dimension of families $V$ of (not necessarily reduced) curves in $S$ with prescribed homology class and genus, together with a dimensional characterization of logarithmic Severi varieties: it turns out that the irreducible components of logarithmic Severi varieties are essentially those families for which the upper bound is reached. The content of this section is basically [1, Cor. 2.7]. The main result, Proposition (4.3), is a generalization of Zariski's dimensional characterization of (classical) Severi varieties [18].

**(4.1)** Let $V$ be an irreducible locally closed subset of the Hilbert scheme of curves on $S$, and suppose that it parametrizes genus $g$ curves in the following sense: for a general member $X$ of $V$, there exists a smooth curve $C$ of genus $g$, and a morphism $\phi : C \to X$, not constant on any component of $C$, and such that the push-forward in the sense of cycles $\phi_*[C]$ equals the fundamental cycle of $X$.

**(4.2)** Let $\xi \in \mathrm{NS}(S)$ be the homology class of members of $V$. Letting the general member $X$ of $V$ have irreducible components $X_1, \ldots, X_n$, with respective classes $\xi_1, \ldots, \xi_n$, we assume furthermore:
$$\forall i = 1, \ldots, n : \quad -(K_S + R) \cdot \xi_i > 0,$$
so that we may apply (1.7.0) without any restriction on the contact conditions with $R$.

**(4.3) Proposition.** *Consider a finite subset $\Omega \subseteq R$. Under the above assumptions in (4.1) and (4.2), for all general member $X$ of $V$, one has*

(4.3.1) $$\dim(V) \leqslant -(K_S + R) \cdot \xi + g - 1 + \mathrm{Card}\big((X \cap R) \setminus \Omega\big),$$

*where the last number is defined set-theoretically (i.e., we do not count with multiplicities).*

*Moreover, if equality holds in (4.3.1), then there exist $\alpha, \beta$ such that $V$ is a dense subset of a component of a log-Severi variety $V_g^\xi(\alpha, \beta)(\Omega)$ if and only if*

(4.3.2) $$\mathrm{Card}\big(\phi^{-1}(R)\big) = \mathrm{Card}(X \cap R),$$

*with $\phi : C \to X$ a genus $g$ morphism as above.*

Note that $\alpha$ and $\beta$ are of course uniquely determined by the contacts of $X$ and $L$. Before proving this statement we make a few observations, in the hope that they will clarify some of its subtleties.

**(4.4) Remark.** Proposition (4.3) is really a result about families of embedded curves in $S$, not families of maps.

Indeed, if a general member $X$ is not reduced, then genus $g$ maps $\phi : C \to X$ as above involve multiple covers, and there is in general a positive dimensional family of maps giving the same $X$, so that the family of maps corresponding to $V$ has in general dimension larger than $V$.

When equality holds in (4.3.1), the map $\phi$ is necessarily a birational isomorphism on each component of $C$. If moreover $g \geqslant 1$, then the general member of $V$ is actually reduced. See (4.9) for precisions about this.

**(4.5) Remark.** A straightforward corollary of Proposition (4.3) which may be useful for the applications is that the inequality (4.3.1) still holds if we only assume $C$ to be reduced and replace $g$ by the *arithmetic* genus of $C$; this alternative inequality is always strict when $C$ is not smooth.

**(4.6) Remark.** Assumption (4.3.2) ensures that the normalization of $X$ is unibranch over the points in $X \cap R$.

*(4.6.1) Example.* Set $(S, R) = (\mathbf{P}^2, L)$ whith $L$ a line, and $\xi = 3[H]$ with $[H]$ the hyperplane class. The family $V$ of plane cubics with one node on $L$ and otherwise smooth has dimension 7, which equals
$$\underbrace{-(K_S + R) \cdot \xi}_{=2[H] \cdot 3[H] = 6} + \underbrace{g}_{=0} - 1 + \underbrace{\mathrm{Card}\big((X \cap R) \setminus \Omega\big)}_{=2} \quad \text{with } \Omega = \varnothing.$$

It is a divisor in the log-Severi variety $V_0^{3[H]}(0, 3)(\varnothing)$, but it is not a component of the family $V_0^{3[H]}(0, [1, 1])(\varnothing)$ of rational plane cubics with one variable tangency along $L$, which has dimension 7 as well.

**(4.7) Remark.** In the equality case for (4.3.1), if one replaces (4.3.2) with the weaker condition that

$$\text{(4.7.1)} \qquad \text{Card}\big(\phi^{-1}(R \setminus \Omega)\big) = \text{Card}\big((X \cap L) \setminus \Omega\big),$$

the families that we get that are not log-Severi varieties may be considered as "log-Severi varieties with $\Omega$ containing repeated points", see Example (4.7.2) below. (These are not log-Severi varieties according to Definition (1.3)).

Beware that in [1, Cor. 2.7], the weaker (4.7.1) is used instead of (4.3.2), which is slightly inaccurate.

*(4.7.2) Example.* Set $(S, R) = (\mathbf{P}^2, L)$ as above, and $\xi = 4[H]$, and fix $p \in L$. The family of plane quartics with one triple point at $p$ has dimension $\binom{6}{2} - 1 - 6 = 8$, which equals

$$\underbrace{-(K_S + R) \cdot \xi}_{=2[H] \cdot 4[H]=8} + \underbrace{g}_{=0} - 1 + \underbrace{\text{Card}\big((X \cap L) \setminus \Omega\big)}_{=1} \quad \text{with } \Omega = \{p\}.$$

One may wish to consider it as $V_0^{3[H]}(3, 1)(p, p, p)$.

**Proof of Proposition (4.3)**

We divide the proof into several steps, which correspond to paragraphs (4.8)–(4.10). We first treat the case when $X$ is irreducible; the general case is taken care of by induction in (4.10).

**(4.8)** We first prove the proposition under the assumption that $X$ is reduced and irreducible; in this case $\phi$ is the normalization of $X$ and $C \cong \bar{X}$.

If $e := \text{Card}\big((X \cap R) \setminus \Omega\big) = 0$, then the statement is a slight variant of part (1.7.0) of the main Theorem (1.7): we get the required inequality (4.3.1) exactly as in paragraph (3.4), with $D = \phi^* R$ and $D_0 = 0$ in (3.4.1) and (3.4.2) respectively. The only difference with the setting of (3.4) is that here two distinct points of the support of $\phi^* R \subseteq \bar{X}$ may have the same image by $\phi$ in $X$; this makes absolutely no difference in the argument.

Now if (4.3.2) holds, then $\bar{X}$ is unibranch over the points of $X \cap R$. This implies that $V$ is contained in a certain log-Severi variety $V_g^\xi(\alpha, 0)(\Omega)$. If moreover equality holds in (4.3.1), then $V$ is dense in an irreducible component of this same log-Severi variety by (1.7.0).

In the case when $e > 0$, we consider the map

$$\rho : V \to \text{Sym}^e R$$

sending a curve to its reduced intersection scheme with $R \setminus \Omega$; this may not be well-defined everywhere since $e$ may drop along certain closed subschemes of $V$, but it is in a neighbourhood of $[X]$.

Then we can apply the $e = 0$ case of the Proposition to the fibres of $\rho$; for a general $\Sigma \in \text{Sym}^e R$, setting $\Omega' = \Omega \cup \text{Supp}\,\Sigma$ one gets that the fibre $\rho^{-1}(\Sigma)$ has dimension at most $-(K_S + R) \cdot \xi + g - 1$. Inequality (4.3.1) follows, and the equality case of the Proposition as well, again applying the $e = 0$ case to the fibres of $\rho$.

*Remark.* Note that the above reasoning also proves that, if $V$ is dense in a suitable irreducible component of a log-Severi variety, then the map $\rho : V \to \text{Sym}^e R$ is dominant.

**(4.9)** Let us now consider the case when $X$ is non-reduced, but still irreducible; we shall show that inequality (4.3.1) holds, and that it is almost always strict; when equality holds in (4.3.1), condition (4.3.2) does not hold. This will also prove the second part of the proposition, since $V$ is not a dense subset of a component of a log-Severi variety when $X$ is non-reduced.

We have to consider $\phi : C \to X$ where $X$ is non-reduced but irreducible, and $C$ may be reducible. We let $m > 1$ be the degree of $\phi$, i.e. the sum of the degrees of the various $\phi_i : C_i \to X$. The key is to write the appropriate version of the Riemann–Hurwitz Formula: we have

$$2g - 2 = \deg(\omega_C) \geqslant m \deg(\omega_{X_{\mathrm{red}}}) = m(2q - 2),$$

with $q$ the arithmetic genus of $X_{\mathrm{red}}$. Let $h$ be the geometric genus of $X_{\mathrm{red}}$: we have $h \leqslant q$, so that eventually $g - 1 \geqslant m(h - 1)$.

Then we apply Proposition (4.3) in the reduced case, which we have already proven above, to the reduced curve $X_{\mathrm{red}}$ underlying $X$: it has class $\frac{1}{m}\xi$, so we get

$$\dim(V) \leqslant -(K_S + R) \cdot \frac{1}{m}\xi + h - 1 + \mathrm{Card}\big((X \cap R) \setminus \Omega\big)$$
$$\leqslant -(K_S + R) \cdot \frac{1}{m}\xi + \frac{g-1}{m} + \mathrm{Card}\big((X \cap R) \setminus \Omega\big);$$

recall that $X \cap R$ is considered as a set, so that

$$\mathrm{Card}\big((X \cap R) \setminus \Omega\big) = \mathrm{Card}\big((X_{\mathrm{red}} \cap R) \setminus \Omega\big).$$

Now the required inequality (4.3.1) follows from the basic inequality:

$$(\star) \qquad -(K_S + R) \cdot \frac{1}{m}\xi + \frac{g-1}{m} \leqslant -(K_S + R) \cdot \xi + g - 1.$$

If $g \geqslant 1$, the latter inequality $(\star)$ is straightforward, and always strict, as $-(K_S + R) \cdot \xi > 0$ and $m > 1$. If $g \leqslant 0$ however, one has $-(K_S + R) \cdot \xi \geqslant m$ on the one hand, because $\frac{1}{m}\xi$ is an integral class, and $g \geqslant -m + 1$, with equality holding if and only if $C$ is the disjoint union of $m$ smooth rational curves; thus inequality $(\star)$ holds also in this case, and it is an equality if and only if $-(K_S + R) \cdot \frac{1}{m}\xi = 1$ and $C$ is the disjoint union of $m$ smooth rational curves. When the latter condition is realized, condition (4.3.2) does not hold.

*(4.9.1) Example.* Let $S$ be the projective plane blown-up at eight general points, with exceptional divisors $E_1, \ldots, E_8$, and take $R = 0$. We consider the class

$$\xi = m(3H - \textstyle\sum_{i=1}^{8} E_i)$$

where $H$ as usual denotes the pull-back to $S$ of the line class on $\mathbf{P}^2$. Let $X_0$ be an irreducible rational curve in the linear system $|3H - \sum_{i=1}^{8} E_i|$, and consider the family $V$ consisting solely of $X = mX_0$ as a 0-dimensional family of curves of genus $-m + 1$, by means of the morphism $\phi : C \to X$ consisting of $m$ disjoint copies of the normalisation of $X_0$. Equality holds for $V$ in (4.3.1).

**(4.10)** It remains to consider the case when $X$ is reducible. Proceeding by induction on the number of irreducible components, we may assume that $X = X_1 \cup X_2$ where $X_1$ and $X_2$ move in two families $V_1$ and $V_2$ of curves of genera $g_1$ and $g_2$ respectively, such that $V \subseteq V_1 \times V_2$ and the Proposition holds for $V_1$ and $V_2$. Adding the two corresponding inequalities gives
(4.10.1)
$$\dim(V) \leqslant \dim(V_1) + \dim(V_2)$$
$$= -(K_S + R) \cdot \xi + (g_1 + g_2 - 1) - 1 + \mathrm{Card}\big((X_1 \cap R) \setminus \Omega\big) + \mathrm{Card}\big((X_2 \cap R) \setminus \Omega\big).$$

If $(X_1 \cap X_2 \cap R) \setminus \Omega$ is empty then this readily gives the required inequality, and the second part of the proposition in the equality case also follows.

So let us consider the case when $(X_1 \cap X_2 \cap R) \setminus \Omega$ is non-empty, and let us assume for simplicity that this set consists of only one point $p$, the general case being strictly similar. In this case condition (4.3.2) does not hold because $p$ has at least two preimages.

If $p$ is a fixed point of the intersection with $R$ for either one of the two families, say $V_1$, then it is also a fixed point for the other family: indeed, otherwise, by the generality of $X$ we could perturb $X_2$ a little so that it does not pass through $p$, and then $p$ is no longer a point in $X_1 \cap X_2 \cap R$. In this case, applying the proposition to $V_1$ and $V_2$ with $\Omega' := \Omega \cup \{p\}$, one gets

$$\dim(V) \leqslant -(K_S + R) \cdot \xi + (g_1 + g_2 - 1) - 1 + \underbrace{\mathrm{Card}\big((X_1 \cap R) \setminus \Omega'\big) + \mathrm{Card}\big((X_2 \cap R) \setminus \Omega'\big)}_{=\mathrm{Card}\big((X \cap R) \setminus \Omega\big) - 1},$$

and the result follows; note that inequality (4.3.1) is strict in this case.

Otherwise $p$ is variable for both $V_1$ and $V_2$; in this case $V$ necessarily has codimension at least 1 in $V_1 \times V_2$ (this may be proved for instance as in (4.8) by applying the Proposition to the fibres of the projection $V \to V_1$), and therefore (4.10.1) gives the required inequality (4.3.1). Equality may hold in (4.3.1), but $V$ cannot be dense in a component of a log-Severi variety, since it already has positive codimension in $V_1 \times V_2$.

*(4.10.2) Example.* Take $S = \mathbf{P}^2$ and $R$ a line, let $V_1 = V_2 = V_0^{[H]}(0,1)(\varnothing)$, and consider $V$ the sub-family of $V_1 \times V_2$ parametrizing the sums of two lines, the intersection point of which lies on $R$. Then $V$ has dimension 3, and parametrizes curves of genus $-1$, so that equality holds in (4.3.1) because $\mathrm{Card}(X \cap R) = 1$ for a general membre $X$ of $V$. It has codimension 1 in $V_1 \times V_2$, which is an irreducible component of $V_{-1}^{[2H]}(0,2)(\varnothing)$. Concretely, on $V$ we are artificially decreasing $\mathrm{Card}(X \cap R)$ by imposing that two contact points with $R$ coincide, and the latter prescription is not of Severi type.

The proof of Proposition (4.3) is now over. $\square$

## 5 – Examples

### 5.1 – An example with contact points in special position

In this short subsection I give an example illustrating (i) the effect of having the points in $\Omega$ in special position (which is also a theme in Subsection 5.3 below), and (ii) the importance of considering the positivity of $-K_S \cdot C - \deg \phi_* D$ (in the notation of Theorem (1.7)) on all irreducible components of $C$ separately.

**(5.1) Example.** Let $S$ be the projective plane, $R$ be a smooth cubic, and $\Omega$ be cut out on $R$ by another cubic $R'$. We consider the Severi variety $V_0^4(9,3)(\Omega)$ of rational quartics passing through the 9 intersection points of $R$ and $R'$ (and thus with 3 unassigned additional contact points with $R$); it has expected dimension

$$-\underbrace{(K_S + R)}_{=0} \cdot \xi + \underbrace{g}_{=0} - 1 + \underbrace{|\beta|}_{=3} = 2.$$

It has a 3-dimensional irreducible component $V_{\mathrm{ni}}$ parametrizing reducible quartics $C_1 + C_2$ where $C_1$ is a cubic in the pencil generated by $R$ and $R'$, and $C_2$ is a line. The component $V_{\mathrm{ni}}$ thus has dimension larger than expected, and indeed (1.7.0) does not apply in this case because, in the notation of Theorem (1.7),
$$-K_S \cdot C_1 - \deg \phi_* D|_{C_1} = 0.$$

The Severi variety $V_0^4(9,3)(\Omega)$ has another component $V_{\text{irr}}$, which parametrizes irreducible quartics $C$. This time,
$$-K_S \cdot C - \deg \phi_* D = 3,$$
so (1.7.0) applies and all irreducible components of $V_{\text{irr}}$ have dimension 2 as expected. It is however not obvious that there indeed exist irreducible rational quartics passing through $\Omega$ (so that $V_{\text{irr}}$ is non-empty): we give an argument in the next paragraph.

**(5.2)** Le me now show that $V_{\text{irr}}$ introduced above is non-empty, following a suggestion of Edoardo Sernesi. Consider the family $W$ of reducible quartics $C_1 + C_2$, where $C_1$ and $C_2$ are two conics, passing respectively through the first five points of $\Omega$, and through the remaining four points. Since no four points of $\Omega$ are aligned, $C_1$ is fixed and $C_2$ moves in a pencil, hence $W$ has dimension 1 as expected. Besides,
$$-K_S \cdot C_1 - \deg \phi_* D|_{C_1} = 1 \quad \text{and} \quad -K_S \cdot C_2 - \deg \phi_* D|_{C_2} = 2,$$
so (1.7.0) indeed applies to $W$.

The fact that the dimension of $W$ equals the expected dimension implies that the nodes of a general member of $W$ may be smoothed independently, i.e., a general member of $W$ may be deformed in such a way that an arbitrary subset of its nodes are preserved while the other are smoothed. Indeed the assumption implies that $W$ is smooth at a general point $[C]$, and the conditions defining its tangent space as a linear subspace of $H^0(C, \mathcal{O}_C(C))$ are independent, hence these conditions may be relaxed independently. This is a standard argument going back to Severi himself, of which a modern account may be found in [16].

The upshot is that one may smooth one and only one node of a general member $C_1 + C_2$ of $W$, while preserving the condition of passing through the points in $\Omega$: the result is a necessarily irreducible quartic curve with three nodes, passing through $\Omega$, as required.

## 5.2 – Superabundant log Severi varieties coming from double covers

In this subsection I give examples of superabundant logarithmic Severi varieties, i.e., such that the dimension exceeds the expected dimension: the inequality in (1.7.00) is strict, and (1.7.0) does not apply. These examples live on the projective plane, and come from linear systems on double covers. This is taken from [8]. The examples are given in (5.5) below, after a few recaps on double covers.

**(5.3)** We shall use the following elementary facts about double covers. Let $d$ be a positive integer, and $B$ be a degree $2d$ curve in $\mathbf{P}^2$. We consider the double cover $\pi : S \to \mathbf{P}^2$ branched over $B$. Let $H$ be the line class on $\mathbf{P}^2$, and $L$ be its pull-back to $S$. For all $k \in \mathbf{N}$ we have
$$H^0(S, kL) = \pi^* H^0(\mathbf{P}^2, kH) \oplus \pi^* H^0(\mathbf{P}^2, kH - \tfrac{1}{2}B),$$
which is the isotypical decomposition of $H^0(S, kL)$ as a representation of $\mathbf{Z}/2\mathbf{Z}$. The first summand corresponds to divisors that are double covers of degree $k$ curves in $\mathbf{P}^2$, and the second to divisors that decompose as $B$ (seen as the ramification divisor in $S$) plus the double cover of a degree $k - d$ curve in $\mathbf{P}^2$.

**(5.4) Proposition.** *For $k \geqslant d$, the general member $C$ of $|kL|$ is not a double cover of some hypersurface in $\mathbf{P}^2$, the restriction $\pi|_C$ is birational on its image, a degree $2k$ hypersurface $C^\flat$ in $\mathbf{P}^2$ everywhere tangent to $B$, with a node at every point of a complete intersection $Z$ of type $(k, k-d)$.*

*Proof.* The divisor $C$ belongs to a unique pencil $\langle A', B+D' \rangle$, with $A'$ and $D'$ the double covers of curves $A$ and $D$ in $\mathbf{P}^2$ of respective degrees $k$ and $k-d$. Thus $C^\flat := \pi(C)$ belongs to the pencil $\langle 2A, B+2D \rangle$, from which it follows that $C^\flat$ is double along $Z := A \cap D$, and touches $B$ doubly along $A \cap B$, which accounts for the whole intersection scheme of $C^\flat$ and $B$. The base locus of this pencil is the scheme defined by the ideal sheaf $\mathcal{I}_Z^2(\mathcal{I}_A^2 + \mathcal{I}_B)$.

The pull-back $\pi^* C^\flat \in |2kL|$ splits as $C + i(C)$, with $i$ the involution on $V$ associated to $\pi$; it has a double singularity along $Z' := \pi^{-1}(Z)$ and $\pi^{-1}(B \cap A)$, with at each point one local sheet belonging to $C$ and another to $i(C)$. The union $Z' \cup \pi^{-1}(B \cap A)$ is the base locus of the pencil $\langle A', B+D' \rangle$. □

**(5.5) Example.** We consider the image $V_{B,k}$ in $|2kH|$ of the linear system $|kL|$ on $S$. It has dimension
$$h^0(S, kL) - 1 = h^0(\mathbf{P}^2, kH) + h^0(\mathbf{P}^2, (k-d)H) - 1,$$
and parametrizes curves of geometric genus
$$g_{k,d} = \tfrac{1}{2}(2k-1)(2k-2) - k(k-d),$$
everywhere tangent to $B$; the number of contact points is thus $2kd$.

The family of curves $V_{B,k}$ is therefore contained in the log-Severi variety $V_{g_{k,d}}^{2k}(0, [0, 2kd])$ of the pair $(\mathbf{P}^2, B)$, for which the expected dimension is
$$-(K_{\mathbf{P}^2} + B) \cdot 2kH + g_{k,d} - 1 + 2kd = k(k+3-d).$$

By (1.7.0) a component of the Severi variety has the expected dimension if it has an irreducible member and
$$-K_{\mathbf{P}^2} \cdot 2kH - 2kd \geqslant 1 \iff 2k(3-d) \geqslant 1;$$
we note that the latter inequality holds if and only if $d \leqslant 2$.

It turns out that the dimension of our family $V_{B,k}$ exceeds the expected dimension of the log-Severi variety. Indeed a direct computation shows that
$$\dim(V_{B,k}) - \mathrm{expdim}\bigl(V_{g_{k,d}}^{2kH}(0, [0, 2kd])\bigr) = \tfrac{(d-1)(d-2)}{2}$$
$$= p_g(S)$$

(cf. [5, V.22 p.237] for the last equality).

### 5.3 – Logarithmic $K3$ surfaces

In this subsection I discuss logarithmic Severi varieties of pairs $(S, R)$ with $K_S + R = 0$, in relation with Severi varieties of $K3$ surfaces. This is a central theme of this volume, and as such it will be revisited in various other chapters.

**(5.6) $K3$ surfaces.** Let us first consider the case of "absolute" $K3$ surfaces: let $S$ be a $K3$ surface, and $R = 0$. In this case Theorem (1.7) is not quite accurate, a (well-known) prominent issue being that the condition in order for (1.7.0) to apply is not verified, and indeed the expected dimension given in (1.7.0) is not the actual dimension.

Suppose $S$ is equipped with a polarisation $L$ of genus $p$ (i.e., $L^2 = 2p - 2$). The expected dimension of $V_g^L$ as defined in (3.3) is $g - 1$, whereas its actual dimension is $g$, if $0 \leqslant g \leqslant p$ and $S$ is very general, say. Technically, the deformations of $[C] \in V_g^L$ are governed by the invertible sheaf $L|_C \cong \omega_C$, hence the obstruction space is $H^1(C, \omega_C)$ which is 1-dimensional; it turns out however that the equigeneric deformations of $C$ are unobstructed, even though the obstruction space is non-trivial. A conceptual explanation for this is that there exist non-algebraic $K3$ surfaces, which do not carry any curve at all; the expected dimension for Severi varieties on $K3$ surfaces thus becomes accurate if considered in family over all $K3$ surfaces, including the non-algebraic ones. I refer to [10, §4.2] for a detailed account; in the present volume, this phenomenon will be considered again in [**??**, subsection **??**], in the context of Gromov–Witten theory.

We shall now describe some analogous phenomena for $K3$ pairs, by looking at typical examples.

**(5.7) Surfaces with canonical curve sections.** Let $S$ be $\mathbf{P}^2$, and $R$ be a smooth cubic; note that in this case one has $K_S + R = 0$. Let $C$ be a smooth curve of degree $d \geqslant 4$ on $S$, and set $\Omega = C \cap R$ (for simplicity we may assume that $C$ and $R$ intersect transversely). The blow-up $S'$ of $\mathbf{P}^2$ at $\Omega$ is a smooth surface having a unique anticanonical divisor, namely the proper transform $R'$ of $R$. The linear system $|C'|$ of the proper transform $C'$ of $C$ gives a birational model of $S'$ in $\mathbf{P}^p$, in which the proper transform of $R$ is contracted, and whose hyperplane sections are the canonical models of the degree $d$ plane curves passing through $\Omega$ (we have let $p$ equal p($d$), the genus of smooth plane $d$-ics).

In many aspects the surface $S'$ behaves like a $K3$ surface. Beware, however, that it may not always be deformed to a $K3$ surface: this is the case if and only if $d \leqslant 6$; if $d > 6$ the elliptic singularity on $S'$ gotten from the contraction of the proper transform of $R$ is not smoothable, see [2].

Surfaces of this kind naturally occur in degenerations of $K3$ surfaces. For instance, surfaces $S'$ as above with $d = 4$ appear in degenerations of smooth quartics in $\mathbf{P}^3$ to the union of a smooth cubic and a plane, see [I, Section 7.2]; there we have seen that the logarithmic Severi variety $V_g^4(12, 0)(\Omega)$ of the pair $(\mathbf{P}^2, R)$ appears in the limit of the Severi varieties of genus $g$ hyperplane sections of smooth quartics.

There are of course many variants of this situation. For instance, Du Val surfaces[5] have been brought in the foreground recently in [1]. In this volume we shall often consider the case when $S$ is a toric surface, and $R$ is a cycle of smooth rational curves, the sum of all toric divisors of $S$. An emblematic instance of this situation is that of the pair consisting of $\mathbf{P}^2$ and a triangle, endowed with the linear system of plane quartics passing through twelve points cut out as above on the triangle by a quartic; it occurs in the degeneration of quartic $K3$ surfaces in $\mathbf{P}^3$ to a tetrahedron, see [I, 8.2] in which it corresponds to the contribution of a face (see ibid.). Other examples are given in [**??**, Section **??**].

---

[5]a genus $p$ Du Val surface is the projective plane blown-up at nine points on a smooth cubic $R$, endowed with the linear system of the proper transforms of plane curves of degree $3g$ with a $g$-tuple ordinary point at the first eight blown-up points, and a $(g - 1)$-tuple ordinary point at the remaining blown-up point; these curves have geometric genus $g$, and the linear system has dimension $g$ (note that it has a tenth base point lying on $R$).

**(5.8) Logarithmic $K3$ surfaces.** Let $S' \subseteq \mathbf{P}^p$ be as in (5.7). For all $g = 0, \ldots, p$, we may view the Severi variety of genus $g$ hyperplane sections of $S' \subseteq \mathbf{P}^p$ as the logarithmic Severi variety $V_p^d(3d, 0)(\Omega)$ of the pair $(\mathbf{P}^2, R)$. Its expected dimension, as defined in (3.3), is $g - 1$, whereas its actual dimension is $g$; note that once again (1.7.0) does not apply, because the required inequality does not hold. In this case the discrepancy between the actual and expected dimensions comes from the fact that the points in $\Omega$ are not general points on $R$. This is illustrated in Examples (5.8.1) and (5.8.2) below.

It is pleasant to think of the pair $(\mathbf{P}^2, R)$ endowed with a set $\Omega$ as a *logarithmic $K3$ surface*, *algebraic* when $\Omega$ is cut out on $R$ by a plane curve as above, and *non-algebraic* otherwise, when the points in $\Omega$ are in general position. One may prefer to consider that the logarithmic $K3$ surface is the equivalent data of the pair $(S', R')$, where $S'$ is the blow-up of $\mathbf{P}^2$ at $\Omega$, and $R'$ is the proper transform of $R$. For a more conceptual definition, as well as a classification, I refer to [13].

In the logarithmic context, the manoeuver of taking non-algebraic $K3$'s into account to adjust the expected dimension takes the following guise: one chooses an arbitrary point $a \in \Omega$, and considers $V_p^d(3d - 1, 1)(\Omega - a)$ instead of $V_p^d(3d, 0)(\Omega)$. Statement (1.7.0) applies to the former, which thus has both expected and actual dimension equal to $p$, and it turns out that the mobile contact point with $R$ remains in fact immobile, so that all members of $V_p^d(3d-1, 1)(\Omega - a)$ automatically pass through the point $a$ as well.

*(5.8.1) Example.* Let us illustrate this in the concrete case $d = 4$ (then, $p = 3$). The linear system of plane quartics has dimension 14. If we take a set $\Omega$ of 12 general points on $R$, then the Severi variety $V_3^4(12, 0)(\Omega)$ is empty: $\Omega$ imposes 12 independent conditions on quartics, the linear system of quartics through $\Omega$ is 2-dimensional, and all its members are made of $R$ plus a line. On the other hand if we take $\Omega$ the complete intersection of $R$ with a smooth quartic $C$, then one sees using the restriction exact sequence

$$0 \to \mathcal{O}_{\mathbf{P}^2}(1) \to \mathcal{O}_{\mathbf{P}^2}(4) \to \mathcal{O}_R(4) \to 0$$

that the linear system of quartics through $\Omega$ is 3-dimensional, generated by $C$ and the net of reducible quartics containing $R$. If we take a set $\Omega$ of 11 general points on $R$, it imposes 11 independent conditions on quartics, and the linear system of quartics through $\Omega$ is 3-dimensional with a 12-th base point on $R$.

*(5.8.2) Example.* We may consider the linear system $|2C'|$ in a similar fashion. Seen on $\mathbf{P}^2$, it is the system of plane $(2d)$-ics with a node at each of the $3d$ points of $\Omega = C \cap R$.

Again we shall work this out in the case $d = 4$. One has $(2C')^2 = 4 \cdot (C')^2 = 16 = 2 \cdot 9 - 2$, so the adjunction formula on $S'$, which is essentially the same as on a $K3$ surface, tells us that curves in $|2C'|$ have genus 9. Moreover the Riemann–Roch Formula, which works as on a $K3$ as well, tells us that $|2C'|$ has dimension 9.

On the other hand plane octics have arithmetic genus 21, so an octic with 12 nodes has geometric genus 9, confirming the above computations carried out on $S'$. The linear system of plane octics has dimension 44. Since a node at a prescribed point is 3 conditions, the expected dimension of a linear system of octics with 12 nodes at prescribed points is $44 - 36 = 8$. When the nodes are imposed at a general set of points $\Omega$ on $R$, this is indeed the correct dimension, and the linear system consists only of curves made of $R$ plus a quintic curve passing (in general simply) through $\Omega$. When $\Omega$ is cut out on $R$ by a quartic curve $C_1$, the linear system has an extra generator (namely, the curve $2C_1$), and thus has one extra dimension.

This may also be verified using a resolution of the ideal sheaf $\mathcal{I}_\Omega^2$, where $\mathcal{I}_\Omega$ is the ideal sheaf of $\Omega \subseteq \mathbf{P}^2$. Let $r$ and $f$ be homogeneous equations of the curves $R$ and $C$ respectively. While for $\mathcal{I}_\Omega$ there is the Koszul resolution, for its square we have the exact sequence

$$0 \to \mathcal{O}(-10) \oplus \mathcal{O}(-11) \xrightarrow{\begin{pmatrix} -f & 0 \\ r & -f \\ 0 & r \end{pmatrix}} \mathcal{O}(-6) \oplus \mathcal{O}(-7) \oplus \mathcal{O}(-8) \xrightarrow{(r^2, rf, f^2)} \mathcal{I}_\Omega^2 \to 0,$$

which gives $h^0(\mathcal{O}_{\mathbf{P}^2}(8) \otimes \mathcal{I}_\Omega^2) = 10$ as required.

This carries over to all systems $|kC'|$, to the effect that the condition of having a $k$-uple point at all $3d$ points of $\Omega = C \cap R$ imposes one less condition on plane $d$-ics than if the $3d$ points of $\Omega$ were in general position. We leave this to the reader.

**(5.9) Severi varieties of curves with full contact.** Logarithmic Severi varieties of logarithmic $K3$ surfaces $(S, R)$ parametrizing curves having one point of contact with high order with $R$ have attracted the attention for several reasons. A first reason is the possibility, if $(S, R)$ comes from a degeneration of $K3$ surfaces, to deform such a logarithmic Severi variety (roughly speaking) to a Severi variety on the nearby $K3$ surfaces, the contact point of high order $m$ deforming to $m-1$ nodes: this is the central theme of the present volume, in particular see [VIII]. Besides, this possibility is not specific to $K3$ surfaces.

Another reason is the search for $\mathbf{A}^1$ curves (i.e., curves isomorphic to the affine line, in other words affine rational curves) on logarithmic $K3$ surfaces, which is the analogue of the search for rational curves on $K3$ surfaces. Let me use an example to clarify this problem. Consider the projective plane together with a smooth cubic curve $R$ (and an empty set of points $\Omega$). The problem is to find $\mathbf{A}^1$ curves lying in the open surface $\mathbf{P}^2 \setminus R$. To this effect it is relevant to consider the logarithmic Severi varieties $V_0^d(0, [0, \ldots, 0, 3d])$ of the pair $(\mathbf{P}^2, R)$. They have expected dimension 0, and of course for any such curve the contact point with $R$ must be one of the $(3d)^2$ order $3d$ torsion points of $R$. The problem of the existence of $\mathbf{A}^1$ curves thus corresponds to that of the non-emptiness of these logarithmic Severi varieties. The characterization of the log-$K3$ surfaces for which there exists infinitely many $\mathbf{A}^1$ curves has been carried out in [7]. The interested reader may also consult [9] for other non-emptiness results, in particular [9, Prop. 3.12].

Logarithmic Severi varieties with one point of full contact have also been related in [11] with the tropical vertex group. Elements of this group are formal families of symplectomorphisms of the 2-dimensional algebraic torus $(\mathbf{C}^*)^2$. The relevant Severi varieties are the following: let $S$ be a toric surface, and let $R$ be the sum of all toric divisors $R_1, \ldots, R_n, R_{\text{out}}$ ($R_{\text{out}}$ is simply one of the toric divisors, to which we want to give a special role); one considers logarithmic Severi varieties of the pair $(S, R)$ of the form

$$V_0^\xi\big((\alpha_1, \ldots, \alpha_n, 0), (0, \ldots, 0, [0, \ldots, 0, m])\big)(\Omega_1, \ldots, \Omega_n, \varnothing),$$

which parametrize rational curves having fixed intersection schemes with $R_1, \ldots, R_n$, and one contact point of maximal order $m = \xi \cdot R_{\text{out}}$ with $R_{\text{out}}$. These logarithmic Severi varieties have expected and actual dimension 0, and the number of points they are made of may be interpreted as Gromov–Witten invariants. The main result of [11] is that these numbers are structure constants of the tropical vertex group; they are in fact determined by ordered product factorizations of commutators in the tropical vertex group, namely they appear as the coefficients of the logarithms of the factors.

## References

**Chapters in this Volume**

[I]   C. Ciliberto and Th. Dedieu, *Limits of nodal curves: generalities and examples.*

# Lecture IV
# Deformations of an $m$-tacnode into $m-1$ nodes

## by Margherita Lelli-Chiesa

This text is dedicated to a somewhat technical-looking statement (Theorem (1.2)) which is the cornerstone of many results in this volume: it is the key to understand the phenomenon of, loosely speaking, $(m-1)$-nodal curves degenerating to $m$-tacnodal curves within a degeneration of smooth surfaces to the transverse union of two surfaces. About this phenomenon see, in particular, [I, VI, B].

## 1 − Statement of results

Let $C$ be a planar curve. A point $P \in C$ is called an $m$-tacnode if the local equation of $C$ at $P$ has the form $y(y + x^m) = 0$, where $(x, y)$ are local coordinates centered at $P$. The Jacobian ideal of $C$ at $P$ is $J = (2y + x^m, yx^{m-1})$ and

$$\mathcal{O}_{\mathbf{A}^2}/J = \langle y, xy, \ldots, x^{m-2}y, 1, x, \ldots, x^{m-1}\rangle.$$

We denote by $\pi : S \to \Delta$ the versal deformation space of $P \in C$, where $\Delta \simeq \mathbf{A}^{2m-1}_{\alpha_0,\ldots,\alpha_{m-2},\beta_0,\ldots,\beta_{m-1}}$, and $S \subseteq \Delta \times \mathbf{A}^2_{x,y}$ is defined by the equation

(1.0.1) $\quad y^2 + yx^m + \alpha_{m-2}yx^{m-2} + \ldots + \alpha_1 yx + \alpha_0 y + \beta_{m-1}x^{m-1} + \ldots + \beta_1 x + \beta_0 = 0,$

see [II, Section 5]. We set

(1.0.2) $\quad\quad\quad \Delta_m \;\; := \;\; \{(\underline{\alpha}, \underline{\beta}) \in \Delta \,|\, \pi^{-1}(\underline{\alpha}, \underline{\beta}) \text{ has } m \text{ nodes}\},$

(1.0.3) $\quad\quad\quad \Delta_{m-1} \;\; := \;\; \{(\underline{\alpha}, \underline{\beta}) \in \Delta \,|\, \pi^{-1}(\underline{\alpha}, \underline{\beta}) \text{ has } m-1 \text{ nodes}\}.$

**(1.1) Proposition.** (1.1.1) *The locus $\Delta_m$ is smooth of dimension $m-1$, and it is defined in $\Delta$ by the equations $\beta_0 = \beta_1 = \ldots = \beta_{m-1} = 0$.*
(1.1.2) *The locus $\Delta_{m-1}$ has dimension $m$, and has $m$ local sheets at a general point of $\Delta_m$.*

*Proof.* The equation of $\pi^{-1}(\underline{\alpha}, \underline{\beta})$ has degree two in $y$ and its discriminant is given by

$$\delta_{\underline{\alpha},\underline{\beta}}(x) = (x^m + \alpha_{m-2}x^{m-2} + \ldots + \alpha_1 x + \alpha_0)^2 - 4(\beta_{m-1}x^{m-1} + \ldots + \beta_1 x + \beta_0).$$

If $V$ denotes the space of monic polynomials of degree $2m$ in $x$ with no degree $2m-1$ term, the regular map $\delta : \Delta \to V$ mapping a point $(\underline{\alpha}, \underline{\beta})$ to its discriminant $\delta_{\underline{\alpha},\underline{\beta}}(x)$ is an isomorphism;

indeed, $\underline{\alpha}$ and $\underline{\beta}$ can be easily expressed in terms of the coefficients of $\delta_{\underline{\alpha},\underline{\beta}}(x)$. We remark that, for fixed $(\underline{\alpha},\underline{\beta})$, the locus $\pi^{-1}(\underline{\alpha},\underline{\beta})$ is a double cover of the $x$-axis branched along the vanishing locus of $\delta_{\underline{\alpha},\underline{\beta}}(x)$, which is a divisor of degree $2m$.

A point $(\underline{\alpha},\underline{\beta})$ lies in $\Delta_m$ whenever $\delta_{\underline{\alpha},\underline{\beta}}(x)$ has $m$ double roots, that is, if and only if $\delta_{\underline{\alpha},\underline{\beta}}(x)$ is a square, equivalently if and only if $\beta_0 = \beta_1 = \ldots = \beta_{m-1} = 0$. This proves (1.1.1).

Analogously, $(\underline{\alpha},\underline{\beta})$ lies in $\Delta_{m-1}$ exactly when $\delta_{\underline{\alpha},\underline{\beta}}(x)$ has $m-1$ double roots, that is, we can rewrite the discriminant as

$$(1.1.3) \qquad \delta_{\underline{\alpha},\underline{\beta}}(x) = \left(\prod_{i=1}^{m-1}(x-a_i)^2\right) \cdot (x^2 + bx + c),$$

where in fact $b = 2\sum_{i=1}^{m-1} a_i$ because the degree $2m-1$ term of $\delta_{\underline{\alpha},\underline{\beta}}(x)$ vanishes. Hence, $\Delta_{m-1}$ has dimension $m$, and the $m$ local sheets at a general point of $\Delta_m$ as in (1.1.2) reflect the possibilities for choosing $m-1$ double roots out of $m$. □

Our main goal is to prove the following theorem.

**(1.2) Theorem** ([2, Lemma 2.8]). *Let $\Lambda \subseteq \Delta$ be a smooth $m$-dimensional variety such that:*
  *(i) $\Lambda$ contains $\Delta_m$;*
  *(ii) the tangent space $T_{\underline{0}}\Lambda$ of $\Lambda$ at $\underline{0}$ is not contained in the hyperplane $H$ defined by the equation $\beta_0 = 0$.*
*Then, we have $\Lambda \cap \Delta_{m-1} = \Delta_m \cup \Psi$, where $\Psi$ is a smooth curve intersecting $\Delta_m$ at $\underline{0}$ with multiplicity $m$.*

Note that $\Delta_{m-1}$ and $\Lambda$ both have codimension $m-1$ in $\Delta$, which has dimension $2m-1$, so that the expected dimension for their intersection is 1. The theorem thus states that if $\Lambda$ contains $\Delta_m$, under a suitable transversality assumption, namely (ii), the intersection $\Lambda \cap \Delta_{m-1}$ residual to the obvious superabundant component $\Delta_m$ has the expected dimension.

The reader may also consult [A] for complements and detailed examples.

## 2 – Proof of the main theorem in a model case

We first prove Theorem (1.2) in the special case where $\Lambda = \Lambda_0 = \{\beta_1 = \ldots = \beta_{m-1} = 0\}$.
Equation (1.0.1) of $S$ restricted to $\Lambda$ is:

$$(2.0.1) \qquad y^2 + yx^m + \alpha_{m-2}yx^{m-2} + \ldots + \alpha_1 yx + \alpha_0 y + \beta_0 = 0,$$

and the discriminant of a point $(\underline{\alpha},\underline{\beta}) \in \Lambda$ has the form:

$$(2.0.2) \qquad \delta_{\underline{\alpha},\underline{\beta}}(x) = (x^m + \alpha_{m-2}x^{m-2} + \ldots + \alpha_1 x + \alpha_0)^2 - 4\beta_0.$$

If $\beta_0 = 0$, then $(\underline{\alpha},\underline{\beta})$ lies in $\Delta_m$. Therefore, $\Lambda \cap \Delta_{m-1}$ contains $\Delta_m$ with multiplicity $m$, since by Proposition (1.1) $\Delta_{m-1}$ has multiplicity $m$ along $\Delta_m$.

From now on, we assume $\beta_0 \neq 0$, set $\nu(x) := x^m + \alpha_{m-2}x^{m-2} + \ldots + \alpha_1 x + \alpha_0$ and rewrite (2.0.2) as:

$$(2.0.3) \qquad \delta_{\underline{\alpha},\underline{\beta}}(x) = (\nu(x) - 2\sqrt{\beta_0})(\nu(x) + 2\sqrt{\beta_0}).$$

Note that the polynomials $\nu(x) - 2\sqrt{b}$ and $\nu(x) + 2\sqrt{b}$ have no common roots since $\beta_0 \neq 0$. We impose that $\delta_{\underline{\alpha},\underline{\beta}}(x)$ has $m-1$ double roots, taking advantage of the following lemma.

**(2.1) Lemma.** *Let $\gamma := 2\sqrt{\beta_0}$ with $\beta_0 \neq 0$. Then, for every integer $m \geqslant 2$ there exists a polynomial*
$$\nu_\gamma(x) := x^m + c_{m-2}(\gamma)x^{m-2} + \ldots + c_1(\gamma)x + c_0(\gamma)$$
*in $x$ with complex coefficients $c_i(\gamma)$ such that $c_{m-2}(\gamma) \neq 0$ and the following are satisfied:*

(a) *if $m = 2l+1$ is odd, then both $\nu_\gamma(x) + \gamma$ and $\nu_\gamma(x) - \gamma$ have $l$ double roots;*

(b) *if $m = 2l$ is even, then $\nu_\gamma(x) + \gamma$ (respectively, $\nu_\gamma(x) - \gamma$) has $l$ (resp., $l-1$) double roots.*

*Furthermore, $\nu_\gamma(x)$ is unique up to replacing it with $\nu_\gamma(\xi x)$ for some $m$-th root $\xi$ of unity, and it is odd (respectively, even) if $m$ is.*

Taking the lemma for granted, we set $m := 2l + \varepsilon$ with $\varepsilon \in \{0, 1\}$ depending on the parity of $m$ and consider the polynomial
$$\nu_1(x) = x^m + c_{m-2}(1)x^{m-2} + c_{m-4}(1)x^{m-4} + \ldots + c_\varepsilon(1)x^\varepsilon$$
obtained for $\gamma = 1$. For any other $\gamma \in \mathbf{C}^*$, by setting $\gamma := u^m$, we get $\nu_\gamma(x) = u^m \nu_1\left(\frac{x}{u}\right)$. The coefficients of $\nu_\gamma(x)$ are thus obtained from those of $\nu_1(x)$ in the following way:
$$c_{2h+\varepsilon}(\gamma) = c_{2h+\varepsilon}(u) = u^{2l-2h} c_{2h+\varepsilon}(1), \text{ for } 0 \leqslant h \leqslant l-1.$$

By the change of variable $t := u^2$, the equations
$$c_{2h+\varepsilon}(t) = t^{l-h} c_{2h+\varepsilon}(1), \text{ for } 0 \leqslant h \leqslant l-1$$
parametrize a curve $\Psi_0$ in $\Lambda \cap \Delta_{m-1}$. This curve is smooth because $c_{m-2}(1) \neq 0$, and its contact order with $\Delta_m$ at $\underline{0}$ is $m$ because $\beta_0 = \frac{t^m}{4}$. Therefore, it only remains to prove the lemma.

*Proof of Lemma (2.1).* We first prove the existence of a degree $m$ covering $f : \mathbf{P}^1 \to \mathbf{P}^1$ mapping $\infty$ to $\infty$ and with $\{\infty, \gamma, -\gamma\}$ as branch locus; we also require $f$ to be totally ramified at $\infty$ and the number of ramification points in the fiber $f^{-1}(-\gamma)$ (respectively, in $f^{-1}(\gamma)$) to coincide with the number of double roots we are requiring for the polynomial $\nu_\gamma(x) + \gamma$ (respectively, $\nu_\gamma(x) - \gamma$).

Assume that $m = 2l+1$ is odd. By Riemann's Existence Theorem, the existence of $f$ is equivalent to the existence of two permutations $\tau, \sigma$ in the symmetric group $S_m$, each one the product of $l$ disjoint transpositions, such that $\tau \cdot \sigma$ is cyclic of order $m$ in $S_m$. The permutations $\tau := (1\,2)(3\,4)\cdots(2l-1\,2l)$ and $\sigma := (2\,3)(4\,5)\cdots(2l\,2l+1)$ satisfy these requirements and they are the only elements of $S_m$ doing that, up to conjugation. As a consequence, $f$ exists and is unique up to automorphisms of the domain fixing $\infty$. The even case is analogous.

Since $f(\infty) = \infty$, the covering $f$ defines a polynomial $\nu_\gamma(x)$ satisfying (a) or (b) respectively. By acting with an automorphism of $\mathbf{P}^1$ fixing $\infty$ (i.e., composing $\nu_\gamma(x)$ with a polynomial of degree 1), we can assume $\nu_\gamma(x)$ to be monic and with no degree $m-1$ term. The polynomial $\nu_\gamma(x)$ is then unique up to replacing it with $\nu_\gamma(\xi x)$ for some $m$-th root $\xi$ of unity. By Proposition (A.1) in the Appendix, we can therefore assume that $\nu_1(x)$ coincides with the Chebyshev polynomial of the first kind $T_m(x)$; in particular, it is even or odd depending on the parity of $m$, and the coefficient $c_{m-2}(1)$ is nonzero, as one can easily check using (A.0.1) in the Appendix. The same properties hold for $\nu_\gamma(x)$ for any $\gamma = u^m \in \mathbf{C}^*$ because $\nu_\gamma(x) = u^m \nu_1\left(\frac{x}{u}\right)$. □

## 3 – Proof of the main theorem in general

We now prove Theorem (1.2) in full generality.

We consider the blow-up $\tau : \widetilde{\Delta} \to \Delta$ of $\Delta$ along $\Delta_m$ and denote by $Z := \tau^{-1}(\Delta_m)$ the exceptional divisor and by $\Phi := \tau^{-1}(\underline{0})$ the inverse image of $\underline{0}$. Then $\Phi \simeq \mathbf{P}^{m-1}_{\beta_0,\ldots,\beta_{m-1}}$ can be identified with the space of polynomials of degree at most $m-1$ in $x$ modulo scalars. The open subset $\Phi \supset \Phi_0 := \{\beta_0 \neq 0\}$ parametrizes polynomials that do not vanish in $x = 0$ and the point $Q = [1 : 0 : \ldots : 0] \in \Phi_0$ corresponds to constants. The proof of Theorem (1.2) in general is based on the description of the strict transform $\widetilde{\Delta_{m-1}}$ of $\Delta_{m-1}$ in $\widetilde{\Delta}$ and proceeds by steps.

**Step 1.** *The point $Q$ lies in $\widetilde{\Delta_{m-1}}$.*

Let $\widetilde{\Lambda_0}$ be the strict transform of the linear space $\Lambda_0 = \{\beta_1 = \ldots = \beta_{m-1} = 0\}$ under $\tau$. We have already shown that $\Lambda_0 \cap \Delta_{m-1} = \Delta_m \cup \Psi_0$ as in Theorem (1.2). Since $\widetilde{\Lambda_0} \cap \Phi = \{Q\}$, then $Q$ lies in the strict transform $\widetilde{\Psi_0}$ of the curve $\Psi_0$ and hence also in $\widetilde{\Delta_{m-1}}$.

**Step 2.** *The locus $\Phi$ is an irreducible component of $\widetilde{\Delta_{m-1}} \cap Z$ and is the only component containing $Q$.*

We identify $\Delta_m$ with the space of monic polynomials of degree $m$ in $x$ with no degree $m-1$ term. Let $\alpha(t) \subseteq \Delta_m$ be any arc such that $\alpha(0) = \underline{0}$. As $t$ goes to 0, all the roots of the polynomial corresponding to $\alpha(t)$ tend to $x = 0$ and hence the limit of $\tau^{-1}(\alpha(t)) \cap \widetilde{\Delta_{m-1}}$ is not contained in $\Phi_0$. As a consequence, any component of $\widetilde{\Delta_{m-1}} \cap Z$ containing $Q$ is contained in $\Phi$, and thus coincides with it by dimensional reason.

**Step 3.** *The strict transform $\widetilde{\Delta_{m-1}}$ is smooth at $Q$.*

The previous step implies that the only irreducible component of $\widetilde{\Lambda_0} \cap \widetilde{\Delta_{m-1}}$ containing $Q$ is the curve $\widetilde{\Psi_0}$. In particular, the intersection is proper in a neighborhood of $Q$ (indeed, $\dim \widetilde{\Lambda_0} = \dim \widetilde{\Delta_{m-1}} = m$, and $\dim \widetilde{\Delta} = 2m - 1$); smoothness of $\widetilde{\Psi_0}$ thus yields smoothness of $\widetilde{\Delta_{m-1}}$ at $Q$.

**Step 4.** *The strict transform $\widetilde{\Delta_{m-1}}$ is smooth at every point of $\Phi_0$.*

For every $c \in \mathbf{C}^*$, let $\sigma_c$ be the automorphism of $\Delta$ such that $\sigma_c(\alpha_i) = c^{m-i}\alpha_i$ for $0 \leqslant i \leqslant m-2$ and $\sigma_c(\beta_j) = c^{2m-j}\beta_j$ for $0 \leqslant j \leqslant m-1$. This defines an action of the multiplicative group $\mathbf{C}^*$ on $\Delta$, under which $\Delta_m$ is clearly invariant. Since $\delta_{\sigma_c(\underline{\alpha}),\sigma_c(\underline{\beta})}(x) = c^{2m}\delta_{\underline{\alpha},\underline{\beta}}\left(\frac{x}{c}\right)$, then $\Delta_{m-1}$ is invariant, too. As a consequence, the action of $\mathbf{C}^*$ on $\Delta$ lifts to an action on $\widetilde{\Delta}$ preserving $\widetilde{\Delta_{m-1}}$. Furthermore, $\underline{0}$ is a fixed point lying in the closure of any orbit and hence the lifted action also preserves $\Phi$. In particular, given $c \in \mathbf{C}^*$, we obtain an automorphism $\widetilde{\sigma_c}$ of $\Phi$ mapping a point $P$ of homogeneous coordinates $[\beta_0 : \beta_1 : \ldots : \beta_{m-1}]$ to

$$\widetilde{\sigma_c}(P) = [c^{2m}\beta_0 : c^{2m-1}\beta_1 : \ldots : c^{m+1}\beta_{m-1}] = [\beta_0 : c^{-1}\beta_1 : \ldots : c^{-(m-1)}\beta_{m-1}].$$

Since the limit for $c^{-1}$ going to 0 of $\widetilde{\sigma_c}(P)$ is $Q$ as soon as $\beta_0 \neq 0$, this shows that $Q$ lies in the closure of any orbit of $\widetilde{\Delta_{m-1}}$ intersecting $\Phi_0$. The statement thus follows from Step 3.

**Step 5.** *The intersection multiplicity of $\widetilde{\Delta_{m-1}}$ and $Z$ along $\Phi$ is $m$.*

Since the intersection multiplicity of $\Psi_0$ and $\Delta_m$ at $\underline{0}$ is $m$, then the same holds for that of $\widetilde{\Psi_0}$ and $Z$ at $Q$. On the other hand, the intersection $\widetilde{\Psi_0} \cap Z$ is contained in $\widetilde{\Delta_{m-1}} \cap Z \cap \widetilde{\Lambda_0}$, and the intersection multiplicity at $Q$ of $\Phi$ (which is the only component of $\widetilde{\Delta_{m-1}} \cap Z$ containing $Q$ by Step 2) with $\widetilde{\Lambda_0}$ is 1, because $\widetilde{\Lambda_0} \cap Z$ is a section of $\tau|_Z : Z \to \Delta_m$. The statement follows.

**Step 6.** *Conclusion.*

Let $\Delta_m \subseteq \Lambda \subseteq \Delta$ be as in the hypotheses. Since $\widetilde{\Lambda} \cap Z$ is a section of $\tau|_Z : Z \to \Delta_m$, then the intersection $\widetilde{\Lambda} \cap \Phi$ is a point $R$ with multiplicity one. The assumption $T_{\underline{0}}\Lambda \not\subseteq H = \{\beta_0 = 0\}$ ensures that $R$ lies in $\Phi_0$ and hence $\widetilde{\Delta_{m-1}}$ is smooth at $R$ by Step 4. The $m$-dimensional tangent space $T_R\widetilde{\Delta_{m-1}}$ contains $T_R\Phi$, that has dimension $m-1$. Since $T_R\widetilde{\Lambda} \cap T_R\Phi = \{0\}$ and $\dim T_R\widetilde{\Lambda} = m$, then $T_R\widetilde{\Delta_{m-1}}$ and $T_R\widetilde{\Lambda}$ generate the $(2m-1)$-dimensional tangent space $T_R\widetilde{\Delta}$. Equivalently, $\widetilde{\Delta_{m-1}}$ and $\widetilde{\Lambda}$ intersect transversally along a smooth curve $\widetilde{\Psi}$ in a neighborhood of $R$. Moreover, the curve $\widetilde{\Psi}$ is not tangent at $\Phi$ in $R$ (again because $T_R\widetilde{\Lambda} \cap T_R\Phi = \{0\}$) and its intersection multiplicity with $Z$ at $Q$ is $m$ by Step 5. Therefore, $\Psi := \tau(\widetilde{\Psi})$ is a smooth curve having intersection multiplicity $m$ with $\Delta_m$ at $\underline{0}$.

# A – Appendix: Chebyshev polynomials

There exist four families of Chebyshev polynomials, $T_m, U_m, V_m, W_m$ — called of the first, second, third and fourth kind — each satisfying

(A.0.1) $\qquad P_0(x) = 1, \quad \text{and} \quad \forall m \geqslant 1: \quad P_{m+1}(x) = 2xP_m(x) - P_{m-1}(x).$

They may be defined by the following choice of initial conditions:

$$\begin{aligned}
(A.0.2) \qquad T_1(x) &= x, \\
(A.0.3) \qquad U_1(x) &= 2x, \\
(A.0.4) \qquad V_1(x) &= 2x - 1, \\
(A.0.5) \qquad W_1(x) &= 2x + 1.
\end{aligned}$$

They are the unique polynomials satisfying

$$\begin{aligned}
(A.0.6) \qquad T_m(\cos\theta) &= \cos(m\theta), \\
(A.0.7) \qquad U_m(\cos\theta) &= \frac{\sin((m+1)\theta)}{\sin\theta}, \\
(A.0.8) \qquad V_m(\cos\theta) &= \frac{\cos((m+1/2)\theta)}{\cos(\theta/2)}, \\
(A.0.9) \qquad W_m(\cos\theta) &= \frac{\sin((m+1/2)\theta)}{\sin(\theta/2)}
\end{aligned}$$

for all $\theta \in \mathbf{R}$. The Chebyshev polynomials of the first kind can be explicitly written down as

$$T_m(x) = \frac{1}{2}\Big[\big(x + \sqrt{x^2-1}\big)^m + \big(x - \sqrt{x^2-1}\big)^m\Big],$$

and are even (respectively, odd) whenever $m$ is even (respectively, odd). Polynomials of the third and fourth kind satisfy $W_m(x) = (-1)^m V_m(x)$.

These are the first few values of Chebyshev polynomials:

$$\begin{aligned} T_0(x) &= 1, \\ T_1(x) &= x, \\ T_2(x) &= 2x^2 - 1, \\ T_3(x) &= 4x^3 - 3x, \\ T_4(x) &= 8x^4 - 8x^2 + 1, \ldots \end{aligned}$$

$$\begin{aligned} U_0(x) &= 1, \\ U_1(x) &= 2x, \\ U_2(x) &= 4x^2 - 1, \\ U_3(x) &= 8x^3 - 4x, \\ U_4(x) &= 16x^4 - 12x^2 + 1, \ldots \end{aligned}$$

$$\begin{aligned} V_0(x) &= 1, \\ V_1(x) &= 2x - 1, \\ V_2(x) &= 4x^2 - 2x - 1, \\ V_3(x) &= 8x^3 - 4x^2 - 4x + 1, \\ V_4(x) &= 16x^4 - 8x^3 - 12x^2 + 4x + 1, \ldots \end{aligned}$$

$$\begin{aligned} W_0(x) &= 1, \\ W_1(x) &= 2x + 1, \\ W_2(x) &= 4x^2 + 2x - 1, \\ W_3(x) &= 8x^3 + 4x^2 - 4x - 1, \\ W_4(x) &= 16x^4 + 8x^3 - 12x^2 - 4x + 1 \ldots \end{aligned}$$

The next result summarizes some properties of Chebyshev polynomials.

**(A.1) Proposition.** *If $m = 2l$ is even, then*

*(a1)* $T_m(x) + 1 = 2(T_l(x))^2$,

*(a2)* $T_m(x) - 1 = 2(x^2 - 1)(U_l(x))^2$;

*in particular, the polynomial $T_m(x) + 1$ is a perfect square, while $T_m(x) - 1$ is the product of a perfect square by a degree two polynomial.*

*If instead $m = 2l + 1$ is odd, then:*

*(b1)* $T_m(x) + 1 = (x + 1)V_l(x)^2$,

*(b2)* $T_m(x) - 1 = (x - 1)W_l(x)^2$;

*in particular, both $T_m(x) + 1$ and $T_m(x) - 1$ are the product of a perfect square by a linear polynomial.*

*Proof.* Set $x = \cos\theta$. If $m = 2l$ is even, then $T_m(x) + 1 = \cos(2l\theta) + 1 = 2\cos^2(l\theta) = 2(T_l(x))^2$. Similarly,

$$T_m(x) - 1 = \cos(2l\theta) - 1 = -2\sin^2(l\theta) = -2\sin^2\theta \cdot \frac{\sin^2(l\theta)}{\sin^2\theta} = 2(x^2 - 1)(U_l(x))^2.$$

If instead $m = 2l + 1$ is odd, one has

$$T_m(x) + 1 = \cos(m\theta) + 1 = 2\cos^2(m\theta/2) = 2\cos^2(\theta/2) \cdot \frac{\cos^2(m\theta/2)}{\cos^2(\theta/2)} = (x+1)(V_l(x))^2,$$

and

$$T_m(x) - 1 = \cos(m\theta) - 1 = -2\sin^2(m\theta/2) = -2\sin^2(\theta/2) \cdot \frac{\sin^2(m\theta/2)}{\sin^2(\theta/2)} = (x-1)(W_l(x))^2.$$

□

# Lecture V
# Products of deformation spaces of tacnodes

## by Francesco Bastianelli

## 1 − Introduction

This lecture retraces [1, Section 4.2], and concerns Caporaso and Harris' description of products of deformation spaces of higher-order tacnodes. Of course, this description relies on the analysis included in [IV] about deformations of a single tacnode. The main result of this lecture shall be involved in the proof of [VI, **??**], which is [1, Theorem 1.3], in order to describe the local geometry of hyperplane sections of generalized Severi varieties $V^{d,\delta}(\alpha,\beta)$, around a point $[X_0]$ parameterizing a reducible curve $X_0 = X \cup L$ endowed with some tacnodes of orders $m_j$, each deforming to $m_j - 1$ nodes on the curves parameterized nearby.

We consider a sequence $\underline{m} := (m_1, m_2, \ldots, m_n)$ of integers $m_j \geqslant 2$, and we set

$$\lambda := \operatorname{lcm}\{m_j\}, \qquad \mu := \prod_j m_j, \qquad \kappa := \frac{\mu}{\lambda},$$

where lcm denotes the least common multiple.

For any $1 \leqslant j \leqslant n$, let $(C_j, p_j)$ be a tacnode of order $m_j$, and let $\pi_j \colon \mathcal{S}_j \longrightarrow \Delta_j$ be the corresponding versal deformation family described in [IV]. Therefore, $\Delta_j \cong \mathbb{A}^{2m_j - 1}$ is an affine space with coordinates $(\underline{a}_j, \underline{b}_j) = (a_{j,m_j-2}, \ldots, a_{j,0}, b_{j,m_j-1}, \ldots, b_{j,0})$, the subscheme $\mathcal{S}_j \subseteq \mathbb{A}^2 \times \Delta_j$ has equation

$$y^2 + \left(x^{m_j} + a_{j,m_j-2}x^{m_j-2} \cdots + a_{j,0}\right) y + b_{j,m_j-1}x^{m_j-1} + \cdots + b_{j,0} = 0,$$

and $\pi_j \colon \mathcal{S}_j \longrightarrow \Delta_j$ is the second projection. Moreover, we introduce the subloci $\Delta_{j,m_j}$ and $\Delta_{j,m_j-1}$ of $\Delta_j$ as

$$\Delta_{j,m_j} := \overline{\left\{(\underline{a}_j, \underline{b}_j) \in \Delta_j \,\middle|\, \pi_j^{-1}(\underline{a}_j, \underline{b}_j) \text{ is a curve having } m_j \text{ nodes}\right\}}$$

and

$$\Delta_{j,m_j-1} := \overline{\left\{(\underline{a}_j, \underline{b}_j) \in \Delta_j \,\middle|\, \pi_j^{-1}(\underline{a}_j, \underline{b}_j) \text{ is a curve having } m_j - 1 \text{ nodes}\right\}}.$$

We recall that $\Delta_{j,m_j} \cong \mathbb{A}^{m_j-1}$ is the linear $(m_j-1)$-dimensional subspace given by the equations $b_{j,m_j-1} = \cdots = b_{j,0} = 0$, and it is the locus over which the fibers $\pi_j^{-1}(\underline{a}_j, \underline{b}_j)$ are reducible plane curves, whereas $\Delta_{j,m_j-1}$ is an $m_j$-dimensional subvariety containing $\Delta_{m_j}$.

Then we define

$$\Delta := \Delta_1 \times \cdots \times \Delta_n,$$

$$\Delta_{\underline{m}} := \Delta_{1,m_1} \times \cdots \times \Delta_{n,m_n},$$

and

$$\Delta_{\underline{m}-1} := \Delta_{1,m_1-1} \times \cdots \times \Delta_{n,m_n-1},$$

so that $\Delta \cong \mathbb{A}^{\Sigma(2m_j-1)}$ is an affine space endowed with coordinates $(\underline{a}_1, \underline{b}_1, \ldots, \underline{a}_n, \underline{b}_n)$, the linear subspace $\Delta_{\underline{m}} \cong \mathbb{A}^{\Sigma(m_j-1)}$ is given by equations $\{b_{j,i} = 0,\ 1 \leqslant j \leqslant n,\ 0 \leqslant i \leqslant m_j - 1\}$, and $\Delta_{\underline{m}-1}$ is a subvariety of dimension $\sum m_j$ containing $\Delta_{\underline{m}}$. Finally, we set

$$H := \bigcup_j \{b_{j,0} = 0\},$$

which is a union of hyperplanes of $\Delta$.

We now state the main result of this lecture (cf. [1, Lemma 4.3]). Along the same lines as [IV, (1.2)] (cf. [1, Lemma 4.1]), it describes the local intersection at the origin $\underline{0} \in \Delta$ between $\Delta_{\underline{m}-1}$ and varieties $W \subseteq \Delta$ containing $\Delta_{\underline{m}}$ as a subvariety of codimension 1 and satisfying a suitable transversality assumption.

**(1.1) Lemma.** *Let $W \subseteq \Delta$ be a smooth subvariety of dimension $\sum(m_j - 1) + 1$ containing $\Delta_{\underline{m}}$, and such that the tangent space $T_{\underline{0}}W$ at the origin $\underline{0} \in \Delta$ is not contained in $H$. Then, in an étale neighborhood of the origin,*

$$W \cap \Delta_{\underline{m}-1} = \Delta_{\underline{m}} \cup \Gamma_1 \cup \Gamma_2 \cup \cdots \cup \Gamma_\kappa,$$

*where $\Gamma_1, \ldots, \Gamma_\kappa \subseteq W$ are distinct reduced unibranched curves such that each $\Gamma_\alpha$ has intersection multiplicity $\big(\Gamma_\alpha \cdot \Delta_{\underline{m}}\big)_{\underline{0}} = \lambda$ with $\Delta_{\underline{m}}$ at the origin, and the origin is a point of multiplicity $\mathrm{mult}_{\underline{0}}(\Gamma_\alpha) = \frac{\lambda}{\max_j\{m_j\}}$ of $\Gamma_\alpha$.*

Note that $\Delta_{\underline{m}-1}$ has codimension $\sum(m_j - 1)$ in $\Delta$, so the expected dimension of the intersection $W \cap \Delta_{\underline{m}-1}$ is 1.

**(1.2) Remark.** (i) The intersection multiplicity at $\underline{0} \in \Delta$ between $\Delta_{\underline{m}}$ and the reducible curve $\Gamma := \Gamma_1 \cup \cdots \cup \Gamma_\kappa$ is $\big(\Gamma \cdot \Delta_{\underline{m}}\big)_{\underline{0}} = \kappa\lambda = \mu$, in analogy with the case of a single tacnode.

(ii) If $\lambda = m_j$ for some $j$, then all the curves $\Gamma_\alpha$ are smooth.

Before proving Lemma (1.1), we state another result which shall turn out to be equivalent to it. We denote by

$$\widetilde{\Delta} := \mathrm{Bl}_{\Delta_{\underline{m}}} \Delta \xrightarrow{\pi} \Delta$$

the blow-up of $\Delta$ along $\Delta_{\underline{m}}$. Let $E \subseteq \widetilde{\Delta}$ be the exceptional divisor, and let $F := \pi^{-1}(\underline{0}) \subseteq E$ be the fibre over the origin $\underline{0} \in \Delta$. Finally, let $\widetilde{\Delta}_{\underline{m}-1}$ and $\widetilde{H}$ denote the proper transforms of $\Delta_{\underline{m}-1}$ and $H$, respectively. Then the following holds (cf. [1, Lemma 4.4]).

**(1.3) Lemma.** *The intersection $\widetilde{\Delta}_{\underline{m}-1} \cap E$ contains $F$ as a component of multiplicity $\mu$. Moreover, in an étale neighborhood of any point $p \in F$ not contained in $\widetilde{H}$, the variety $\widetilde{\Delta}_{\underline{m}-1}$ consists of $\kappa$ reduced branches, each having multiplicity $\frac{\lambda}{\max_j\{m_j\}}$ along $F$, intersection number $\lambda$ with $E$ along $F$, and tangent cone at $p$ supported on a linear space contained in $E$.*

## 2 – Proofs

In analogy with the argument used to deduce [1, Lemma 4.1], the proof of Lemma (1.1) proceeds in three steps. The first one consists of proving the assertion when $W$ is assumed to be a linear space. Then we use it to achieve Lemma (1.3). Finally, we conclude by deducing Lemma (1.1) for an arbitrary subvariety $W$ satisfying the hypothesis of the statement.

**(2.1) Proof of Lemma (1.1) with linear $W \subseteq \Delta$.** We assume that $W \cong \mathbb{A}^{\Sigma(m_j-1)+1}$ is a linear subspace of $\Delta$ such that $\Delta_{\underline{m}} \subseteq W$ and $W = T_{\underline{0}}W \not\subseteq H$. Thus $\Delta_{\underline{m}} \subseteq H$ is a hyperplane in $W$, and we may choose a point $v \in \Delta \smallsetminus H$ such that $W = \langle \Delta_{\underline{m}}, v \rangle$, i.e., $W$ is the linear span of $\Delta_{\underline{m}}$ and $v$.

For any $1 \leqslant j \leqslant n$, let $\rho_j \colon \Delta \longrightarrow \Delta_j$ be the $j^{\text{th}}$ projection, and let $t \colon W \longrightarrow \mathbb{C}$ be a non-zero linear function vanishing along the hyperplane $\Delta_{\underline{m}}$. We note that $\rho_j(\Delta_{\underline{m}-1}) = \Delta_{m_j-1}$ and $\rho_j(\Delta_{\underline{m}}) = \Delta_{m_j}$. Moreover, setting $W_j := \rho_j(W)$ and $v_j := \rho_j(v)$, we have that $v_j \in W_j \smallsetminus \Delta_{m_j}$ as $v \notin H$. Thus $W_j = \langle \Delta_{m_j}, v_j \rangle \subseteq \Delta_j$ is an $m_j$-dimensional plane not contained in the hyperplane $\{b_{j,0} = 0\}$. In particular $W_j$ satisfies the hypothesis of [IV, (1.2)], so the intersection between $\Delta_{m_j-1}$ and $W_j$ in $\Delta_j$ consists of the union of $\Delta_{m_j}$ and a smooth curve $\Gamma_j$ having intersection multiplicity $(\Gamma_j \cdot \Delta_{m_j})_{\underline{0}} = m_j$ with $\Delta_{m_j}$ at the origin $\underline{0} \in \Delta_j$. Therefore, in some étale neighborhood of the origin $\underline{0} \in W_j$, we may choose local coordinates $(x_{j,0}, \ldots, x_{j,m_j-2}, t_j)$ such that $\rho_j^*(t_j) = t$,[1] $\Delta_{m_j}$ is the hyperplane with equation $t_j = 0$, and the curve $\Gamma_j$ is defined by the equations
$$x_{j,1} = \cdots = x_{j,m_j-2} = 0 \quad \text{and} \quad t_j = (x_{j,0})^{m_j}.$$

Then the functions $t = \rho_j^*(t_j)$ and $y_{j,i} := \rho_j^*(x_{j,i})$, with $1 \leqslant j \leqslant n$ and $1 \leqslant i \leqslant m_j - 2$, define a system of local coordinates on some étale neighborhood of $W$ at the origin. It follows from their very definition that the intersection $W \cap \Delta_{\underline{m}-1}$ consists of the plane $\Delta_{\underline{m}} = \{t = 0\}$ and of the subvariety $\Gamma$ defined by the equations

(2.1.1) $\qquad t = (y_{j,0})^{m_j} \quad \text{and} \quad y_{j,i} = 0, \qquad \text{for } 1 \leqslant j \leqslant n \text{ and } 1 \leqslant i \leqslant m_j - 2.$

We point out that $\Gamma$ is a curve, as for any fixed value of $t \in \mathbb{C}$, there exist finitely many $((\sum m_j - 1) + 1)$-tuples $(y_{1,0}, \ldots, y_{n,m_n-2}, t)$ satisfying (2.1.1), which depend only on the choice of an $m_j^{\text{th}}$ root of $t$ for each coordinate $y_{j,0}$.

In order to parameterize any branch of $\Gamma$, let us consider a sequence $\underline{\zeta} := (\zeta_1, \zeta_2, \ldots, \zeta_n)$ such that each $\zeta_j$ is a $m_j^{\text{th}}$ root of unity. Then we define a local parameterization $\varphi_{\underline{\zeta}} \colon \mathbb{C} \longrightarrow W$ as

(2.1.2) $\qquad z \longmapsto \begin{cases} t &= z^\lambda \\ y_{j,0} &= \dfrac{z^{\frac{\lambda}{m_j}}}{\zeta_j} & \text{for } 1 \leqslant j \leqslant n \\ y_{j,i} &= 0 & \text{for } 1 \leqslant j \leqslant n \text{ and } 1 \leqslant i \leqslant m_j - 2, \end{cases}$

and we denote by $\Gamma_{\underline{\zeta}}$ its image. Notice that the map $\varphi_{\underline{\zeta}}$ is injective.[2] Furthermore, every branch of $\Gamma$ can be parameterized in this way: for any point $\underline{y} \in \Gamma$ with coordinates $(y_{1,0}, \ldots, y_{n,m_n-2}, t)$, it is enough to choose a $\lambda^{\text{th}}$ root $z$ of $t$ and set $\zeta_j = z^{\frac{\lambda}{m_j}}/y_{j,0}$, which is an $m_j^{\text{th}}$ root of unity by (2.1.1). We also point out that the multiplicity of $\Gamma_{\underline{\zeta}}$ at the origin is the least power of $z$ appearing in the parameterization, i.e., $\text{mult}_{\underline{0}}(\Gamma_{\underline{\zeta}}) = \dfrac{\lambda}{\max_j\{m_j\}}$. On the other hand, the intersection multiplicity at the origin between $\Gamma_{\underline{\zeta}}$ and the hyperplane $\Delta_{\underline{m}} = \{t = 0\}$ is $(\Gamma_{\underline{\zeta}} \cdot \Delta_{\underline{m}})_{\underline{0}} = \lambda$, as $t = z^\lambda$ on $\Gamma_{\underline{\zeta}}$. Therefore the behavior of the branches of $\Gamma$ satisfies the statement, and it remains to enumerate them.

To this aim, we note that two distinct sequences $\underline{\zeta} := (\zeta_1, \zeta_2, \ldots, \zeta_n)$ and $\underline{\eta} := (\eta_1, \eta_2, \ldots, \eta_n)$ parameterize the same branch if and only if $\zeta_j = \varepsilon^{\frac{\lambda}{m_j}} \eta_j$ for all $j$, where $\varepsilon$ is some $\lambda^{\text{th}}$ root of

---

[1]It suffices to choose $t_j \colon W_j \longrightarrow \mathbb{C}$ to be a linear function such that $t_j(\Delta_{m_j}) = 0$ and $t_j(v_j) = t(v)$.

[2]Suppose that $\varphi_{\underline{\zeta}}(z) = \varphi_{\underline{\zeta}}(w)$ for some $z, w \in \mathbb{C}$. Then $(z/w)^{\frac{\lambda}{m_j}} = 1$ for any $j$, hence $z/w$ is a $\nu^{\text{th}}$ root of unity, where $\nu$ divides any $\frac{\lambda}{m_j}$. Thus $\nu = 1$ as $\lambda$ is the least common multiple of the integers $m_j$, hence $z = w$.

unity.[3]

Moreover, $\zeta_j = \varepsilon^{\frac{\lambda}{m_j}} \zeta_j$ for any $j$ if and only if $\varepsilon = 1$ [4]. We conclude that the number of branches of $\Gamma$ equals the number of sequences $\underline{\zeta} := (\zeta_1, \zeta_2, \ldots, \zeta_n)$ divided by the number of $\lambda^{\text{th}}$ roots of unity, that is $\kappa = \frac{\mu}{\lambda}$. Thus we achieved Lemma (1.1) when $W$ is assumed to be a linear subspace of $\Delta$. □

We now turn to prove Lemma (1.3).

**(2.2) Proof of Lemma (1.3).** Let $G$ be the Grassmannian of vector subspaces of $\Delta \cong \mathbb{A}^{\Sigma(2m_j-1)}$ of dimension $\sum(m_j - 1) + 1$, and let $B \cong \mathbb{P}^{\Sigma m_j - 1}$ be the subvariety of $G$ parameterizing those subspaces containing $\Delta_{\underline{m}} \cong \mathbb{A}^{\Sigma(m_j-1)}$. Let $W_{\underline{b}}$ denote the subspace parameterized by $\underline{b} \in B$, and consider the incidence variety $\Lambda := \left\{ (\underline{y}, \underline{b}) \in \Delta \times B \,\middle|\, \underline{y} \in W_{\underline{b}} \right\} \subseteq \Delta \times B$. Then $\Lambda$ is the restriction to $B$ of the universal family over $G$, and it fits in the following diagram

(2.2.1)
$$\begin{array}{ccc} & \Lambda & \\ \text{pr}_1 \swarrow & & \searrow \text{pr}_2 \\ \Delta & \dashrightarrow^{\psi} & B, \end{array}$$

where $\text{pr}_1$ and $\text{pr}_2$ are the natural projections, while $\psi$ is the map sending any $\underline{y} \in \Delta \smallsetminus \Delta_{\underline{m}}$ to the point parameterizing the linear span $\langle \underline{y}, \Delta_{\underline{m}} \rangle \subseteq \Delta$ [5]. We point out that $\overline{\text{pr}}_1 \colon \Lambda \longrightarrow \Delta$ is an isomorphism outside $(\text{pr}_1)^{-1}(\Delta_{\underline{m}})$, whereas the fibre $(\text{pr}_1)^{-1}(\underline{y})$ over any point $\underline{y} \in \Delta_{\underline{m}}$ is a section of $\text{pr}_2$, that is a projective space of dimension $\sum m_j - 1$. Therefore $\text{pr}_1 \colon \Lambda \longrightarrow \Delta$ coincides with the blow-up $\widetilde{\Delta} := \text{Bl}_{\Delta_{\underline{m}}} \Delta \xrightarrow{\pi} \Delta$ of $\Delta$ along $\Delta_{\underline{m}}$ (see, e.g., [2, p.604]). [6].

Consider the proper transforms $\widetilde{\Delta}_{\underline{m}-1}, \widetilde{H} \subseteq \widetilde{\Delta}$ of $\Delta_{\underline{m}-1}$ and $H$. Let $\widetilde{\Delta}^0_{\underline{m}-1} \subseteq \widetilde{\Delta}_{\underline{m}-1}$ and $\Delta^0_{\underline{m}-1} \subseteq \Delta_{\underline{m}-1}$ denote the open subsets of points not lying on $\widetilde{H}$ and $H$, respectively, and let $B_0 \subseteq B$ be the open subvariety parameterizing subspaces not contained in $H$. Then the diagram (2.2.1) restricts to

(2.2.2)
$$\begin{array}{ccc} & \widetilde{\Delta}^0_{\underline{m}-1} & \\ \pi_0 \swarrow & & \searrow \tau_0 \\ \Delta^0_{\underline{m}-1} & \xrightarrow{\psi_0} & B_0. \end{array}$$

We note that the closure $\overline{(\psi_0)^{-1}(\underline{b})} \subseteq \Delta_{\underline{m}-1}$ of the fibre over a point $\underline{b} \in B_0$ is the residual intersection of $\Delta_{\underline{m}-1}$ and the subspace $W_{\underline{b}}$ when we remove $\Delta_{\underline{m}}$. By Lemma (1.1) applied to $W_{\underline{b}} \cong \mathbb{A}^{\Sigma(m_j-1)+1}$, in some étale neighborhood of the origin $\underline{0} \in \Delta$, the latter intersection consists of $\kappa$ distinct reduced unibranched curves $\Gamma_1, \ldots, \Gamma_\kappa$, each satisfying $\left(\Gamma_\alpha \cdot \Delta_{\underline{m}}\right)_{\underline{0}} = \lambda$ and $\text{mult}_{\underline{0}}(\Gamma_\alpha) = \frac{\lambda}{\max_j\{m_j\}}$. Notice that the proper transform $\widetilde{W}_b = (\text{pr}_2)^{-1}(\underline{b})$ of $W_b$ meets the fibre $F := \pi^{-1}(\underline{0})$ at a single point $p \in \widetilde{\Delta} \smallsetminus \widetilde{H}$ (recall that $\pi : \widetilde{\Delta} \to \Delta$ is the blow-up is the blow-up along $\Delta_{\underline{m}}$, which identifies with $\text{pr}_1 : \Lambda \to \Delta$). It follows that in some étale neighborhood of $p$, the closure of the fibre $(\tau_0)^{-1}(\underline{b}) \subseteq \widetilde{\Delta}_{\underline{m}-1}$ consists of the proper transforms $\widetilde{\Gamma}_1, \ldots, \widetilde{\Gamma}_\kappa$ of the curves $\Gamma_\alpha \subseteq \Delta_{\underline{m}-1}$. Moreover, any $\widetilde{\Gamma}_\alpha$ is a reduced unibranched curve such that

(2.2.3) $$\left(\widetilde{\Gamma}_\alpha \cdot E\right)_p = \left(\Gamma_\alpha \cdot \Delta_{\underline{m}}\right)_{\underline{0}} = \lambda$$

and

(2.2.4) $$\left(\widetilde{\Gamma}_\alpha \cdot F\right)_p = \text{mult}_{\underline{0}}(\Gamma_\alpha) = \frac{\lambda}{\max_j\{m_j\}},$$

by the Projection Formula and [3, Lemma 1.40 p.28].

---

[3]Indeed, $\varphi_{\underline{\zeta}}(z) = \varphi_{\underline{\eta}}(w)$ for some $z, w \in \mathbb{C}$ if and only if $(z/w)^\lambda = 1$ and $\eta_j z^{\frac{\lambda}{m_j}} = \zeta_j w^{\frac{\lambda}{m_j}}$, that is $z = \varepsilon w$ for some $\lambda^{\text{th}}$ root $\varepsilon$ of unity and $\zeta_j = \varepsilon^{\frac{\lambda}{m_j}} \eta_j$.

[4]Notice that $\varepsilon^{\frac{\lambda}{m_j}} = 1$ for any $j$ if and only if $\varepsilon$ is a $\nu^{\text{th}}$ root of unity, where $\nu$ divides any $\frac{\lambda}{m_j}$. Hence $\nu = 1$.

[5]We note that $B$ is naturally identified with the projectivization $\mathbb{P}\left(\Delta/\Delta_{\underline{m}}\right) \cong \mathbb{P}^{\Sigma m_j - 1}$ of the quotient $\Delta/\Delta_{\underline{m}}$. Under this identification, we may view $\psi$ as the map sending a point $\underline{y} = \left(\underline{a}_1, \underline{b}_1, \ldots \underline{a}_n, \underline{b}_n\right) \in \Delta \smallsetminus \Delta_{\underline{m}}$ to the point $\underline{b} \in \mathbb{P}\left(\Delta/\Delta_{\underline{m}}\right)$ having homogeneous coordinates $\underline{b} = \left[\underline{b}_1, \ldots, \underline{b}_n\right]$.

[6]In other terms, the dominant rational map $\psi \colon \Delta \dashrightarrow \mathbb{P}\left(\Delta/\Delta_{\underline{m}}\right)$ can be resolved to a morphism by blowing up $\Delta$ along the indeterminacy locus $\Delta_{\underline{m}}$.

Finally, as we vary $\underline{b} \in B_0$, each $\widetilde{\Gamma}_\alpha$ describes a reduced branch $\widetilde{\Delta}^\alpha$ of $\widetilde{\Delta}_{\underline{m}-1}$ having multiplicity $\frac{\lambda}{\max_j\{m_j\}}$ along $F$, intersection number $\lambda$ with $E$ along $F$, and since the curves $\widetilde{\Gamma}_\alpha$ are unibranched, the tangent cone to $\widetilde{\Delta}^\alpha$ at $p$ is supported on a linear space contained in $E$. In particular, the intersection $\widetilde{\Delta}_{\underline{m}-1} \cup E$ contains the whole fibre $F$, with multiplicity given by the sum of the intersection multiplicities along $F$ between $E$ and each branch of $\widetilde{\Delta}_{\underline{m}-1}$, that is $\kappa\lambda = \mu$. □

Finally, we conclude this lecture by proving Lemma (1.1) for an arbitrary subvariety $W \subseteq \Delta$ satisfying the hypothesis of the statement.

**(2.3) Proof of Lemma (1.1).** Let $W \subseteq \Delta$ be a smooth $(\sum(m_j-1)+1)$-dimensional variety containing $\Delta_{\underline{m}}$ and such that $T_{\underline{0}}W \not\subseteq H$. By the smoothness of $W$, its proper transform $\widetilde{W} \subseteq \widetilde{\Delta}$ intersects the exceptional divisor $E$ transversally, and cuts out on the latter a section over $\widetilde{\Delta}_{\underline{m}}$. Moreover, since $T_{\underline{0}}W \not\subseteq H$, we have that $\widetilde{W}$ meets $F := \pi^{-1}(\underline{0})$ at a single point $p \in \widetilde{\Delta}$ which lies off $\widetilde{H}$, hence $T_p\widetilde{W} \cap T_pF = \{0\}$. By Lemma (1.3), in some étale neighborhood of $p$, the variety $\widetilde{\Delta}_{\underline{m}-1}$ consists of $\kappa$ reduced sheets $\widetilde{\Delta}^1, \ldots, \widetilde{\Delta}^n$, and the tangent cone to a branch $\widetilde{\Delta}^\alpha$ at $p$ is supported on a linear space $L \subseteq E$. In particular, $T_pF \subseteq L \subseteq E$. We notice that $\dim L = \dim \widetilde{\Delta}^\alpha = \sum m_j$, and $\dim T_pF = \dim F = \sum m_j - 1 = \dim \widetilde{\Delta}^\alpha - 1$. Thus $\dim(T_p\widetilde{W} \cap L) \leqslant 1$, so that the linear span of $T_p\widetilde{W}$ and $L$ satisfies

$$\dim\langle T_p\widetilde{W}, L\rangle \geqslant \dim T_p\widetilde{W} + \dim L - 1 = \sum(m_j-1) + 1 + \sum m_j - 1$$
$$= \sum(2m_j - 1) = \dim \widetilde{\Delta},$$

that is $\langle T_p\widetilde{W}, L\rangle = T_p\widetilde{\Delta}$. Therefore, the intersection around $p$ between $\widetilde{W}$ and $\widetilde{\Delta}^\alpha$ is transverse, and it consists of a unibranched reduced curve $\widetilde{\Gamma}_\alpha$ such that $\left(\widetilde{\Gamma}_\alpha \cdot F\right)_p$ equals the multiplicity $\frac{\lambda}{\max_j\{m_j\}}$ of $\widetilde{\Delta}^\alpha$ along $F$, and $\left(\widetilde{\Gamma}_\alpha \cdot E\right)_p$ is the intersection number $\lambda$ between $\widetilde{\Delta}^\alpha$ and $E$ along $F$.

For any $1 \leqslant \alpha \leqslant \kappa$, we consider the reduced curve $\Gamma_\alpha := \pi(\widetilde{\Gamma}_\alpha)$, and in some étale neighborhood of the origin, we have

$$\Delta_{\underline{m}-1} \cap W = \Delta_{\underline{m}} \cup \Gamma_1 \cup \cdots \cup \Gamma_\kappa.$$

Since each $\widetilde{\Gamma}_\alpha$ meets $E$ only at the point $p$, the curve $\Gamma_\alpha$ is unibranched as well. Furthermore, by arguing as in (2.2.3) and (2.2.4), we deduce that each $\Gamma_\alpha$ satisfies $\text{mult}_{\underline{0}}(\Gamma_\alpha) = \left(\widetilde{\Gamma}_\alpha \cdot F\right)_p = \frac{\lambda}{\max_j\{m_j\}}$ and $\left(\Gamma_\alpha \cdot \Delta_{\underline{m}}\right)_{\underline{0}} = \left(\widetilde{\Gamma}_\alpha \cdot E\right)_p = \lambda$. Thus Lemma (1.1) follows. □

# Appendix A
# Geometry of the deformation space of an $m$-tacnode into $m-1$ nodes

## by Thomas Dedieu

## 1 – Introduction

This appendix complements the previous lecture [IV], of which we shall keep the notation: let $\pi : S \to \Delta$ be the versal deformation space of the $m$-tacnode defined by the equation $y(y + x^m) = 0$ in the affine plane; we have $\Delta \cong \mathbf{A}^{2m-1}_{\alpha_0,\ldots,\alpha_{m-2},\beta_0,\ldots,\beta_{m-1}}$, and $S$ is defined in $\mathbf{A}^2_{x,y} \times \Delta$ by the equation

$$y^2 + y(x^m + \alpha_{m-2}x^{m-2} + \cdots + \alpha_0) + \beta_{m-1}x^{m-1} + \cdots + \beta_0.$$

Inside $\Delta$ we consider the two subvarieties $\Delta_{m-1}$ and $\Delta_m$, the closures of the two loci of those $t \in \Delta$ such that the curves $\pi^{-1}(t)$ have $m-1$ nodes and $m$ nodes, respectively.

The goal of [IV] was to establish some properties of $\Delta_{m-1}$ around the origin in $\Delta$, of fundamental importance in order to carry out the program of [1], see [VI]. Here we shall provide complementary results, and work out detailed examples for $m = 2$ and 3. The first motivation is that this will hopefully give more insight into the main result of [IV], as well as help to follow the arguments given there. On the other hand, the variety $\Delta_{m-1}$ for $m = 2$ will show up under another guise in this volume, see [??], namely as the local model for the dual variety $X^\vee \subseteq \check{\mathbf{P}}^3$ of a surface $X \subseteq \mathbf{P}^3$, around a point corresponding to a tacnodal hyperplane section; this makes $\Delta_{m-1}$ one of the main characters of this volume. We are therefore interested in more properties of $\Delta_{m-1}$ than those strictly necessary for the Caporaso–Harris Program.

Finally, I would like to point out before I start that the singularity of $\Delta_{m-1}$ at the origin is known for $m = 2$ as a swallowtail singularity, *queue d'aronde* in french, one of seven elementary situations in René Thom's *théorie des catastrophes*.

## 2 – General results

Let $\tau : \tilde{\Delta} \to \Delta$ be the blow-up along $\Delta_m$; recall that $\Delta_m \subseteq \Delta$ is the affine plane with equations

$$\beta_0 = \cdots = \beta_{m-1} = 0,$$

and it is contained in $\Delta_{m-1}$. Let $Z \subseteq \tilde{\Delta}$ be the exceptional divisor, a $\mathbf{P}^{m-1}$-bundle over $\Delta_m$, and let $\Phi$ be the fibre of $Z$ over the origin. Let $\tilde{\Delta}_{m-1}$ be the proper transform of $\Delta_{m-1}$. Let $H$ be the hyperplane $\beta_0 = 0$ in $\Delta$, and let $\tilde{H}$ be its proper transform in $\tilde{\Delta}$.

The main result of [IV] may be reformulated as follows, cf. [1, Lemma 4.2], where the equivalence of the two statements is proved (all the ingredients for this proof are present in [IV]).

**(2.1) Proposition.** *Let $Z_{m-1}$ be the intersection $\tilde{\Delta}_{m-1} \cap Z$ (this is the exceptional locus of $\tau : \tilde{\Delta}_{m-1} \to \Delta_{m-1}$). Then $Z_{m-1}$ contains $\Phi$ as a component of multiplicity $m$. Moreover, $\tilde{\Delta}_{m-1}$ is smooth at all points of $\Phi$ not in $\tilde{H}$.*

We shall now prove the result below, following [2, §2.4.4]. It is interesting in itself, but also as a pretext to give a local description of $\Delta_{m-1}$ at an arbitrary point of $\Delta_m$.

**(2.2) Proposition.** *The fibres of $Z_{m-1}$ over $\Delta_m$ are unions of linear spaces.*

Before we prove this, we introduce a natural stratification of $\Delta_m$. Each stratum is the locus in which the $m$-tacnode deforms into $k$ non-infinitely near tacnodes of orders $m_1, \ldots, m_k$ respectively, with $m_1 + \cdots + m_k = m$.

The space $\Delta_m$ identifies with the space of monic degree $m$ polynomials in $x$ with no $x^{m-1}$-term, and we shall consider in this space the loci of polynomials having roots with given multiplicities. For all partition $m = m_1 + \cdots + m_k$, we let

$$\Delta\{m_1, \ldots, m_k\} \subseteq \Delta_m \subseteq \Delta$$

be the locus corresponding to monic degree $m$ polynomials with no $x^{m-1}$-term of the form $(x - \lambda_1)^{m_1} \cdots (x - \lambda_k)^{m_k}$, with $\lambda_1, \ldots, \lambda_k$ pairwise distinct. Note that $\Delta\{m_1, \ldots, m_k\}$ has dimension $k - 1$ (because we impose that there is no $x^{m-1}$-term, equivalently the $\lambda_i$'s sum up to zero), so its codimension in $\Delta_m$ is $m - k = \sum_{i=1}^k (m_i - 1)$.

For $\alpha \in \Delta\{m_1, \ldots, m_k\}$, the curve $S_\alpha = \pi^{-1}(\alpha) \subseteq \mathbf{A}^2_{x,y}$ is a reducible curve with two branches, defined by the two equations $y = 0$ and $y = (x - \lambda_1)^{m_1} \cdots (x - \lambda_k)^{m_k}$ respectively; these two branches intersect at the $k$ points $r_1 = (\lambda_1, 0), \ldots, r_k = (\lambda_k, 0)$, with multiplicities $m_1, \ldots, m_k$ respectively.

For all $i = 1, \ldots, k$ we consider the versal deformation space $\Delta(S_\alpha, r_i)$ of the singular point $r_i \in S_\alpha$; note that this is an $m_i$-tacnode. There is then a natural map

$$\sigma : U \to \prod_{i=1}^k \Delta(S_\alpha, r_i),$$

where $U$ is an open neighbourhood $U$ of $\alpha$ in $\Delta$, and by the openness of versality it has surjective differential at $\alpha$; the fibre of $\sigma$ over the origin parametrizes the equisingular deformations of $S_\alpha$, for which only the location of the points $r_i$ varies along the $x$-axis.

For all $i = 1, \ldots, k$ we consider the subvarieties $\Delta_{m_i-1}(r_i)$ and $\Delta_{m_i}(r_i)$ of $\Delta(S_\alpha, r_i)$, defined as the closures of the loci where $r_i$ deforms to $m_i - 1$ nodes and $m_i$ nodes respectively. Then, in the neighbourhood $U$ of $\alpha$ in $\Delta$, we have

$$\Delta_m = \sigma^{-1}\bigl(\Delta_{m_1}(r_1) \times \cdots \Delta_{m_k}(r_k)\bigr)$$

and

$$\Delta_{m-1} = \bigcup_{l=1}^k \sigma^{-1}\bigl(\Delta_{m_1}(r_1) \times \cdots \times \Delta_{m_l-1}(r_l) \times \cdots \Delta_{m_k}(r_k)\bigr).$$

From this we see that $\Delta_{m-1}$ has $k$ local sheets in a neighbourhood of $\alpha$, each containing $\Delta_m$. The $l$-th sheet looks like the product of $\Delta_{m_l-1} \subseteq \Delta_{m_l}$ with some smooth factor, and along this sheet the fibres of $\pi : S \to \Delta$ have $m_i$ nodes tending to $r_i$ for $i \neq l$, and $m_l - 1$ nodes tending to $r_l$.

We are now ready to prove the statement.

*Proof of Proposition (2.2).* In the above setup, applying Proposition (2.1) to the spaces $\Delta(S_\alpha, r_i)$, we will be able to describe the fibre of $Z_{m-1}$ over $\alpha \in \Delta_m$. Note that $\Phi$, the fibre of $Z$ over the origin, identifies with the projectivization of the space of degree $m-1$ polynomials in $x$ (with coefficients $\beta_0, \ldots, \beta_{m-1}$), equivalently with the projectivization of the space of polynomials modulo those vanishing to order $m$ at the origin (i.e., $\mathbf{C}[x]/(x^m)$).

Then, for all $l = 1, \ldots, k$, by Proposition (2.1) applied to $\Delta_{m_l-1}(r_l)$, the proper transform of the $l$-th sheet of $\Delta_{m-1}$ intersects the fibre $\Phi_\alpha$ of $Z$ over $\alpha$ in the linear subspace of $\Phi_\alpha$ corresponding to polynomials in $x$ vanishing to order $m_h$ at $r_h$ for all $h \neq l$; in other words, if we identify $\Phi_\alpha$ with the projectivization of

$$\mathbf{C}[x]/\bigl((x-\lambda_1)^{m_1} \cdots (x-\lambda_k)^{m_k}\bigr) \cong \prod_{h=1}^{k} \mathbf{C}[x]/\bigl((x-\lambda_h)^{m_h}\bigr),$$

then the proper transform of the $l$-th sheet of $\Delta_{m-1}$ intersects $\Phi_\alpha$ in its linear subspace which is the projectivization of $\mathbf{C}[x]/\bigl((x-\lambda_l)^{m_l}\bigr)$ (with multiplicity $m_l$).

In turn, the intersection of $\tilde{\Delta}_{m-1}$ with the fibre $\Phi_\alpha$ of $Z$ over $\alpha$ is the union of these linear spaces for $l$ ranging from 1 to $k$. □

**(2.3) A complement.** Consider a small arc $t \in \mathbf{D} \mapsto \alpha(t) \in \Delta_m$ tending to the origin as $t$ tends to 0, where $\mathbf{D}$ stands for the unit disc in $\mathbf{C}$. If this arc is small enough, then all points $\alpha(t)$ with $t \neq 0$ lie in the same stratum $\Delta\{m_1, \ldots, m_k\}$. Moreover, as $t$ goes to $0 \in \mathbf{D}$ the singular points $r_i(t)$ of the curve $S_{\alpha(t)} \subseteq \mathbf{A}^2_{x,y}$ tend to the origin. It follows that the limiting position in $\Phi$ of the intersection with $\Phi_{\alpha(t)}$ of the $l$-th sheet of $\Delta_{m-1}$ along $\Delta\{m_1, \ldots, m_k\}$ is the space of polynomials vanishing to order $m - m_l$ at 0 (in the identification of $\Phi$ with polynomials in $x$ modulo those vanishing to order $m$ at the origin).

This implies, as in [IV, 2], that $\Phi$ is an irreducible component of $Z_{m-1}$, the general point of which corresponds to the tangent direction of an arc in $\Delta_{m-1}$ which intersects $\Delta_m$ only in the origin.

## 3 – Complete description in the case $m = 2$

When $m = 2$ we are considering the versal deformation space of the ordinary tacnode, defined in $\mathbf{A}^2_{x,y}$ by the equation $y(y+x^2) = 0$; this is $S \to \Delta \cong \mathbf{A}^3_{\alpha_0, \beta_0, \beta_1}$, with $S$ defined in $\mathbf{A}^2_{x,y} \times \Delta$ by the equation

$$y^2 + y(x^2 + \alpha_0) + \beta_1 x + \beta_0 = 0.$$

We are interested in the surface $\Delta_1 \subseteq \Delta$ of deformations of cogenus at least one of the tacnode, especially along the curve $\Delta_2 \subseteq \Delta$ of deformations of cogenus two, which is the line $\beta_1 = \beta_0 = 0$, along wich the tacnode deforms to two nodes. The picture looks as follows, as we shall see.

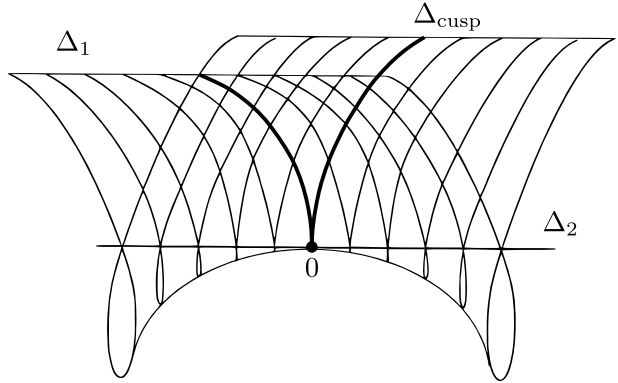

Figure 1: Geometry of $\Delta_1$ along $\Delta_2$ (beware the singularity along $\Delta_{\text{cusp}}$)
(reproduction from [2, p.181])

The surface $\Delta_1$ has two local sheets along the line $\Delta_2$; they intersect transversely along the complement of the origin in $\Delta_1$, but they have the same tangent plane at the origin. Also, it is important to keep in mind that the surface $\Delta_1$ is singular along the curve $\Delta_{\text{cusp}}$, which is the locus along which the tacnode deforms to an ordinary cusp: the surface $\Delta_1$ has a cuspidal singularity at the generic point of $\Delta_{\text{cusp}}$, although this not pictured on the above figure.

The proper transform of $\Delta_1$ under the blow-up of $\Delta$ along the line $\Delta_2$ looks as follows.

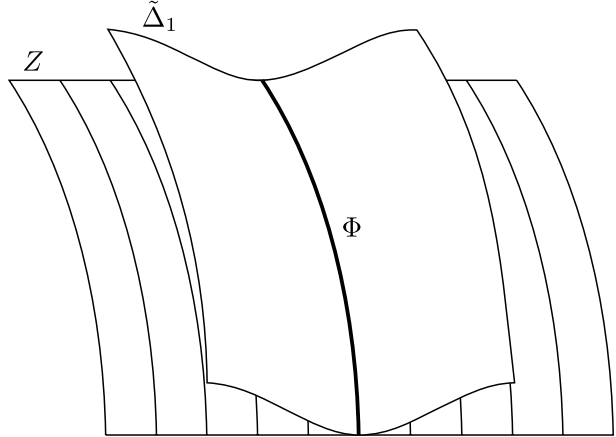

Figure 2: Intersection of $\tilde{\Delta}_1$ with the exceptional divisor
(reproduction from [1, p.379])

**(3.1) Equation of $\Delta_1$.** The surface $\Delta_1$ is defined in $\mathbf{A}^3_{\alpha_0,\beta_0,\beta_1}$ by the condition that the polynomial

$$\delta_{\alpha_0,\beta_0,\beta_1}(x) = (x^2 + \alpha_0)^2 - 4(\beta_1 x + \beta_0) = x^4 + 2\alpha_0 x^2 - 4\beta_1 x + \alpha_0^2 - 4\beta_0$$

has at least one non-simple root (see [IV, Section 1]), equivalently its discriminant vanishes. Using the formula for the discriminant of monic degree 4 polynomials with no $x^3$-term, namely

$$\mathrm{Disc}(x^4 + ax^2 + bx + c) = -4a^3b^2 + 16a^4c - 27b^4 + 144ab^2c - 128a^2c^2 + 256c^3,$$

we find the equation for $\Delta_1$ in $\mathbf{A}^3_{\alpha_0,\beta_0,\beta_1}$,

(3.1.1) $$256 \cdot (16\alpha_0^3\beta_1^2 + 16\alpha_0^2\beta_0^2 - 72\alpha_0\beta_0\beta_1^2 - 27\beta_1^4 - 64\beta_0^3) = 0.$$

Note that it is quasi-homogeneous of degree 12 for the weights $2, 4, 3$ of $\alpha_0, \beta_0, \beta_1$ respectively.

**(3.2) Equations of the curve $\Delta_{\mathrm{cusp}}$.** We consider the locus $\Delta_{\mathrm{cusp}}$ of those $t \in \Delta$ such that the curve $S_t$ has (at least) a cuspidal singularity. The locus $\Delta_{\mathrm{cusp}}$ is defined by the condition that the polynomial $\delta_{\alpha_0, \beta_0, \beta_1}$ has at least one triple root, in other words that there exists $u \in \mathbf{A}^1$ such that

$$\delta_{\alpha_0, \beta_0, \beta_1}(x) = (x - u)^3 \cdot (x + 3u)$$
$$= x^4 - 6u^2 x^2 + 8u^3 x - 3u^4$$

(recall that the roots of $\delta_{\alpha_0, \beta_0, \beta_1}$ sum up to zero as it has no $x^3$-term). Therefore, $\Delta_{\mathrm{cusp}}$ is the projection on the second factor of the complete intersection in $\mathbf{A}^1_u \times \mathbf{A}^3_{\alpha_0, \beta_0, \beta_1}$ defined by the equations

$$\begin{cases} 2\alpha_0 = -6u^2 \\ -4\beta_1 = 8u^3 \\ \alpha_0^2 - 4\beta_0 = -3u^4. \end{cases}$$

Eliminating $u$ from these equations (note that this can be done algorithmically), one finds the following equations defining $\Delta_{\mathrm{cusp}}$ in $\mathbf{A}^3_{\alpha_0, \beta_0, \beta_1}$,

$$\begin{cases} \alpha_0^2 - 3\beta_0 = 0 \\ 4\alpha_0^3 + 27\beta_1^2 = 0. \end{cases}$$

**(3.3) Singularities of the surface $\Delta_1$.** To compute the equations of the singularities of $\Delta_1$, we consider the Jacobian ideal $(\partial_{\alpha_0} F, \partial_{\beta_0} F, \partial_{\beta_1} F)$, where $F$ is the equation of $\Delta_1$ given in (3.1.1). It readily defines the singular subscheme of $\Delta_1$, since $F$ is quasi-homogeneous hence liable to the generalized Euler formula, cf. [II, Section 1]. Using a software for symbolic computation with polynomials, e.g., Macaulay2 [3] and its command `primaryDecomposition`, one finds the ideal defining the singular subscheme of $\Delta_1$ is the intersection of the two following primary ideals,

$$(\beta_0, \beta_1) \quad \text{and} \quad (9\alpha_0 \beta_0^2 - 4\alpha_0^2 \beta_1 + 24\beta_1^2,\ 8\alpha_0^3 - 27\beta_0^2 - 36\alpha_0 \beta_1,\ 27\beta_0^4 + 16\alpha_0^2 \beta_1^2 - 64\beta_1^3)$$

The former is the defining ideal of $\Delta_2$, while the latter is not prime, with radical the defining ideal of the curve $\Delta_{\mathrm{cusp}}$.

Also, one sees directly from the equation of $\Delta_1$, (3.1.1), that it has a triple point at the origin, with tangent cone the triple plane at the origin $\beta_0^3 = 0$.

**(3.4) Blow-up along the line $\Delta_2$.** The blow-up of $\mathbf{A}^3_{\alpha_0, \beta_0, \beta_1}$ along the line $\Delta_2$ is the hypersurface in $\mathbf{P}^1_{(u_0:u_1)} \times \mathbf{A}^3_{\alpha_0, \beta_0, \beta_1}$ defined by the equation $u_0 \beta_1 - u_1 \beta_0 = 0$. Pulling back the equation of $\Delta_1$, (3.1.1), we find

$$\beta_1^2 \cdot \left(16\alpha_0^2 u_0^2 - 64\beta_1 u_0^3 + 16\alpha_0^3 - 72\alpha_0 \beta_1 u_0 - 27\beta_1^2\right) = 0$$

in the affine chart $u_1 \neq 0$, and

$$\beta_0^2 \cdot \left(16\alpha_0^2 - 64\beta_0 + 16\alpha_0^3 u_1^2 - 72\alpha_0 \beta_0 u_1^2 - 27\beta_0^2 u_1^4\right) = 0$$

in the affine chart $u_0 \neq 0$; the factors between parentheses are thus the defining equations of the proper transform $\tilde{\Delta}_1$ in $\mathbf{A}^1_{u_0} \times \mathbf{A}^3_{\alpha_0, \beta_1}$ and $\mathbf{A}^1_{u_1} \times \mathbf{A}^3_{\alpha_0, \beta_0}$ respectively. From these equations we see that $\tilde{\Delta}_1$ is smooth in the latter affine chart, and has a double point at the origin of the former, with tangent cone there given by $\beta_1^2 = 0$.

To obtain the intersection with the exceptional divisor $Z$, we plug $\beta_1 = 0$ and $\beta_0 = 0$ respectively in these two equations; we find

$$16\alpha_0^2(u_0^2 + \alpha_0) = 0 \quad \text{and} \quad 16\alpha_0^2(1 + \alpha_0 u_1^2) = 0$$

for the two equations of $Z_1 = \tilde{\Delta}_1 \cap Z$ in the two charts $\mathbf{A}^1_{u_0} \times \mathbf{A}^3_{\alpha_0, \beta_1}$ and $\mathbf{A}^1_{u_1} \times \mathbf{A}^3_{\alpha_0, \beta_0}$ respectively. The exceptional divisor $Z_1$ thus has two irreducible component, the fibre $\Phi_0$, and the hypersurface $u_0^2 + \alpha_0 u_1^2) = 0$; the latter maps $2:1$ onto the line $\Delta_2$, and corresponds to the normal directions of the two sheets of $\Delta_1$ along $\Delta_2$.

## 4 – The case $m = 3$ and beyond

From this point on it is no longer advisable to do the computations by hand, and I have used Macaulay2 [3]. In fact, even with a computer it is not realistic to go much further. Still we can make some interesting observations, that are summed up in Subsection 4.2.

### 4.1 – The case $m = 3$

When $m = 3$, we are considering the versal deformation space of the 3-tacnode, defined in $\mathbf{A}^2_{x,y}$ by the equation $y(y + x^3) = 0$; this is $S \to \Delta \cong \mathbf{A}^5_{\alpha_0, \alpha_1, \beta_0, \beta_1, \beta_2}$, with $S$ defined in $\mathbf{A}^2_{x,y} \times \Delta$ by the equation

$$y^2 + y(x^3 + \alpha_1 x + \alpha_0) + \beta_2 x^2 + \beta_1 x + \beta_0 = 0.$$

We are interested in the solid $\Delta_2 \subseteq \Delta$ of deformations of cogenus at least two of the 3-tacnode, especially along the locus $\Delta_3 \subseteq \Delta$ of deformations of cogenus three, which is the plane $\beta_2 = \beta_1 = \beta_0 = 0$ along wich the 3-tacnode deforms to three (possibly infinitely near) nodes.

The singularities of $\Delta_2$ are more complicated than those of $\Delta_1$ in the case $m = 2$, and we will content ourselves with the description of the exceptional divisor of $\tilde{\Delta}_2$.

**(4.1) Equations of $\Delta_2$.** We shall proceed by elimination, as in (3.2). The solid $\Delta_2$ is defined by the condition that the polynomial

$$(4.1.1) \quad \begin{aligned} \delta_{\alpha_0, \alpha_1, \beta_0, \beta_1, \beta_2}(x) &= (x^3 + \alpha_1 x + \alpha_0)^2 - 4(\beta_2 x^2 + \beta_1 x + \beta_0) \\ &= x^6 + 2\alpha_1 x^4 + 2\alpha_0 x^3 + (\alpha_1^2 - 4\beta_2)x^2 + (2\alpha_0\alpha_1 - 4\beta_1)x + (\alpha_0^2 - 4\beta_0) \end{aligned}$$

has at least two non-simple roots, equivalently that there exist $(u, v, w) \in \mathbf{A}^3$ such that

$$(4.1.2) \qquad \delta_{\alpha_0, \alpha_1, \beta_0, \beta_1, \beta_2}(x) = (x-u)^2(x-v)^2(x-w)(x + 2u + 2v + w).$$

It is thus the projection to the second factor of the complete intersection in $\mathbf{A}^3_{u,v,w} \times \mathbf{A}^5_{\alpha_0, \alpha_1, \beta_0, \beta_1, \beta_2}$ defined by the four non-trivial equations obtained by identifying the $x$-terms in (4.1.1) and (4.1.2). Eliminating $u, v, w$ from these equations, for instance by using the function `eliminate` from Macaulay2 [3], one finds the ideal of $\Delta_2$ in $\Delta$ is the following:

$$\begin{aligned}
\big(&48\alpha_0\alpha_1^2\beta_2^2 + 3\alpha_1^2\beta_0\beta_1 + 45\alpha_0\alpha_1\beta_1^2 - 72\alpha_0\alpha_1\beta_0\beta_2 + 27\alpha_0^2\beta_1\beta_2 \\
&- 128\alpha_1\beta_1\beta_2^2 + 64\alpha_0\beta_2^3 + 27\alpha_0\beta_0^2 - 125\beta_1^3 + 180\beta_0\beta_1\beta_2, \\
9\alpha_1^3&\beta_0\beta_1 - 9\alpha_0\alpha_1^2\beta_1^2 + 63\alpha_0\alpha_1^2\beta_0\beta_2 - 54\alpha_0^2\alpha_1\beta_1\beta_2 - 81\alpha_1^3\beta_2^2 + 16\alpha_1^2\beta_1\beta_2^2 - 128\alpha_0\alpha_1\beta_2^3 - 54\alpha_0\alpha_1\beta_0^2 \\
&- 81\alpha_0^2\beta_0\beta_1 + 25\alpha_1\beta_1^3 - 252\alpha_1\beta_0\beta_1\beta_2 + 135\alpha_0\beta_1^2\beta_2 + 324\alpha_0\beta_0\beta_2^2 + 192\beta_1\beta_2^3 + 405\beta_0^2\beta_1, \\
16\alpha_1^4&\beta_2^2 + 16\alpha_1^3\beta_1^2 - 16\alpha_1^3\beta_0\beta_2 - 216\alpha_0^2\alpha_1\beta_2^2 - 128\alpha_1^2\beta_2^3 - 81\alpha_0\alpha_1\beta_0\beta_1 - 135\alpha_0^2\beta_1^2 + 135\alpha_0^2\beta_0\beta_2 \\
&- 120\alpha_1\beta_1^2\beta_2 + 576\alpha_1\beta_0\beta_2^2 + 288\alpha_0\beta_1\beta_2^2 + 256\beta_2^4 + 675\beta_0\beta_1^2 - 432\beta_0^2\beta_2,
\end{aligned}$$

$$16\alpha_1^4\beta_0\beta_2 - 16\alpha_0\alpha_1^3\beta_1\beta_2 - 16\alpha_1^3\beta_0^2 - 81\alpha_0\alpha_1^2\beta_0\beta_1 + 81\alpha_0^2\alpha_1\beta_1^2 - 216\alpha_0^2\alpha_1\beta_0\beta_2 + 135\alpha_0^3\beta_1\beta_2$$
$$+ 16\alpha_1^2\beta_1^2\beta_2 - 128\alpha_1^2\beta_0\beta_2^2 + 192\alpha_0\alpha_1\beta_1\beta_2^2 + 135\alpha_0^2\beta_0^2 + 240\alpha_1\beta_0\beta_1^2 - 225\alpha_0\beta_1^3 + 576\alpha_1\beta_0^2\beta_2$$
$$- 252\alpha_0\beta_0\beta_1\beta_2 - 320\beta_1^2\beta_2^2 + 256\beta_0\beta_2^3 - 432\beta_0^3,$$
$$81\alpha_0\alpha_1^3\beta_0\beta_2 - 81\alpha_0^2\alpha_1^2\beta_1\beta_2 - 1215\alpha_0^3\alpha_1\beta_2^2 + 240\alpha_1^3\beta_1^2\beta_2 - 81\alpha_0\alpha_1^2\beta_0^2 - 486\alpha_0^2\alpha_1\beta_0\beta_1 - 729\alpha_0^3\beta_1^2$$
$$+ 240\alpha_1^2\beta_1^3 + 729\alpha_0^3\beta_0\beta_2 - 420\alpha_1^2\beta_0\beta_1\beta_2 + 261\alpha_0\alpha_1\beta_1^2\beta_2 + 1980\alpha_0\alpha_1\beta_0\beta_2^2 + 2052\alpha_0^2\beta_1\beta_2^2$$
$$- 2240\alpha_1\beta_1\beta_2^3 + 2560\alpha_0\beta_2^4 + 243\alpha_1\beta_0^2\beta_1 + 3645\alpha_0\beta_0\beta_1^2 - 1836\alpha_0\beta_0^2\beta_2 - 2300\beta_1^3\beta_2 + 3312\beta_0\beta_1\beta_2^2\bigr).$$

It has five generators, quasi-homogeneous of degrees 15, 17, 16, 18, 19 respectively, for the weights 3, 2 for $\alpha_0, \alpha_1$ and 6, 5, 4 for $\beta_0, \beta_1, \beta_2$.

**(4.2) Tangent cone at the origin.** Using the ideal of $\Delta_2$ in $\Delta$, one finds (e.g., with the function `tangentCone` of Macaulay2) its tangent cone at the origin is defined by the primary ideal
$$\bigl(25\beta_0\beta_1^2 - 16\beta_0^2\beta_2, \beta_0^2\beta_1, \beta_0^3, 27\alpha_0\beta_0^2 - 125\beta_1^3 + 180\beta_0\beta_1\beta_2\bigr).$$
It defines a degree 6 projective scheme, supported over the 3-plane $\beta_1 = \beta_2 = 0$.

**(4.3) Blow-up along the plane $\Delta_3$.** We consider the blow-up of $\Delta$ along $\Delta_3$ as the subvariety of $\mathbf{P}^2_{(u_0:u_1:u_2)} \times \mathbf{A}^5_{\alpha_0,\alpha_1,\beta_0,\beta_1,\beta_2}$ defined by the equations $u_j\beta_i - u_i\beta_j = 0$ for all distinct $i$ and $j$ in $\{0, 1, 2\}$.

Let us consider the affine chart $u_2 = 1$ to fix ideas; it may be seen in $\mathbf{A}^2_{u_0,u_1} \times \mathbf{A}^3_{\alpha_0,\alpha_1,\beta_2}$ by eliminating $\beta_0$ and $\beta_1$, since $\beta_0 = u_0\beta_2$ and $\beta_1 = u_1\beta_2$ in this chart. The pull-back of the ideal of $\Delta_2$ in $\Delta$ is the intersection of the two primary ideals $(\beta_2^2)$ and

$$\bigl(3\alpha_1^2 u_0 u_1 + 45\alpha_0\alpha_1 u_1^2 - 125\beta_2 u_1^3 + 48\alpha_0\alpha_1^2 - 72\alpha_0\alpha_1 u_0 + 27\alpha_0 u_0^2 + 27\alpha_0^2 u_1 - 128\alpha_1\beta_2 u_1 + 180\beta_2 u_0 u_1 + 64\alpha_0\beta_2,$$
$$16\alpha_1^3 u_1^2 + 16\alpha_1^4 - 16\alpha_1^3 u_0 - 81\alpha_0\alpha_1 u_0 u_1 - 135\alpha_0^2 u_1^2 - 120\alpha_1\beta_2 u_1^2 + 675\beta_2 u_0 u_1^2 - 216\alpha_0^2 - 128\alpha_1^2\beta_2$$
$$+ 135\alpha_0^2 u_0 + 576\alpha_1\beta_2 u_0 - 432\beta_2 u_0^2 + 288\alpha_0\beta_2 u_1 + 256\beta_2^2,$$
$$144\alpha_0\alpha_1^2 u_1^2 - 400\alpha_1\beta_2 u_1^3 + 144\alpha_0\alpha_1^3 - 279\alpha_0\alpha_1^2 u_0 + 135\alpha_0\alpha_1 u_0^2 + 135\alpha_0^2\alpha_1 u_1 - 400\alpha_1^2\beta_2 u_1 + 81\alpha_0^2 u_0 u_1$$
$$+ 792\alpha_1\beta_2 u_0 u_1 - 405\beta_2 u_0^2 u_1 - 135\alpha_0\beta_2 u_1^2 + 81\alpha_0^3 + 320\alpha_0\alpha_1\beta_2 - 324\alpha_0\beta_2 u_0 - 192\beta_2^2 u_1,$$
$$\alpha_1^4 u_0 - \alpha_1^3 u_0^2 - \alpha_0\alpha_1^3 u_1 + 81\alpha_0^2\alpha_1 u_1^2 + \alpha_1^2\beta_2 u_1^2 + 15\alpha_1\beta_2 u_0 u_1^2 - 225\alpha_0\beta_2 u_1^3 + 81\alpha_0^2\alpha_1^2 - 135\alpha_0^2\alpha_1 u_0$$
$$- 8\alpha_1^2\beta_2 u_0 + 54\alpha_0^2 u_0^2 + 36\alpha_1\beta_2 u_0^2 - 27\beta_2 u_0^3 + 54\alpha_0^3 u_1 - 204\alpha_0\alpha_1\beta_2 u_1$$
$$+ 288\alpha_0\beta_2 u_0 u_1 - 20\beta_2^2 u_1^2 + 108\alpha_0^2\beta_2 + 16\beta_2^2 u_0\bigr).$$

The latter has four generators, and defines the proper transform $\tilde{\Delta}_2$ of $\Delta_2$ in this affine chart. To obtain the intersection $Z_2$ of $\tilde{\Delta}_2$ with the exceptional divisor, we plug in the equation $\beta_2 = 0$. We find the ideal of $Z_2$ in $\mathbf{A}^2_{u_0,u_1} \times \mathbf{A}^3_{\alpha_0,\alpha_1,\beta_2}$ is the intersection of the three primary ideals

$$\bigl(\alpha_0\alpha_1, \alpha_0^2, \alpha_1^2 u_1 + 9\alpha_0 u_0, \alpha_1^3\bigr),$$
$$\bigl(u_1^2 + \alpha_1 - u_0, u_0 u_1 + \alpha_0, \alpha_1 u_0 - u_0^2 - \alpha_0 u_1\bigr), \text{ and}$$
$$\bigl(\alpha_1 u_0 u_1^2 - \alpha_0 u_1^3 + \alpha_1^2 u_0 - 2\alpha_1 u_0^2 + u_0^3 - \alpha_0\alpha_1 u_1 + 3\alpha_0 u_0 u_1 + \alpha_0^2,$$
$$16\alpha_1^2 u_1^2 + 16\alpha_1^3 - 31\alpha_1^2 u_0 + 15\alpha_1 u_0^2 + 15\alpha_0\alpha_1 u_1 + 9\alpha_0 u_0 u_1 + 9\alpha_0^2,$$
$$\alpha_1^2 u_0 u_1 + 15\alpha_0\alpha_1 u_1^2 + 16\alpha_0\alpha_1^2 - 24\alpha_0\alpha_1 u_0 + 9\alpha_0 u_0^2 + 9\alpha_0^2 u_1,$$
$$\alpha_1^3 u_0 - \alpha_1^2 u_0^2 - \alpha_0\alpha_1^2 u_1 - 6\alpha_0\alpha_1 u_0 u_1 - 9\alpha_0^2 u_1^2 - 15\alpha_0^2\alpha_1 + 9\alpha_0^2 u_0,$$
$$16\alpha_1^4 - 15\alpha_1^2 u_0^2 - 96\alpha_0\alpha_1^2 u_1 - 90\alpha_0\alpha_1 u_0 u_1 - 135\alpha_0^2 u_1^2 - 360\alpha_0^2\alpha_1 + 216\alpha_0^2 u_0\bigr).$$

The first one has radical $(\alpha_0, \alpha_1)$, and defines the fibre $\Phi_0$ of $Z$ over the origin with multiplicity 3. The second one is radical, and defines a locus that dominates $\Delta_3$ (this can be seen by eliminating the variables $u_0$ and $u_1$, which gives the zero ideal); moreover, one finds there are three points in this locus over the generic point of $\Delta_3$. As for the third one, if we eliminate $u_0$ and $u_1$, we find the ideal
$$\big((4\alpha_1^3 + 27\alpha_0^2)^2\big),$$
which is generated by the square of the discriminant of the polynomial $x^3 + \alpha_1 x + \alpha_0$; over the generic point of the locus defined by this discriminant in $\Delta_3$, we find the the third ideal defines a line with multiplicity 2.

Computations in the two other affine charts are similar, although they take a little more time to the computer. In conclusion we find that the exceptional divisor $Z_2$ of $\tilde\Delta_2 \to \Delta_2$ has three irreducible components: (i) the fibre $\Phi_0$ of $Z$ over the origin; (ii) a surface dominating $\Delta_3$, and corresponding to the normal directions of the three sheets of $\Delta_2$ along $\Delta_3$; (iii) a surface fibered in lines over the cuspidal curve $\Delta\{2,1\} \subseteq \Delta_3$ (in the notation of Section 2).

## 4.2 – General conclusions

**(4.4) Some computations in the case $m=4$.** For $m=4$ it is still possible to compute the defining ideal of the fourfold $\Delta_3$ in $\mathbf{A}^7_{\alpha_0,\ldots,\alpha_2,\beta_0,\ldots,\beta_3}$, although it took about half an hour on my personal computer.

This ideal has 20 generators, all quasihomogenous if we assign the weights $4,3,2$ to $\alpha_2,\alpha_1,\alpha_0$ and $8,\ldots,5$ to $\beta_3,\ldots,\beta_0$ respectively. For what it's worth, the degrees of the generators are $18, 19, 20, 20, 21, 21, 21, 22, 22, 22, 23, 23, 24, 24, 28, 28, 29, 30, 30, 31$.

The tangent cone to $\Delta_2$ at the origin is supported on the 4-plane $\beta_0 = \beta_1 = \beta_2 = 0$, and has degree 10.

I have not been able to compute a primary decomposition of the intersection of $\tilde\Delta_3$ with the exceptional divisor $Z$ of the blow-up of $\Delta$ along $\Delta_4$.

**(4.5) General conjectures.** From the above computations one can conjecture that the defining ideal of $\Delta_{m-1}$ in $\Delta$ will always be quasihomogeneous for the weights $m,\ldots,2$ for $\alpha_0, \alpha_{m-2}$ and $2m, \ldots, m+1$ for $\beta_0, \ldots, \beta_{m-1}$, respectively.

The origin will be a point of multiplicity $\sum_{k=1}^m k$, with tangent cone there supported on the $m$-plane $\beta_0 = \cdots = \beta_{m-2} = 0$.

The exceptional divisor of $\tilde\Delta_{m-1} \to \Delta_{m-1}$ will have $m-1$ irreducible components, respectively dominating the following loci in $\Delta_m$ : $\Delta\{1,\ldots,1\} = \Delta_m$, $\Delta\{2,1,\ldots,1\}$, …, $\Delta\{m-1,1\}$.

The latter two facts are probably not very hard to prove using the setup of Section 2. It is plausible that the first conjecture may be approached using reduced discriminant theory (see [**??**]).

## References

**Chapters in this Volume**

# Lecture VI
# Degenerations and deformations of generalized Severi varieties

by Thomas Dedieu and Concettina Galati

This text is dedicated to the recursive formula obtained by Caporaso and Harris [1] as a solution to the problem of "counting plane curves of any genus", see Corollary (1.5) and Theorem (5.4). Fix a degree $d \in \mathbf{N}^*$, and denote by $p_a(d) = \frac{1}{2}(d-1)(d-2)$ the genus of smooth plane curves of degree $d$. For all $g = 0, \ldots, p_a(d)$, (reduced) plane curves of degree $d$ and geometric genus $g$ form a family of pure dimension $\nu_g^d = 3d + g - 1$. The problem is to find the number of plane curves of degree $d$ and geometric genus $g$ passing through $\nu_g^d$ general points in the plane. It turns out that this number equals the number of $\delta$-nodal plane curves of degree $d$ which pass through $\nu_g^d$ general points in the plane, with $\delta = p_a(d) - g$ (a curve is $\delta$-nodal if its singularities are exactly $\delta$ nodes, i.e., $\delta$ ordinary double points).

Caporaso and Harris obtain their inductive formula by specializing one after another the aforementioned $\nu_g^d$ imposed crossing points to general points on a fixed line $L$. This may be conveniently formulated in terms of a degeneration of the projective plane $\mathbf{P}^2$ to the transverse union along a line of $\mathbf{P}^2$ and a minimal rational ruled surface $\mathbf{F}_1$, using the formalism presented in [I]. We will do it on an example in [VII], but here we stick to the presentation given in [1]. A similar result for curves on minimal rational ruled surfaces has been obtained by Vakil [6], following a similar degeneration procedure, and around the same period of the time. We shall nevertheless concentrate on the case of plane curves.

The above described degeneration procedure makes it necessary to generalize the problem. Thus, one is led to count curves of given degree and genus with, in addition, a prescribed intersection pattern with the line $L$, both at assigned and unassigned points. We shall call *generalized Severi varieties* the families formed by these curves: they are a particular case of the logarithmic Severi varieties studied in [III].

We give precise definitions and statements in Section 1 below. There, we shall also provide more details about the general strategy of the proof, which will explain the title of this text.

# 1 – The counting problem and its solution

Let us fix once and for all a line $L \subseteq \mathbf{P}^2$. We shall use freely the notation and definitions introduced in [III, Section 1], in the special case when the pair $(S, R)$ is $(\mathbf{P}^2, L)$; every class $\xi \in \mathrm{NS}(\mathbf{P}^2)$ is $d$ times the linear equivalence class of a line for some integer $d$, and therefore it shall be denoted simply by $d$.

Let $d \in \mathbf{N}$ and $g \in \mathbf{Z}$, let $\alpha = (\alpha_i)_{i \geqslant 1}$ and $\beta = (\beta_i)_{i \geqslant 1}$ be two sequences in $\underline{\mathbf{N}}$ such that $I\alpha + I\beta = d$, and let $\Omega = \left(\{p_{i,j}\}_{1 \leqslant j \leqslant \alpha_i}\right)_{i \geqslant 1}$ be a set of $\alpha$ points on the line $L$. The *logarithmic Severi variety* $V_g^d(\alpha, \beta)(\Omega)$ is the locally closed subset of $|\mathcal{O}_{\mathbf{P}^2}(d)|$ (the projective space of dimension $\frac{1}{2}d(d+3)$ parametrizing all plane curves of degree $d$) parametrizing those reduced plane curves of degree $d$ and geometric genus $g$ that intersect $L$ at the points of $\Omega$ with the multiplicities prescribed by $\alpha$ and at further unassigned points with the multiplicities prescribed by $\beta$; for the rigorous definition, see [III, Definition (1.5)]. The points of $\Omega$ shall be referred to as the *fixed*, or *assigned*, contact points, and the other intersection points with $L$ as *mobile*, or *unassigned*, contact points (as it turns out, in the present situation, the unassigned contact points indeed move as the curves moves in the Severi variety, see Proposition (2.2) below). The relevant notation for these points and their relatives is explained in Paragraph (2.1). Forgetting about the location of assigned intersection points at $\Omega$, we may also loosely say that the curves in $V_g^d(\alpha, \beta)(\Omega)$ have intersection pattern with $L$ prescribed by $\alpha$ and $\beta$.

As it turns out (see Proposition (2.2) below), the general element of every irreducible component of $V_g^d(\alpha, \beta)(\Omega)$ has $\delta$ nodes and no further singularities, with

$$\delta = p_a(d) - g = \frac{1}{2}(d-1)(d-2) - g.$$

The upshot is that it is indifferent to require either that the geometric genus be $g$ or that the singularities be exactly $\delta$ nodes.

**(1.1) Definition.** *For all natural integers $d$ and $\delta$, and all $\alpha, \beta \in \underline{\mathbf{N}}$ and $\Omega$ as above, the generalized Severi variety $V^{d,\delta}(\alpha, \beta)(\Omega)$ is the Zariski closure in the projective space $|\mathcal{O}_{\mathbf{P}^2}(d)|$ of the logarithmic Severi variety $V_{p_a(d)-\delta}^d(\alpha, \beta)(\Omega)$ of the pair $(\mathbf{P}^2, L)$.*

In the particular case when $\alpha = 0$, and thus $\Omega = \varnothing$, and $\beta = (d, 0, \ldots)$ (i.e., when no conditions are imposed on the intersection with $L$), the Severi varieties $V_g^d(\alpha, \beta)(\Omega)$ and $V^{d,\delta}(\alpha, \beta)(\Omega)$ will simply be denoted by $V_g^d$ and $V^{d,\delta}$, respectively, and called *plain* Severi varieties.

Although this should not create any problem, beware that in our notation $V_g^d$ and $V^{d,\delta}$ (with any additional decorations) are respectively locally closed and closed in $|\mathcal{O}_{\mathbf{P}^2}(d)|$.

For any point $p \in \mathbf{P}^2$, the degree $d$ plane curves passing through the point $p$ form a hyperplane $p^\perp$ in the projective space $|\mathcal{O}_{\mathbf{P}^2}(d)|$. Therefore, the problem of counting curves of degree $d$ and genus $g$ (and no further conditions) mentioned in the introduction amounts to the problem of computing the degree of the plain Severi variety $V^{d,\delta} = \overline{V}_g^d$, with $\delta = p_a(d) - g$. We shall denote this number by $N^{d,\delta}$. More generally, for all $\alpha, \beta \in \underline{\mathbf{N}}$, we shall denote by $N^{d,\delta}(\alpha, \beta)$ the common degree of all $V_g^d(\alpha, \beta)(\Omega)$, where $\Omega$ is a general, cardinality $\alpha$, set of points of $L$.

The recursion procedure of Caporaso and Harris consists in computing the number $N^{d,\delta}$ by cutting $V^{d,\delta}$ by hyperplanes $p^\perp$ with $p \in L$. Their main result is Theorem (1.4) below. We need some additional notation to state it.

**(1.2) Notation.** For all $k \in \mathbf{N}^*$, we let $e_k = (0, \ldots, 0, 1) \in \underline{\mathbf{N}}$: it is the sequence with only one non-zero entry, in $k$-th position and equal to 1 (recall the convention made in [III] that we may omit the infinite sequence of zeros at the end of elements of $\underline{\mathbf{N}}$). For all $\beta, \beta' \in \underline{\mathbf{N}}$, we let

$$I^\beta = \prod_{k \geqslant 1} k^{\beta_k} \qquad \text{and} \qquad \binom{\beta'}{\beta} = \prod_{k \geqslant 1} \binom{\beta'_k}{\beta_k}.$$

If $\Omega = (\Omega_1, \Omega_2, \ldots)$ is a set of cardinality $\alpha \in \underline{\mathbf{N}}$, a subset of $\Omega$ of cardinality $\alpha' \leqslant \alpha$ is simply a set $\Omega' = (\Omega'_1, \Omega'_2, \ldots)$ of cardinality $\alpha' \in \underline{\mathbf{N}}$ such that $\Omega'_k \subseteq \Omega_k$ for all $k \geqslant 1$.

The following notation is not needed to state the main theorem below, but we shall use it for the proof and in the more precise Theorem (4.1). For all $\beta \in \underline{\mathbf{N}}$, we let

$$\mathrm{lcm}(\beta) = \mathrm{lcm}\{i : \beta_i \neq 0\} \qquad \text{and} \qquad \max(\beta) = \max\{i : \beta_i \neq 0\}.$$

**(1.3) Convention.** In the statement below, generalized Severi varieties $V^{d-1,\delta'}(\alpha',\beta')(\Omega') \subseteq |\mathcal{O}_{\mathbf{P}^2}(d-1)|$ are considered as the projective subvarieties of $|\mathcal{O}_{\mathbf{P}^2}(d)|$ parametrizing all curves of the form $C + L$ with $C \in V^{d-1,\delta'}(\alpha',\beta')(\Omega')$. We will make this identification freely all along the text, where needed.

**(1.4) Theorem** (Caporaso–Harris)**.** *Let $d$ and $\delta$ be natural integers, and let $\alpha, \beta \in \underline{\mathbf{N}}$ be sequences such that $I\alpha + I\beta = d$, and $\Omega$ be a set of $\alpha$ general points on the line $L$. Let $p$ be a general point of $L$. Then one has the following equality of cycles:*

$$(1.4.1) \quad V^{d,\delta}(\alpha,\beta)(\Omega) \cap p^\perp = \sum_{k \geqslant 1:\, \beta_k > 0} k \cdot V^{d,\delta}(\alpha + e_k, \beta - e_k)(\Omega \cup \{p\})$$

$$+ \sum_{\substack{\Omega' \subseteq \Omega,\, \mathrm{Card}(\Omega') = \alpha' \\ \beta' \geqslant \beta:\, I\alpha' + I\beta' = d-1 \\ \delta' = \delta + |\beta' - \beta| - d + 1}} I^{\beta' - \beta} \binom{\beta'}{\beta} \cdot V^{d-1,\delta'}(\alpha',\beta')(\Omega')$$

*where, in the first sum, $\Omega \cup \{p\}$ means that we add the point $p$ to the $k$-th term of $\Omega = (\Omega_1, \Omega_2, \ldots)$.*

Note that, in the above theorem, if we start with a plain Severi variety for $V^{d,\delta}(\alpha,\beta)(\Omega)$ (i.e., if $\alpha = 0$ and thus $\Omega = \varnothing$, and $\beta = (d)$), then the intersection with $p^\perp$ decomposes as a sum of genuinely generalized (i.e., non-plain) Severi varieties. This is why one needs to generalize the problem of counting curves of fixed degree and genus for the inductive solution proposed by Caporaso and Harris.

It follows from Proposition (2.2) below that if $p$ is outside of $\Omega$, then no irreducible component of $V^{d,\delta}(\alpha,\beta)(\Omega)$ is contained in the hyperplane $p^\perp$. Therefore, all irreducible varieties appearing on the right-hand-side of (1.4.1) above have codimension one in $V^{d,\delta}(\alpha,\beta)(\Omega)$, and indeed we shall verify this by hand later in the text. The number $N^{d,\delta}(\alpha,\beta)$ we are looking for is the degree of $V^{d,\delta}(\alpha,\beta)(\Omega)$ as well as that of the cycle on the right-hand-side of (1.4.1). Thus, Theorem (1.4) has the following direct corollary.

**(1.5) Corollary.** *Let $d$ and $\delta$ be natural integers, and let $\alpha, \beta \in \underline{\mathbf{N}}$ be sequences such that $I\alpha + I\beta = d$. Then one has the following equality:*

$$N^{d,\delta}(\alpha,\beta) = \sum_{k \geqslant 1:\, \beta_k > 0} k \cdot N^{d,\delta}(\alpha + e_k, \beta - e_k) + \sum_{\substack{\alpha' \leqslant \alpha \\ \beta' \geqslant \beta:\, I\alpha' + I\beta' = d-1 \\ \delta' = \delta + |\beta' - \beta| - d + 1}} I^{\beta' - \beta} \binom{\beta'}{\beta} \binom{\alpha'}{\alpha} \cdot N^{d-1,\delta'}(\alpha',\beta').$$

Since the numbers $N^{1,\delta}(\alpha,\beta)$ are readily computed, the above recursive formula enables one to compute all numbers $N^{d,\delta}(\alpha,\beta)$, and thus to completely solve our counting problem. Note that the recursion indeed terminates, as those terms appearing in the first sum are all of the form $N^{d,\delta}(\alpha'',\beta'')$ with $\beta'' < \beta$, so that unassigned intersection with $L$ are exhausted after finitely many applications of the formula, and then only terms of the form $N^{d-1,\delta'}(\alpha',\beta')$ appear. It is of course helpful to watch the inductive process work on a concrete example. For this we refer to [VII], or [1, p. 349]; see also Example (5.1).

Although, arguably, the relevant enumerative problem is to count irreducible plane curves of any genus, there is no irreducibility requirement in our definition of Severi varieties, and therefore the numbers $N^{d,\delta}(\alpha,\beta)$ count curves regardless of their being irreducible or not. It is yet possible to derive a recursion formula counting irreducible curves using the inductive process of Caporaso and Harris (and so do they in the original article). We postpone this until Section 5, as the formula is yet more complicated than the above one, and having the geometry of the inductive process in mind helps in understanding it.

The proof of Theorem (1.4) is in two parts. The first one consists in using stable reduction in order to limit the possible irreducible components of the intersection of $V^{d,\delta}(\alpha,\beta)(\Omega)$ and $p^\perp$. This is the degeneration part; it is carried out in Section 3 and the result is Theorem (3.1). The dimensional characterization of Severi varieties (Proposition (2.3)) is fundamental in identifying the components of the intersection as generalized Severi varieties.

The second part consists in showing that all the generalized Severi varieties found in the first part indeed occur as irreducible components of the intersection of $V^{d,\delta}(\alpha,\beta)(\Omega)$ and $p^\perp$, and to compute their multiplicites in this intersection. This is the deformation part; it is carried out in Section 4 and the result is Theorem (4.1). A fundamental result for this part is the study of the deformation space of an $m$-tacnode into $m-1$ nodes, also due to Caporaso and Harris, and described in [IV]. The latter result is certainly the central result of this whole volume of notes.

Theorems (3.1) and (4.1) together directly imply Theorem (1.4).

## 2 – Recap on Severi varieties

**(2.1)** We shall use consistent notation to denote assigned and unassigned contact points with the line $L$. Consider a generalized Severi variety $V^{d,\delta}(\alpha,\beta)(\Omega)$ as in Definition (1.1). The assigned contact points with $L$ are the elements of $\Omega$, which we write as $\Omega = \big(\{p_{i,j}\}_{1 \leqslant j \leqslant \alpha_i}\big)_{i \geqslant 1}$. Let $[C]$ be a general member of any irreducible component of $V^{d,\delta}(\alpha,\beta)(\Omega)$. Denote by $\bar{C} \to C$ the normalization of $C$, and by $\phi : \bar{C} \to \mathbf{P}^2$ its composition with the inclusion $C \subseteq \mathbf{P}^2$. By definition of the Severi varieties, there exist $\alpha$ points $q_{i,j} \in \bar{C}$, for all $i \geqslant 1$ and $1 \leqslant j \leqslant \alpha_i$, and $\beta$ points $r_{i,j} \in \bar{C}$, for all $i \geqslant 1$ and $1 \leqslant j \leqslant \beta_i$, such that

$$\forall\, i \geqslant 1,\ \forall\, j = 1,\ldots,\alpha_i : \quad \phi(q_{i,j}) = p_{i,j}$$

and

$$\phi^* L = \sum_{i \geqslant 1} \sum_{1 \leqslant j \leqslant \alpha_i} i\, q_{i,j} + \sum_{i \geqslant 1} \sum_{1 \leqslant j \leqslant \beta_i} i\, r_{i,j}.$$

Moreover, we let $s_{i,j} = \phi(r_{i,j})$ for all $i \geqslant 1$ and $1 \leqslant j \leqslant \beta_i$.

Thus, the assigned and unassigned contact points of $C$ with $L$ are, respectively, the $p_{i,j}$'s and the $s_{i,j}$'s, and their respective counterparts on the normalization $\bar{C}$ are the $q_{i,j}$'s and the $r_{i,j}$'s.[1]

The two following results are key ingredients in the proof of both Theorems (3.1) and (4.1). They are direct applications of Theorem (1.7) and Proposition (4.3) in [III], respectively.

---

[1] The reader may wish to keep in mind that $\phi(q_{i,j}) = p_{i,j}$ whereas $\phi(r_{i,j}) = s_{i,j}$; thus, beware, $\phi$ goes backwards with the alphabetical order for the assigned contact points, but forward for the unassigned ones.

**(2.2) Proposition.** *Let the notation be as in (2.1) above. The generalized Severi variety $V^{d,\delta}(\alpha,\beta)(\Omega)$ has pure dimension $2d+g-1+|\beta|$, with $g = p_a(d)-\delta$. Morever:*
— *the curve $C \subseteq \mathbf{P}^2$ is $\delta$-nodal, and smooth at its intersection points with the line $L$;*
— *the contact points of $C$ with $L$, i.e., all the points $p_{i,j}$ and $s_{i,j}$, are altogether pairwise distinct;*
— *the counterparts on the normalization $\bar{C}$ of the contact points of $C$ with $L$, i.e., all the points $q_{i,j}$ and $r_{i,j}$, are altogether pairwise distinct;*
— *for any curve $G$ and finite set $\Gamma$ in $\mathbf{P}^2$, if $(G \cup \Gamma) \cap \Omega = \varnothing$ and $[C]$ is general with respect to $G$ and $\Gamma$, then $C$ intersects $G$ transversely and does not intersect $\Gamma$.*

**(2.3) Proposition.** *Let $V \subseteq |\mathcal{O}_{\mathbf{P}^2}(d)|$ be an irreducible variety, parametrizing genus $g$ curves in the following sense: for a general member $[X]$ of $V$, there exists a smooth curve $\tilde{X}$ of genus $g$, and a morphism $f : \tilde{X} \to X$, not constant on any component of $\tilde{X}$, and such that the push-forward in the sense of cycles $f_*\tilde{X}$ equals the fundamental cycle of $X$.*

*Let $\Omega \subseteq L$ be a finite set of points. For any general member $[X]$ of $V$, one has*

$$\dim(V) \leqslant 2d + g - 1 + \mathrm{Card}\big((X \cap L) \setminus \Omega\big),$$

*where the last number is defined set-theoretically (i.e., multiplicities do not count).*

*Moreover, if equality holds, then there exist $\alpha, \beta \in \underline{\mathbf{N}}$ such that $V$ is a dense subset of a component of the generalized Severi variety $V^{d,\delta}(\alpha,\beta)(\Omega)$, $\delta = p_a(d) - g$, if and only if*

$$\mathrm{Card}\big(f^{-1}(X \cap L)\big) = \mathrm{Card}(X \cap L),$$

*with $f : \tilde{X} \to X$ a genus $g$ morphism as above.*

Lastly, let us mention that, although we pretend not to know it in this text, all generalized Severi varieties $V^{d,\delta}(\alpha,\beta)(\Omega)$ have their irreducible components in bijection with the possible splittings of degree $d$, $\delta$-nodal curves, into several irreducible components. Thus, for instance, there is at most one irreducible component of $V^{d,\delta}(\alpha,\beta)(\Omega)$ parametrizing irreducible curves. This is proved in [14] for plain Severi varieties, and in [7] for generalized Severi varieties.

## 3 – Degenerations of generalized Severi varieties

This section is dedicated to the analysis of how curves in a generalized Severi variety $V^{d,\delta}(\alpha,\beta)(\Omega)$ degenerate when one imposes the passing through a general point of the line $L$. The main result of this section is Theorem (3.1) below, but we will also provide in Section 3.4 a geometric decription of the degeneration at play, of fundamental importance for the proof of the theorem in the other direction, i.e., Theorem (4.1) on deformations of generalized Severi varieties, which is the object of the next section.

In the following statement, and as in Convention (1.3), for all $\delta', \alpha', \beta', \Omega'$, the generalized Severi variety $V^{d-1,\delta'}(\alpha',\beta')(\Omega')$ is identified with a variety parametrizing degree $d$ curves containing the line $L$, namely

$$\big\{[C \cup L] : \ [C] \in V^{d-1,\delta'}(\alpha',\beta')(\Omega')\big\} \ \subseteq \ |\mathcal{O}_{\mathbf{P}^2}(d)|.$$

**(3.1) Theorem** ([1, Theorem 1.2]). *Let $p \in L$ be a general point, and let $V$ be an irreducible component of the intersection $V^{d,\delta}(\alpha,\beta)(\Omega) \cap p^{\perp}$. Then $V$ is an irreducible component of one of the following varieties:*
*a) for all $k$ such that $\beta_k > 0$, the variety*

$$V^{d,\delta}(\alpha+e_k, \beta-e_k)(\Omega \cup \{p_{k,\alpha_k+1}\})$$

*with $p_{k,\alpha_k+1} = p$;*

b) *for all $\alpha' \leqslant \alpha$, $\beta' \geqslant \beta$, and $\delta' \leqslant \delta$ such that (i) $I\alpha' + I\beta' = d-1$, and (ii) $\delta - \delta' + |\beta' - \beta| = d-1$, the variety*

$$V^{d-1,\delta'}(\alpha', \beta')(\Omega'),$$

*where $\Omega'$ is an arbitrary cardinality $\alpha'$ subset of $\Omega$.*

**(3.2) Remark.** The condition $\delta - \delta' + |\beta' - \beta| = d - 1$ in case b), is equivalent to the condition $\dim(V^{d-1,\delta'}(\alpha', \beta')(\Omega')) = \dim(V^{d,\delta}(\alpha, \beta)(\Omega)) - 1$. Indeed,

$$\dim(V^{d,\delta}(\alpha, \beta)) - \dim(V^{d-1,\delta'}(\alpha', \beta')) = (2d + g - 1 + |\beta|) - (2(d-1) + g' - 1 + |\beta'|)$$
$$= g - g' + |\beta - \beta'| - 2,$$

where we let $g$ and $g'$ be the genera of the members of $V^{d,\delta}(\alpha, \beta)$ and $V^{d-1,\delta'}(\alpha', \beta')$ respectively, and thus

$$g - g' = \binom{d-1}{2} - \delta - \binom{d-2}{2} + \delta' = d - 2 - \delta + \delta'$$

by the Pascal Formula.

Later on we will be able to explain geometrically where the $\delta - \delta'$ nodes degenerate, see Subsection 3.4.

**(3.3) Remark.** Note moreover that, still in case b), the conditions $\alpha' \leqslant \alpha$, $\beta' \geqslant \beta$, and $I\alpha' + I\beta' = d-1$ altogether imply that $\alpha' < \alpha$, i.e., there exists an $i_0$ such that $\alpha'_{i_0} < \alpha_{i_0}$.

It is however possible that $\beta' = \beta$, but in that case all the irreducible components of $V^{d,\delta}(\alpha, \beta)$ having $V$ in their intersection with $p^\perp$ parametrize reducible curves, see Remark (4.21).

*Proof of Theorem (3.1).* Let $[X_0]$ be a general point of $V$. There are two cases to be considered. The first case is when the curve $X_0$ does not contain $L$. Then, in order for $[X_0]$ to sit in $p^\perp$, we must have one of the mobile contact points (which we may loosely refer to as points of type $\beta$) be at the prescribed point $p$ (note that by generality, $p$ is off $\Omega$). Thus one mobile contact point must be turned into an assigned contact point (i.e., a point of type $\beta$ is turned into one of type $\alpha$), hence $V$ must be contained in a variety

$$V^{d,\delta}(\alpha + e_k, \beta - e_k)(\Omega \cup \{p_{k,\alpha_k+1} = p\}).$$

Since the latter varieties all have pure dimension $\dim(V^{d,\delta}(\alpha, \beta)(\Omega)) - 1 = \dim(V)$, this implies that $V$ is an irreducible component of a generalized Severi variety as in case a) of the theorem.

The remaining case is when $X_0$ has $L$ as an irreducible component; then we shall loosely say that $X_0 = L \cup C$. In this case the content of Theorem (3.1) is Proposition (3.15) below. □

The proof of Proposition (3.15) occupies most of this section. The general idea is to apply (a suitable version of) semistable reduction to a family of general members of $V^{d,\delta}(\alpha, \beta)(\Omega)$ degenerating to $X_0$, in order to derive some necessary conditions on $X_0$ which will limit the number of dimensions in which it can move. Comparing this with the fact that $X_0$ has to move in dimension $\dim(V^{d,\delta}(\alpha, \beta)(\Omega)) - 1$, and taking advantage of the dimensional characterization of generalized Severi varieties (Corollary (2.3)), we will then be able to conclude.

**(3.4) Setup.** Let $V$ be an irreducible component of $V^{d,\delta}(\alpha, \beta) \cap p^\perp$, the general member of which has $L$ as an irreducible component. Let $X_0$ be a general member of $V$, and let $C$ be the sum of the components of $X_0$ supported on curves different from $L$; we shall loosely write $X_0 = L \cup C$.

**(3.5) Semistable reduction.** Let $\Gamma$ be a general analytic curve in $V^{d,\delta}(\alpha,\beta)$ passing through the point $[X_0]$, and let $\mathcal{X} \to \Gamma$ be the corresponding family of curves. By the generality of $\Gamma$, its general member $[X_\gamma]$ is general in an irreducible component of $V^{d,\delta}(\alpha,\beta)$, hence it corresponds to a reduced curve $X_\gamma$ with $\delta$ nodes as its only singularities.

Let $\nu : \Gamma^\nu \to \Gamma$ be the normalization of $\Gamma$, and let $b_0 \in \Gamma^\nu$ be a point with $\nu(b_0) = [X_0]$. We let $\mathcal{X}^\nu$ be the normalization of the total space $\mathcal{X} \times_\Gamma \Gamma^\nu$, obtained from $\mathcal{X}$ by the base change $\Gamma^\nu \to \Gamma$. The general fibre of $\mathcal{X}^\nu \to \Gamma^\nu$ is smooth of genus $g = p_a(d) - \delta$, by the Teissier Simultaneous Resolution Theorem, see [III].

Applying a suitable version of stable reduction to the family $\mathcal{X}^\nu \to \Gamma^\nu$, we arrive at a family of curves $f: \mathcal{Y} \to B$, with special fibre $Y_0 = f^{-1}(b_0)$, satisfying the following conditions:
- $\mathcal{Y} \to B$ is a nodal reduction of $\mathcal{X}^\nu \to \Gamma^\nu$, i.e., all fibres of $f$ are reduced, nodal curves, and $B$ is a branched cover of $\Gamma^\nu$;
- the total space $\mathcal{Y}$ is smooth;
- there are $|\alpha|$ sections $Q_{i,j}$ ($1 \leqslant j \leqslant \alpha_i$) and $|\beta|$ sections $R_{i,j}$ ($1 \leqslant j \leqslant \beta_i$) of $\mathcal{Y} \to B$, altogether pairwise disjoint, such that $\pi(Q_{i,j}) = p_{i,j}$ (the points of $\Omega$), and for all $b \neq b_0$ one has
$$(\pi|_{Y_b})^* L = \sum i \, Q_{i,j}|_{Y_b} + \sum i \, R_{i,j}|_{Y_b};$$
- the map sending the general fibre $Y_b$ of $\mathcal{Y}$ to the corresponding plane curve $X_\gamma$ extends to a regular map $\pi: \mathcal{Y} \to \mathbb{P}^2$ (here $\gamma = \nu(\eta(b))$, see diagram (3.5.1) below);
- the family $\mathcal{Y} \to B$ is minimal with respect to the above properties.

Note that, since $\mathcal{Y}$ is smooth and the $Q_{i,j}$'s and $R_{i,j}$'s are sections, they do not pass through any singular point of the central fibre $Y_0$.

The family $\mathcal{Y} \to B$ is obtained from $\mathcal{X}^\nu \to \Gamma^\nu$ by applying a finite number of base changes and blow-ups. A number of additional such operations may be required with respect to the usual nodal reduction in order to first organize the points $q_{i,j}$ and $p_{i,j}$ on the fibres into sections (a priori they may be permuted by the monodromy), and then separate these sections. The situation is summarized in the diagram below.

(3.5.1) 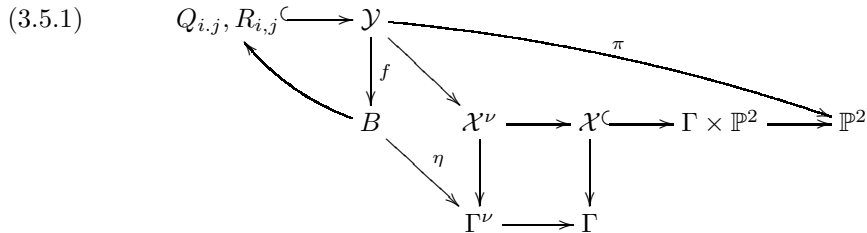

## 3.1 – General picture of the degeneration

**(3.6) Components of the central fibre of $\mathcal{Y}$.** The irreducible components of the central fibre $Y_0$ of $\mathcal{Y}$ come in three types. Namely, we let:
- $\tilde{L}$ be the union of the components of $Y_0$ on which $\pi$ is nonconstant and maps to $L$;
- $\tilde{C}$ be the union of the components of $Y_0$ on which $\pi$ is nonconstant and maps to $C$;
- $Z$ be the union of the components of $Y_0$ on which $\pi$ is constant.

**(3.7) Degenerations of assigned and mobile contact points.** We shall now relabel the sections $(Q_{i,j})_{1\leqslant j\leqslant \alpha_i}$ and $(R_{i,j})_{1\leqslant j\leqslant \beta_i}$ according to the type of the irreducible component of the central fibre they pass through. We advise the reader to read this paragraph with Figure 1 at hand.

For all $i$, we let $\alpha_i^C$ and $\alpha_i^L$ be the number of the sections $Q_{i,j}$, $1 \leqslant j \leqslant \alpha_i$, passing through $\tilde{C}$ and $\tilde{L}$, respectively. We package all these numbers as sequences $\alpha^C = (\alpha_i^C)_{i \geqslant 1} \in \underline{\mathbf{N}}$ and $\alpha^L = (\alpha_i^L)_{i \geqslant 1} \in \underline{\mathbf{N}}$, and label the corresponding sections as $(Q_{i,j}^C)_{1 \leqslant j \leqslant \alpha_i^C}$ and $(Q_{i,j}^L)_{1 \leqslant j \leqslant \alpha_i^L}$, their respective intersection points with the central fibre as $(q_{i,j}^C)_{1 \leqslant j \leqslant \alpha_i^C}$ and $(q_{i,j}^L)_{1 \leqslant j \leqslant \alpha_i^L}$, and the images by $\pi$ in $L$ of the latter as $(p_{i,j}^C)_{1 \leqslant j \leqslant \alpha_i^C}$ and $(p_{i,j}^L)_{1 \leqslant j \leqslant \alpha_i^L}$. (Note that, by definition, also the $q_{i,j}^C \in \tilde{C}$ are mapped to $L$ by $\pi$).

Similarly, for all $i$ we let $\beta_i^C$ be the number of the sections $R_{i,j}$, $1 \leqslant j \leqslant \beta_i$, passing through $\tilde{C}$. We also let $\beta^C = (\beta_i^C)_{i \geqslant 1} \in \underline{\mathbf{N}}$, label the corresponding sections $(R_{i,j}^C)_{1 \leqslant j \leqslant \beta_i^C}$, and name $(r_{i,j}^C)_{1 \leqslant j \leqslant \beta_i^C}$ their respective points on the central fibre $Y_0$, and $(s_{i,j}^C)_{1 \leqslant j \leqslant \beta_i^C}$ their images on $L$ by $\pi$.

The aware reader will have noticed that we have not considered all possibilities, as there are neither $\alpha^Z$, nor $\beta^Z$ or $\beta^L$. The reason for this is that these turn out not to exist, as we will see in Claims (3.9) and (3.10) below.

**(3.8) New mobile contact points.** We now identify some new points on $\tilde{C}$ which, as we shall see, induce new mobile contact points with $L$ as $C$ moves in the family $V$. These are the intersection points of $\tilde{C}$ with the vertical part of $\pi^*L$, i.e., the part of the divisor $\pi^*L \subseteq \mathcal{Y}$ with support contained in the central fibre $Y_0$, which is

$$(\pi^*L)^{\mathrm{vert}} = \pi^*L - \sum i Q_{i,j} - \sum i R_{i,j}.$$

Thus, for all $i$, we let $\beta_i^{C \cap L}$ be the number of points in $\tilde{C} \cap (\pi^*L)^{\mathrm{vert}}$ appearing with multiplicity $i$ in the divisor $(\pi^*L)^{\mathrm{vert}}$, and label these points as $(r_{i,j}^{C \cap L})_{1 \leqslant j \leqslant \beta_i^{C \cap L}}$. We also let $(s_{i,j}^{C \cap L})_{1 \leqslant j \leqslant \beta_i^{C \cap L}}$ be the images of these points by $\pi$, and $\beta^{C \cap L} = (\beta_i^{C \cap L})_{i \geqslant 1} \in \underline{\mathbf{N}}$.

The reason for the notation is that it will turn out that the support of $(\pi^*L)^{\mathrm{vert}}$ is the union of $\tilde{L}$ and those connected components of $Z$ which intersect $\tilde{L}$ (this follows from Claim (3.9) below and the argument given in (3.19)). Therefore, on the surface $\mathcal{Y}^\flat$ obtained from $\mathcal{Y}$ by contracting $Z$, we can see the points $r_{i,j}^{C \cap L}$ as the intersection points of (the images in $\mathcal{Y}^\flat$ of) $\tilde{C}$ and $\tilde{L}$.

The situation is summarized in Figure 1 below. In fact, we want this figure to picture the situation as it occurs in reality, and therefore we shall make a number of claims before we draw the picture. The proof of these claims occupies a substantial part of the remainder of this section.

The final interpretation in terms of degenerations of Severi varieties of the data introduced above is given in Proposition (3.15) below.

**(3.9) Claim.** The sections $(Q_{i,j})_{1 \leqslant j \leqslant \alpha_i}$ and $(R_{i,j})_{1 \leqslant j \leqslant \beta_i}$ are disjoint from $Z$.

**(3.10) Claim.** All sections $(R_{i,j})_{1 \leqslant j \leqslant \beta_i}$ pass through $\tilde{C}$.

Note that Claims (3.9) and (3.10) imply that

$$\{Q_{i,j}\}_{1 \leqslant i \leqslant \alpha_i} = \{Q_{i,j}^C\}_{1 \leqslant i \leqslant \alpha'_i} \cup \{Q_{i,j}^L\}_{1 \leqslant i \leqslant \alpha''_i} \quad \text{and} \quad \{R_{i,j}\}_{1 \leqslant i \leqslant \alpha_i} = \{R_{i,j}^C\}_{1 \leqslant i \leqslant \alpha'_i}.$$

**(3.11) Claim.** The curve $\tilde{L}$ consists of a unique irreducible component, on which $\pi$ is an isomorphism.

This claim implies that the line $L$ appears with multiplicity 1 in $X_0$, so $X_0 = C \cup L$ with $C$ a degree $d-1$ curve.

**(3.12) Claim.** *The curve $Z$ consists of a disjoint union of chains of rational curves joining $\tilde{L}$ to $\tilde{C}$.*

This Claim implies that the surface $\mathcal{Y}^\flat$ obtained from the smooth surface $\mathcal{Y}$ by contracting $Z$ has only rational double points as singularities.

**(3.13) Claim.** *The map $\pi$ is a birational isomorphism on each irreducible component of $\tilde{C}$.*

This Claim implies that the curve $C$ is reduced.

**(3.14) Claim.** *The curve $\tilde{C}$ is smooth (though not connected, in general).*

We may now draw a picture of the situation.

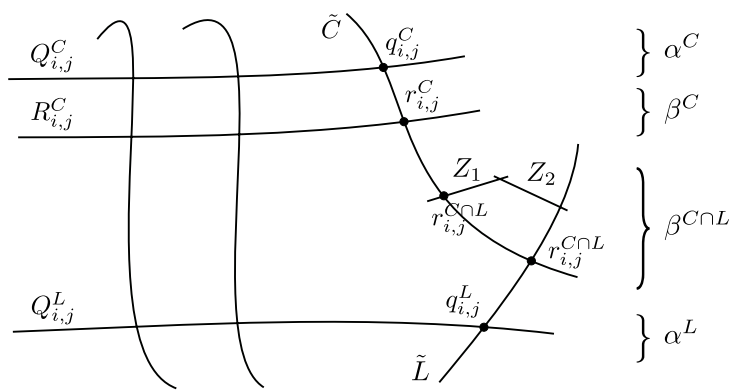

Figure 1: General picture of the degeneration

We shall prove the following result, which encapsulates Theorem (3.1) in the case of an irreducible component having general member of the form $X_0 = C \cup L$.

**(3.15) Proposition.** *In the situation of (3.1), the curve $C$ is a general member of the generalized Severi variety $V^{d-1,\delta'}(\alpha',\beta')(\Omega')$, where:*
— $\alpha' = \alpha^C$, and $\Omega' = (p_{i,j}^C)_{1 \leqslant j \leqslant \alpha_i^C}$;
— $\beta' = \beta^C + \beta^{C \cap L}$;
— $\delta' = \delta - d + 1 + |\beta' - \beta|$.

In other words, we will prove that the assigned contact points for the curves $X_\gamma$ (recall the notation from (3.5)) that end up on $\tilde{C}$ transfer to assigned contact points for $C$, whereas those that end up on $\tilde{L}$ vanish, being accounted for by the component $L$ of $X_0 = C \cup L$; moreover, all mobile contact points of $X_\gamma$ end up on $\tilde{C}$ and thus transfer to mobile contact points on $C$; the curve $C$ has additional mobile contact points with $L$, which come out of the intersection points of $\tilde{C}$ and $\tilde{L}$ in $\mathcal{Y}^\flat$.

We will provide an explicit geometric description of the degeneration of the curves $X_\gamma$ in a neigbourhood of all these points in Subsection 3.4 below. Then, we will be able to track the limits of the $\delta$ nodes of the curves $X_\gamma$, and thus provide a geometric explanation for the formula relating $\delta'$ and $\delta$.

## 3.2 – Simplified proof of Proposition (3.15)

We follow the presentation of [1], and first give a proof of Proposition (3.15) under some simplifying assumptions, in order to describe as clearly as possible what is actually going on in reality. We shall prove later on that these simplifying assumptions do indeed hold, in Subsection 3.3.

We have made six claims in Subsection 3.1 above, with numbers (3.9)–(3.14). Our simplifying assumptions consist in taking the three Claims (3.9), (3.11), and (3.12) for granted, for the moment.

**(3.16) Proof of Proposition (3.15) assuming Claims (3.9), (3.11), and (3.12).** Claim (3.12) ensures that the surface $\mathcal{Y}^\flat$ obtained from $\mathcal{Y}$ by contracting $Z$ has at worst rational double points as singularities, and is still a nodal reduction of the family $\mathcal{X} \to \Gamma$. From now on we work on this new family, and denote any curve in $\mathcal{Y}^\flat$ image of a curve in $\mathcal{Y}$ by the same symbol (with the exception of the central fibre of $\mathcal{Y}^\flat$, which we will denote by $Y_0^\flat$).

Since $Y_0$ and $Y_0^\flat$ are semistable limits of the curves $X_\gamma$, which are smooth of genus $g$, it holds that
$$p_a(Y_0) = p_a(Y_0^\flat) = g.$$

Claim (3.11) ensures that $\tilde{L} \simeq L$ is rational. On the other hand, $Y_0^\flat = \tilde{L} \cup \tilde{C}$ is nodal, and the total number of points in $\tilde{C} \cap \tilde{L}$ is $|\beta^{C \cap L}|$. Thus one finds

$$(3.16.1) \qquad g = p_a(Y_0^\flat) = p_a(\tilde{L} \cup \tilde{C}) = p_a(\tilde{C}) + |\beta^{C \cap L}| - 1.$$

We shall now apply Corollary (2.3), which gives a dimension bound on families of genus $g$ curves and a dimensional characterization of generalized Severi varieties. Let $V$ be the irreducible component of $V^{d,\delta}(\alpha, \beta) \cap p^\perp$ of which $X_0 = C \cup L$ is a general member. The dimension bound tells us that

$$(3.16.2) \qquad \dim(V) \leqslant 2(d-1) + g(\tilde{C}) - 1 + |\beta^C + \beta^{C \cap L}|,$$

since the number of unassigned contact points of $C$ with $L$ is at most $|\beta^C + \beta^{C \cap L}|$, and $C$ has degree at most $d-1$ (indeed, the unassigned contact points of $C$ with $L$ are the points in the support of $\pi_*\tilde{C} \cdot L$ off $\Omega$, one has $\pi_*\tilde{C} \cdot L = \pi_*(C \cdot \pi^*L)$ by the projection formula, and the intersection points of $\tilde{C}$ with the sections $Q_{ij}$ are mapped to $\Omega$; besides, Claim (3.9) ensures that the number $|\beta^C + \beta^{C \cap L}|$ does not count any point in $\tilde{C} \cap \tilde{L} \subseteq Y_0^\flat$ several times).

But $p^\perp$ is a hyperplane which does not contain any irreducible component of $V^{d,\delta}(\alpha, \beta)(\Omega)$, by Proposition (2.2), and since $p \notin \Omega$. Hence

$$(3.16.3) \qquad \dim(V) = \dim\bigl(V^{d,\delta}(\alpha, \beta)\bigr) - 1 = 2d + g - 2 + |\beta|.$$

In the upshot, we have the following sequence of equalities and inequalities:

$$(3.16.4) \quad \begin{aligned} 2d + g + |\beta| - 2 &\leqslant 2(d-1) + g(\tilde{C}) - 1 + |\beta^C + \beta^{C \cap L}| &&\text{by (3.16.2) and (3.16.3)} \\ &\leqslant 2(d-1) + p_a(\tilde{C}) - 1 + |\beta^C + \beta^{C \cap L}| &&\text{since } g(\tilde{C}) \leqslant p_a(\tilde{C}) \\ &= 2(d-1) + g - |\beta^{C \cap L}| + |\beta^C + \beta^{C \cap L}| &&\text{by (3.16.1)} \\ &= 2(d-1) + g + |\beta^C| \\ &\leqslant 2(d-1) + g + |\beta|. \end{aligned}$$

We conclude that equality holds throughout. In particular, (3.16.2) is an equality and therefore, by Corollary (2.3), $V$ is an irreducible component of the Severi variety $V^{d-1,\delta'}(\alpha', \beta')(\Omega')$ where $\alpha'$, $\Omega'$, and $\beta'$ are as stated in Proposition (3.15), and $\delta'$ is determined by the geometric genus of $C$ (note in particular that the orders of contacts of $C$ with $L$ are given by the multiplicities in $\pi^*L \cdot \tilde{C}$, and therefore they are indeed $\alpha'$ and $\beta'$). Yet a couple of words are in order as to why the condition $\mathrm{Card}(C \cap L) = \mathrm{Card}(\pi^{-1}(C \cap L))$ of Corollary (2.3) is indeed verified: on the one hand, it follows from Claim (3.9) that the sections $Q_{i,j}$ are still mutually disjoint in $\mathcal{Y}^\flat$, so that there is a unique point in $\tilde{C}$ over each point of $\Omega'$; on the other hand, the equality in (3.16.2) implies (by Corollary (2.3) again) that $C$ indeed intersects $L$ at $|\beta^C + \beta^{C \cap L}|$ pairwise distinct points, so that eventually all points in $\pi^{-1}(C \cap L)$ have distinct images in $L$, as required.

To finish, let us explicitly derive the formula for $\delta'$: it follows from the series of equalities in (3.16.4) that the geometric genus of $C$ is

$$g(C) = g(\tilde{C}) = p_a(\tilde{C}) = g - |\beta^{C \cap L}| + 1,$$

so that

$$\begin{aligned}\delta' &= \binom{d-2}{2} - g(C) = \binom{d-2}{2} - g + |\beta^{C \cap L}| - 1 \\ &= \binom{d-2}{2} - \binom{d-1}{2} + \delta + |\beta^{C \cap L}| - 1 = \delta - (d-1) + |\beta^{C \cap L}|,\end{aligned}$$

which is the value asserted in Proposition (3.15), as $\beta^{C \cap L} = \beta' - \beta$, since $\beta^C = \beta$ by the series of equalities in (3.16.4). □

**(3.17) Proof of Claims (3.10), (3.13), and (3.14) assuming Claims (3.9), (3.11), and (3.12).** As we have already observed, it follows from the series of equalities in (3.16.4) that $\beta = \beta^C$, which proves Claim (3.10).

On the other hand, the fact that equality holds in (3.16.2) implies by Corollary (2.3) that $C$ is reduced. Thus, Claim (3.13) holds. Eventually, we have already observed that the series of equalities in (3.16.4) implies the equality $g(\tilde{C}) = p_a(\tilde{C})$, which proves that $\tilde{C}$ is smooth, i.e., Claim (3.14) holds. □

## 3.3 – Complete proof of Proposition (3.15)

**(3.18) Labelling of the remaining sections.** In order to prove Claims (3.9), (3.11), and (3.12), we need to consider the situation as complicated as it could a priori be, and therefore to introduce some more notation with respect to that already given in Subsection 3.1.

Thus, for all $i$ we define (in the surface $\mathcal{Y}$) $\beta_i^{C,Z}$ as the number of the sections $(R_{i,j})_{1 \leqslant j \leqslant \beta_i}$ passing through a connected component of $Z$ that intersects $\tilde{C}$, and package these numbers in the sequence $\beta^{C,Z} = (\beta_i^{C,Z})$.

On the other hand we redefine $|\beta^{C \cap L}|$ as the sum of the number of points of $\tilde{C}$ meeting $\tilde{L}$ plus the number of connected components of $Z$ meeting both $\tilde{C}$ and $\tilde{L}$. Thus it is the number of intersection points between $\tilde{C}$ and $\tilde{L}$ after we contract $Z$. Note that $|\beta^{C \cap L}|$ no longer corresponds to an actual sequence $\beta^{C \cap L} = (\beta_i^{C \cap L})$; the point is, without Claim (3.12), there could be a connected component of $Z$ meeting $\tilde{C}$ at several distinct points with different multiplicities in $(\pi^* L)^{\text{vert}}$.

**(3.19) Bounding the number of mobile contact points of $C$ and $L$.** To this end, we note that every connected component of $Z$ mapping to a point of $L$ by $\pi$ necessarily intersects either $\tilde{L}$, or a section among $Q_{i,j}$ and $R_{i,j}$. Indeed, otherwise there would exist a connected component of $\pi^{-1}(L)$ entirely contained in $Z$, and thus after contracting $Z$ we would obtain a surface with an isolated point in the preimage of $L$, which is impossible.

Therefore the number of mobile contact points of $C$ and $L$ (in $\mathbf{P}^2$) is bounded from above as follows:

$$(3.19.1) \qquad \operatorname{Card}\bigl(C \cap (L - \Omega)\bigr) \leqslant |\beta^C + \beta^{C,Z}| + |\beta^{C \cap L}|.$$

**(3.20) Bounding the genus of** $C$**.** As in (3.16), $Y_0$ is a semistable limit of the curves $X_\gamma$, and therefore $p_a(Y_0) = g$. Since the number of intersection points between $\tilde{C}$ and $\tilde{L}$ after we contract $Z$ is $|\beta^{C \cap L}|$, we find that

$$(3.20.1) \qquad g = p_a(Y_0) \geqslant p_a(\tilde{L}) + p_a(\tilde{C}) + |\beta^{C \cap L}| - 1.$$

Note that this inequality is strict if some connected component of $Z$ has positive arithmetic genus, or intersects $\tilde{C}$ or $\tilde{L}$ at more than one point.

**(3.21) Proof of Claim (3.11).** Let $e$ be the multiplicity of $L$ in $X_0 = C \cup L$. Then $C$ has degree $d - e$. On the other hand, $\tilde{L}$ is mapped $e : 1$ to $L$ by $\pi$, so by Riemann–Hurwitz we have the inequality $p_a(\tilde{L}) \geqslant -e + 1$. Together with (3.20.1), this gives the bound:

$$(3.21.1) \qquad p_a(\tilde{C}) \leqslant g - |\beta^{C \cap L}| + e.$$

Now, we argue exactly as in (3.16): by Corollary (2.3), we find that the irreducible component $V$ of $V^{d,\delta}(\alpha, \beta) \cap p^\perp$ that $C$ moves in has dimension at most

$$(3.21.2) \qquad \begin{aligned} 2(d-e) + p_a(\tilde{C}) - 1 + |\beta^C + \beta^{C,Z}| + |\beta^{C \cap L}| &\leqslant 2(d-e) + g + |\beta^C + \beta^{C,Z}| + e - 1 \\ &= 2d + g + |\beta^C + \beta^{C,Z}| - e - 1. \end{aligned}$$

But since

$$\dim(V) = \dim\bigl(V^{d,\delta}(\alpha, \beta)\bigr) - 1 = 2d + g + |\beta| - 2,$$

and $|\beta| \leqslant |\beta^C + \beta^{C,Z}|$, one necessarily has $e = 1$. $\square$

**(3.22) Proof of Claim (3.12).** By the fact that (3.21.2) is in fact an equality, we see that also (3.16.1) must be an equality. Thus, as already observed in (3.20), every connected component of $Z$ has arithmetic genus $0$, i.e., it is a tree of rational curves, and intersects each of $\tilde{L}$ and $\tilde{C}$ at most once.

In order to prove that such a connected component is a chain of rational curves connecting $\tilde{L}$ and $\tilde{C}$, it is therefore sufficient to prove that any "end component" of a connected component of $Z$ meets either $\tilde{C}$ or $\tilde{L}$. By an end component, we mean an irreducible component $Z_1$ that meets at most one other irreducible component of $Z$ (one if the corresponding connected component of $Z$ has several irreducible components, zero if $Z_1$ is itself a connected component of $Z$).

The key argument is that if an end component $Z_1$ of $Z$ intersects neither $\tilde{C}$ nor $\tilde{L}$, then by the minimality of $\mathcal{Y}$ there must be at least two of the sections $(Q_{i,j})_{1 \leqslant j \leqslant \alpha_i}$ and $(R_{i,j})_{1 \leqslant j \leqslant \beta_i}$ meeting $Z_1$, for otherwise we could blow down $Z_1$ in $\mathcal{Y}$ and still satisfy all the conditions imposed on $\mathcal{Y}$. (Note that $Z_1$, being a rational curve intersecting $Y_0 - Z_1$ in exactly one point since $Y_0$ is connected, is a $(-1)$-curve).

These two sections cannot be both of type $\alpha$, for otherwise they would be contracted by $\pi$ to the same point of $L$, which is excluded by the definition of generalized Severi varieties. To rule out the other possibilities, we use the fact that since (3.21.2) is an equality, then also (3.19.1) must be an equality. But if two sections either both of type $\beta$, or of types $\alpha$ and $\beta$ respectively, meet the same connected component of $Z$, then the two corresponding points of $C$ in $\mathbf{P}^2$ are the same, hence (3.19.1) is a strict inequality. This ends the proof of Claim (3.12). $\square$

**(3.23) Proof of Claim (3.9).** Now we know that Claim (3.12) holds, i.e., the connected components of $Z$ are chains of rational curves connecting $\tilde{C}$ and $\tilde{L}$. Therefore we can indeed consider the points $(r_{i,j}^{C \cap L})_{1 \leqslant j \leqslant \beta_i^{C \cap L}}$ introduced in (3.8), and their respective images by $\pi$, the points $(s_{i,j}^{C \cap L})_{1 \leqslant j \leqslant \beta_i^{C \cap L}}$ (see also (3.18)).

Our main argument here is that (3.19.1) is an equality, as we have already observed in (3.22) above. If a section of type $\alpha$ intersects a connected component of $Z$, then the image in $\mathbf{P}^2$ of the point $r_{i,j}^{C \cap L}$ corresponding to this component must lie in $\Omega$, which implies that (3.19.1) is a strict inequality, a contradiction.

Similarly, if a section of type $\beta$ intersects a connected component of $Z$, then the point $r_{i,j}^{C,Z}$ corresponding to this section and the point $r_{i,j}^{C \cap L}$ corresponding to the component of $Z$ have the same image in $\mathbf{P}^2$, hence (3.19.1) is a strict inequality, a contradiction. □

We have now proved all three Claims (3.9), (3.11), and (3.12). In Subsection 3.2 we had proved Proposition (3.15) assuming these three claims. Thus, we now have a complete proof of Proposition (3.15), and the proof of Theorem (3.1) is over.

### 3.4 – Local picture of the degeneration of curves

We consider the situation set up in (3.1) and (3.5), when a general member $X_\gamma$ of $V^{d,\delta}(\alpha,\beta)$ degenerates to a curve $X_0 = L \cup C$. We shall now give a complete description of this degeneration, in its concrete incarnation as a degeneration of plane curves. This will be useful in the next section, when we prove the theorem going the other direction than Theorem (3.1) above.

We carry on with the notation introduced above for Proposition (3.15). Something nontrivial happens only at the points where $C$ and $L$ intersect, which come in the following four types: (3.24)–(3.27).

**(3.24) Vanishing assigned contact points.** These are the points of $\Omega \setminus \Omega'$, i.e., the points $p_{i,j}^L$. Since $C$ is general in the component $V$ of $V^{d-1,\delta'}(\alpha',\beta')(\Omega')$, by Proposition (2.2) it does not pass through any point $p_{i,j}^L$. Thus the curve $X_0$ is smooth at those points, i.e., it has only one local branch there, which is an open subset of the line $L$ itself.

Therefore the degeneration happens like this: locally at a point $p_{i,j}^L$, the curves $X_\gamma$, $\gamma \neq 0$, have one smooth local branch tangent to the order $i$ with $L$, which degenerates in $X_0$ to an open neighbourhood of $p_{i,j}^L$ in $L$.

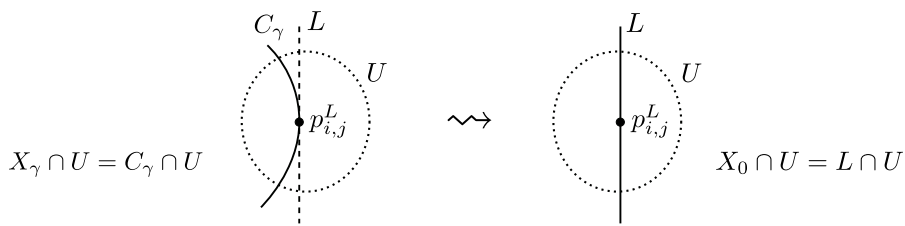

Figure 2: Degeneration in a neighbourhood $U$ of a vanishing assigned contact point

**(3.25) Non-vanishing assigned contact points.** These are the points of $\Omega'$, i.e., the points $p_{i,j}^C$. Such a point is by definition the image of the point $q_{i,j}^C \in \tilde{C}$, and $C$ is smooth at $p_{i,j}^C$ by Proposition (2.2). On the other hand $p_{i,j}^C$ also has a preimage in $\tilde{L}$, and exactly one as $\tilde{L} \to L$ is an isomorphism: let's call it $(q_{i,j}^C)'$. It follows that $X_0$ has two smooth local branches at $p_{i,j}^C$, tangent to $i$-th order, so that $p_{i,j}^C$ is an $i$-tacnode of $X_0$.

As we have just seen, the map $Y_0 \to X_0$ factors through the normalization of $p_{i,j}^C$ in $X_0$ (recall that $Y_0$ is the central fibre of the semistable reduction introduced in (3.5)), and thus the family $\mathcal{Y}$ is smooth in a neighbourhood of the two preimages of $p_{i,j}^C$ in $Y_0$. This implies that the curves $X_\gamma$, as they approach $X_0$, have two local branches analytically locally around $p_{i,j}^C \in \mathbf{P}^2$.

The difference in arithmetic genus between $Y_b$ and $X_\gamma$ over a neighbourhood of $p_{i,j}^C$ in $\mathbf{P}^2$ is constant as $b$ moves (with $\gamma = \nu(\eta(b))$, see diagram (3.5.1)), because both $\mathcal{Y}$ and $\mathcal{X}$ are flat families. Moreover, $\mathcal{X}$ is equigeneric over $\Gamma \setminus \{0\}$, and $\mathcal{Y}$ is a nodal reduction of $\mathcal{X}$, so for $b \neq 0$, the difference in arithmetic genus between $Y_b$ and $X_\gamma$ is just the sum of the $\delta$-invariants of the singularities of $X_\gamma$; since the curves $X_\gamma$ are nodal for $\gamma \neq 0$, this is therefore the number of nodes of $X_\gamma$ tending to the point $p_{i,j}^C$ as $\gamma$ tends to $0$. On the other hand, over a neighbourhood of $p_{i,j}^C$, $Y_0 \to X_0$ is the normalization of an $i$-tacnode, so the difference in arithmetic genus between $Y_0$ and $X_0$ over a neighbourhood of $p_{i,j}^C$ equals $i$.

The upshot is that locally at $p_{i,j}^C$, the curves $X_\gamma$, $\gamma \neq 0$, have two local branches intersecting transversely in $i$ points tending to the $i$-tacnode of $X_0$ at $p_{i,j}^C$; one of the two branches of $X_\gamma$ is tangent to order $i$ with $L$ at the fixed point $p_{i,j}^C \in \Omega' \subseteq \Omega$ (this is the branch corresponding to the branch of $Y_b$ intersecting the section $Q_{i,j}^C$, namely that tending to the branch of $Y_0$ contained in $\tilde{C}$).

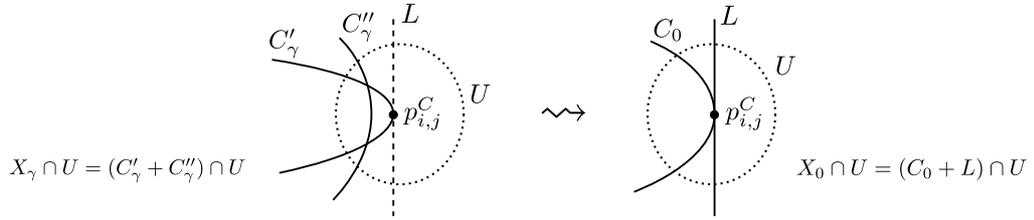

Figure 3: Degeneration in a neighbourhood $U$ of a non-vanishing assigned contact point

**(3.26) Limits of unassigned contact points.** These are the points $s_{i,j}^C$. At those points the situation is almost exactly the same as that at the points $p_{i,j}^C$ described above, the key point being that the points $s_{i,j}^C \in L$ have two points in their preimage in $Y_0$, namely $r_{i,j}^C \in \tilde{C}$ and an additional point on $\tilde{L}$, which we may call $(r_{i,j}^C)'$ so that, as in (3.25), the map $Y_0 \to X_0$ factors through the normalization of the $i$-tacnode of $X_0$ at $s_{i,j}^C$. The only difference is that here the branch of $X_\gamma$ which is tangent to order $i$ with $L$ touches $L$ at a mobile point, different from $s_{i,j}^C$ by Proposition (2.2) (this is the branch corresponding to the branch of $Y_b$ intersecting the section $R_{i,j}^C$, which tends to a local branch of $\tilde{C}$).

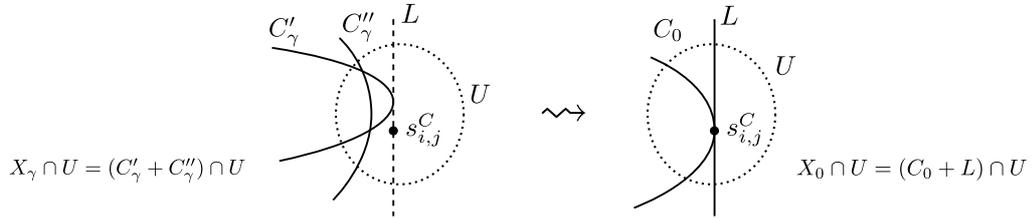

Figure 4: Degeneration in a neighbourhood $U$ of a limit of unassigned contact points

**(3.27) New unassigned contact point.** These are the points $s_{i,j}^{C \cap L}$. Such a point is the image of a point $r_{i,j}^{C \cap L}$ which, in the surface $\mathcal{Y}^\flat$ obtained from $\mathcal{Y}$ by contracting $Z$, lies at the intersection of $\tilde{C}$ and $\tilde{L}$; it is thus an ordinary node of the central fibre $Y_0^\flat$ (in fact, seen in $\mathcal{Y}$ it is a node as well, although there it may lie on $\tilde{C} \cap Z$ rather than on $\tilde{C} \cap \tilde{L}$), and the point $s_{i,j}^{C \cap L}$ therefore has only one point in its preimage in $Y_0$. Since the nearby fibres $Y_b$ are smooth (because $\mathcal{X}$ is equigeneric over $\Gamma \setminus \{0\}$, and $\mathcal{Y}$ is a nodal reduction of $\mathcal{X}$), they have only one branch in an analytic neighbourhood of the preimage of $s_{i,j}^{C\cap L}$ in $\mathcal{Y}$, and therefore the curves $X_\gamma$, $\gamma \neq 0$, have only one local branch in a neighbourhood of $s_{i,j}^{C\cap L}$ in $\mathbf{P}^2$.

The points $s_{i,j}^{C\cap L}$ are again $i$-tacnodes of the curve $X_0$, this time because both $L$ and $C$ are smooth at $s_{i,j}^{C\cap L}$ by Proposition (2.2), and intersect there with multiplicity $i$ because by definition $s_{i,j}^{C\cap L}$ lies on a component of multiplicity $i$ of $(\pi^*L)^{\mathrm{vert}}$. In particular the difference in arithmetic genus between $Y_0$ and $X_0$ over $s_{i,j}^{C\cap L}$ is $i-1$. Since this difference for the families $\mathcal{Y}$ and $\mathcal{X}$ is constant over a neighbourhood of $s_{i,j}^{C\cap L}$, and the curves $X_\gamma$ are nodal for $\gamma \neq 0$, we see that the latter will have $i-1$ nodes tending to $s_{i,j}^{C\cap L}$.

To sum up, analytically locally around $s_{i,j}^{C\cap L}$, the curves $X_\gamma$, $\gamma \neq 0$, have one irreducible, $(i-1)$-nodal local branch, transverse to $L$ and not passing by $s_{i,j}^{C\cap L}$, degenerating to the $i$-tacnode at $s_{i,j}^{C\cap L}$ formed by $C$ and $L$.

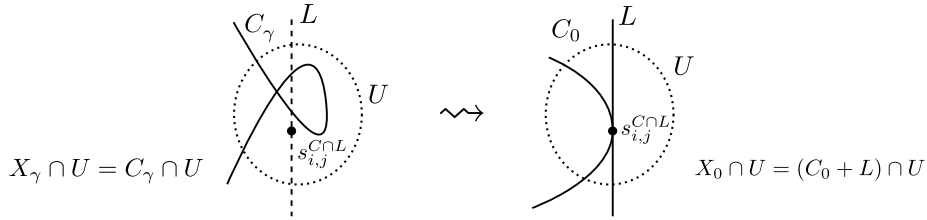

Figure 5: Degeneration in a neighbourhood $U$ of a new unassigned contact point

We emphasize that this phenomenon is arguably the central topic of this whole volume; its local model has been studied in detail in the previous chapters [IV] and [A].

Whereas the above description of the degeneration will be needed in the next section, the remainder of this subsection is not necessary for the proof of the main theorem. I believe it is nevertheless helpful.

First of all, having the above descriptions at hand, we may now give the promised geometric interpretation for the relation between $\delta$ and $\delta'$ in Theorem (3.1) and Proposition (3.15).

**(3.28) Book-keeping of the nodes in the degeneration.** As a result of (3.24)–(3.27) above, we see that as $X_\gamma$ degenerates to $X_0 = C \cup L$ there are $i$ nodes of $X_\gamma$ tending to each point of types $\alpha_i^C$ and $\beta_i^C$ (respectively the preserved tangencies at assigned points, and the limits of tangencies at unassigned points, which are all preserved), and $i-1$ nodes of $X_\gamma$ tending to each point of type $\beta_i^{C\cap L}$ (the new tangency conditions, all at unassigned points); the degeneration is otherwise equisingular. Therefore, we have:
$$\delta - \delta' = I\alpha^C + I\beta^C + \left(I\beta^{C\cap L} - |\beta^{C\cap L}|\right).$$
On the other hand we know that, in the notation of Proposition (3.15), $\alpha' = \alpha^C$ and $\beta' = \beta^C + \beta^{C\cap L}$, so that
$$I\alpha^C + I\beta^C + I\beta^{C\cap L} = I\alpha' + I\beta' = C \cdot L = d-1;$$
since moreover $\beta = \beta^C$, we also have $|\beta^{C\cap L}| = |\beta' - \beta|$. Finally, we find
$$\delta - \delta' = (d-1) - |\beta' - \beta|$$
as required.

Lastly, let us analyze the multiplicities of the irreducible components of $Y_0$ in the divisor $\pi^*L$. This will be helpful in understanding the multiplicities of the components of the type $V^{d-1,\delta'}(\alpha', \beta')$ in the intersection of $V^{d,\delta}(\alpha, \beta)$ with the hyperplane $p^\perp$, cf. Theorem (4.1).

**(3.29) Multiplicities of the divisor $\pi^*L$ in $\mathcal{Y}$.** Recall that the central fibre $Y_0$ of $\mathcal{Y}$ is reduced and consists of (i) $\tilde{L}$, an irreducible curve mapped isomorphically to $L$, (ii) $\tilde{C}$, the normalizations of the various irreducible components of $C$, and (iii) $Z$, the union of disjoint chains of rational curves joining $\tilde{L}$ to $\tilde{C}$.

Let $m$ be the multiplicity of $\tilde{L}$ in the divisor $\pi^*L$ on $\mathcal{Y}$. Let us remind that the integer $m$ records the ramification of $\pi : \mathcal{Y} \to \mathbf{P}^2$ along $\tilde{L}$ in the direction normal to $\tilde{L}$, hence $m > 1$ does not contradict $\pi|_{\tilde{L}}$ being an isomorphism onto $L$; instead, if $\pi|_{\tilde{L}} : \tilde{L} \to L$ had degree $k$, this would translate into $\pi_*\tilde{L} = kL$.

The main output of this paragraph will be the divisibility of $m$ by $\mathrm{lcm}(\beta^{C\cap L}) = \mathrm{lcm}(\beta'-\beta)$. To see this, consider a point $r_{i,j}^{C\cap L} \in \tilde{C}$. By definition, see Paragraph (3.8), $r_{i,j}^{C\cap L}$ sits both on $\tilde{C}$ and on a multiplicity $i$ component of $(\pi^*L)^{\mathrm{vert}}$, and we have seen that $(\pi^*L)^{\mathrm{vert}}$ is supported on $\tilde{L}+Z$. Thus, if $i < m$, then $r_{i,j}^{C\cap L}$ sits on a smooth irreducible rational curve $Z_1 \subseteq Z$ which has multiplicity $i$ in $\pi^*L$. Let $Z_2$ be the next curve in the chain of rational curves joining $\tilde{C}$ to $\tilde{L}$ and containing $Z_1$ (here, 'next' is intended relative to the direction pointing to $\tilde{L}$; if $Z_1$ directly joins $\tilde{C}$ to $\tilde{L}$, then we take $Z_2 = \tilde{L}$), and let $i_2$ be its multiplicity in $\pi^*L$. Since $\pi$ contracts $Z_1$, we have $\pi^*L \cdot Z_1 = 0$, hence

$$(iZ_1 + i_2Z_2) \cdot Z_1 = 0.$$

Since $Z_1$ is a $(-2)$-curve[2] and $Z_1 \cdot Z_2 = 1$, we find that $i_2 = 2i$. We continue by induction to compute the successive multiplicities: the connected component of $Z$ containing $Z_1$ is a chain of $(-2)$-curves; let us denote them by $Z_1,\ldots,Z_l$, in that order, and let $Z_{l+1} = \tilde{L}$. Let $i_1,\ldots,i_{l+1}$ be their multiplicities in $\pi^*L$. Once we know, for some $k \leqslant l$, that $i_s = si$ for all $s = 1,\ldots,k$, we have

$$\pi^*L \cdot Z_k = (i_{k-1}Z_1 + i_kZ_k + i_{k+1}Z_{k+1}) \cdot Z_k = 0,$$

hence $i_{k+1} = 2i_k - i_{k-1} = (k+1)i$, as required. Finally, the multiplicity $m$ of $\tilde{L} = Z_{l+1}$ in $\pi^*L$ equals $(l+1)i$; thus, $m$ is divisible by $i$, and the length of the chain of $(-2)$-curves is $l = m/i - 1$.

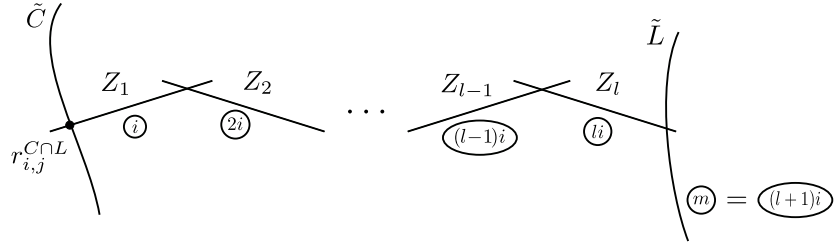

Figure 6: Multiplicities of the components of $(\pi^*L)^{\mathrm{vert}}$

The same argument shows that there are no points $r_{i,j}^{C\cap L}$ with $i > m$, and all points $r_{m,j}^{C\cap L}$ appear directly on $\tilde{L}$.

In conclusion, for all $i$ such that $\beta_i^{C\cap L} = (\beta'-\beta)_i > 0$, the multiplicity $m$ is divisible by $i$; equivalently, $m$ is divisible by $\mathrm{lcm}(\beta'-\beta)$.

## 4 − Deformations of generalized Severi varieties

In this section we prove the result going the other direction than Theorem (3.1). It asserts that all the generalized Severi varieties which, according to Theorem (3.1), may appear in the hyperplane section $V^{d,\delta}(\alpha,\beta) \cap p^\perp$, do indeed appear.

---

[2]this is well-known and may be seen as follows: all fibres of $\mathcal{Y}$ are linearly equivalent, and have intersection number 0 with $Z_1$; therefore, $Y_0 \cdot Z_1 = (\tilde{C} + Z_1 + Z_2) \cdot Z_1 = 0$, hence $Z_1^2 = -(\tilde{C}+\tilde{L}) \cdot Z_1 = -2$.

**(4.1) Theorem** ([1, Theorem 1.3]). *Let $p$ be a general point of $L$.*

*a) Let $V'$ be an irreducible component of $V^{d,\delta}(\alpha + e_k, \beta - e_k)(\Omega \cup \{p\})$ as in part a) of Theorem (3.1). Then $V'$ is a component of the intersection $V^{d,\delta}(\alpha,\beta)(\Omega) \cap p^\perp$, and at a general point of $V'$, the variety $V^{d,\delta}(\alpha,\beta)$ is smooth and has intersection multiplicity $k$ with $p^\perp$ along $V'$.*

*b) Let $V'$ be an irreducible component of $V^{d-1,\delta'}(\alpha',\beta')(\Omega')$ as in part b) of Theorem (3.1). Then $V'$ is a component of the intersection $V^{d,\delta}(\alpha,\beta)(\Omega) \cap p^\perp$ and, at a general point of $V'$, the variety $V^{d,\delta}(\alpha,\beta)$ has $\binom{\beta'}{\beta}I^{\beta'-\beta}/\mathrm{lcm}(\beta'-\beta)$ local sheets, each of which has intersection multiplicity $\mathrm{lcm}(\beta'-\beta)$ with $p^\perp$ along $V'$. Moreover, each of these sheets has the generic point of $V'$ as a point of multiplicity $\mathrm{lcm}(\beta'-\beta)/\max(\beta'-\beta)$,*

In both cases a) and b) the general strategy will be to show that the general member of $V'$ may be deformed to general members of $V^{d,\delta}(\alpha,\beta)(\Omega)$. We treat the two cases separately, in Subsections 4.1 and 4.2 respectively.

The multiplicity in case a) displays a phenomenon reminiscent of the fundamental principle of projective duality. The deformation argument for case b) is much more demanding, and eventually relies on the local material on deformations of tacnodes gathered in Chapters [IV] and [V], which is arguably the keystone of the present volume.

## 4.1 – Fixing an unassigned contact point

This section is devoted to the proof of part a) of Theorem (4.1), which is the case when $V'$ parametrizes curves that do not contain the line $L$, and for which the condition that they lie on the hyperplane $p^\perp$ is accounted for by the fact that $p$ is one of the assigned contact points.

It is fairly clear that $V^{d,\delta}(\alpha + e_k, \beta - e_k)(\Omega \cup \{p\})$ is contained in the hyperplane section of $V^{d,\delta}(\alpha,\beta)(\Omega)$ by $p^\perp$, so the point here is to establish the assertion on the smoothness and intersection multiplicity.

The multiplicity is as one would expect it to be from the point of view of projective duality. Indeed, the latter tells us the following. Consider a variety $X \subseteq \mathbf{P}^N$, and its dual $X^\vee \subseteq \check{\mathbf{P}}^N$ in the dual projective space, which parametrizes hyperplanes in $\mathbf{P}^N$ tangent to $X$. Then for a general point $p \in X$, the hyperplane $p^\perp \subseteq \check{\mathbf{P}}^N$ (consisting of those hyperplanes of $\mathbf{P}^N$ that contain the point $p$) is tangent to $X^\vee$ at the points $T^\perp \in \check{\mathbf{P}}^N$ corresponding to hyperplanes that are tangent to $X$ at $p$. More generally, one may consider "osc-dual" varieties $X^{\vee_\beta} \subseteq \check{\mathbf{P}}^N$ of $X$, parametrizing hyperplanes of $\mathbf{P}^N$ osculating $X$ to the order $\beta$, and then $p^\perp$ osculates $X^{\vee_\beta}$ to the order $\beta$ at points $O^\perp$ corresponding to osculating hyperplanes to $X$ at $p$. Some of these are discussed in Chapter [**??**].

In the situation under consideration, this has to be applied to the case when $X$ is the degree $d - I\alpha$ rational normal curve image of $L$ by the linear system cut out by degree $d$ plane curves with the contact conditions at assigned points determined by $\alpha$ and $\Omega$. In fact this is exactly what we do in practice, following the approach of [1, §4.3], and establishing along the way the appropriate biduality principle for osculating hyperplanes to rational normal curves, which is the following.

**(4.2) Proposition** ([1, Lemma 4.7]). *Consider the line $L \subseteq \mathbf{P}^2$, and a general point $p$ of $L$. We assume the data of $d, \alpha, \beta$, and $\Omega = (p_{i,j})_{1 \leqslant j \leqslant \alpha_i}$ is given as in Theorem (4.1); in particular, $I\alpha + I\beta = d$. We define the following three loci in $|\mathcal{O}_L(d)|$:*

- $p^\perp = \{D : p \in D\}$;
- $\Phi_{\alpha,\beta}(\Omega) = \{D = \sum_{1 \leqslant j \leqslant \alpha_i} i \cdot p_{i,j} + \sum_{1 \leqslant j \leqslant \beta_i} i \cdot s_{i,j} \quad \text{for some } s_{i,j}\text{'s in } L\}$;
- $\Psi_{\alpha,\beta}(\Omega) = \{D = k \cdot p + \sum_{1 \leqslant j \leqslant \alpha_i} i \cdot p_{i,j} + \sum_{1 \leqslant j \leqslant (\beta - e_k)_i} i \cdot s_{i,j} \quad \text{for some } s_{i,j}\text{'s in } L\}$.

*In a neighbourhood of a general point $[D_0] \in \Psi_{\alpha,\beta}(\Omega)$, the variety $\Phi_{\alpha,\beta}(\Omega)$ is smooth and has intersection multiplicity $k$ with the hyperplane $p^\perp$ along $\Psi_{\alpha,\beta}(\Omega)$.*

*Proof.* We will proceed by an explicit local computation. Let $x$ be an affine coordinate on $L$ centered at the point $p$. Let the points $p_{i,j}$ have respective coordinates $\lambda_{i,j}$, and let the points $s_{i,j}^0$ such that
$$D_0 = k \cdot p + \sum_{1 \leq j \leq \alpha_i} i \cdot p_{i,j} + \sum_{1 \leq j \leq (\beta - e_k)_i} i \cdot s_{i,j}^0$$
have respective coordinates $\mu_{i,j}$. By generality of $D_0$, we may assume that these coordinates are altogether mutually distinct.

Then we can parametrize a neighbourhood of $[D_0]$ in $\Phi_{\alpha,\beta}(\Omega)$ by
$$(\varepsilon, \varepsilon_{i,j}) \mapsto [f(x)] = \left[ (x - \varepsilon)^k \cdot \prod_{1 \leq j \leq \alpha_i} (x - \lambda_{i,j})^i \cdot \prod_{1 \leq j \leq (\beta - e_k)_i} (x - \mu_{i,j} - \varepsilon_{i,j})^i \right] \in \mathbf{C}_d[x],$$
from which we see first of all that $\Phi_{\alpha,\beta}(\Omega)$ is smooth at $[D_0]$. On the other hand, writing
$$f(x) = x^d + b_{d-1} x^{d-1} + \cdots + b_1 x + b_0 \in \mathbf{C}_d[x],$$
we get a system of affine coordinates on some chart of $|\mathcal{O}_L(d)|$, in which the hyperplane $p^\perp$ is defined by the linear equation $b_0 = 0$. Then, in the local coordinates $(\varepsilon, \varepsilon_{i,j})$ on $\Phi_{\alpha,\beta}(\Omega)$ around $[D_0]$, the intersection with $p^\perp$ is defined by the equation
$$\prod \lambda_{i,j}^i \cdot \varepsilon^k \cdot \prod (\varepsilon_{i,j} + \mu_{i,j})^i = 0 \iff \varepsilon^k \cdot \prod (\varepsilon_{i,j} + \mu_{i,j})^i = 0,$$
hence, locally around $[D_0]$, it is $k$ times the divisor defined by the equation $\varepsilon = 0$, and the latter divisor is exactly $\Psi_{\alpha,\beta}(\Omega)$. □

In order to relate the local geometry of generalized Severi varieties to the model situation of Proposition (4.2), we consider the restriction map $\pi : |\mathcal{O}_{\mathbf{P}^2}(d)| \dashrightarrow |\mathcal{O}_L(d)|$. Note that this is merely the projection from the codimension $d+1$ linear subspace $L^\perp$ of $|\mathcal{O}_{\mathbf{P}^2}(d)|$ parametrizing curves wich contain $L$. Finally, observe that the hyperplane $p^\perp$ of $|\mathcal{O}_{\mathbf{P}^2}(d)|$ is the preimage by $\pi$ of the hyperplane $p^\perp$ of $|\mathcal{O}_L(d)|$.

Let $W$ be the plain Severi variety $V^{d,\delta}$, and denote by $\pi_W$ the restriction of $\pi$ to $W$. By definition, using the notation of Proposition (4.2), we have
$$V^{d,\delta}(\alpha, \beta)(\Omega) = \pi_W^{-1}\bigl(\Phi_{\alpha,\beta}(\Omega)\bigr) \quad \text{and} \quad V^{d,\delta}(\alpha + e_k, \beta - e_k)(\Omega \cup \{p\}) = \pi_W^{-1}\bigl(\Psi_{\alpha,\beta}(\Omega)\bigr).$$

Therefore, part a) of Theorem (4.1) follows directly from the following proposition together with Proposition (4.2) above.

**(4.3) Proposition.** *Let $X_0$ be a general member of (any irreducible component of) $V^{d,\delta}(\alpha + e_k, \beta - e_k)(\Omega \cup \{p\})$. The differential of the map $\pi_W : W \dashrightarrow |\mathcal{O}_L(d)|$ is surjective at the point $[X_0]$.*

*Proof.* First note that, being general, $X_0$ does not contain $L$ and therefore $\pi$ is well defined at $[X_0]$. Since the map $\pi$ is a linear projection, it suffices to show that the projective tangent space $\mathbf{T}_{[X_0]} W$ to $W = V^{d,\delta}$ at the point $[X_0]$ intersects the center of the projection transversely, i.e., that the intersection $\mathbf{T}_{[X_0]} W \cap L^\perp$ has codimension $d+1$ in $\mathbf{T}_{[X_0]} W$.

Now, the projective tangent space $\mathbf{T}_{[X_0]} W$ identifies with the linear series of degree $d$ curves passing through the nodes of $X_0$, and thus its intersection with $L^\perp$ identifies with the linear series of degree $d - 1$ curves passing through the nodes of $X_0$. By Lemma (4.4), these two linear series have dimensions $\frac{d(d+3)}{2} - \delta$ and $\frac{(d-1)(d+2)}{2} - \delta$ respectively, so the codimension of $\mathbf{T}_{[X_0]} W \cap L^\perp$ in $\mathbf{T}_{[X_0]} W$ equals that of $|\mathcal{O}_{\mathbf{P}^2}(d-1)|$ in $|\mathcal{O}_{\mathbf{P}^2}(d)|$, which is $d+1$, as we wanted. □

The following lemma is the classical fact that, on a regular surface, the adjoint series cut out the complete canonical series. We include the proof for completeness.

**(4.4) Lemma.** *Let $C \subseteq \mathbf{P}^2$ be a $\delta$-nodal curve of degree $d$. For all $e \geqslant d - 3$, the nodes of $C$ impose independent conditions on curves of degree $e$.*

*Proof.* It is sufficient to prove that the nodes of $C$ impose independent conditions on curves of degree $d - 3$. Let $\mathcal{I} \subseteq \mathcal{O}_{\mathbf{P}^2}$ be the ideal sheaf of the nodes of $C$ as a subscheme of $\mathbf{P}^2$, and let $\mathcal{A} \subseteq \mathcal{O}_C$ be the conductor ideal sheaf of $C$ (with respect to its normalization $\bar{C}$). The ideal $\mathcal{A}$ is the restriction of $\mathcal{I}$ to $C$ (see, e.g., [II, Lemma (2.5)]), and we have an exact sequence

$$0 \longrightarrow \mathcal{O}_{\mathbf{P}^2}(-3) \longrightarrow \mathcal{O}_{\mathbf{P}^2}(d-3) \otimes \mathcal{I} \longrightarrow \mathcal{O}_C(d-3) \otimes \mathcal{A} \longrightarrow 0.$$

It follows that there is an isomorphism $H^0(\mathcal{O}_{\mathbf{P}^2}(d-3) \otimes \mathcal{I}) \cong H^0(\mathcal{O}_C(d-3) \otimes \mathcal{A})$. Now, by [II, Corollary (2.4)], $H^0(\mathcal{O}_C(d-3) \otimes \mathcal{A})$ is isomorphic to $H^0(\bar{C}, \omega_{\bar{C}})$, where $\bar{C}$ still denotes the normalization of $C$, hence it has dimension $g$, the geometric genus of $C$, so that finally

$$h^0\big(\mathcal{O}_{\mathbf{P}^2}(d-3) \otimes \mathcal{I}\big) = g = p_a(d) - \delta$$
$$= h^0\big(\mathcal{O}_{\mathbf{P}^2}(d-3)\big) - \delta,$$

as we wanted to prove. □

## 4.2 – Merging a line

This section is devoted to the proof of part b) of Theorem (4.1), which is the case when $V'$ parametrizes curves containing the line $L$. The proof will be divided in two parts (Sections 4.2.1 and 4.2.2) which we shall describe shortly, but we first need to introduce some notation. We advise the reader that he should have Sections 3.1 and 3.4 clear in mind in order to understand what is going on.

**(4.5) The situation.** We consider a curve $X_0 = C \cup L$, where $[C]$ is a general member of the component $V'$ of $V^{d-1,\delta'}(\alpha', \beta')(\Omega')$ we have started with.

We have $\alpha' \leqslant \alpha$, and correspondingly $\Omega'$ is a subset of $\Omega$; by this we mean that $\Omega' = (p'_{i,j})_{1 \leqslant j \leqslant \alpha'_i}$ such that for all $i$, $(p'_{i,j})_{1 \leqslant j \leqslant \alpha'_i}$ is a subsequence of $(p_{i,j})_{1 \leqslant j \leqslant \alpha_i}$.

On the other hand we have $\beta' \geqslant \beta$, and we choose a subset $\Lambda = (s_{i,j})_{1 \leqslant j \leqslant \beta_i}$ of cardinality $\beta$ of the sequence $(s'_{i,j})_{1 \leqslant j \leqslant \beta'_i}$ of unassigned contact points of $C$ with $L$.

**(4.6) General strategy.** The first part of the proof will consist in deforming the curve $X_0$ to degree $d$ curves locally equigenerically around the points of $\Omega'$ and $\Lambda$, while maintaining the contact conditions encoded in $\alpha, \beta$, and $\Omega$. At a point $p'_{i,j}$ or $s_{i,j}$, the curve $C \cup L$ has an order $i$ tacnode; the requirement that the deformation be locally equigeneric around such a point amounts in pratice to the requirement that the $i$-th order tacnode deforms into $i$ nodes.

In the second part of the proof we will take care of the $\beta' - \beta$ remaining mobile contact points, i.e., those points $s'_{i,j}$ that are not in $\Lambda$. At these points, the curve $C \cup L$ has order $i$ tacnodes as well, but these ones will be deformed into $i - 1$ nodes only. This will ensure in particular that the line $L$ be merged with $C$ along the deformation, so that for instance if $C$ is irreducible, then the curves deformed from $C \cup L$ will be irreducible as well.

The reader should consult Section 3.4 to see that these are indeed the appropriate deformations. The reason for this two steps deformation process is that we want to isolate the deformations of $i$-tacnodes to $i - 1$ nodes from the rest, in order to apply the results from the two previous lectures [IV] and [V].

**4.2.1 The relaxed Severi variety**

Here we perform the first step of the proof, as described in (4.6) above. The deformation procedure described there will amount to defining a relaxed Severi variety, seeing $X_0 = C \cup L$ as one of its members, and showing that our relaxed Severi variety is smooth at $[X_0]$.

We use the notation introduced in (4.5) above, and in particular we fix the choice of a subsequence $\Lambda$. Moreover we label the points in $\Omega \setminus \Omega'$ as $(p^L_{i,j})_{1 \leqslant j \leqslant \alpha^L_i}$; then, letting $\alpha^L = (\alpha^L_i)_{i \geqslant 1}$, we have $\alpha = \alpha' + \alpha^L$. Note that around a point $p^L_{i,j}$, the curve $X_0$ has only one local branch, which is an open subset of the line $L$.

On the other hand, we label the mobile contact points of $C$ with $L$ that have not been included in $\Lambda$ as $(s^{C \cap L}_{i,j})_{1 \leqslant j \leqslant \beta^{C \cap L}_i}$; then, letting $\beta^{C \cap L} = (\beta^{C \cap L}_i)_{i \geqslant 1}$, we have $\beta' = \beta + \beta^{C \cap L}$.

**(4.7) The relaxed Severi variety.** In an analytic neighbourhood of the point $[X_0]$ in $|\mathcal{O}_{\mathbf{P}^2}(d)|$, we define the *relaxed Severi variety* $W_\Lambda$ to be the closure of the locus of reduced, degree $d$ curves $X_t$ satisfying the following six conditions:
 (i) $X_t$ preserves the $\delta'$ nodes of $C$, i.e., for every node of $X_0$ off $L$, $X_t$ will have a node nearby;
 (ii) at each point $p^L_{i,j} \in \Omega \setminus \Omega'$, $X_t$ has contact of order $i$ with $L$;
 (iii) in a neighbourhood of each point $p'_{i,j} \in \Omega'$, $X_t$ has singularities with total cogenus $i$, i.e.,
  if we let $\bar{X}_t$ be the normalization of exactly these singularities, then $p_a(\bar{X}_t) = p_a(X_t) - i$;
 (iv) in a neighbourhood of each point $s_{i,j} \in \Lambda$, $X_t$ has singularities with total cogenus $i$;
 (v) in a neighbourhood of each point $p'_{i,j} \in \Omega'$, it follows from condition (iii) that $X_t$ has two local branches; we require that the one branch that is a deformation of an open subset of $C$ has contact of order $i$ with $L$ at $p'_{i,j}$;
 (vi) in a neighbourhood of each point $s_{i,j} \in \Lambda$, the curve $X_t$ has two local branches; we require that the one branch that is a deformation of an open subset of $C$ has contact of order $i$ with $L$ at a point near $s_{i,j}$.

Just to be sure, we recall that the curve $X_0 = C \cup L$ has an $i$-th order tacnode at each point $p'_{i,j} \in \Omega'$ or $s_{i,j} \in \Lambda$, with two local branches that are open subsets of $C$ and $L$ respectively. Conditions (iii) and (iv) imply that in a neighbourhood of such a point, the curves $X_t$ have two local branches as well, which are deformations of open subsets of $C$ and $L$ respectively.

With this definition for the relaxed Severi variety, it is clear that $[X_0] \in W_\Lambda$. We could have defined $W_\Lambda$ alternatively, by requiring in (iii) and (iv) that in a neighbourhood of each point $p'_{i,j}$ or $s_{i,j}$, the curve $X_t$ has $i$ nodes. The reason why $W_\Lambda$ will be smooth at $[X_0]$ is because the only special feature of $X_0$ with respect to a general member of $W_\Lambda$ is that it has several groups of $i$ nodes coming together as $i$-th order tacnodes (for various $i$'s), and this is perfectly harmless from the point of view of the deformation theory of nodal curves: somehow, this is the classical view that an $i$-th order tacnode is just $i$ infinitely near ordinary nodes. This will be made precise when we compute the Zariski tangent space to $W_\Lambda$ at the point $[X_0]$.

In order to show the smoothness of $W_\Lambda$ at $[X_0]$ we will compare the dimension of $W_\Lambda$ with that of its tangent space at $[X_0]$. We proceed to compute the dimension of $W_\Lambda$ in the next paragraph.

**(4.8) The relaxed Severi variety as a local sheet of a Severi variety.** The basic observation in this paragraph is that the six conditions defining the relaxed Severi variety $W_\Lambda$ in fact essentially define a generalized Severi variety in the plain sense. Namely, we are requiring that the curves $X_t$, (i) have a certain genus (or, equivalently, cogenus), defined by the $\delta'$ nodes deformations of those of $C$ and the various groups of $i$ nodes deformations of the points $p'_{i,j} \in \Omega'$ and $s_{i,j} \in \Lambda$, and (ii) satisfy certain contact conditions with $L$. Specifically, the cogenus, i.e., the difference in arithmetic genus between the curves $X_t$ and their normalizations, is required to be:

$$(4.8.1) \quad \begin{aligned} \delta'' &:= \delta' + I\alpha' + I\beta = \delta' + I\alpha' + I\beta' - I(\beta' - \beta) \\ &= \delta' + (d-1) - I(\beta' - \beta) = \delta - \big(I(\beta' - \beta) - |\beta' - \beta|\big) \end{aligned}$$

(for the last equality, recall the relation between $\delta$ and $\delta'$, given for instance in Theorem (3.1)). On the other hand the contact conditions with $L$ are exactly those prescribed by $\alpha, \beta$ and $\Omega$. Altogether, these are the conditions defining the (generalized) Severi variety $V^{d,\delta''}(\alpha,\beta)(\Omega)$.

The additional requirements in the definition of $W_\Lambda$ have the effect of selecting an open subset of some local sheets of $V^{d,\delta''}(\alpha,\beta)(\Omega)$ at $[X_0]$. We are first of all requiring that the $i$ nodes corresponding to a point $p'_{i,j}$ or $s_{i,j}$ be in a neighbourhood of this same point; in other words we prescribe among the tacnodes of $X_0$ which deform to $i$ nodes, and which are ignored. This already selects local sheets of $V^{d,\delta''}(\alpha,\beta)(\Omega)$ at $[X_0]$, determined by $\Lambda$: indeed, if we consider $W_\Lambda$ as a space of maps from curves of cogenus $\delta''$ to $\mathbf{P}^2$, then the member of $W_\Lambda$ corresponding to $X_0$ is the partial normalization of $X_0 = C \cup L$ at the points of $\Omega'$ and $\Lambda$; if one considers another partial normalization of $X_0$ at the points of $\Omega'$ and some other $\Lambda'$, one gets a curve member of the space of maps corresponding to $V^{d,\delta''}(\alpha,\beta)(\Omega)$, which however does not belong to $W_\Lambda$.

We are moreover requiring that, locally around a point $p'_{i,j}$ or $s_{i,j}$, the order $i$ contact with $L$ be supported on the one local branch that is a deformation of the local branch of $C$. This again, and in the same way as above, has the effect of operating a selection among the local sheets of $V^{d,\delta''}(\alpha,\beta)(\Omega)$ at $[X_0]$.

The upshot is that the relaxed Severi variety $W_\Lambda$ is an open subset in some local sheets of the generalized Severi variety $V^{d,\delta''}(\alpha,\beta)(\Omega)$. Since the latter is equidimensional, we have

$$\dim(W_\Lambda) = \dim\big(V^{d,\delta''}(\alpha,\beta)\big) = 2d + g'' + |\beta| - 1,$$

with $g'' = p_a(d) - \delta'' = \binom{d-1}{2} - \delta''$.

**(4.9) The ideal defining the tangent space to the relaxed Severi variety.** We shall prove, in Proposition (4.10) below, that the tangent space at $[X_0]$ of the relaxed Severi variety $W_\Lambda$ identifies with the subspace $H^0(X_0, \mathcal{I}(d))$ of $H^0(X_0, \mathcal{O}(d))$ defined by an ideal sheaf $\mathcal{I}$ of $\mathcal{O}_{X_0}$ which we shall now describe, based on the description of the tangent space of generalized Severi varieties given in [III].

If $X$ is a $\delta$-nodal curve member of the ordinary Severi variety $V^{d,\delta}$, then the tangent space $T_{[X]}V^{d,\delta}$ is $H^0(X, \mathcal{A}(d))$, with $\mathcal{A}$ the ideal sheaf of the points supporting the nodes of $X$. When $i$ of the $\delta$ nodes become infinitely near so as to form an $i$-th order tacnode, the ideal sheaf $\mathcal{A}$ has to be modified in the most straightforward way, i.e., it is merely the ideal sheaf of the possibly infinitely near points supporting the nodes of $X$. In other words, it is cut out locally at an order $i$ tacnode by curves osculating to the order $i$ the two local branches of $X$ at the tacnode.

On the other hand we have seen in a previous Lecture, see [III, Lemma (2.11)], that an order $i$ contact condition with $L$ at an assigned (resp. unassigned) point translates into the vanishing to the order $i$ (resp. $i-1$) of the sections in $H^0(X, N_X)$ corresponding to vectors tangent to such deformations.

These considerations motivate the definition of the ideal sheaf $\mathcal{I}$.

*(4.9.1) Definition.* We consider the curve $X_0 = C \cup L$ introduced in (4.5) above. For all $i, j$ with $1 \leqslant j \leqslant \alpha'_i$, we let $^\alpha\mathcal{I}_{i,j}$ be the subsheaf of $\mathcal{O}_{X_0}$ of those regular functions $f$, such that the restrictions $f|_C$ and $f|_L$ vanish at $p'_{i,j}$ to the orders $i$ and $2i$ respectively.

Similarly, for all $i,j$ with $1 \leqslant j \leqslant \beta_i$, we let ${}^\beta\mathcal{I}_{i,j}$ be the sheaf of regular functions $f$ on $X_0$, such that $f|_C$ and $f|_L$ vanish at $s_{i,j}$ to the orders $i$ and $2i-1$ respectively.

For all points $x \in X_0$, we denote by $\mathfrak{m}_x$ the ideal sheaf defining $x$ in $X_0$. Lastly, we call $u_1, \ldots, u_{\delta'}$ the $\delta'$ nodes of $C$. Finally, we set:

$$\mathcal{I} = \prod_{1 \leqslant i \leqslant \delta'} \mathfrak{m}_{u_i} \cdot \prod_{1 \leqslant j \leqslant \alpha_i^L} \mathfrak{m}_{p_{i,j}^L}^i \cdot \prod_{1 \leqslant j \leqslant \alpha_i'} {}^\alpha\mathcal{I}_{i,j} \cdot \prod_{1 \leqslant j \leqslant \beta_i} {}^\beta\mathcal{I}_{i,j}.$$

**(4.10) Proposition.** *The relaxed Severi variety $W_\Lambda$ is smooth at $[X_0]$, with tangent space*

$$T_{[X_0]}W_\Lambda = H^0(X_0, \mathcal{I}(d)) \ \subseteq \ T_{[X_0]}\bigl(|\mathcal{O}_{\mathbf{P}^2}(d)|\bigr) = H^0(X_0, \mathcal{O}_{X_0}(d)),$$

*where $\mathcal{I}$ is the sheaf of ideals defined in (4.9.1) above.*

*Proof.* We will compare the deformations of $X_0$ along $W_\Lambda$ with those of the map given by a certain partial normalization of $X_0$. Let $\mu : \tilde{X}_0 \to X_0$ be the normalization of $X_0 = C \cup L$ at the nodes $u_i$ of $C$ and at the tacnodes $p'_{i,j} \in \Omega'$ and $s_{i,j} \in \Lambda$ (at which $C$ and $L$ touch), but not at the tacnodes $s'_{i,j}$, because the latter are off $\Lambda$. We will consider the deformations of the map $\phi : \tilde{X}_0 \to \mathbf{P}^2$ obtained as the composition of $\mu$ with the inclusion $X_0 \subseteq \mathbf{P}^2$.

The map $\phi$ is unramified, hence its normal sheaf $N_\phi$ is locally free, equal to $\omega_{\tilde{X}_0} \otimes \phi^*\omega_{\mathbf{P}^2}^{-1}$, see [II, Lemma (8.3)]. By Lemma (4.11) below, $H^1(\tilde{X}_0, N_\phi)$ vanishes, hence the deformation space of $\phi$ is smooth of dimension

$$h^0(\tilde{X}_0, N_\phi) = h^0(\tilde{X}_0, \omega_{\tilde{X}_0} \otimes \phi^*\omega_{\mathbf{P}^2}^{-1}).$$

The smoothness of this space ensures in particular the existence of a map $F$ from a neighbourhood of $[\phi]$ to $|\mathcal{O}_{\mathbf{P}^2}(d)|$, mapping a deformation of $\phi$ to its image in $\mathbf{P}^2$, as in [III, (2.6)].

The image of $F$ is what we will call the *superrelaxed Severi variety* $\tilde{W}_\Lambda$, which we define as $W_\Lambda$ but only with the three requirements (i), (iii), and (iv) out of the six defining $W_\Lambda$; simply put, we just drop the contact conditions with $L$. Arguing as in (4.8) above, we see that $\tilde{W}_\Lambda$ is a local sheet of the plain Severi variety $V^{d,\delta''}$ at $[X_0]$, which implies that it has dimension

$$\dim(\tilde{W}_\Lambda) = \dim(V^{d,\delta''}) = 3d + g'' - 1.$$

That $F$ surjects onto a neighbourhood of $[X_0]$ in $\tilde{W}_\Lambda$ is a consequence of Teissier's Résolution Simultanée Theorem, as in [III, (2.6)].

In order to analyze the tangent space to $\tilde{W}_\Lambda$ at $[X_0]$, we let $\mathcal{A} = (\mathcal{O}_{X_0} : \mu_*\mathcal{O}_{\tilde{X}_0})$ be the conductor of $\tilde{X}_0$ in $X_0$, see [II, §2]. It is a sheaf of ideals of $\mathcal{O}_{X_0}$, trivial on the open subset onto which $\mu$ is an isomorphism, equal to $\mathfrak{m}_{u_i}$ locally at the nodes $u_i$ of $C$, while locally at the points $p'_{i,j} \in \Omega'$ and $s_{i,j} \in \Lambda$ it is the sheaf of regular functions vanishing to the order $i$ on both $C$ and $L$; this can be seen for instance with [II, Lemma (2.5)].

By [II, Corollary (8.4)], one has $\mu_*N_\phi = \mathcal{A} \otimes N_{X_0/\mathbf{P}^2}$. Therefore, there is a canonical isomorphism

$$H^0(\tilde{X}_0, N_\phi) \cong H^0(X_0, \mathcal{A} \otimes N_{X_0/\mathbf{P}^2}).$$

This isomorphism identifies with the differential of $F$ at $[\phi]$, and thus we see that $F$ is an isomorphism onto a neighbourhood of $[X_0]$ in $\tilde{W}_\Lambda$; in particular, the superrelaxed Severi variety $\tilde{W}_\Lambda$ is smooth at $[X_0]$, with tangent space $H^0(X_0, \mathcal{A} \otimes N_{X_0/\mathbf{P}^2})$.

Now, it follows from [III, (2.11)] that the Zariski tangent space at $[\phi]$ to the image by $F^{-1}$ of the relaxed Severi variety $W_\Lambda \subseteq \tilde{W}_\Lambda$ is contained in

$$H^0\biggl(\tilde{X}_0, N_\phi\Bigl(-\sum_{1 \leqslant j \leqslant \alpha_i^L} i \cdot q_{i,j}^L - \sum_{1 \leqslant j \leqslant \alpha_i'} i \cdot q'_{i,j} - \sum_{1 \leqslant j \leqslant \beta_i} (i-1) \cdot r_{i,j}\Bigr)\biggr),$$

where (i) $q_{i,j}^L$ is the unique point of $\tilde{X}_0$ over $p_{i,j}^L$, and (ii) $q'_{i,j}$ and $r_{i,j}$ are the points on the proper transform of $C$ in $\tilde{X}_0$ lying over $p'_{i,j}$ and $s_{i,j}$ respectively. Moreover, it follows from the description of $\mathcal{A}$ given above that

$$\mathcal{A} \otimes \phi_* \mathcal{O}_{\tilde{X}_0}\Big(-\sum_{1 \leqslant j \leqslant \alpha_i^L} i \cdot q_{i,j}^L - \sum_{1 \leqslant j \leqslant \alpha'_i} i \cdot q'_{i,j} - \sum_{1 \leqslant j \leqslant \beta_i} (i-1) \cdot r_{i,j}\Big) = \mathcal{I},$$

hence

$$\phi_* N_\phi\Big(-\sum_{1 \leqslant j \leqslant \alpha_i^L} i \cdot q_{i,j}^L - \sum_{1 \leqslant j \leqslant \alpha'_i} i \cdot q'_{i,j} - \sum_{1 \leqslant j \leqslant \beta_i} (i-1) \cdot r_{i,j}\Big) = \mathcal{I} \otimes N_{X_0/\mathbf{P}^2}.$$

Therefore, the Zariski tangent space to the relaxed Severi variety $W_\Lambda$ at $[X_0]$ is contained in $H^0(X_0, \mathcal{I} \otimes N_{X_0/\mathbf{P}^2})$. Then, the proposition follows from the equality

$$h^0(X_0, \mathcal{I} \otimes N_{X_0/\mathbf{P}^2}) = \dim(W_\Lambda),$$

which is given by (4.12.2) below. $\square$

The rest of the present Section 4.2.1 is devoted to some elementary considerations leading to the non-speciality of the invertible sheaves on $\tilde{X}_0$ involved in the above proof, and the computation of the dimensions of their spaces of global sections with the Riemann–Roch Formula.

**(4.11) Lemma.** *Let $p$ be any point of the partial normalization $\tilde{X}_0$, and*

$$D = \sum_{1 \leqslant j \leqslant \alpha_i^L} i \cdot q_{i,j}^L + \sum_{1 \leqslant j \leqslant \alpha'_i} i \cdot q'_{i,j} + \sum_{1 \leqslant j \leqslant \beta_i} (i-1) \cdot r_{i,j},$$

*with the notation as in the proof of Proposition (4.10) above. The two line bundles $N_\phi(-D)$ and $N_\phi(-D-p)$ on $\tilde{X}_0$ are non-special.*

*Proof.* Let $\tilde{C}$ and $\tilde{L}$ be the proper transforms in $\tilde{X}_0$ of $C$ and $L$ respectively. Since $N_\phi = \omega_{\tilde{X}_0} \otimes \phi^* \omega_{\mathbf{P}^2}^{-1}$ as we have seen above, one has

$$\omega_{\tilde{X}_0} \otimes \big(N_\phi(-D)\big)^{-1} = \phi^* \mathcal{O}_{\mathbf{P}^2}(-3) \otimes \mathcal{O}_{\tilde{X}_0}(D).$$

The degree of the restriction of this line bundle to $\tilde{C}$ has degree

$$-3(d-1) + I\alpha' + \big(I\beta - |\beta|\big) \leqslant -3(d-1) + I\alpha' + I\beta'$$
$$= -2(d-1) < 0,$$

hence every global section of $\omega_{\tilde{X}_0} \otimes N_\phi(-D)^{-1}$ vanishes identically on $\tilde{C}$. One proves similarly that the same holds for $\omega_{\tilde{X}_0} \otimes N_\phi(-D-p)^{-1}$, noting that it restricts on $\tilde{C}$ to a line bundle of degree $-2(d-1) + \varepsilon$, where $\varepsilon$ equals 0 or 1 depending on whether $p$ sits on $\tilde{C}$ or not.

Thus, every global section of $\omega_{\tilde{X}_0} \otimes N_\phi(-D)^{-1}$ restricts on $\tilde{L}$ to a global section of the line bundle

$$\phi^* \mathcal{O}_{\mathbf{P}^2}(-3) \otimes \mathcal{O}_{\tilde{L}}(D|_{\tilde{L}}) \otimes \mathcal{O}_{\tilde{L}}(-\tilde{C}|_{\tilde{L}}),$$

which has degree

$$-3 + I\alpha^L - I\beta^{C \cap L} = -3 + I(\alpha - \alpha') - I(\beta' - \beta)$$
$$= -3 + (I\alpha + I\beta) - (I\alpha' + I\beta')$$
$$= -3 + d - (d-1) = -2.$$

Therefore every global section of $\omega_{\tilde{X}_0} \otimes N_\phi(-D)^{-1}$ vanishes identically on $\tilde{L}$ as well, hence

$$h^0(\tilde{X}_0, \omega_{\tilde{X}_0} \otimes N_\phi(-D)^{-1}) = h^1(\tilde{X}_0, N_\phi(-D)) = 0.$$

One concludes in the same way that every global section of $\omega_{\tilde{X}_0} \otimes N_\phi(-D-p)^{-1}$ vanishes identically, hence $N_\phi(-D-p)$ is non-special. $\square$

**(4.12) Lemma.** *The following equalities hold:*

(4.12.1) $$h^0(\tilde{X}_0, N_\phi) = 3d + g'' - 1;$$
(4.12.2) $$h^0(\tilde{X}_0, N_\phi(-D)) = 2d + g'' - 1 + |\beta|.$$

*Moreover, the line bundle $N_\phi(-D)$ has no base point.*

*Proof.* Both equalities follow from the application of the Riemann–Roch Formula on the curve $\tilde{X}_0$, which has arithmetic genus $g''$, see (4.8). The line bundle $N_\phi$ equals $\omega_{\tilde{X}_0} \otimes \phi^* \omega_{\mathbf{P}^2}^{-1}$ which is obviously non-special, hence

$$\begin{aligned} h^0(N_\phi) &= 1 - g'' + \deg\bigl(\omega_{\tilde{X}_0} \otimes \phi^* \omega_{\mathbf{P}^2}^{-1}\bigr) \\ &= 1 - g'' + (2g'' - 2) + 3d \\ &= 3d + g'' - 1. \end{aligned}$$

The line bundle $N_\phi(-D)$ is non-special by Lemma (4.11) above, hence

$$\begin{aligned} h^0\bigl(N_\phi(-D)\bigr) &= 1 - g'' + \deg(N_\phi) - \deg(D) \\ &= 3d + g'' - 1 - \bigl(I\alpha^L + I\alpha' + I\beta - |\beta|\bigr) \\ &= 2d + g'' - 1 + |\beta|, \end{aligned}$$

because

$$I\alpha^L + I\alpha' + I\beta = I\alpha + I\beta = d.$$

Eventually, since for all $p \in \tilde{X}_0$ the line bundle $N_\phi(-D-p)$ is non-special as well, we have

$$h^0\bigl(\tilde{X}_0, N_\phi(-D-p)\bigr) = h^0\bigl(\tilde{X}_0, N_\phi(-D)\bigr) - 1,$$

hence not all global sections of $N_\phi(-D)$ vanish at $p$. $\square$

### 4.2.2   Deformation of order $i$ tacnodes to $i-1$ nodes

Here we perform the second step of the proof of part b) of Theorem (4.1), as described in (4.6) above. We assume that $\beta' > \beta$, for otherwise there is nothing left to do in the second step (see Remark (4.21) below for what happens in this case).

We consider the map
$$F_\Lambda : W_\Lambda \to \Delta,$$
from the relaxed Severi variety $W_\Lambda$, introduced and studied in the previous Section 4.2.1, to the product $\Delta$ of the versal deformation spaces of the tacnodes $s_{i,j}^{C \cap L}$ of $X_0 = C \cup L$, i.e., those tacnodes of $X_0$ not included in $\Lambda$.

We refer to [IV] and [V] for background and the necessary material on deformation spaces of tacnodes and their product. The result we shall need here is [V, (1.1)], which we now recall in a form adapted to our context.

**(4.13)** The tacnodes $s_{i,j}^{C\cap L}$ are indexed by $\beta' - \beta$, and thus

$$\Delta = \prod_{i\in \mathbf{N}^*} \prod_{1\leqslant j\leqslant \beta'_i-\beta_i} \Delta_{i,j},$$

where each $\Delta_{i,j}$ is the versal deformation space of the $i$-th order tacnode $s_{i,j}^{C\cap L}$. For all $i$ and $j$, we let $\Delta_{i,j;i}$ and $\Delta_{i,j;i-1}$ be the subvarieties of $\Delta_{i,j}$ of those deformations of $s_{i,j}^{C\cap L}$ into $i$ and $i-1$ nodes respectively. We let

$$\mathbf{m} = (m_l)_{1\leqslant l\leqslant |\beta'-\beta|} = (\underbrace{1,\ldots,1}_{\beta'_1-\beta_1}, \quad \ldots, \quad \underbrace{n,\ldots,n}_{\beta'_n-\beta_n})$$

with $n = \max(\beta' - \beta)$, and consider

$$\Delta_{\mathbf{m}} = \prod_{i\in \mathbf{N}^*} \prod_{1\leqslant j\leqslant \beta'_i-\beta_i} \Delta_{i,j;i}, \quad \text{and} \quad \Delta_{\mathbf{m}-\mathbf{1}} = \prod_{i\in \mathbf{N}^*} \prod_{1\leqslant j\leqslant \beta'_i-\beta_i} \Delta_{i,j;i-1}.$$

Each $\Delta_{i,j}$ is a $(2i-1)$-dimensional affine space, with coordinates $a_0,\ldots,a_{i-2},b_0,\ldots,b_{i-1}$, such that the corresponding versal deformation is the hypersurface

(4.13.1) $$y^2 + y(x^i + a_{i-2}x^{i-2} + \cdots + a_0) + b_{i-1}x^{i-1} + \cdots + b_0$$

in $\mathbf{A}_{x,y}^2 \times \mathbf{A}_{\mathbf{a},\mathbf{b}}^{2i-1}$. In these terms, $\Delta_{i,j;i}$ is the dimension $i-1$ affine subspace defined by the equations $b_0 = \cdots = b_{i-1} = 0$. We consider in addition the hyperplane $\mathcal{H}_{i,j} \subseteq \Delta_{i,j}$ defined by the equation $b_0 = 0$ in this system of coordinates. Observe that $\mathcal{H}_{i,j}$ parametrizes those deformations of the tacnode defined by $y(y+x^m) = 0$ that still pass through the origin in $\mathbf{A}_{x,y}^2$; in other words, these are the deformations maintaining one assigned contact point of order 1 with the line $y = 0$. The general member of $\mathcal{H}_{i,j}$ is an irreducible curve. Eventually we let

$$\mathcal{H} = \bigcup_{i,j} \mathcal{H}_{i,j} \subseteq \Delta.$$

Now, the result is the following. Let $W \subseteq \Delta$ be a smooth subvariety of dimension $\dim(\Delta_{\mathbf{m}-\mathbf{1}})+1$, containing $\Delta_{\mathbf{m}-\mathbf{1}}$, and such that the tangent space at the origin $T_0W$ is not contained in $\mathcal{H}$. Then, in an étale neighbourhood of the origin $0 \in \Delta$,

$$W \cap \Delta_{\mathbf{m}-\mathbf{1}} = \Delta_{\mathbf{m}} \cup \Gamma_1 \cup \cdots \cup \Gamma_\kappa,$$

where

$$\kappa = \frac{\prod_{1\leqslant l\leqslant |\beta'-\beta|} m_l}{\mathrm{lcm}_{1\leqslant l\leqslant |\beta'-\beta|} m_l} = \frac{I^{\beta'-\beta}}{\mathrm{lcm}(\beta'-\beta)},$$

and $\Gamma_1,\ldots,\Gamma_\kappa$ are distinct reduced unibranched curves, each of which has the origin as a point of multiplicity

$$\frac{\mathrm{lcm}_{1\leqslant l\leqslant |\beta'-\beta|} m_l}{\max_{1\leqslant l\leqslant |\beta'-\beta|} m_l} = \frac{\mathrm{lcm}(\beta'-\beta)}{\max(\beta'-\beta)},$$

and has intersection multiplicity $\mathrm{lcm}_{1\leqslant l\leqslant |\beta'-\beta|} m_l = \mathrm{lcm}(\beta'-\beta)$ with $\Delta_{\mathbf{m}}$ at the origin.

We now proceed to prove a series of claims in order to see that the hypotheses of the above statement are verified in our situation, taking $W = F_\Lambda(W_\Lambda)$.

**(4.14) Claim.** The inverse image $W_{\Lambda,\mathbf{m}} = F_\Lambda^{-1}(\Delta_\mathbf{m})$ is the locus of those members of $W_\Lambda$ that contain $L$. It is also the intersection of $W_\Lambda$ with the hyperplane $p^\perp \subseteq |\mathcal{O}_{\mathbf{P}^2}(d)|$, for general $p \in L$.

*Proof.* Let us first prove the first assertion. The inverse image $F_\Lambda^{-1}(\Delta_\mathbf{m})$ is the locus in $W_\Lambda$ along which all the tacnodes $s_{i,j}^{C \cap L}$ of $X_0$ deform equigenerically. Taking in addition into account the requirements (i), (iii), and (iv) in the definition of the relaxed Severi variety $W_\Lambda$, we see that the curves in $F_\Lambda^{-1}(\Delta_\mathbf{m})$ are equigeneric deformations of $X_0$. Therefore they all come from deformations of the map $\bar\phi : \bar X_0 \to \mathbf{P}^2$ obtained from the total normalization $\nu : \bar X_0 \to X_0$, by the same argument that has been used in the proof of Proposition (4.10), which is essentially Teissier's Résolution Simultanée. Now since the normalization $\bar X_0$ is disconnected, the domain of any deformation of $\bar\phi$ is disconnected as well by the Principle of Connectedness (see, e.g., [16, Exercise III.11.4]). The upshot is that all curves in $F_\Lambda^{-1}(\Delta_\mathbf{m})$ consist of a deformation of $C$ plus a deformation of $L$.

On the other hand we are assuming here that $\beta' > \beta$, for otherwise the second step of the proof of part b) of Theorem (4.1) is pointless. Since
$$I\alpha + I\beta = d \quad \text{and} \quad I\alpha' + I\beta' = d-1,$$
it follows that $I(\alpha - \alpha') \geqslant 2$, i.e., $I\alpha^L \geqslant 2$. Now by requirement (ii) in the definition of the relaxed Severi variety, for each member of $F^{-1}(\Delta_\mathbf{m})$, the component that is a deformation of $L$ must satisfy the contact conditions with $L$ encoded in $\alpha^L$. Therefore it intersects $L$ in at least two points (possibly infinitely near), hence it is $L$ itself, and our first assertion about $F^{-1}(\Delta_\mathbf{m})$ is proved.

As for the second assertion, the first assertion gives us of course the inclusion $F_\Lambda^{-1}(\Delta_\mathbf{m}) \subseteq W_\Lambda \cap p^\perp$ for all $p \in L$. To prove the other inclusion, we note that for all $p \in L$ neither in $\Omega$ nor in $\Lambda$, any curve in $W_\Lambda \cap p^\perp$ in a neighbourhood of $[X_0]$ has a total of
$$I\alpha + I\beta + 1 = d+1$$
imposed points of intersection with $L$ whereas it has degree $d$, so that it must contain $L$. □

**(4.15) Claim.** The image of $F_\Lambda$ contains the locus $\Delta_\mathbf{m}$.

*Proof.* The projection of $F_\Lambda(W_\Lambda)$ to the factor $\Delta_{i,j}$ of $\Delta$ contains $\Delta_{i,j;i}$ if and only if we can find a deformation $X_t$ of $X_0$ along $W_\Lambda$ in which a neighbourhood of the tacnode $s_{i,j}^{C \cap L} \in X_0$ deforms to a reducible open curve $X_t^{\circ\prime} + X_t^{\circ\prime\prime}$, such that $X_t^{\circ\prime}$ and $X_t^{\circ\prime\prime}$ intersect transversely in $i$ points that can be chosen arbitrarily on $X_t^{\circ\prime}$.

To see that this holds simultaneously for all tacnodes $s_{i,j}^{C \cap L} \in X_0$, we note that $F_\Lambda^{-1}(\Delta_\mathbf{m})$ identifies as in (4.8) with a local sheet at $[X_0]$ of the Severi variety
$$V^{d-1,\delta'}\bigl(\alpha', \beta + I(\beta'-\beta)e_1\bigr)(\Omega') :$$
we have seen that $F_\Lambda^{-1}(\Delta_\mathbf{m})$ consists of those curve in $W_\Lambda$ having $L$ as an irreducible component, and these curves correspond to equigeneric deformations of $C$ maintaining the contact conditions with $L$ in the points $p_{i,j}' \in \Omega'$ and $s_{i,j} \in \Lambda$, but dropping those in the points encoded in the points $s_{i,j} \notin \Lambda$; the latter dropping so to speak amounts to turn $\beta' - \beta$ into $I(\beta'-\beta)e_1$.

Then the condition we need to verify essentially follows from Proposition (2.2) applied to $V^{d-1,\delta'}\bigl(\alpha', \beta + I(\beta'-\beta)e_1\bigr)(\Omega')$. More explicitly, the latter Severi variety contains locally around $[X_0]$ the Severi varieties
$$V^{d-1,\delta'}\bigl(\alpha' + I(\beta'-\beta)e_1, \beta\bigr)(\Omega' \cup \Xi_t),$$
where $\Xi_t$ is an arbitrary choice of $i$ points on $L$ around each of the points $s_{i,j} \notin \Lambda$, which is exactly the condition we needed to verify. □

**(4.16) Claim.** The inverse image $W_{\Lambda,\mathbf{m}} = F_\Lambda^{-1}(\Delta_\mathbf{m})$ has codimension one in $W_\Lambda$.

*Proof.* This follows from the identification given above of $F_\Lambda^{-1}(\Delta_\mathbf{m})$ with a local sheet at $[X_0]$ of the Severi variety $V^{d-1,\delta'}(\alpha', \beta + I(\beta' - \beta)e_1)$. Explicitly, this identifications tells us the dimension of $F_\Lambda^{-1}(\Delta_\mathbf{m})$ as in (4.8), then

$$\dim(W_\Lambda) - \dim\bigl(F_\Lambda^{-1}(\Delta_\mathbf{m})\bigr)$$
$$= \left[2d + \binom{d-1}{2} - \delta'' - 1 + |\beta|\right] - \left[2(d-1) + \binom{d-2}{2} - \delta' - 1 + |\beta| + I(\beta' - \beta)\right]$$
$$= 2 + (d-2) - (\delta'' - \delta') - I(\beta' - \beta),$$

and this equals one by (4.8.1). □

**(4.17) Remark.** It is instructive, although not strictly necessary for our purposes, to identify the fibre of $F_\Lambda$ over a general point of $\Delta_\mathbf{m}$. We shall obtain such an identification from the considerations in the proof of Claim (4.15) above, after we observe that a point of $\Delta_{i,j;i} \subseteq \Delta_{i,j}$ corresponds to the choice of $i$ points on one branch of the tacnode at $s_{i,j}$, but only up to reparametrization. This amounts to the fact that, in (4.13.1), the polynomial

$$x^i + a_{i-2}x^{i-2} + \cdots + a_0$$

defining the unions of these $i$ points has no $x^{i-1}$ term; any degree $i$ polynomial in $x$ can be put in this form by a change of variable, which amounts to translating the roots so that they sum up to zero.

Then the fibre of $F_\Lambda$ over a general point of $\Delta_\mathbf{m}$ identifies with an open subset of the Severi variety

(4.17.1) $$V^{d-1,\delta'}\bigl(\alpha' + \bigl[I(\beta' - \beta) - |\beta' - \beta|\bigr]e_1,\ \beta + |\beta' - \beta|e_1\bigr)(\Omega' \cup \Xi_t^\circ),$$

with $\Xi_t^\circ$ an arbitrary choice of $i-1$ points on $L$ around each of the points $s_{i,j} \notin \Lambda$. The fibre of $F_\Lambda$ over the origin of $\Delta$ on the other hand identifies with an open subset of

(4.17.2) $$V^{d-1,\delta'}(\alpha', \beta')(\Omega'),$$

which was indeed our starting point. It may be useful for the conclusion in (4.20) below to recall that it has dimension one less than $V^{d,\delta}(\alpha,\beta)$.

The fibre over the general point of $\Delta_\mathbf{m}$ and the fibre over the origin have the same dimension: indeed the Severi varieties (4.17.1) and (4.17.2) have the same dimension

$$2(d-1) + g' - 1 + |\beta| + |\beta' - \beta| = 2(d-1) + g' - 1 + |\beta'|,$$

with $g' = p_a(d-1) - \delta' = \binom{d-2}{2} - \delta'$.

We leave it to the reader to verify that the latter common dimension of the fibre over the general point of $\Delta_\mathbf{m}$ and the fibre over the origin also equals

$$\dim(W_\Lambda) - \dim(\Delta_\mathbf{m}) - 1;$$

this is yet another straightforward computation after one notes that

$$\dim(\Delta_\mathbf{m}) = \sum\nolimits_{1 \leqslant l \leqslant |\beta' - \beta|} (m_l - 1) = I(\beta' - \beta) - |\beta' - \beta|.$$

**(4.18) Claim.** Let $dF_\Lambda : T_{[X_0]}W_\Lambda \to T_0\Delta$ be the differential map of $F_\Lambda$ at the point $[X_0]$. One has

$$\dim\bigl(dF_\Lambda^{-1}(T_0\Delta_\mathbf{m})\bigr) = \dim\bigl(F_\Lambda^{-1}(\Delta_\mathbf{m})\bigr) = \dim(W_\Lambda) - 1.$$

*Proof.* The second equality is Claim (4.16) above. For the first one, the key observation is that, under the isomorphism
$$T_{[X_0]}W_\Lambda \cong H^0(X_0, \mathcal{I}(d)),$$
the inverse image of $T_0\Delta_\mathbf{m} \subseteq T_0\Delta$ by the differential $dF_\Lambda$ at $[X_0]$ is the subspace of sections vanishing identically along $L \subseteq X_0$.

The restriction $\mathcal{I}(d)|_L$ has degree
$$d - I\alpha - I\beta = 0,$$
and is therefore trivial. Now, the restriction map on global sections,
$$H^0(X_0, \mathcal{I}(d)) \to H^0(L, \mathcal{I}(d)|_L) \cong \mathbf{C},$$
is non-zero by Lemma (4.12), hence
$$\dim(dF_\Lambda^{-1}(T_0\Delta_\mathbf{m})) = h^0(X_0, \mathcal{I}(d)) - 1 = \dim(W_\Lambda) - 1,$$
where the last equality comes from the smoothness of $W_\Lambda$ at $[X_0]$. □

**(4.19) Claim.** The image of the differential map $dF_\Lambda : T_{[X_0]}W_\Lambda \to T_0\Delta$ is not contained in $\mathcal{H}$.

*Proof.* From the description of $\mathcal{H}$ in (4.13), the inverse image $(dF_\Lambda)^{-1}(\mathcal{H})$ is the union of the subspaces in $H^0(X_0, \mathcal{I}(d))$ of sections vanishing at one point $s_{i,j}^{C \cap L}$. By Lemma (4.12) these are all proper subspaces (in fact, hyperplanes) of $H^0(X_0, \mathcal{I}(d))$, and the claim is proved. □

**(4.20) Conclusion.** We may now finish the proof of part b) of Theorem (4.1). It follows from the description given in Section 3.4 that, locally around the point $[X_0]$, the sum of the local sheets of the Severi variety $V^{d,\delta}(\alpha,\beta)(\Omega)$ corresponding to the choice of $\Lambda$ is given as the closure in the relaxed Severi variety $W_\Lambda$ of the inverse image
$$F_\Lambda^{-1}(\Delta_{\mathbf{m-1}} \setminus \Delta_\mathbf{m}).$$

Now, it follows from Claim (4.18) that the differential map
$$dF_\Lambda : T_{[X_0]}W_\Lambda \to T_0\Delta$$
has rank $\dim(\Delta_\mathbf{m}) + 1$, equal to the dimension of $W = F_\Lambda(W_\Lambda)$. Thus $W$ is smooth at $[X_0]$, of dimension $\dim(\Delta_\mathbf{m}) + 1$, and the map $F_\Lambda$ is a submersion over $W$, locally around the origin in $\Delta$.

The subvariety $W \subseteq \Delta$ is smooth and, by Claim (4.19), its tangent space is not contained in $\mathcal{H}$. We may thus apply the main result of [V], to the effect that
$$W \cap \Delta_{\mathbf{m-1}} = \Delta_\mathbf{m} \cup \Gamma_1 \cup \cdots \cup \Gamma_\kappa,$$
with $\Gamma_1, \ldots, \Gamma_\kappa$ curves as in (4.13) above. The upshot is thus that the local sheets of $V^{d,\delta}(\alpha,\beta)(\Omega)$ corresponding to the choice of $\Lambda$ are in the number $\kappa = I^{\beta'-\beta}/\mathrm{lcm}(\beta'-\beta)$, each a fibration over a curve $\Gamma_l$, $l = 1, \ldots, \kappa$, with the fibre $V^{d-1,\delta'}(\alpha',\beta')(\Omega')$ over the origin. The local sheet corresponding to the curve $\Gamma_l$ has multiplicity $\mathrm{lcm}(\beta'-\beta)/\max(\beta'-\beta)$ at the generic point of the fibre over the origin, i.e., the generic point of $V^{d-1,\delta'}(\alpha',\beta')(\Omega')$, and intersection multiplicity $\mathrm{lcm}(\beta'-\beta)$ with the inverse image $F_\Lambda^{-1}(\Delta_\mathbf{m})$, i.e., with the hyperplane $p^\perp \subseteq |\mathcal{O}_{\mathbf{P}^2}(d)|$ for general $p \in L$ by Claim (4.14), along the fibre over the origin.

Finally, note that $\binom{\beta'}{\beta}$ is the number of possible choices for the subsequence $\Lambda$ of the sequence of unassigned contact points of $C$ with $L$. □

The proof of Theorem (4.1) is now over. We end this section with a remark about what happens when $\beta' = \beta$.

As we have already observed, in this case the first step performed in Subsubsection 4.2.1 above is sufficient to prove Theorem (4.1). Indeed, if $\beta = \beta'$, then the relaxed Severi variety $W_\Lambda$ itself is the unique local sheet of the Severi variety $V^{d,\delta}(\alpha,\beta)(\Omega)$ at $[X_0]$, we already know that it is smooth at $[X_0]$, and there is nothing else to prove.

**(4.21) Remark.** If $\beta' = \beta$, then all the irreducible components of $V^{d,\delta}(\alpha,\beta)$ having an irreducible component of $V^{d-1,\delta'}(\alpha',\beta')$ in their intersection with the hyperplane $p^\perp \subseteq |\mathcal{O}_{\mathbf{P}^2}(d)|$ parametrize reducible curves.

Indeed, the members of such a component $W$ in a neighbourhood of $[X_0] = [C \cup L]$ with $[C] \in V^{d-1,\delta'}(\alpha',\beta')$ correspond to deformations of the total normalization $\bar{X}_0 \to X_0 \subseteq \mathbf{P}^2$, and therefore they are reducible curves because $\bar{X}_0$ is disconnected. This is the same argument we have used in the proof of Claim (4.14) above.

Conversely, if $\beta' > \beta$, and if $V$ is an irreducible component of $V^{d-1,\delta'}(\alpha',\beta')$ that parametrizes irreducible degree $d-1$ curves, then all the irreducible components of $V^{d,\delta}(\alpha,\beta)$ having $V$ in their intersection with $p^\perp$ parametrize irreducible curves.

Indeed, if $\beta' > \beta$, then there is at least one order $i$ tacnode of $X_0 = C \cup L$ that is deformed to $i-1$ nodes only, so that locally around such a point the two local branches of $X_0$, corresponding respectively to an irreducible component $C_1$ of $C$ and to $L$, deform to only one, irreducible, local branch, and thus $C_1$ and $L$ merge into the same irreducible component of the curves $X_t$ neighbouring $X_0$.

## 5 – The formula for irreducible curves

In this section we give and explain the version of Caporaso and Harris' recursion formula enabling the counting of irreducible curves. The formula is given in Theorem (5.4) below. The key observation to pass from the standard formula to the formula for irreducible curves is Remark (4.21) at the end of the previous section. Let us start with an example.

**(5.1) Example.** We consider the enumeration of quartics with three nodes and two assigned simple contact points with the line $L$, i.e., curves parametrized by $V^{4,3}(2,2)(\Omega)$. By the recursion formula, we have

(5.1.1) $$N^{4,3}(2,2) = N^{4,3}(3,1) + 3N^{3,1}(0,3) + 2N^{3,0}(1,2).$$

The irreducible components of $V^{4,3}(2,2)(\Omega)$ parametrize either irreducible quartics, or quartics that are the sum of a cubic and a line. We let $N^{4,3}(2,2)_{\text{irr}}$ and $N^{4,3}(2,2)_{3+1}$ be the respective contributions to $N^{4,3}(2,2)$ of irreducible and decomposable quartics. Thus,

$$N^{4,3}(2,2) = N^{4,3}(2,2)_{\text{irr}} + N^{4,3}(2,2)_{3+1}.$$

Formula (5.1.1) splits accordingly into

(5.1.2) $$\begin{aligned} N^{4,3}_{\text{irr}}(2,2) &= N^{4,3}_{\text{irr}}(3,1) + 3N^{3,1}_{\text{irr}}(0,3) + 2N^{3,0}_{\text{irr}}(1,2) \quad \text{and} \\ N^{4,3}_{3+1}(2,2) &= N^{4,3}_{3+1}(3,1) + 3N^{3,1}_{3+1}(0,3) + 2N^{3,0}_{3+1}(1,2). \end{aligned}$$

Note that all components of $V^{3,1}(0,3)$ and $V^{3,0}(1,2)$ parametrize irreducible cubics, hence the distribution of $N^{3,1}(0,3)$ into $N^{3,1}_{\text{irr}}(0,3)$ and $N^{3,1}_{3+1}(0,3)$ (and similarly that of $N^{3,0}(1,2)$) depends on whether the merging of these irreducible cubics with the line $L$ produces an irreducible or a decomposable quartic. The upshot of Remark (4.21) is that

$$N^{3,0}(1,2) = N^{3,0}(1,2)_{3+1} \quad \text{and} \quad N^{3,1}(0,3) = N^{3,1}(0,3)_{\text{irr}},$$

because $\beta' = \beta$ for the former and $\beta' > \beta$ for the latter, in the notation of Remark (4.21). Therefore, (5.1.2) reduces to

$$N^{4,3}_{\text{irr}}(2,2) = N^{4,3}_{\text{irr}}(3,1) + 3N^{3,1}(0,3) \quad \text{and}$$
$$N^{4,3}_{3+1}(2,2) = N^{4,3}_{3+1}(3,1) + 2N^{3,0}(1,2).$$

Now let us show these formulas in action. First note that

$$N^{3,1}(0,3) = 12 \quad \text{and} \quad N^{3,0}(1,2) = 1$$

(the former is the plain number of rational cubics, while the latter is merely the number of cubics passing through one assigned point on $L$, and passing through eight general crossing points in both cases).

Let us compute the contributions $N^{4,3}(2,2)_{3+1}$ and $N^{4,3}(3,1)_{3+1}$. Since $V^{4,3}(2,2)$ has dimension 9, $N^{4,3}(2,2)_{3+1}$ is the number of quartics decomposed as a cubic plus a line passing through two assigned points $p_1, p_2$ on $L$, and through nine general points $a_1, \ldots, a_9$ in the plane. There are only the two following possibilities (recall that in the definition of Severi varieties it is requested that its members should not contain $L$ itself):

(i) the line passes through one of $p_1$ and $p_2$ and through one of the $a_i$'s, and the cubic passes through the remaining nine points among $p_1, p_2$, and $a_1, \ldots, a_9$;
(ii) the line passes through two of the $a_i$'s, and the cubic passes through $p_1, p_2$ and the seven points remaining among $a_1, \ldots, a_9$.

We thus find

$$N^{4,3}(2,2)_{3+1} = 2 \times 9 + \binom{9}{2} = 54.$$

Similarly, $N^{4,3}(3,1)_{3+1}$ counts decomposed quartics passing through three assigned points $p_1, p_2, p_3$ on $L$ and eight general points $a_1, \ldots, a_8$ on the plane. We have the same possibilities as in the previous case, and thus

$$N^{4,3}(3,1)_{3+1} = 3 \times 8 + \binom{8}{2} = 52.$$

We can now observe that, indeed,

$$N^{4,3}(2,2)_{3+1} = 54 = N^{4,3}(3,1)_{3+1} + 2N^{3,0}(1,2) = 52 + 2 \cdot 1.$$

Besides, one has $N^{4,3}(2,2)_{\text{irr}} = N^{4,3}(0,4)_{\text{irr}}$ because putting two crossing points on $L$ doesn't put the crossing points in special position, and therefore $N^{4,3}(2,2)_{\text{irr}} = 620$, see [I, Section 7.2.2], or [1, p. 349]. As a sanity check, note that

$$N^{4,3}(2,2) = 674 = N^{4,3}(2,2)_{\text{irr}} + N^{4,3}(2,2)_{3+1} = 620 + 54$$

(see [1, p. 349] for $N^{4,3}(2,2) = 674$). On the other hand, $N^{4,3}(3,1) = 636$ by [1, p. 349], hence $N^{4,3}(3,1)_{\text{irr}} = 636 - 52 = 584$. We can now observe that, indeed,

$$N^{4,3}_{\text{irr}}(2,2) = 620 = N^{4,3}_{\text{irr}}(3,1) + 3N^{3,1}(0,3) = 584 + 3 \cdot 12.$$

The main point in passing from the standard formula to the formula for irreducible curves is to reconsider Remark (4.21) carefully in order to make it cover all possible cases. The upshot is the following.

**(5.2)** Let the notation be as in Section 4.2 above. In particular, we assume that $V^{d-1,\delta'}(\alpha',\beta')$ appears in the intersection of $V^{d,\delta}(\alpha,\beta)$ with the hyperplane $p^\perp$, and let $V'$ be an irreducible component of $V^{d-1,\delta'}(\alpha',\beta')$. Let $C$ be a general member of $V'$, and assume that it decomposes as a union of irreducible curves $C_1,\ldots,C_k$.

Then, a local sheet of $V^{d,\delta}(\alpha,\beta)$ at the point $[C \cup L]$ parametrizes irreducible curves if and only if the corresponding choice of $\Lambda$ as in Paragraph (4.5) is such that, for all $l = 1,\ldots,k$, there is at least one unassigned contact point of $C$ with $L$ on the component $C_l$ that is off $\Lambda$.

Indeed, $\Lambda$ is the subset of cardinality $\beta$ of the set of unassigned contact points of $C$ with $L$ around which two local branches are maintained in the deformation of $C \cup L$. Moreover, two local branches are maintained at all assigned contact points of $C$ with $L$. Therefore, the contact points of $C$ with $L$ at which the two local branches deform to a unique irreducible branch are exactly the points $s_{i,j}^{C \cap L}$ in the notation introduced at the beginning of Section 4.2.1, which deform to $i-1$ nodes only; those form the complement of $\Lambda$ in the set of unassigned contact points of $C$ with $L$.

The other ingredient of the proof of the recursive formula for irreducible curves is the following basic enumerative theoretic computation.

**(5.3)** Let $V$ be the sum of irreducible components of $V^{d,\delta}(\alpha,\beta)$ that parametrizes all the curves that decompose as unions of irreducible curves $C_1,\ldots,C_k$, such that for all $l = 1,\ldots,k$, the curve $C_l$ is a general member of an irreducible component of $V^{d_l,\delta_l}(\alpha^l,\beta^l)$. For all $l$, let $V_l = V_{\mathrm{irr}}^{d_l,\delta_l}(\alpha^l,\beta^l)$ be the sum of irreducible components of $V^{d_l,\delta_l}(\alpha^l,\beta^l)$ that parametrizes all irreducible members of $V^{d_l,\delta_l}(\alpha^l,\beta^l)$. Then, $V$ is the image of the product $\prod_{l=1}^k V_l$ by the Segre map

$$\Sigma : \prod_{l=1}^k |\mathcal{O}_{\mathbf{P}^2}(d_l)| \longrightarrow |\mathcal{O}_{\mathbf{P}^2}(d)|.$$

The degree of $V$ is the intersection number $V \cdot H^{\dim(V)}$, where $H$ denotes the hyperplane class in the Chow ring of $|\mathcal{O}_{\mathbf{P}^2}(d)|$. It may be computed by pulling-back by $\Sigma$. One has $\Sigma^* H = H_1 + \cdots + H_k$, where $H_l$ is the hyperplane class of $|\mathcal{O}_{\mathbf{P}^2}(d_l)|$ for all $l$, and then

$$\left(\prod_{l=1}^k V_l\right) \cdot (\Sigma^* H)^{\dim(V)} = \left(\prod_{l=1}^k V_l\right) \cdot (H_1 + \cdots + H_k)^{\dim(V)}$$

$$= \left(\prod_{l=1}^k V_l\right) \cdot \sum_{\nu_1 + \cdots + \nu_k = \dim(V)} \binom{\dim V}{\nu_1, \ldots, \nu_k} H_1^{\nu_1} \cdots H_k^{\nu_k}$$

$$= \binom{\dim V}{\dim V_1, \ldots, \dim V_k} \prod_{l=1}^k V_l \cdot H_l^{\dim(V_l)},$$

by computing in the Chow ring of $\prod_{l=1}^k |\mathcal{O}_{\mathbf{P}^2}(d_l)|$, as explained in [3, Section 2.1.4].

The above number equals $V \cdot H^{\dim(V)}$ if the restriction of $\Sigma$ to $\prod_{l=1}^k V_l$ is birational, which happens if and only if the $V_l$'s are pairwise distinct; otherwise the restriction of $\Sigma$ has finite degree, equal to the order of the group of permutations of the factors that leaves the product unchanged. The upshot is that, denoting by $\sigma$ the order of this group (we shall in fact use a refined version of this in the recursion formula),

$$\deg(V) = \frac{1}{\sigma} \binom{\dim V}{\dim V_1, \ldots, \dim V_k} \prod_{l=1}^k \deg(V_l).$$

For the application below, recall that the dimensions of generalized Severi varieties is given in Proposition (2.2). Now, fasten your seatbelts, the formula is the following.

**(5.4) Theorem.** *For all $d, \delta, \alpha, \beta$ as in Corollary (1.5), let $N_{\mathrm{irr}}^{d,\delta}(\alpha, \beta)$ be the number of irreducible members of the generalized Severi variety $V^{d,\delta}(\alpha,\beta)(\Omega)$ passing through $2d+g-1+|\beta|$ general points of the plane ($g = p_a(d) - \delta$), for any general $\Omega$. Then,*

$$N_{\mathrm{irr}}^{d,\delta}(\alpha, \beta) = \sum_{k \geqslant 1:\, \beta_k > 0} k \cdot N_{\mathrm{irr}}^{d,\delta}(\alpha + e_k, \beta - e_k)$$
$$+ \sum \left[ \frac{1}{\sigma} \binom{2d+g-1+|\beta|}{2d_1+g_1-1+|\beta^1|,\, \ldots,\, 2d_k+g_k-1+|\beta^k|} \cdot \binom{\alpha}{\alpha^1, \ldots, \alpha^k, \alpha - \alpha'} \right.$$
$$\left. \cdot \prod_{l=1}^{k} \binom{\beta^l}{\beta^l - \gamma^l} \cdot \prod_{l=1}^{k} I^{\gamma^l} \cdot \prod_{l=1}^{k} N_{\mathrm{irr}}^{d_l,\delta_l}(\alpha^l, \beta^l) \right],$$

*where the second sum is taken over all integers $k > 0$, and over all collections of integers $d_1, \ldots, d_k$, and $\delta_1, \ldots, \delta_k$, and all collections of sequences $\alpha^1, \ldots, \alpha^k, \beta^1, \ldots, \beta^k$, and $\gamma^1, \ldots, \gamma^k$, subject to the relations*

$$\alpha' = \alpha^1 + \cdots + \alpha^k \leqslant \alpha$$
$$\beta' = \beta^1 + \cdots + \beta^k = \beta + \gamma^1 + \cdots + \gamma^k$$
$$\forall l = 1, \ldots, k: \qquad \gamma^l \neq 0$$
$$d_1 + \cdots + d_k = d - 1$$
$$\delta_1 + \cdots + \delta_k = \delta + |\beta' - \beta| - d + 1 - \sum_{h < l} d_h d_l,$$

*and $\sigma$ is defined as follows: define an equivalence relation on $\{1, \ldots, k\}$ by declaring $h \sim l$ if*

$$d_h = d_l, \quad \delta_h = \delta_l, \quad \alpha^h = \alpha^l, \quad \beta^h = \beta^l, \quad \text{and} \quad \gamma^h = \gamma^l;$$

*then $\sigma$ is the product of the factorials of the cardinalities of the equivalence classes.*

In the recursion formula, we have used multinomials for elements of $\underline{\mathbf{N}}$; we define them as

$$\binom{\alpha}{\alpha^1, \ldots, \alpha^k, \alpha^0} = \prod_{i \geqslant 1} \binom{\alpha_i}{\alpha_i^1, \ldots, \alpha_i^k, \alpha_i^0}.$$

The interpretation of the first sum in the formula should be fairly clear; let us decipher the second sum. Each summand corresponds to the contribution of

$$\prod_{l=1}^{k} V_{\mathrm{irr}}^{d_l,\delta_l}(\alpha^l, \beta^l) \subseteq V^{d_1+\cdots+d_k,\, \delta_1+\cdots+\delta_k + \sum_{h<l} d_h d_l}\left(\sum_l \alpha^l, \sum_l \beta^l\right)$$

to $V_{\mathrm{irr}}^{d,\delta}(\alpha, \beta) \cap p^\perp$ determined by the sequences $\gamma^1, \ldots, \gamma^k \in \underline{\mathbf{N}}$ in the following way: for all $l = 1, \ldots, k$, $\gamma^l \in \underline{\mathbf{N}}$ is the number of unassigned contact points of the general member $C_l$ of $V_{\mathrm{irr}}^{d_l,\delta_l}(\alpha^l, \beta^l)$ at which an open neighbourhood of $C_l \cup L$ deforms to only one irreducible branch; the deformation of $C_1 \cup \cdots \cup C_k \cup L$ is irreducible if and only if there is at least one such point for all $l$ (see Paragraph (5.2) above), which corresponds to the requirement that $\gamma^l \neq 0$ for all $l$. Note that the number of nodes of $C_1 \cup \cdots \cup C_k$ is indeed $\delta_1 + \cdots + \delta_k + \sum_{h<l} d_h d_l$.

This contribution in fact means the sum of the contributions of all $\prod_{l=1}^{k} V_{\mathrm{irr}}^{d_l,\delta_l}(\alpha^l, \beta^l)(\Omega^l)$, for all sets $\Omega^1, \ldots, \Omega^k$ of cardinalities $\alpha^1, \ldots, \alpha^k \in \underline{\mathbf{N}}$ such that $\Omega = \bigcup_{l=1}^{k} \Omega^l$; the multinomial $\binom{\alpha}{\alpha^1, \ldots, \alpha^k, \alpha - \alpha'}$ corresponds to the number of possible $\Omega^1, \ldots, \Omega^k$. On the other hand, the factor $\frac{1}{\sigma}$ times the first multinomial is the degree of $\prod_{l=1}^{k} V_{\mathrm{irr}}^{d_l,\delta_l}(\alpha^l,\beta^l,\gamma^l)$ (see Paragraph (5.3)); note that we have added the decoration $\gamma^l$ (with the meaning explained above), and the "stabilizer term" $\sigma$ has been changed accordingly in the theorem with respect to Paragraph (5.3).

We leave to the reader the excellent exercise of filling out the details of the proof of Theorem (5.4).

# References

**Chapters in this Volume**

[I]     C. Ciliberto and Th. Dedieu, *Limits of nodal curves: generalities and examples.*

[II]    Th. Dedieu and E. Sernesi, *Deformations of curves on surfaces.*

[III]   Th. Dedieu, *Geometry of logarithmic Severi varieties at a general point.*

[IV]    M. Lelli-Chiesa, *Deformations of an m-tacnode into m − 1 nodes.*

[V]     F. Bastianelli, *Products of deformation spaces of tacnodes.*

[A]     Th. Dedieu, *Geometry of the deformation space of an m-tacnode into m − 1 nodes.*

[VII]   Th. Dedieu, *The Caporaso–Harris recursion within a degeneration of the plane.*

[**??**]  Th. Dedieu, *Some classical formulæ for curves and surfaces.*

**External References**

# Lecture VII
# The Caporaso–Harris recursion within a degeneration of the plane

## by Thomas Dedieu

In this text we show on examples how the Caporaso–Harris recursion procedure [1], explained in details in [VI], can be formulated in terms of a degeneration of the projective plane to the union of a projective plane and a rational ruled surface $\mathbf{F}_1$ intersecting transversely along a line.

We formulate the degeneration procedure in a general framework in Section 1, and work it out explicitly on two examples in Sections 2 and 3, namely for the enumerations of 2-nodal plane quartics, and of 1-nodal cubics tangent to a line, respectively. In Section 4 we follow a different degeneration procedure, by way of comparison, in order to enumerate 2-nodal plane quartics.

We assume the reader is familiar with the two chapters [I] and [VI], and will freely use terminology and notation from them.

## 1 – The degeneration procedure

**(1.1)** Consider the projective plane $\mathbf{P}^2$ and a line $R$ in it. Let $d, \delta$ be integers and $\alpha, \beta \in \underline{\mathbf{N}}$ be sequences such that $d > 0$, $\delta \geqslant 0$, and $I\alpha + I\beta = d$. Let $\Omega$ be a (sequence of) set(s) of $\alpha$ points on the line $R$. Let $Z$ be the union of $2d + p_a(d) - \delta - 1 + |\beta|$ general points in the plane, where $p_a(d) = \frac{1}{2}(d-1)(d-2)$. We want to count the (finitely many) plane curves of degree $d$ with $\delta$ nodes, intersecting $R$ at the points of $\Omega$ with the multiplicities prescribed by $\alpha$ and at further unassigned points with the multiplicities prescribed by $\beta$, and passing through all points of $Z$. The result is the number $N^{d,\delta}(\alpha, \beta)$ in the notation of [VI] and [1].

**(1.2)** Let $\mathcal{S}$ be the blow-up of $\mathbf{P}^2 \times \mathbf{A}^1$ along the line $R \times \{0\}$, and let $\pi : \mathcal{S} \to \mathbf{A}^1$ be the composition of the blow-up map with the second projection. For $t \neq 0$, the fibre $S_t$ of $\pi$ over $t$ is a projective plane, whereas the central fibre $S_0$ is the union of a projective plane and a rational ruled surface $\mathbf{F}_1$, which we shall call $\mathbf{P}$ and $\mathbf{F}$ respectively, intersecting transversely along a curve $E$ which is a line in $\mathbf{P}$, and the section of the ruling with self-intersection $-1$ in $\mathbf{F}$; we shall write $S_0 = \mathbf{P} \cup_E \mathbf{F}$.

Let $\mathcal{L}$ be the line bundle on $\mathcal{S}$ obtained by pulling-back $\mathcal{O}_{\mathbf{P}^2}(1)$ successively by the first projection $\mathbf{P}^2 \times \mathbf{A}^1 \to \mathbf{P}^2$ and by the blow-up map. The line bundle $\mathcal{L}$ restricts to $\mathcal{O}_{\mathbf{P}}(H)$ and $\mathcal{O}_{\mathbf{F}}(F)$ on $\mathbf{P}$ and $\mathbf{F}$ respectively, where $H$ is the line class and $F$ is the class of the ruling.

Let $\mathcal{R}$ be the proper transform of $R \times \mathbf{A}^1$ in $\mathcal{S}$. For $t \neq 0$, the fibre $R_t$ of $\mathcal{R}$ over $t \in \mathbf{A}^1$ is the line $R$ in $\mathbf{P}^2$, whereas the central fibre $R_0$ is a section of the ruling of $\mathbf{F}$ that has self-intersection $+1$. Let $\mathcal{W}$ be the proper transform of $\Omega \times \mathbf{A}^1$, and denote by $\Omega_0$ its central fibre: this is merely $\Omega$ seen on the section $R_0 \cong \mathbf{P}^1$ of $\mathbf{F}$.

In order to perform the Caporaso–Harris degeneration procedure, we let the set of points $Z$ degenerate as follows. Let $p$ be one point of $Z$, and consider a general section $\sigma_p$ of $\mathbf{P}^2 \times \mathbf{A}^1 \to \mathbf{A}^1$ which intersects $R \times \{0\}$, and whose fibre over $1 \in \mathbf{A}^1$ coincides with $p$. Then we consider the proper transform $\mathcal{Z}$ in $\mathcal{S}$ of $(Z - p) \times \mathbf{A}^1 + \sigma_p$; its fibre over $1 \in \mathbf{A}^1$ is $Z \subseteq \mathbf{P}^2$, whereas its central fibre is the sum of $Z - p \subseteq \mathbf{P}$ and a general point $p_0$ of $\mathbf{F}$.

**(1.3)** We consider the family of surfaces $\pi : \mathcal{S} \to \mathbf{A}^1$, equipped with the sheaf $\mathcal{L}^{\otimes d}$. It provides a degeneration to $S_0 = \mathbf{P} \cup_E \mathbf{F}$ of the linear system $|dH|$ on $\mathbf{P}^2$. Moreover, $\mathcal{R}$ provides a degeneration of the line $R \subseteq \mathbf{P}^2$, and $\mathcal{W}$ and $\mathcal{Z}$ provide degenerations of $\Omega$ and $Z$ respectively. We want to solve the enumerative problem stated in (1.1) by considering its limit within this degeneration.

Let $\mathcal{I}_{\mathcal{Z}/\mathcal{S}}$ be the ideal sheaf defining $\mathcal{Z}$ in $\mathcal{S}$, and let $\mathcal{I}_{\mathcal{W}}$ be the ideal sheaf encoding the contact conditions with $\mathcal{R}$ at $\mathcal{W}$ corresponding to $\alpha \in \underline{\mathbf{N}}$. It follows from the results in [VI], which are taken from [1], that the pair $(\mathcal{S}/\mathbf{A}^1, \mathcal{L}^{\otimes d} \otimes \mathcal{I}_{\mathcal{Z}/\mathcal{S}} \otimes \mathcal{I}_{\mathcal{W}})$ is $\delta$-well behaved, in the terminology of [I, Section 5]. Thus, our enumerative problem stated in (1.1) can be solved by computing the regular part of the limit Severi variety

$$\mathfrak{V}_\delta(\mathcal{S}/\mathbf{A}^1, \mathcal{L}^{\otimes d} \otimes \mathcal{I}_{\mathcal{Z}/\mathcal{S}} \otimes \mathcal{I}_{\mathcal{W}});$$

in other words, the number we are looking for equals the number of curves in the linear systems on $\mathbf{P} \cup \mathbf{F}$ corresponding to all possible twists of $\mathcal{L}^{\otimes d}$, passing through $Z - p$ on $\mathbf{P}$ and through $p_0$ on $\mathbf{F}$, with contact conditions with $R_0$ determined by $\alpha \in \underline{\mathbf{N}}$ at the points $\Omega_0$ and by $\beta \in \underline{\mathbf{N}}$ at unassigned points, and with $\delta_{\mathbf{P}}$ nodes on $\mathbf{P}$, $\delta_{\mathbf{F}}$ nodes on $\mathbf{F}$, and $\tau$ tacnodes along $E = \mathbf{P} \cap \mathbf{F}$ ($\tau \in \underline{\mathbf{N}}$), such that

(1.3.1) $$\delta = \delta_{\mathbf{P}} + \delta_{\mathbf{F}} + \nu(\tau)$$

with $\nu(\tau) = \sum_{m \geqslant 2} \tau_m (m-1)$; it is understood that these curves are counted with the multiplicity $\mu(\tau) = \prod_{m \geqslant 2} m^{\tau_m}$.[1]

The framework of [I] does not allow to formally take in consideration the contact conditions at unprescribed points; however this appendix is intended to be illustrative, so I will remain informal. Moreover, the framework of [I] would require to have an invertible sheaf instead of $\mathcal{L}^{\otimes d} \otimes \mathcal{I}_{\mathcal{Z}/\mathcal{S}} \otimes \mathcal{I}_{\mathcal{W}}$; this can be easily remedied by considering a suitable blow-up of $\mathcal{S}$, but I prefer not to do it here.

## 2 – Two-nodal quartics

In this section we enumerate 2-nodal plane quartics by repeatedly applying the procedure described in Section 1.

---

[1] Beware the difference in notation between $\alpha = (\alpha_m)_{m \geqslant 1}$ and $\beta = (\beta_m)_{m \geqslant 1}$ on the one hand, and $\tau = (\tau_m)_{m \geqslant 2}$ on the other hand: while the former encode all contact points with $R_0$, including the transverse ones, the latter only encodes $m$-tacnodes with $m > 1$. One can extend $\tau$ to $\tilde{\tau}$ to put it in $\alpha/\beta$ style, by setting $\tau_1 = \deg(\mathcal{L}^{\otimes d}(-W)\big|_E) - \sum_{m \geqslant 2} m \tau_m$ and $\tilde{\tau} = (\tau_1, \tau)$, where $\mathcal{L}^{\otimes d}(-W)$ is the twist of $\mathcal{L}^{\otimes d}$ under consideration; then, $\nu(\tau) = I\tilde{\tau} - |\tilde{\tau}|$, and $\mu(\tau) = I^{\tilde{\tau}}$.

**(2.1)** Our initial problem is to count 2-nodal plane quartics passing through 12 general points. In the notation of (1.1), we consider $d = 4$, $\delta = 2$, $\beta = 4$, and $\alpha = 0$. The result is the number $N^{4,2} = N^{4,2}(0,4)$.

**(2.2)** We will let our 12 points degenerate one by one to general points on $\mathbf{F}$, and we will need to consider the restriction to $\mathbf{P} \cup \mathbf{F}$ of all possible twists of $\mathcal{L}^{\otimes 4}$. The line bundle $\mathcal{L}^{\otimes 4}$ itself restricts to $\mathcal{O}_\mathbf{P}(4H)$ and $\mathcal{O}_\mathbf{F}(4F)$ on $\mathbf{P}$ and $\mathbf{F}$ respectively. Thus, members of the corresponding linear system on $S_0$ consist of a plane quartic on $\mathbf{P}$ and four fibres of the ruling of $\mathbf{F}$, which have to match along their intersections with $E = \mathbf{P} \cap \mathbf{F}$, which consist of four points.

The only other relevant twist will be $\mathcal{L}^{\otimes 4}(-\mathbf{F})$, which restricts to $\mathcal{O}_\mathbf{P}(3H)$ and $\mathcal{O}_\mathbf{F}(E + 4F)$ on $\mathbf{P}$ and $\mathbf{F}$ respectively. The linear systems $|3H|_\mathbf{P}$ and $|E + 4F|_\mathbf{F}$ have dimensions 9 and 8 respectively. Members of the linear system on $S_0$ corresponding to $\mathcal{L}^{\otimes 4}(-\mathbf{F})$ consist of a member of $|3H|_\mathbf{P}$ and a member of $|E+4F|_\mathbf{F}$ which match along their intersections with $E$, which consist of three points; the dimension of the linear system on $S_0$ is therefore $9+8-3 = 14$ (both $|3H|_\mathbf{P}$ and $|E + 4F|_\mathbf{F}$ restrict to the complete linear system $|3H|_E$ on $E \simeq \mathbf{P}^1$), equal to that of the linear system of all plane quartics.

To rule out the other twists, we will use (i) that $\mathcal{L}^{\otimes 4}(-a\mathbf{F})$ with $a > 1$ restricts to $\mathcal{O}_{\mathbf{P}^2}(4-a)$ on $\mathbf{P}^2$, hence the corresponding curves on $\mathbf{P}$ have degree at most 2, and (ii) that $\mathcal{L}^{\otimes 4}(a\mathbf{F})$ with $a > 0$ restricts to $\mathcal{O}_{\mathbf{F}_1}(-aE + 4F)$ on $\mathbf{F}$, which is not effective, so that it is already clear that the latter twists will never give any contribution to our enumerations.

Now, let us start the enumeration. The relations we get at each step are all spelled out in Paragraph (2.8).

**(2.3)** $N^{4,2}(0,4)$. Let us apply the procedure described in Section 1 to our initial enumerative problem: we end up with 11 general points on $\mathbf{P}$, and 1 general point on $\mathbf{F}$, as indicated on the figure below.

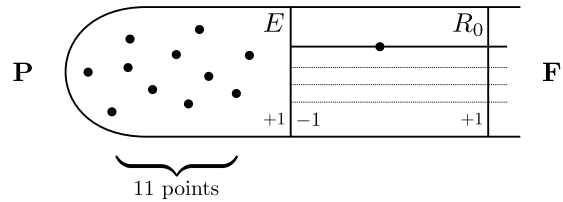

Only the trivial twist is relevant, because there is no curve of degree $d' < 4$ on $\mathbf{P}$ passing through the 11 base points. Thus we are only considering curves consisting of a plane quartic $C_\mathbf{P}$ on $\mathbf{P}$ and of four rulings of $\mathbf{F}$. One of these rulings must pass through the base point on $\mathbf{F}$, and has a fixed intersection point with $E$, which imposes a fixed passing point to $C_\mathbf{P}$ on $E$. On the other hand, the only way to fulfil Relation (1.3.1) with $\delta = 2$ with the curves under consideration is to have $\delta_\mathbf{P} = 2$ (and $\delta_\mathbf{F} = 0$ and $\tau = 0$). The upshot is that $N^{4,2}(0,4) = N^{4,2}(1,3)$.

**(2.4)** $N^{4,2}(1,3)$. We now apply the procedure of Section 1 to the enumerative problem corresponding to the number $N^{4,2}(1,3)$. We end up with 10 general points of $\mathbf{P}$, 1 general point of $R_0$, and 1 general point of $\mathbf{F}$, as pictured below.

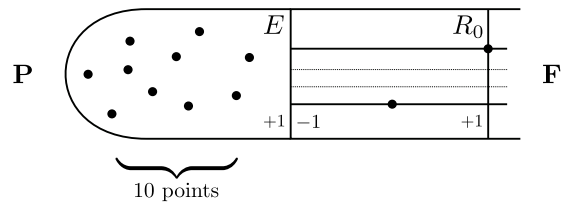

Again, only the trivial twist is relevant, and one must have $\delta_{\mathbf{P}} = 2$. This time, two of the rulings on $\mathbf{F}$ must pass through the two base points on $\mathbf{F}$, and they each impose one fixed passing point on $E$. The upshot is that $N^{4,2}(1,3) = N^{4,2}(2,2)$.

**(2.5)** $N^{4,2}(2,2)$. We repeat the procedure to compute this number. We end up with 9 general points of $\mathbf{P}$, 2 general point of $R_0$, and 1 general point of $\mathbf{F}$, as pictured below.

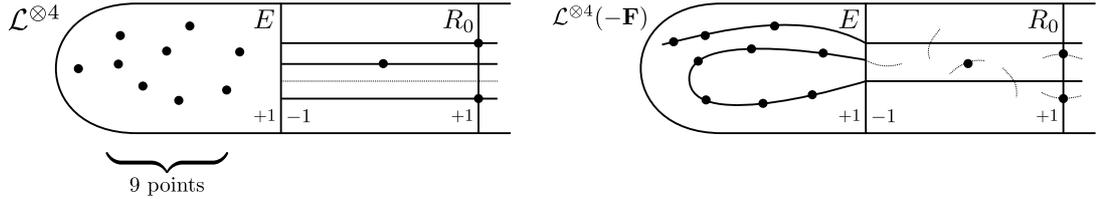

This time the two twists $\mathcal{L}^{\otimes 4}$ and $\mathcal{L}^{\otimes 4}(-\mathbf{F})$ will contribute. The trivial twist gives a contribution of $N^{4,2}(3,1)$ as in the previous cases. The contribution of the other twist comes from curves made of a cubic $C_{\mathbf{P}}$ on $\mathbf{P}$ and a curve $C_{\mathbf{F}}$ on $\mathbf{F}$ linearly equivalent to $E + 4F$. The cubic $C_{\mathbf{P}}$ is uniquely determined by the 9 base points on $\mathbf{P}$ (if one will, $N^{3,0}(0,3) = 1$), it is smooth and fixes 3 pairwise distinct passing points on $E$ for $C_{\mathbf{F}}$. Therefore, the curve $C_{\mathbf{F}}$ must be 2-nodal. Since $C_{\mathbf{F}} \cdot F = 1$, this implies that it must be of the form $F_1 + F_2 + C'$ with $F_1$ and $F_2$ two rulings, and $C' \sim E + 2F$. The two rulings must take up two of the fixed points on $E$, which gives $\binom{3}{2}$ possibilities, and then $C'$ is uniquely determined by the third point on $E$ and the base points on $\mathbf{F}$, as the linear system $|E + 2F|$ has dimension 4.

**(2.6)** $N^{4,2}(3,1)$. We repeat the procedure, and end up with 8 general points of $\mathbf{P}$, 3 general point of $R_0$, and 1 general point of $\mathbf{F}$, as pictured below.

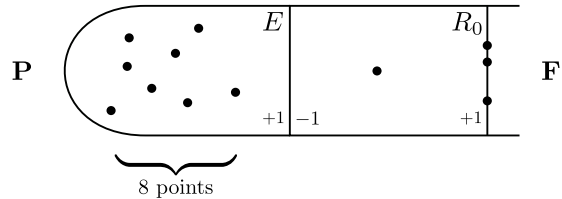

The trivial twist $\mathcal{L}^{\otimes 4}$ contributes by $N^{4,2}(4,0)$ as in the previous cases, and this time the twist $\mathcal{L}^{\otimes 4}(-\mathbf{F})$ contributes in several different ways: on the $\mathbf{P}$ side we have cubics through the 8 base points, which thus move in a pencil, and cut out a $g_3^1$ on $E$, while on the $\mathbf{F}$ side we have curves in $|E + 4F|$ through the $3 + 1$ base points, which impose independent linear conditions, so that we end up with a 4-dimensional linear system. Therefore we have the following three possibilities to fulfil (1.3.1) with $\delta = 2$.

(a) $\delta_{\mathbf{P}} = \delta_{\mathbf{F}} = 1$. The number of 1-nodal cubics through the 8 points on $\mathbf{P}$ is $N^{3,1}(0,3)$, and each of these fixes 3 points on $E$. In turn, the curve on $\mathbf{F}$ has to split as $F_1 + C'$ with $F_1 \sim F$ and $C' \sim E + 3F$ in order to be 1-nodal; moreover, in order to match with the curve on $\mathbf{P}$, the ruling $F_1$ must pass through one of the 3 fixed points on $E$. Since $|E + 3F|$ has dimension 6, the curve $C_1$ is then uniquely determined by the two other fixed points on $E$ and the $3 + 1$ base points on $\mathbf{F}$. The number of curves of this kind is thus $\binom{3}{1} N^{3,1}(0,3)$.

(b) $\tau_2 = 1$ (i.e., one 2-tacnode) and $\delta_{\mathbf{F}} = 1$. The number of cubics through the 8 points on $\mathbf{P}$ that are tangent to $E$ is $N^{3,0}(0, [1, 1])$ (this number is the degree of the dual surface of a twisted cubic, which is 4). The curve on $\mathbf{F}$ has to split as $F_1 + C'$ as in the previous case, and this time the ruling $F_1$ can only pass through the simple point of the divisor fixed on $E$ by the curve on $\mathbf{P}$. Then, the curve $C_1$ is uniquely determined by its passing through the double point of the divisor on $E$ and the $3+1$ base points on $\mathbf{F}$. The number of curves of this kind is thus $N^{3,0}(0,[1,1])$, to be counted with multiplicity 2 because of the tacnodal contribution in (1.3.1).

(c) $\delta_{\mathbf{F}} = 2$. In this case the curve on $\mathbf{F}$ has to split as $F_1 + F_2 + C'$ with $F_1, F_2 \sim F$ and $C' \sim E + 2F$. One of $F_1$ and $F_2$, say $F_1$, must pass through one of the three base points on $R_0$, and $C'$ must pass through the other two: this gives $\binom{3}{1}$ choices. Then $F_1$ fixes a point on $E$, which determines a unique cubic through the 8 base points on $\mathbf{P}$ (which we can formulate as $N^{3,0}(1,2) = 1$), which fixes two other points on $E$. The ruling $F_2$ must pass through one of them, which gives $\binom{2}{1}$ possibilities, after what the curve $C'$ is uniquely determined by the third point on $E$, the two remaining points on $R_0$, and the general point on $\mathbf{F}$. The number of curves of this kind is thus $\binom{3}{1}\binom{2}{1}N^{3,0}(1,2)$.

**(2.7)** $N^{4,2}(4,0)$. We repeat the procedure, and end up with 7 general points of $\mathbf{P}$, 4 general point of $R_0$, and 1 general point of $\mathbf{F}$, as pictured below.

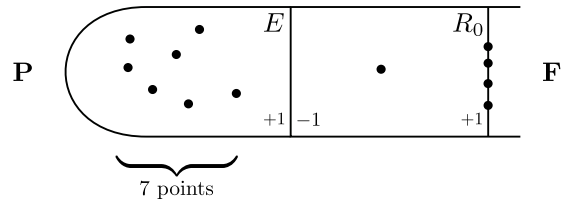

The trivial twist $\mathcal{L}^{\otimes 4}$ no longer contributes, because there is no curve in $\mathbf{F}$ linearly equivalent to $4F$ and passing through the $4+1$ base points on $F$. This is thus the point in the recursion where "contributions of degree 4 curves have been entirely exhausted.". For the contribution of the twist $\mathcal{L}^{\otimes 4}(-\mathbf{F})$, we have the following possibilities.

(a) $\delta_{\mathbf{P}} = 2$. The number of 2-nodal cubics through the 7 points on $\mathbf{P}$ is $N^{3,2}(0,3)$ (which equals $\binom{7}{2}$ because all such cubics are decomposed as a conic plus a line). For each such cubic, there is a unique matching curve in $|E + 4F|$ on $\mathbf{F}$ passing through the $4+1$ base points.

(b) $\delta_{\mathbf{P}} = \delta_{\mathbf{F}} = 1$. The curve on $\mathbf{F}$ must split as $F_1 + C'$, with $F_1 \sim F$ and $C' \sim E + 3F$. The ruling $F_1$ must pass through one of the base points on $R_0$ (which gives 4 possibilities), and therefore fixes a point on $E$. Then, there is a finite number $N^{3,1}(1,2)$ of 1-nodal cubics on $\mathbf{P}$ passing through this point and the 7 base points on $\mathbf{P}$ (the number $N^{3,1}(1,2)$ is easily seen to equal $N^{3,1}(0,3)$, which is 12). In turn, each of these fixes two new points on $E$, and the curve $C'$ is uniquely determined by these two points and the remaining $3+1$ base points on $\mathbf{F}$. Curves of this kind are thus in the number $\binom{4}{1}N^{3,1}(1,2)$.

(c) $\delta_{\mathbf{F}} = 2$. The curve on $\mathbf{F}$ must split as $F_1 + F_2 + C'$, with $F_1, F_2 \sim F$ and $C' \sim E + 2F$, and the two rulings $F_1$ and $F_2$ each must pass through one of the 4 base points on $R_0$, which gives $\binom{4}{2}$ possibilities. These two rulings then fix a unique cubic on $\mathbf{P}$ through the 7 points (which we can phrase as $N^{3,0}(2,1) = 1$), which in turn fixes a third point on $E$. The curve $C'$ is then uniquely determined by the latter point and the remaining $2+1$ base points on $\mathbf{F}$. We thus find $\binom{4}{2}N^{3,0}(2,1)$ curves.

(d) $\delta_{\mathbf{P}} = 1$ and $\tau_2 = 1$ (i.e., one 2-tacnode). There is a finite number of possible curves on $\mathbf{P}$, in the number $N^{3,1}(0,[1,1])$, and each of these determines a unique matching curve in $|E + 4F|$ on $\mathbf{F}$. These curves count with multiplicity 2.

(e) $\delta_{\mathbf{F}} = 1$ and $\tau_2 = 1$. The curve on $\mathbf{F}$ must split as $F_1 + C'$, with $F_1 \sim F$ and $C' \sim E + 3F$. The ruling $F_1$ must pass through one of the base points on $R_0$ (which gives 4 possibilities), and therefore fixes a point on $E$. Then, there is a pencil of cubics on $\mathbf{P}$ passing through this point and the 7 base points on $\mathbf{P}$; it cuts out (residually) a $g_2^1$ on $E$, which has two members made of a double point, meaning that there are two cubics in the pencil that are tangent to $E$: in other words, $N^{3,0}(1,[0,1]) = 2$. Finally, for each such cubic there is a unique matching curve $C'$ which also pass through the $3 + 1$ remaining base points on **F**. Thus we find $\binom{4}{1} N^{3,0}(1, [0, 1])$ curves; they count with multiplicity 2.

(f) $\tau_3 = 1$. There is a finite number of possible curves on **P**, which is $N^{3,0}(0, [0, 0, 1])$. This number equals the number of divisors of the form $3q$ in the $g_3^2$ cut out on $E$ by the 2-dimensional linear system of cubics through the 7 base points on **P**, which is 3 (this is the number of flexes of a 1-nodal cubic, see [**??**, Section **??**]). For each such cubic there is a unique matching curve in $|E + 4F|$ on **F** which also passes through the $4 + 1$ base points. We thus find $N^{3,0}(0, [0, 0, 1])$ curves, and they count with multiplicity 3.

**(2.8) Summary and conclusion.** We have applied the degeneration procedure of Section 1 five times, and after that we are reduced to enumerations of curves of degrees $d' < 4$. The relations we have obtained are the following:

$$\begin{aligned}
N^{4,2}(0,4) &= N^{4,2}(1,3) \\
&= N^{4,2}(2,2) \\
&= N^{4,2}(3,1) + \binom{3}{2}\underbrace{N^{3,0}(0,3)}_{=1} \\
N^{4,2}(3,1) &= N^{4,2}(4,0) + \binom{3}{1}\underbrace{N^{3,1}(0,3)}_{=12} + 2\cdot\underbrace{N^{3,0}(0,[1,1])}_{=4} + \binom{3}{1}\binom{2}{1}\underbrace{N^{3,0}(1,2)}_{=1} \\
N^{4,2}(4,0) &= \underbrace{N^{3,2}(0,3)}_{=21} + \binom{4}{1}\underbrace{N^{3,1}(1,2)}_{=12} + \binom{4}{2}\underbrace{N^{3,0}(2,1)}_{=1} + 2\cdot\underbrace{N^{3,1}(0,[1,1])}_{=36} \\
&\quad + 2\cdot\binom{4}{1}\underbrace{N^{3,0}(1,[0,1])}_{=2} + 3\cdot\underbrace{N^{3,0}(0,[0,0,1])}_{=3}.
\end{aligned}$$

With the (possible) exception of $N^{3,1}(0, [1, 1])$, which is computed in Section 3 below, the numbers of cubic curves appearing in the relations above can all be more or less elementarily computed; hints for these computations have been given in the course of the recursive application of the degeneration procedure. In the upshot, we get

$$\begin{aligned}
N^{4,2}(4,0) &= 21 + 48 + 6 + 72 + 16 + 9 = 172 \\
N^{4,2}(3,1) &= 172 + 36 + 8 + 6 = 222 \\
N^{4,2}(0,4) = N^{4,2}(1,3) = N^{4,2}(2,2) &= 222 + 3 = 225.
\end{aligned}$$

## 3 – An ancillary enumeration of cubics

In this section we compute the number of 1-nodal plane cubics tangent to a line; the result is the number $N^{3,1}(0, [1, 1])$. This is necessary for the enumeration of 2-nodal quartics in the previous section, and will serve as an additional illustration of the recursive degeneration procedure. We will need to consider the line bundles $\mathcal{L}^{\otimes 3}$, which restricts to $\mathcal{O}_\mathbf{P}(3H)$ and $\mathcal{O}_\mathbf{F}(3F)$ on **P** and **F** respectively, and $\mathcal{L}^{\otimes 3}(-\mathbf{F})$, which restricts to $\mathcal{O}_\mathbf{P}(2H)$ and $\mathcal{O}_\mathbf{F}(E+3F)$ on **P** and **F** respectively. The generalized Severi variety $V^{3,1}(0, [1, 1])$ has dimension 7, so we start with 7 base points.

**(3.1)** $N^{3,1}(0, [1, 1])$**.** We apply the degeneration procedure of Section 1, which gives 6 general base points on **P**, and one general base point on **F**.

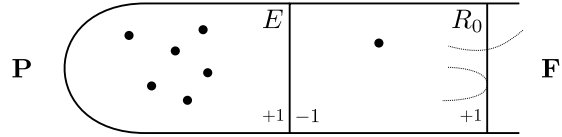

Only the trivial twist contributes, because there is no conic on **P** through the 6 base points, and we must have $\delta_\mathbf{P} = 1$ in order to fulfil (1.3.1) with $\delta = 1$. To have the correct intersection pattern with $R_0$, the curve on **F** must be of the form $F_1 + 2F_2$ with $F_1, F_2 \sim F$. Then, either $F_1$ passes through the base point on **F**, after what there are $N^{3,1}(1, [0,1])$ matching 1-nodal cubics through the 6 base points on **P**, or $F_2$ passes through the base point on **F**, and then there are $N^{3,1}([0,1], 1)$ matching 1-nodal cubics through the 6 base points on **P**; it follows from the analysis in [1], reported on in [VI], that the latter curves count with multiplicity 2.

**(3.2)** $N^{3,1}(1, [0,1])$. Applying the degeneration procedure, we get 5 general base points on **P**, one general base point on $R_0$, and one general point on **F**.

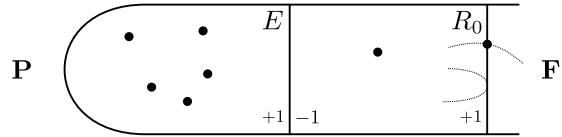

A priori, the two twists $\mathcal{L}^{\otimes 3}$ and $\mathcal{L}^{\otimes 3}(-\mathbf{F})$ contribute in this case. In the contribution of the trivial twist $\mathcal{L}^{\otimes 3}$, the curve on **F** must once again be of the form $F_1 + 2F_2$ with $F_1, F_2 \sim F$ in order to meet the contact condition with $R_0$. Then, $F_1$ must past through the base point on $R_0$, and $F_2$ through the general base point on **F**. There are then $N^{3,1}([1,1], 0)$ matching curves on **P**, and the curves we find in this way count with multiplicity 2 because of the assignement of the tangency point with $R_0$.

In the contribution of the twist $\mathcal{L}^{\otimes 3}(-\mathbf{F})$, we find a unique conic through the five base points on **P**. It is smooth and fixes two distinct points on $E$. The curve on **F** must be decomposed as $F_1 + C'$ with $F_1 \sim F$ and $C' \sim E + 2F$ in order to be 1-nodal. The ruling $F_1$ must pass through one of the two points on $E$, and therefore does not pass through the base point on $R_0$. It follows that it is impossible to choose $C'$ so that the contact conditions with $R_0$ are fulfilled (it is possible to choose $C'$ tangent to $R_0$ at the base point, but this does not give an admissible curve, see the comment on the definition of logarithmic Severi varieties in [III, (1.6)]; this curve will contribute in (3.3) however).

**(3.3)** $N^{3,1}([0,1], 1)$. We apply the degeneration procedure, and end up with 5 general base points on **P**, one general base point on **F**, and a fixed tangency point with $R_0$.

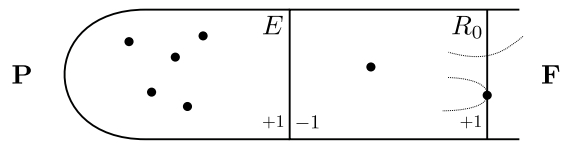

We will have contributions from the two twists $\mathcal{L}^{\otimes 3}$ and $\mathcal{L}^{\otimes 3}(-\mathbf{F})$. For the trivial twist $\mathcal{L}^{\otimes 3}$, the curve on **F** must be of the form $F_1 + 2F_2$ with $F_1$ and $F_2$ two rulings, with $F_2$ passing through the fixed point on $R_0$ and $F_1$ through the general base point on **F**. This curve fixes a tangency point and a simple passing point on $E$, and the possible cubic curves on **P** are in the number $N^{3,1}([1,1], 0)$.

In the situation defined by the twist $\mathcal{L}^{\otimes 3}(-\mathbf{F})$ the curve on **P** is a conic, and there is a unique one through the five base points (i.e., $N^{2,0}(0, 2) = 1$), which fixes two points on $E$. In order to be 1-nodal, the curve on **F** must be decomposed as $F_1 + C'$ with $F_1 \sim F$ and $C' \sim E + 2F$. The curve $F_1$ must pass through one of two points on $E$, which gives two possibilities, and then $C'$ must pass through the second point on $E$, through the base point which is general on **F**, and be tangent to $R_0$ at the prescribed point: these are 4 independent linear conditions, hence the curve $C'$ is uniquely determined. The contribution of $\mathcal{L}^{\otimes 3}(-\mathbf{F})$ is thus by $\binom{2}{1} N^{2,0}(0,2)$.

**(3.4)** $N^{3,1}([1,1],0)$**.** We apply the degeneration procedure and get 4 general base points on **P**, one general base point on **F**, and a fixed tangency plus a fixed simple intersection point with $R_0$.

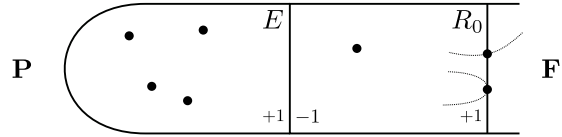

This time only the twist $\mathcal{L}^{\otimes 3}(-\mathbf{F})$ may contribute, because it is impossible for a curve in $|3F|$ to fulfil all the conditions on **F**. This is therefore the point in the recursion procedure where there only remain numbers $N^{d',\delta'}(\alpha',\beta')$ with $d' < 3$. Thus, we have a conic on **P**, which moves in a pencil determined by the four base points. We have the following possibilities to fulfil (1.3.1) with $\delta = 1$.

(a) $\delta_{\mathbf{P}} = 1$. There are $N^{2,1}(0,2)$ 1-nodal conics through the 4 base points on **P**, and each fixes two points on $E$ ($N^{2,1}(0,2) = 3$ because all such conics are decomposed in two lines). Then the curve on **F** must satisfy 6 independent linear conditions, and is thus uniquely determined.

(b) $\delta_{\mathbf{F}} = 1$. Then the curve on **F** must be of the form $F_1 + C'$ with $F_1 \sim F$ and $C' \sim E + 2F$. The ruling $F_1$ must pass through the fixed simple intersection point with $R_0$. This fixes a point on $E$, which determines a unique conic on **P** through the 4 base points ($N^{2,0}(1,1) = 1$), which in turn fixes a second point on $E$. Then the curve $C'$ must pass through this point, through the general base point on **F**, and be tangent with $R_0$ at the prescribed point: it is thus uniquely determined.

(c) $\tau_2 = 1$. The curve on **P** must be a conic tangent to $E$: these are in the number $N^{2,0}(0,[0,1])$ (which equals 2, the number of double points in the $g_2^1$ cut out on $E$). Then the curve on **F** has an imposed tangency with $E$ at a prescribed point, hence it must in total satisfy 6 independent linear conditions, and thus it is uniquely determined. These curve counts with the multiplicity 2.

**(3.5) Summary and conclusion.** In the four above applications of the degeneration procedure, we have obtained the following relations:

$$
\begin{aligned}
N^{3,1}(0,[1,1]) &= N^{3,1}(1,[0,1]) + 2 \cdot N^{3,1}([0,1],1) \\
N^{3,1}(1,[0,1]) &= 2 \cdot N^{3,1}([1,1],0) \\
N^{3,1}([0,1],1) &= N^{3,1}([1,1],0) + \binom{2}{1} \underbrace{N^{2,0}(0,2)}_{=1} \\
N^{3,1}([1,1],0) &= \underbrace{N^{2,1}(0,2)}_{=3} + \underbrace{N^{2,0}(1,1)}_{=1} + 2 \cdot \underbrace{N^{2,0}(0,[0,1])}_{=2}
\end{aligned}
$$
(3.5.1)

Putting these identities all together, we obtain $N^{3,1}([1,1],0) = 8$, $N^{3,1}([0,1],1) = 10$, $N^{3,1}(1,[0,1]) = 16$ and, finally, $N^{3,1}(0,[1,1]) = 36$.

**(3.6) Remark.** The numbers $N^{3,1}([1,1],0)$ and $N^{3,1}([0,1],1)$ may also be computed by appropriately enumerating singular cubics in a pencil, following [**??**, Section **??**], and doing so sheds some light on our degeneration procedure.

For the first number, we consider the pencil of plane cubics defined by 8 base points as follows: 3 of base points on a line $R$, 2 of which are infinitely near, and 5 general points of $\mathbf{P}^2$. There are 12 singular members (counted with multiplicities) in this pencil, as in any pencil of cubics, and the curves counted by $N^{3,1}([1,1],0)$ are all among them. Some singular members of the pencil must be excluded however, namely (a) the unique member of the pencil with a double point at the point where there are two infinitely near base points, and (b) the member of the pencil made of the line $R$ itself plus the unique conic through the 5 general base points.

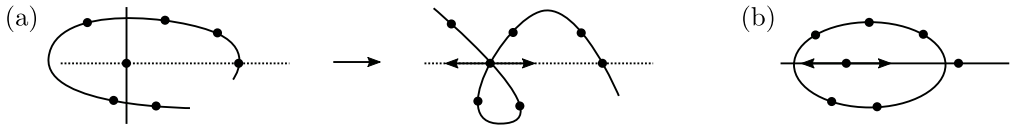

Both these curves have two nodes after we blow-up the base points of the pencil, so they contribute each by 2 to the number 12 of singular members in the pencil, and we thus find $N^{3,1}([1,1],0) = 12 - 2 - 2 = 8$.

For the second number, we put the 8 base points as follows: 2 infinitely near base points on the line $R$, and 6 general points of $\mathbf{P}^2$. This time there is no member of the pencil made of the line $R$ and a conic, and only the member as in case (a) above must be excluded. We thus find $N^{3,1}([1,1],0) = 12 - 2 = 10$.

The difference between these two enumerations explains Relation (3.5.1) above.

## 4 – Alternative degeneration procedure

In this section I present, by way of comparison, an alternative way of enumerating 2-nodal plane quartics by degeneration, which somehow packages all five degenerations in Paragraphs (2.3)–(2.7) in just one degeneration. This is the same degeneration procedure as in [I, Section 7.2.2], where 3-nodal quartics are enumerated, so I will be brief.

We consider always the same degeneration of the projective plane to the transverse union $\mathbf{P} \cup_E \mathbf{F}$, and this time we let 5 of the 12 base points on $\mathbf{P}^2$ degenerate to general points of $\mathbf{F}$.

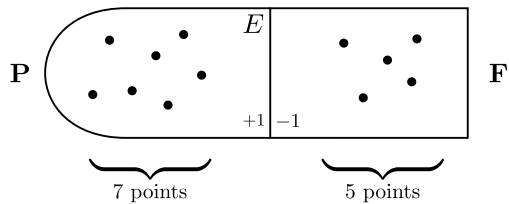

This prevents the trivial twist $\mathcal{L}^{\otimes 4}$ to give any contribution, and therefore the total number $N^{4,2}$ will be accounted for by the contributions from the twist $\mathcal{L}^{\otimes 4}(-\mathbf{F})$. In this twist, we have a net of plane cubics on $\mathbf{P}$ defined by the 7 general base points on $\mathbf{P}$, which cuts out a $g_3^2$ on $E$; on the other hand, the 5 general base points on $\mathbf{F}$ define a 3-dimensional linear subsystem of $|E + 4F|$, which cuts out a $g_3^3$ on $E$. We have the following possibilities to fulfil (1.3.1) with $\delta = 2$.

(a) $\delta_{\mathbf{P}} = 2$. There are $N^{3,2}(0,3) = \binom{7}{2}$ two-nodal cubics through the 7 base points on $\mathbf{P}$, and for each of those there is a unique matching member of $|E + 4F|$ through the 5 base points on $\mathbf{F}$.

(b) $\delta_{\mathbf{P}} = \delta_{\mathbf{F}} = 1$. The 1-nodal curve on $\mathbf{F}$ necessarily decompose as $F_1 + C'$ with $F_1 \sim F$ and $C' \sim E + 3F$. There are the two following combinatorial possibilities.

(i) The ruling $F_1$ passes through 1 base point and $C'$ through the other 4 ($\binom{5}{1} = 5$ possibilities). There are then $N^{3,1}(1,2) = 12$ 1-nodal cubics on $\mathbf{P}$ matching with $F_1$ along $E$, and for each of those there is a unique matching curve $C'$.

(ii) The ruling $F_1$ moves freely in the pencil $|F|$, while $C'$ passes through the 5 base points and thus moves in a pencil. The divisors cut out on $E$ by curves of this kind form a quadric surface in the $g_3^3$ cut out by the complete system $|E + 4F|$. It follows that the number of curves made of a curve of this kind and a nodal cubic on $\mathbf{P}$ through the 7 base points is $2 \times N^{3,1}(0,3) = 24$.

(c) $\delta_{\mathbf{F}} = 2$. The 2-nodal curve on $\mathbf{F}$ must decompose as $F_1 + F_2 + C'$ with $F_1, F_2 \sim F$ and $C' \sim E + 2F$, and there are the two following combinatorial possibilities.

(i) The two rulings $F_1$ and $F_2$ pass through 1 base point each, and $C'$ moves in the pencil defined by the 3 remaining base points ($\binom{5}{2} = 10$ possibilities). There is a cubic through the 7 base points and matching with $F_1$ and $F_2$ ($N^{3,0}(2,1) = 1$), and it cuts out a third point on $E$. Then, there is a unique curve $C'$ passing through this third point and the remaining 3 base points on $\mathbf{F}$.

(ii) The curves $F_1$ and $C'$ pass through 1 and 4 base points respectively ($\binom{5}{1}$ possibilities), thus both are fixed, and $F_2$ moves freely in $|F|$. As in case (i), there is a unique cubic on $\mathbf{P}$ matching with $F_1 + C'$, and in turn a unique ruling $F_2$ through the third intersection point of this cubic with $E$.

(d) $\delta_{\mathbf{P}} = 1$ and $\tau_2 = 1$. There are $N^{3,1}(0,[1,1]) = 36$ one-nodal cubics tangent to $E$, and for each of those there is a unique matching curve on $\mathbf{F}$. These curves count with multiplicity 2.

(e) $\delta_{\mathbf{F}} = 1$ and $\tau_2 = 1$ (with multiplicity 2). The curve on $\mathbf{F}$ decomposes as $F_1 + C'$ with $F_1 \sim F$, and there are the two following possibilities.

(i) The ruling $F_1$ passes through 1 base point ($\binom{5}{1}$ possibilities). There are then $N^{3,0}(1,[0,1]) = 2$ cubics on $\mathbf{P}$ tangent to $E$ and matching with $F_1$. Each of those determines a unique curve $C'$.

(ii) The ruling $F_1$ moves freely in $|F|$. There are finitely many curves $C'$ in $|E + 3F|$ passing through the 5 base points and tangent to $E$ (in the number 2, which is the number of double points in a $g_2^1$ on $E$). For each of those there is a unique matching cubic on $\mathbf{P}$ ($N^{3,0}([0,1],1) = 1$), which cuts out a third point on $E$ and thus fixes $F_1$.

(f) $\tau_3 = 1$ (with multiplicity 3). There are $N^{3,0}(0,[0,0,1]) = 3$ cubics on $\mathbf{P}$ triply tangent with $E$ at some point, and each of those determines a unique matching curve on $\mathbf{F}$.

These eventually add up as :

$$
\begin{array}{lr}
N^{3,2}(0,3) & = 21 \\
+ \binom{5}{1}N^{3,1}(1,2) + \binom{5}{0} \times 2 \times N^{3,1}(0,3) & + 5 \cdot 12 + 2 \cdot 12 \\
+ \binom{5}{2}N^{3,0}(2,1) + \binom{5}{1}N^{3,0}(2,1) & + 10 + 5 \\
+ 2 \cdot N^{3,1}(0,[1,1]) & + 2 \cdot 36 \\
+ 2 \cdot \binom{5}{1}N^{3,0}(1,[0,1]) + 2 \cdot \binom{5}{5} \times 2 \times N^{3,0}([0,1],1) & + 2 \cdot 5 \cdot 2 + 2 \cdot 2 \cdot 1 \\
+ 3 \cdot N^{3,0}(0,[0,0,1]) & + 3 \cdot 3 \\
\hline
N^{4,2}(0,4) & = 225.
\end{array}
$$

### References

**Chapters in this Volume**

# Lecture VIII
# Proving the existence of curves on smooth surfaces by smoothing tacnodes and other singularities on reducible surfaces

by Concettina Galati

**Abstract.** Let $\mathcal{X} \to \mathbf{D}$ be any flat projective family of complex surfaces parametrized by the disc $\mathbf{D}$, with smooth total space $\mathcal{X}$, smooth general fiber $\mathcal{X}_t$, and reducible special fiber $\mathcal{X}_0 = \cup X_i$ with normal crossing singularities. Let $C = \cup_i C_i \subseteq \mathcal{X}_0$ be a reducible Cartier divisor with only nodes on the smooth locus of $\mathcal{X}_0$ and tacnodes and nodes on the singular locus $E$ of $\mathcal{X}_0$. We will show how to study deformations of $C$ to nodal curves $C_t \subseteq \mathcal{X}_t$ by using Caporaso and Harris' local analysis of the versal deformation space of an $m$-tacnode in [1] and [2]. Results of this paper are known. We especially refer to [3] and [10]. But they are usually proved in literature with some special assumptions on the surfaces $\mathcal{X}_t$. Here we work in full generality. This paper is in particular a generalization of [10, Section 3].

## 1 – Introduction

Let $X$ be a smooth complex projective surface and let $D$ be a smooth Cartier divisor of genus $p_a(D)$. Describing the singular elements in the linear system $|D|$ defined by $D$ is not in general an easy problem, even if we restrict our interest to the existence of nodal divisors. There are various deformation arguments which may be helpful case by case. For example, if $X = \mathbf{P}^2$, we know that, for any $d \geqslant 3$, there exist degree $d$ irreducible nodal curves of any genus $0 \leqslant g \leqslant \frac{(d-1)(d-2)}{2}$, cf. [19] and [14]. These curves are deformations of a rational nodal plane curve of degree $d$, obtained as a general projection of a rational normal curve $C_d \subseteq \mathbf{P}^d$ of degree $d$. Notice that $\frac{(d-1)(d-2)}{2}$ is the number of nodes of a rational nodal plane curve and so it is the maximal number of nodes of an irreducible plane curve of degree $d$. If $X$ is any smooth complex projective surface and $D$ is any smooth Cartier divisor, it is not in general obvious how to produce irreducible or even reducible curves $|D|$ with "the maximum allowed number of nodes". One very classical strategy is to degenerate the surface $X \simeq \mathcal{X}_t$ to a reducible normal crossing surface $\mathcal{X}_0$ and to prove the existence of divisors on $X = \mathcal{X}_t$ with certain singularities as deformations of suitable Cartier divisors on $\mathcal{X}_0$. This method was used in particular by Chen in [3] to prove the existence of irreducible rational nodal curves on a very general K3 surface, in the primitive linear system and in all its multiples. The key idea of Chen's theorem is to construct, on a suitable reducible K3 surface $\mathcal{X}_0 = A \cup B$, a curve with tacnodes at points of the singular locus of $\mathcal{X}_0$ and to study its deformations on a very general K3 surface $\mathcal{X}_t$, by using the analysis of the locus of $(m-1)$-nodal curves in the versal deformation space of the $m$-tacnode, carried out in [1] and [2], see [IV]. The local analysis in [1] and [2] of versal deformations of tacnodes has later been extended in [10], where the authors provide a method to construct curves with $A_k$ singularities on smooth surfaces, for every $k \geqslant 1$. In particular [10] contains a new proof of the aforementioned Chen's theorem. Results in [10] have later been applied in [4] and [5].

In [10] the authors work in a family of regular surfaces. An output of the present paper is that this regularity assumption is unnecessary. Here we provide a generalization of [10, Section 3] using a simpler notation. We hope it will be clear to the reader how all results of this paper can be generalized to the study of deformations of reduced curves in a smooth total space with complete intersection singularities, see for example [11].

This paper is divided in four sections. In Section 2 we describe our degeneration argument and we state our deformation problem (Problem (2.2)). We moreover collect several local results from which, in Section 3, we will deduce our main result, Theorem (3.1). The latter proves in particular [I, Proposition (5.5)]; besides, [I] contains many applications of Theorem (3.1).

The degeneration argument which we are going to explain in the next section is useful to the study of many problems. For example, in [6], it has been applied to the study of the number of moduli of Severi varieties of nodal curves on K3 surfaces. Here we will restrict to the existence problem of nodal curves on smooth surfaces.

**(1.1) Terminology and notation.** We will work over $\mathbf{C}$. A curve will be a reduced separated scheme of finite type and dimension 1.

**(1.2) Acknowledgments.** I would like to thank Thomas Dedieu for reading the preliminar version of this paper.

## 2 – Our deformation problem

**(2.1)** Let $\mathcal{X} \to \mathbf{D}$ be a flat projective family of complex surfaces parametrized by the disc $\mathbf{D}$, with smooth total space $\mathcal{X}$, smooth general fiber $\mathcal{X}_t$, and reducible special fiber $\mathcal{X}_0$ with normal crossing singularities, i.e., locally at every point in the analytic topology, $\mathcal{X}_0$ is isomorphic to the closed subset of an open polydisk of $\mathbf{C}^3$ defined by $\{(x_1, x_2, x_3) | \Pi_{i=1}^k x_i = 0\}$, where $k$ may be $1, 2$, or $3$, cf. [8]. Thus every irreducible component $X_i$ of $\mathcal{X}_0 = \cup_{i=1}^r X_i$ is smooth, and the singular locus of $\mathcal{X}_0$ consists of the intersection locus of its irreducible components, where $\mathcal{X}_0$ has singularities of multiplicity 2 or 3. We denote by $E = \cup_{1 \leqslant i < j \leqslant r} E_{ij}$ the singular locus of $\mathcal{X}_0$, where $E_{ij} = X_i \cap X_j$ (possibly empty), and by $E^o \subseteq E \subseteq \mathcal{X}_0$ the set of points of multiplicity exactly 2.

Let $p \in E^o$ be a double point of $\mathcal{X}_0$. Working locally at $p$, we may assume that $\mathcal{X}$ is isomorphic to a closed subset of an open polydisk of $\mathbf{C}^3 \times \mathbf{D}$, with coordinates $(x, y, z, t)$, in such a way the local analytic equation of $\mathcal{X}$ at $p$ is given by

$$xy = t,$$

and $xy = t = 0$ are the analytic equations of $\mathcal{X}_0$ at $p$.

Now we consider a reduced Cartier divisor
$$C = \cup_i C_i \subseteq \mathcal{X}_0 = \cup_{i=1}^r X_i,$$
where every $C_i = C \cap X_i \subseteq X_i$ is a reduced $\delta_i$-nodal curve, that is smooth at every intersection point with $E$ and does not contain any triple point of $\mathcal{X}_0$. The local equations of $C$ at every intersection point $p \in E_{ij} \cap C = C_i \cap C_j$ are

(2.1.1) $$\begin{cases} y + x - z^m = 0 \\ xy = t \\ t = 0, \end{cases}$$

for some $m \geq 1$ (depending on $p$). Thus the singularity of $C$ at $p \in C \cap E_{ij}$ is analytically isomorphic to the planar curve singularity of local equation

(2.1.2) $$f(y,z) = (y - z^m)y = 0,$$

which is called a *tacnode of order m* or an *m-tacnode*; observe that a 1-tacnode is a node.

Thus the curve
$$C = \cup_i C_i \subseteq \mathcal{X}_0 = \cup_{i=1}^r X_i,$$
has only planar singularities that are nodes and tacnodes.

**(2.2) Problem.** *We want to find sufficient conditions for $C \subseteq \mathcal{X}_0$ to be deformed to a nodal curve $C_t \subseteq \mathcal{X}_t$. More precisely we want to understand which is the maximal number of nodes of a curve $C_t \subseteq \mathcal{X}_t$ that is deformation of $C \subseteq \mathcal{X}_0$.*

**(2.3) Remark.** Every tacnode of $C$ at a moving point $p \in E_{ij} = X_i \cap X_j \subseteq E$ imposes at most $m - 1$ independent conditions to the linear system $|\mathcal{O}_{\mathcal{X}_0}(C)|$, being $m - 1$ the expected number of conditions imposed to the linear system $|\mathcal{O}_{X_i}(C)|$ by a tangency of order $m$ with $E_{ij} \subseteq E$ at an unprescribed point. Thus the naive expectation is that one may deform $C$ to a curve $C_t \subseteq \mathcal{X}_t$ in such a way every $m$-tacnode of $C$ on $E$ deforms to $m - 1$ nodes of $C_t$ and all nodes of $C$ on $\mathcal{X}_0 \setminus E$ are preserved. In this case we will say that *an m-tacnode of $C$ is smoothed to $m - 1$ nodes of $C_t$* because $m - 1$ nodes impose one less condition to the arithmetic genus than an $m$-tacnode.

**(2.4) Remark.** Let us recall that an $m$-tacnode is an $A_{2m-1}$ planar curve singularity. An $A_k$ planar curve singularity has analytic equation $y^2 - x^{k+1}$. Every plane curve singularity of multiplicity 2 is an $A_k$-singularity, for some $k$ and two plane curve singularities of multiplicity two are topologically equivalent if and only if they are analytically equivalent. We refer to [13] and [II] for the notion of equisingularity from the topological and the analytic point of view for a family of reduced curves on a smooth surface. In particular, for an $A_k$-singularity the equisingular ideal coincides with the Jacobian ideal, given by $I = J = (y, x^k)$ [13, Proposition 5.6]. This implies that an $A_k$-singularity at a moving point of a smooth surface imposes as most $k = \dim(\mathbf{C}[x,y]/J)$ independent conditions to a linear system. Now, an $A_k$-singularity locally deforms to an $A_{k'}$-singularity for every $k' \leq k$. Thus, in the same vein as in Problem (2.2), one could ask about deformations $C_t \subseteq \mathcal{X}_t$ of $C \subseteq \mathcal{X}_0$ with $A_k$-singularities. As in Remark (2.3), one expects that the curve singularity (2.1.1) deforms on $\mathcal{X}_t$ to $d_i$ singularity of type $A_i$, for every multi-index $(d_1, ..., d_r)$ such that $\sum_i i d_i = m - 1$.

Deformations of $C$ in $\mathcal{X}$ are parametrized by suitable subschemes of the relative Hilbert scheme $\mathcal{H}^{\mathcal{X}|\mathbf{D}}$ of $\mathcal{X} \to \mathbf{D}$ [18]. We denote by $\mathcal{H} \subseteq \mathcal{H}^{\mathcal{X}|\mathbf{D}}$ the union of irreducible components containing the point $[C]$ corresponding to the curve $C$, and by $\mathcal{H}_t$ the fibre of $\mathcal{H}$ over $t \in \mathbf{D}$. We say that *a Cartier divisor $C_t \subseteq \mathcal{X}_t$ is a deformation of $C \subseteq \mathcal{X}_0$ if the point $[C_t] \in \mathcal{H}_t$ belongs to a curve $\gamma \subseteq \mathcal{H}$, which is a multi-section of order $r \geq 1$ of $\mathcal{H} \to \mathbf{D}$, totally ramified at the point $[C]$, i.e., intersecting $\mathcal{H}_0$ at $r[C]$.*

### 2.1 – Some local geometry: the versal deformation space of the singularities of $C$

Every deformations of $C$ in $\mathcal{X}$ induces a deformation of each singularity of $C$. Isolated singularities, i.e., germs of isolated singularities in the analytic topology, have a moduli space with the property of versality, cf. [13, Section 3], [18, Chapter 2], [12, Chapter 2] and [II].

For a plane curve singularity $(D, q)$, a versal deformation family is given by

$$(2.4.1) \quad \mathcal{C}_p = \{g(x,y) + t_1 g_1(x,y) + ... + t_m g_m(x,y) = 0\} \hookrightarrow \mathbf{A}^2 \times \Delta_{D,q}$$

$$\Delta_{D,q} = \mathrm{Spec}(\mathbf{C}[t_1, ..., t_m])$$

where, denoting by $J_p = (g, \frac{\partial g}{\partial y}, \frac{\partial g}{\partial z})$ the Jacobian ideal of the local equation $g(y,z)$ of $D$ at $p$, the polynomials $g_1(x,y), ..., g_m(x,y)$ are such that their images in $\mathbf{C}(y,z)/J_q$ form a basis of $\mathbf{C}(y,z)/J_q$, as a $\mathbf{C}$-vector space. The parameter space

$$\Delta_{D,q} = \mathrm{Spec}(\mathbf{C}[t_1, ..., t_m]),$$

is a (mini)-versal deformation space of $(D,q)$. Using the terminology of [13], we will refer to (2.4.1) as *the* étale versal deformation family of the plane curve singularity $(D,q)$, and to $\Delta_{D,q}$ as *the* étale (mini)-versal deformation space of $(D,q)$.

We now provide an explicit description of the étale versal deformation family of every singularity of our curve $C \subseteq \mathcal{X}$. At every node $p \in \mathcal{X}_0 \setminus E$, the curve $C$ is analytically equivalent to the plane curve $f(y,z) = yz = 0$, thus

$$\mathbf{C}[y,z]/J_p \simeq \mathbf{C}[y,z]/(y,z) \simeq \mathbf{C},$$

and a versal family for $(C,p)$ is $\mathcal{C}_p \to \Delta_{(C,p)} = \mathrm{Spec}(\mathbf{C}[t])$, with

$$\mathcal{C}_p = \{xy + t = 0\} \subseteq \mathbf{A}^2 \times \mathrm{Spec}(\mathbf{C}[t]).$$

The only singular fibre is that over $t = 0$, implying that a deformation of a nodal curve is either a nodal curve or a smooth curve.

At every tacnodal singularity $p \in E_{ij}$, a local equation $f$ of $C$ at $p$ is given by (2.1.2), so that $J_f = (2y - z^m, mz^{m-1}y)$ and

$$\mathbf{C}[y,z]/J_f \simeq \mathbf{C}^{2m-1}.$$

In particular, at every tacnodal singularity $p \in E_{ij}$ of $C$, choosing

$$\{1, z, z^2, \ldots, z^{m-1}, y, yz, yz^2, \ldots, yz^{m-2}\}$$

as a basis for $\mathbf{C}[y,z]/J_f$, a versal deformation family of $(C,p)$ is given by $\mathcal{C}_p \to \Delta_{(C,p)}$, where

$$\Delta_{(C,p)} = \mathrm{Spec}(\mathbf{C}[\alpha_0, ..., \alpha_{m-2}, \beta_0, ..., \beta_{m-1}])$$

and $\mathcal{C}_p \subseteq \mathbf{A}^2 \times \Delta_{(C,p)}$ has equation

$$(2.4.2) \quad \mathcal{C}_p : F(y, z; \underline{\alpha}, \underline{\beta}) = y^2 + \Big(\sum_{i=0}^{m-2} \alpha_i z^i + z^m\Big) y + \sum_{i=0}^{m-1} \beta_i z^i = 0,$$

where we are using the same notation as in [1] and [13], see [IV, V].

We now want to introduce the versal property of the versal family $\mathcal{C}_p \to \Delta_{(C,p)}$. Keeping in mind [16, Appendix B], we will do this in the analytic topology instead of the étale topology as in [13, Section 3], to simplify notation.

Let us denote by $\mathcal{D} \to \mathcal{H}$ the universal family parametrized by the Hilbert scheme $\mathcal{H}$. By versality, for every singular point $p$ of $C$, there exist analytic neighborhoods $U_p$ of $[C]$ in $\mathcal{H}$, $U'_p$ of $p$ in $\mathcal{D}$, and $V_p$ of $\underline{0}$ in $\Delta_{(C,p)}$, and a map $\phi_p : U_p \to V_p$, such that the family $\mathcal{D}|_{U_p} \cap U'_p$ is isomorphic to the pull-back of $\mathcal{C}_p|_{V_p}$ by the map $\phi_p$, making the following diagram

(2.4.3)
$$\begin{array}{ccccccc}
\mathcal{C}_p & \longleftarrow & \mathcal{C}_p|_{V_p} & \longleftarrow & U_p \times_{V_p} \mathcal{C}_p|_{V_p} & \xrightarrow{\simeq} & \mathcal{D}|_{U_p} \cap U'_p \hookrightarrow \mathcal{D} \\
\downarrow & & \downarrow & & & \searrow & \downarrow \\
T^1_{C,p} & \longleftarrow & V_p & & \xleftarrow{\phi_p} & & U_p \hookrightarrow \mathcal{H}
\end{array}$$

commutative. Thus, answering Problem (2.2) is the same as describing the image of the versal map

$$(2.4.4) \qquad \phi = \prod_{p \in \mathrm{Sing}(C)} \phi_p : \bigcap_{p \in \mathrm{Sing}(C)} U_p \to \prod_{p \in \mathrm{Sing}(C)} \Delta_{(C,p)},$$

where $\mathrm{Sing}(C)$ is the singular locus of $C$.

## 2.2 – Equisingular deformations of $C$

In order to understand the image of the versal map $\phi$ of (2.4.4), we describe its differential $d\phi_{[C]}$ at the point $[C] \in \bigcap_p U_p \subseteq \mathcal{H}$, which splits in the direct sum $d\phi_{[C]} = \bigoplus_{p \in \mathrm{Sing}(C)} d\phi_{p,[C]}$, where $d\phi_{p,[C]}$ is the differential of $\phi_p$ at the point $[C]$.

The differential map $d\phi_{[C]}$ can be identified with the map

$$(2.4.5) \quad H^0(\beta) = d\phi_{[C]} = \bigoplus_{p \in \mathrm{Sing}(C)} d\phi_{p,[C]} : H^0(C, \mathcal{N}_{C|\mathcal{X}}) \to H^0(C, T^1_C) = \bigoplus_{p \in \mathrm{Sing}(C)} T^1_{C,p}$$

induced by the exact sequence

$$(2.4.6) \qquad 0 \longrightarrow \Theta_C \longrightarrow \Theta_\mathcal{X}|_C \xrightarrow{\alpha} \mathcal{N}_{C|\mathcal{X}} \xrightarrow{\beta} T^1_C \longrightarrow 0,$$

where $\Theta_C$ and $\Theta_\mathcal{X}$ are the tangent sheaves of $C$ and $\mathcal{X}$, respectively, $\mathcal{N}_{C|\mathcal{X}}$ is the normal bundle of $C$ in $\mathcal{X}$, and $T^1_C$ is the first cotangent sheaf of $C$, cf. [18, Prop. 1.1.9]. We are in particular using standard identifications of the tangent space $T_{[C]}\mathcal{H}$ at $[C]$ with $H^0(C, \mathcal{N}_{C|\mathcal{X}})$ and of $T_0 \Delta_{(C,p)}$ with $T^1_{C,p}$.

The kernel of the sheaf map $\beta$ in (2.4.6) is usually denoted by $\mathcal{N}'_{C|\mathcal{X}}$, and is called the *equisingular normal sheaf of $C$ in $\mathcal{X}$*. This terminology is due to the fact that the origin $\underline{0} \in \Delta_{(C,p)}$ is the only point in the mini-versal deformation space corresponding to a singularity analytically equivalent to $(C, p)$. Thus, by the short exact sequence

$$(2.4.7) \qquad 0 \longrightarrow \mathcal{N}'_{C|\mathcal{X}} \longrightarrow \mathcal{N}_{C|\mathcal{X}} \xrightarrow{\beta} T^1_C \longrightarrow 0,$$

we have that $H^0(C, \mathcal{N}'_{C|\mathcal{X}})$ can be identified with the kernel of the differential $d\phi_{[C]}$ in (2.4.5), and it parametrizes the tangent space to the analytically equisingular deformation locus of $C$ in $\mathcal{X}$. One says that $H^0(C, \mathcal{N}'_{C|\mathcal{X}})$ parametrizes first order infinitesimal deformations of $C$ preserving all singularities of $C$ from the analytic point of view. Similarly, the kernel

$$\ker(d\phi_{p,[C]})$$

of the differential $d\phi_{p,[C]} : H^0(C, \mathcal{N}_{C|\mathcal{X}}) \to T^1_{C,p}$ in (2.4.5) parametrizes first order infinitesimal deformations of $C$ in $\mathcal{X}$ which are analytically equisingular at $p$, and there is an obvious inclusion $H^0(C, \mathcal{N}'_{C|\mathcal{X}}) \subseteq \ker(d\phi_{p,[C]})$.

**(2.5) Definition.** *We denote by*

$$\mathcal{ES}(C) = \phi^{-1}(\underline{0}, ..., \underline{0})$$

*the (scheme theoretic) inverse image of the origin with respect to the versal map $\phi$ in (2.4.4). We will call it the* equisingular deformation locus *of $C$.*

**(2.6) Remark.** As already observed, the inverse image $\phi^{-1}(\underline{0}, ..., \underline{0})$ of $(\underline{0}, ..., \underline{0})$ by the versal map (2.4.4) is, in general, the locus of analytically equisingular deformations of $C$ in $\mathcal{X}$. Since we assume $C$ with only nodes and tacnodes, analytically equisingular deformations of $C$ in $\mathcal{X}$ coincide with topologically equisingular deformations of $C$ in $\mathcal{X}$, cf. Remark (2.4). For this reason we call $\mathcal{ES}(C) = \phi^{-1}(\underline{0}, ..., \underline{0})$ the equisingular deformation locus of $C$, with no further specification. Similarly, we refer to its tangent space $H^0(C, \mathcal{N}'_{C|\mathcal{X}}) \simeq T_{[C]}(\mathcal{ES}(C))$ at $[C]$ just as the space of equisingular first order infinitesimal deformations of $C$ in $\mathcal{X}$.

**(2.7) Remark.** The versal map $\phi$ in (2.4.4), as well as the description of its differential $d\phi_{[C]}$ in (2.4.5), and the exact sequences (2.4.6) and (2.4.7), are defined for any reduced Cartier divisor $C$ on $\mathcal{X}_0$, whose singularities $(C, p)$ are all reduced complete intersection singularities in $\mathbf{C}^3$. Moreover, the description of the étale versal deformation family of $(C, p)$ provided in (2.4.1) under the assumption that $(C, p)$ is a planar curve singularity, extends to complete intersection curve singularities in $\mathbf{C}^n$ (see [12], or [11, Section 3] for an example of explicit computation when $(C, p)$ is a complete intersection curve singularity in $\mathbf{C}^3$).

**(2.8) Lemma.** *Let $\mathcal{X} \to \mathbf{D}$ be a flat projective family of complex surfaces as in (3.1).*

*(a) For* any *reduced Cartier divisor $C \subseteq \mathcal{X}_0$ with non-empty intersection with $E$, one has*

$$(2.8.1) \qquad H^0(C, \mathcal{N}'_{C|\mathcal{X}}) = H^0(C, \mathcal{N}'_{C|\mathcal{X}_0}) \subseteq H^0(C, \mathcal{N}_{C|\mathcal{X}_0}).$$

*More generally, for every intersection point $p \in C \cap E$, the kernel of the differential $d\phi_{p,[C]} : H^0(C, \mathcal{N}_{C|\mathcal{X}}) \to T^1_{C,p}$ in (2.4.5) is such that*

$$(2.8.2) \qquad H^0(C, \mathcal{N}'_{C|\mathcal{X}}) = H^0(C, \mathcal{N}'_{C|\mathcal{X}_0}) \subseteq \ker(d\phi_{p,[C]}) \subseteq H^0(C, \mathcal{N}_{C|\mathcal{X}_0}).$$

*(b) If $C \subseteq \mathcal{X}_0$ is a reduced Cartier divisor as above with only $\delta = \sum_i \delta_i$ nodes on $\mathcal{X}_0 \setminus E$ and tacnodes on $E$, then $H^0(C, \mathcal{N}'_{C|\mathcal{X}}) = H^0(C, \mathcal{N}'_{C|\mathcal{X}_0})$ is contained the linear subspace of sections of $H^0(C, \mathcal{N}_{C|\mathcal{X}_0})$ vanishing at every node $p$ of $C$ on $\mathcal{X}_0 \setminus E$ and vanishing at every m-tacnode $p \in E$ of $C$ with multiplicity $m - 1$.*

*Proof.* We first prove part (a) of the lemma and, in particular, equality (2.8.1) for any reduced Cartier divisor $C$ on $\mathcal{X}_0$ with non-empty intersection with the singular locus $E$ of $\mathcal{X}_0$. Fix any intersection point $p \in E \cap C$. Keeping in mind Remark (2.7), we want to write equisingular infinitesimal first order deformations of $C$ in $\mathcal{X}$ locally at the point $p$. We consider the localized exact sequence

$$(2.8.3) \qquad 0 \longrightarrow \mathcal{N}'_{C|\mathcal{X},p} \longrightarrow \mathcal{N}_{C|\mathcal{X},p} \longrightarrow T^1_{C,p} \longrightarrow 0.$$

Fix local analytic coordinates $x$, $y$, $z$, $t$ at $p$ in such a way that $C$ has equations

$$t = 0, \quad f_1(x,y,z) = 0, \quad \text{and} \quad f_2(x,y,z) = xy = 0$$

at $p$. Then we may identify:

- the local ring $\mathcal{O}_{C,p} = \mathcal{O}_{\mathcal{X},p}/\mathcal{I}_{C|\mathcal{X},p}$ of $C$ at $p$ with the localization at the origin of $\mathbf{C}[x,y,z]/(f_1,f_2)$,

- the $\mathcal{O}_{C,p}$-module $\mathcal{N}_{C|\mathcal{X},p}$ with the free $\mathcal{O}_{\mathcal{X},p}$-module $\mathfrak{hom}_{\mathcal{O}_{\mathcal{X},p}}(\mathcal{I}_{C|\mathcal{X},p}, \mathcal{O}_{C,p})$ generated by the morphisms $f_1^*$ and $f_2^*$, defined by

$$f_i^*\big(s_1(x,y,z)f_1(x,y,z) + s_2(x,y,z)f_2(x,y,z)\big) = s_i(x,y,z), \quad \text{for } i = 1, 2,$$

and, finally,

- the $\mathcal{O}_{C,p}$-module

$$\begin{aligned}(\Theta_{\mathcal{X}}|_C)_p &\simeq \Theta_{\mathcal{X},p}/(I_{C,p} \otimes \Theta_{\mathcal{X},p}) \\ &\simeq \langle \partial/\partial x, \partial/\partial y, \partial/\partial z, \partial/\partial t \rangle_{\mathcal{O}_{C,p}} \Big/ \langle \partial/\partial t - x\partial/\partial y - y\partial/\partial x \rangle\end{aligned}$$

with the free $\mathcal{O}_{\mathcal{X},p}$-module generated by the derivatives $\partial/\partial x, \partial/\partial y, \partial/\partial z$.

With these identifications, the localization $\alpha_p : (\Theta_{\mathcal{X}}|_C)_p \to \mathcal{N}_{C|\mathcal{X},p}$ of the sheaf map $\alpha$ from (2.4.6) is defined by

$$\begin{aligned}\alpha_p(\partial/\partial x) &= \big(s = s_1 f_1 + s_2 f_2 \longmapsto \partial s/\partial x =_{\mathcal{O}_{C,p}} s_1 \partial f_1/\partial x + s_2 \partial f_2/\partial x\big) \\ &= \partial f_1/\partial x f_1^* + y f_2^*, \\ \alpha_p(\partial/\partial y) &= \partial f_1/\partial y f_1^* + x f_2^*, \text{ and} \\ \alpha_p(\partial/\partial z) &= \partial f_1/\partial z f_1^*.\end{aligned}$$

By definition of $\mathcal{N}'_{C|\mathcal{X}}$, a germ $s \in \mathcal{N}_{C|\mathcal{X},p}$ is equisingular at $p$, i.e., $s \in \mathcal{N}'_{C|\mathcal{X},p}$, if and only if there exists a germ

$$u = u_x(x,y,z)\partial/\partial x + u_y(x,y,z)\partial/\partial y + u_z(x,y,z)\partial/\partial z \in (\Theta_{\mathcal{X}}|_C)_p$$

such that $s = \alpha_p(u)$. Hence, locally at $p$, an equisingular first order infinitesimal deformation of $C$ in $\mathcal{X}$ has equation

(2.8.4) $$\begin{cases} f_1(x,y,z) + \varepsilon(u_x \partial f_1/\partial x + u_y \partial f_1/\partial z + u_z \partial f_1/\partial z) = 0 \\ xy + \varepsilon(y u_x + x u_y) = 0. \end{cases}$$

We claim that the second equation of (2.8.4) is the local equation at $p$ of an equisingular deformation of $\mathcal{X}_0$ in $\mathcal{X}$. Indeed, by the exact sequence

$$0 \longrightarrow \Theta_{\mathcal{X}_0} \longrightarrow \Theta_{\mathcal{X}}|_{\mathcal{X}_0} \longrightarrow \mathcal{N}_{\mathcal{X}_0|\mathcal{X}} \cong \mathcal{O}_{\mathcal{X}_0} \longrightarrow T^1_{\mathcal{X}_0} \cong \mathcal{O}_E \longrightarrow 0,$$

one finds that $\mathcal{N}'_{\mathcal{X}_0|\mathcal{X}} \simeq \mathcal{I}_{E|\mathcal{X}_0}$.

Now if $s = \alpha_p(u)$ is the localization of a global section of $\mathcal{N}'_{C|\mathcal{X}}$, by the equality $H^0(\mathcal{X}_0, \mathcal{N}'_{\mathcal{X}_0|\mathcal{X}}) = H^0(\mathcal{X}_0, \mathcal{I}_{E|\mathcal{X}_0}) = 0$, we find that the analytic function

$$y u_x(x,y,z) + x u_y(x,y,z)$$

must be identically zero in the second equation of (2.8.4). This proves inclusion (2.8.1) and (2.8.2).

We now prove part (b) of the lemma. Assume that $C$ is a Cartier divisor with only $\delta = \sum_i \delta_i$ nodes on $\mathcal{X}_0 \setminus E$ and tacnodes on $E$, as in (3.1), and consider $p \in E \cap C$. In this case we may assume
$$f_1(x,y,z) = x + y + z^m$$
so that the equations in (2.8.4) become

(2.8.5) $$\begin{cases} x + y + z^m + \varepsilon(u_x + u_y + mz^m u_z) = 0 \\ xy + \varepsilon(yu_x + xu_y) = 0. \end{cases}$$

If these equations are the equations of the localization of a global section in $H^0(C, \mathcal{N}_{C|\mathcal{X}})$, then, as above, we find that $yu_x + xu_y$ is identically zero. We deduce in particular that

- $u_x(0,0,0) = u_y(0,0,0) = 0$, and
- for every $n \geqslant 1$, no $z^n$-terms appear in $u_x(x,y,z)$ and $u_y(x,y,z)$, no $y^n$-terms and $y^n z^m$-terms appear in $u_x(x,y,z)$ and, finally, no $x^n$-terms and $x^n z^m$-terms appear in $u_y(x,y,z)$.

In particular, if $X_i$, with local equations $x = t = 0$, and $X_j$, with local equations $y = t = 0$, are the two irreducible components of $\mathcal{X}_0$ containing $p$, then the local equations at $p$ on $X_i$ of equisingular infinitesimal deformations of $C$ are given by

(2.8.6) $$\begin{cases} y + z^m + \varepsilon(mz^{m-1}u_z(x,y,z) + yq(y,z)) = 0 \\ x = 0, \end{cases}$$

where $q(y,z)$ is an analytic function in the variables $y$ and $z$, and similarly on $X_j$. This proves that all sections in $H^0(C, \mathcal{N}'_{C|\mathcal{X}})$ vanish with multiplicity $m-1$ at every $m$-tacnode of $C$ on $E$.

The fact that all sections of $H^0(C, \mathcal{N}'_{C|\mathcal{X}_0})$ vanish at every node $p$ of $C$ on $\mathcal{X}_0 \setminus E$ is a very classical result, which can be proved here with the same approach as above, by using that at a node on $\mathcal{X}_0 \setminus E$ the local equations of $C$ are $x - t = t = 0 = yz$ (see also [II]).

Thus the lemma is proved. □

## 2.3 – Image of the versal map

We now consider:

- $\mathcal{X} \to \mathbf{D}$ be a flat projective family of complex surfaces with smooth total space $\mathcal{X}$, smooth general fibre $\mathcal{X}_t$ and central fibre $\mathcal{X}_0$ with normal crossing singularities,
- $C \subseteq \mathcal{X}_0 = \cup X_i$ a reduced Cartier divisor with only tacnodes at the intersection points with the singular locus $E = \cup_{ij} E_{ij}$ of $\mathcal{X}_0$ and nodes on $\mathcal{X}_0 \setminus E$, as in (3.1).

**(2.9) Lemma.** *In the above setting, let $p \in C \cap E_{ij} \subseteq \mathcal{X}_0 \subseteq \mathcal{X}$ be an $m$-tacnode, with $m \geqslant 1$. Then, using the parametrization (2.4.2) for the mini-versal family of $(C,p)$, we have that the image of the differential map*
$$d\phi_{p,[C]} : H^0(C, \mathcal{N}_{C|\mathcal{X}}) \to T^1_{C,p} = T_{\underline{0}}\Delta_{(C,p)}$$
*is always contained in the linear subspace $H_p : \beta_1 = \cdots = \beta_{m-1} = 0$, while the image of the differential map*
$$d\phi_{p,[C]} : H^0(C, \mathcal{N}_{C|\mathcal{X}_0}) \to T^1_{C,p} = T_{\underline{0}}\Delta_{(C,p)}$$

*is always contained in the tangent space* $\Gamma_p = T_{\underline{0}}EG_{(C,p)} : \beta_0 = \cdots = \beta_{m-1} = 0$ *at* $\underline{0}$ *to the equigeneric locus* $EG_{(C,p)}$ *in* $\Delta_{(C,p)}$ *(cf. [13] and [II]). In particular, the subspace*

(2.9.1) $$\ker(d\phi_{p,[C]}) \subseteq H^0(C, \mathcal{N}'_{C|\mathcal{X}_0})$$

*has dimension*

$$\dim(\ker(d\phi_{p,[C]})) \geqslant h^0(C, \mathcal{N}_{C|\mathcal{X}_0}) - m + 1 \geqslant h^0(C, \mathcal{N}_{C|\mathcal{X}}) - m.$$

*The differential map* $d\phi_{p,[C]}$ *is surjective on* $H_p$ *if and only if the subspace* $\ker(d\phi_{p,[C]})$ *of infinitesimal first order deformations of* $C$ *in* $\mathcal{X}_0$ *that are equisingular at* $p$ *has codimension exactly* $m - 1$ *in* $H^0(C, \mathcal{N}_{C|\mathcal{X}_0})$.

*Proof.* We first observe that $H_p$ contains in codimension 1 the tangent space

$$\Gamma_p = T_{\underline{0}}EG_{(C,p)} : \beta_0 = \cdots = \beta_{m-2} = 0$$

at $\underline{0}$ to the equigeneric locus $EG_{(C,p)} \subseteq \Delta_{(C,p)}$.

The image of $U_p \cap \mathcal{H}_0$ with respect to $\phi_p$ (recall the notation from Diagram (2.4.3)) is contained in the equigeneric locus $EG_{(C,p)} \subseteq \Delta_{(C,p)}$. Indeed, no matter how we deform $C$ in $\mathcal{X}_0$, the smooth tangency of the irreducible components $C_i$ and $C_j$ with $E_{ij}$ at $p$ deforms to smooth tangencies with $E_{ij}$ at $r \leqslant m$ points $p_l$, with $1 \leqslant l \leqslant r$, of multiplicity $m_l$ with $\sum_l m_l = m$. This corresponds locally to a deformation of the $m$-tacnode of $C$ at $p$ into $r$ tacnodes of orders $m_1, \ldots, m_r$, each decreasing the geometric genus by $m_l$ with respect to the arithmetic genus. Thus, $d\phi_{p,[C]}(H^0(C, \mathcal{N}_{C|\mathcal{X}_0})) \subseteq \Gamma_p = T_{\underline{0}}EG_{(C,p)}$.

Now, by the exact sequence

$$0 \to \mathcal{N}_{C|\mathcal{X}_0} \to \mathcal{N}_{C|\mathcal{X}} \to \mathcal{N}_{\mathcal{X}_0|\mathcal{X}} = \mathcal{O}_{\mathcal{X}_0} \to 0,$$

we have that $H^0(C, \mathcal{N}_{C|\mathcal{X}_0})$ has codimension less or equal to 1 in $H^0(C, \mathcal{N}_{C|\mathcal{X}})$. If $h^0(C, \mathcal{N}_{C|\mathcal{X}_0}) = h^0(C, \mathcal{N}_{C|\mathcal{X}})$, the lemma is true because $T_{\underline{0}}EG_{(C,p)} \subseteq H_p$. If $h^0(C, \mathcal{N}_{C|\mathcal{X}}) = h^0(C, \mathcal{N}_{C|\mathcal{X}_0}) + 1$, it is enough to find a section

$$\sigma \in H^0(C, \mathcal{N}_{C|\mathcal{X}}) \setminus H^0(C, \mathcal{N}_{C|\mathcal{X}_0})$$

such that $d\phi_{p,[C]}(\sigma) \in H_p \setminus T_{\underline{0}}EG_{(C,p)}$. In this case, one has

$$H^0(C, \mathcal{N}_{C|\mathcal{X}}) \simeq H^0(C, \mathcal{N}_{C|\mathcal{X}_0}) \oplus H^0(C, \mathcal{O}_{\mathcal{X}_0}) \simeq H^0(C, \mathcal{N}_{C|\mathcal{X}_0}) \oplus \mathbf{C}.$$

The section $\sigma \in H^0(C, \mathcal{N}_{C|\mathcal{X}}) \setminus H^0(C, \mathcal{N}_{C|\mathcal{X}_0})$ corresponding to the pair $(\underline{0}, 1)$ via the above isomorphism has local equations

(2.9.2) $$\begin{cases} x + y + z^m = 0 \\ xy = \varepsilon. \end{cases}$$

It follows that the image of $\sigma$ via $d\phi_{p,[C]}$ is the point corresponding to the curve $y(y + z^m) = \beta_0$. This proves the lemma. □

**(2.10) Corollary.** *In the same setting as above and with the same notation and hypotheses as in Lemma (2.9) and Diagram* (2.4.3), *assume that:*

a) $h^0(C, \mathcal{N}_{C|\mathcal{X}_0}) = h^0(C, \mathcal{N}_{C|\mathcal{X}}) - 1;$

b) *the relative Hilbert scheme* $\mathcal{H}^{\mathcal{X}|\mathbf{D}}$ *is smooth at* $[C]$ *of dimension* $h^0(C, \mathcal{N}_{C|\mathcal{X}})$;

c) $\dim(\ker(d\phi_{p,[C]})) = h^0(C, \mathcal{N}_{C|\mathcal{X}_0}) - m + 1$.

*Then the image of the versal map*

$$\phi_p : \mathcal{H} \cap U_p \to \Delta_{(C,p)}$$

*is a smooth variety of dimension $m$ in $\Delta_{(C,p)}$, with tangent space at $\underline{0}$ given by the linear subspace $H_p : \beta_1 = ... = \beta_{m-1} = 0$, while the image of the versal map*

$$\phi_p : \mathcal{H}_0 \cap U_p \to \Delta_{(C,p)}$$

*coincides with the equigeneric deformation locus $EG_{(C,p)} \subseteq \Delta_{(C,p)}$, which is smooth of dimension $m$ at $\underline{0}$, with tangent space given by $\Gamma_p = T_{\underline{0}} EG_{(C,p)} : \beta_0 = ... = \beta_{m-1} = 0$, and whose general element corresponding to an $m$-nodal curve.*

*Moreover, $\phi_p(\mathcal{H} \cap U_p)$ intersects the locus of $(m-1)$-nodal curves in $\Delta_{(C,p)}$ along $\Gamma_p \cup \gamma$, where $\gamma$ is a curve that is smooth at $\underline{0}$ and intersects $\Gamma_p$ with multiplicity $m$. In particular the curve $C$ may be deformed in $\mathcal{H}$ in such a way the tacnode at $p$ is deformed into $m-1$ nodes.*

*Proof.* First observe that, by hypothesis b), there is only one irreducible component $\mathcal{H}$ of $\mathcal{H}^{\mathcal{X}|\mathbf{D}}$ containing $[C]$, which is smooth at $[C]$ of dimension $h^0(C, \mathcal{N}_{C|\mathcal{X}})$. Moreover, by hypothesis a), the central fibre of $\mathcal{H}_0$ of $\mathcal{H}$ is smooth at $[C]$ of dimension $h^0(C, \mathcal{N}_{C|\mathcal{X}_0}) = h^0(C, \mathcal{N}_{C|\mathcal{X}}) - 1$.

In particular, in Diagram (2.4.3), we have that the analytic varieties $\mathcal{H} \cap U_p$ and $\mathcal{H}_0 \cap U_p$ are irreducible. It follows that the image of $\phi_p : \mathcal{H} \cap U_p \to \Delta_{(C,p)}$ is an irreducible variety of dimension $\leqslant m$, since $\phi_p(\mathcal{H}_0 \cap U_p)$ is contained in the equigeneric locus, which has dimension $m - 1$ and $\phi_p(\mathcal{H} \cap U_p)$ has dimension $\leqslant \dim(\phi_p(\mathcal{H}_0 \cap U_p)) + 1$. In particular, $\phi_p$ has general fibre of dimension at least $\dim_{[C]}(\mathcal{H}) - m$. Under the hypotheses of this corollary, the differential $d\phi_{p,[C]}$ has rank $m$. Hence, under hypotheses a), b) and c), the fibre $\phi_p^{-1}(\underline{0})$ has dimension exactly $\dim_{[C]}(\mathcal{H}) - m = h^0(C, \mathcal{N}_{C|\mathcal{X}}) - m$, and $\phi_p(\mathcal{H} \cap U_p)$ is irreducible of dimension $m$ and smooth at $\underline{0}$. The fact that $\phi_p(\mathcal{H} \cap U_p)$ intersects the locus of $(m-1)$-nodal curves in $\Delta_{(C,p)}$ along $\Gamma_p \cup \gamma$, where $\gamma$ is a curve that is smooth at $\underline{0}$ and intersects $\Gamma_p$ with multiplicity $m$ follows from [1, Lemma 4.1], cf. [IV, Theorem (1.2)]. □

**(2.11) Remark.** Under the hypotheses of Lemma (2.10), by [10, Proposition 3.7 and Corollary 3.12], one has more generally that there exist deformations $C_t \subseteq \mathcal{X}_t$ of $C \subseteq \mathcal{X}_0$ such that the $m$-tacnode of $C$ at $p$ deforms to $d_i$ singularity of type $A_i$, for every multi-index $(d_1, ..., d_r)$ such that $\sum_i i d_i = m - 1$.

The following lemma provides a standard way to verify when the hypothesis c) in Corollary (2.10) is satisfied.

Let $W^{|C|}_{p,m-1} \subseteq |\mathcal{O}_{\mathcal{X}_0}(C)|$ be the linear system of Cartier divisors in $|\mathcal{O}_{\mathcal{X}_0}(C)|$ intersecting $E$ at $p$ with multiplicity $m-1$. We have that

$$\dim(W^{|C|}_{p,m-1}) \geqslant \dim(|\mathcal{O}_{\mathcal{X}_0}(C)|) - m + 1.$$

**(2.12) Lemma.** *In the same setting as above, and with the same notation and hypotheses as in Lemma (2.9), assume that:*

(2.12.1) $$\dim(W^{|C|}_{p,m-1}) = \dim(|\mathcal{O}_{\mathcal{X}_0}(C)|) - m + 1.$$

*Then we have*

$$\dim(\ker(d\phi_{p,[C]})) = h^0(C, \mathcal{N}_{C|\mathcal{X}_0}) - m + 1.$$

*Proof.* By the hypothesis that $\dim(W^{|C|}_{p,m-1}) = \dim(|\mathcal{O}_{\mathcal{X}_0}(C)|) - m + 1$, we have

$$\dim(W^{|C|}_{p,r}) = \dim(|\mathcal{O}_{\mathcal{X}_0}(C)|) - r + 1,$$

for every $r \leqslant m - 1$. Via the natural inclusion

$$|\mathcal{O}_{X_0}(C)| = H^0(\mathcal{X}_0, \mathcal{O}_{\mathcal{X}_0}(C))/H^0(\mathcal{X}_0, \mathcal{O}_{\mathcal{X}_0}) \subseteq H^0(C, \mathcal{N}_{C|\mathcal{X}_0})$$

induced by the standard exact sequence

$$0 \to \mathcal{O}_{\mathcal{X}_0} \to \mathcal{O}_{\mathcal{X}_0}(C) \to \mathcal{O}_C(C) \to 0,$$

we find a chain a codimension 1 linear subsystems

$$W^{|C|}_{p,m-1} \subsetneq W^{|C|}_{p,m-2} \subsetneq \ldots \subsetneq W^{|C|}_{p,1} \subsetneq |\mathcal{O}_{X_0}(C)| \subseteq H^0(C, \mathcal{N}_{C|\mathcal{X}_0}).$$

By Lemma (2.8), we know that $W^{|C|}_{p,r} \cap \ker(d\phi_{p,[C]}) \subseteq W^{|C|}_{p,m-1}$, for every $r \leqslant m - 2$. Since $H^0(C, \mathcal{N}_{C|\mathcal{X}_0})$ is mapped by $d\phi_{p,[C]}$ to the $(m-1)$-dimensional linear subspace $T_{\underline{0}}EG \subseteq T_{\underline{0}}\Delta_{(C,p)}$, we find that $d\phi_{p,[C]} : H^0(C, \mathcal{N}_{C|\mathcal{X}_0}) \to T_{\underline{0}}EG$ is surjective by standard linear algebra, and thus $\dim(\ker(d\phi_{p,[C]})) = h^0(C, \mathcal{N}_{C|\mathcal{X}_0}) - m + 1$. □

## 3 – Global deformations of $C$

Let $\mathcal{X} \to \mathbf{D}$ and $C = \cup_i C_i \subseteq \mathcal{X}_0 = \cup_i X_i$ be respectively the flat projective family of surfaces and the Cartier divisor $C$ on $\mathcal{X}_0$ as in (3.1). We recall that $\mathcal{X}$ has smooth general fibre $\mathcal{X}_t$ and reducible central fibre $\mathcal{X}_0 = \cup_i X_i$ with normal crossing singularities. We denoted by $X_i$ the irreducible components of $\mathcal{X}_0$, and by $E = \cup_{i,j} E_{ij}$ the singular locus of $\mathcal{X}_0$, where $E_{ij} = X_i \cap X_j$. Moreover, we assume that $C \cap X_i = C_i$ is a reduced $\delta_i$-nodal curve intersecting $E$ at smooth points, in such a way that the only singularities of $C$ are nodes on $\mathcal{X}_0 \setminus E$ and tacnodes on $E$. For every tacnode $p$ of $C$, we denote by $m_p$ its order and by

(3.0.1) $$\lambda = \operatorname{lcm}\{m_p | p \text{ tacnode of } C\}$$

the least common multiple of the $m_p$. We furthermore set

(3.0.2) $$\mu = \prod_{p \text{ tacnode of } C} m_p \quad \text{and} \quad k = \frac{\mu}{\lambda}.$$

Finally, we assume that $C$ has

$$\delta = \sum_i \delta_i \quad \text{nodes on } \mathcal{X}_0 \setminus E,$$

and we denote by

$$\tau_m \quad \text{the number of } m-\text{tacnodes of } C \text{ on } E,$$

for every $m \geqslant 1$.

Here we will provide sufficient conditions for $C \subseteq \mathcal{X}_0$ to be deformed to a nodal curve $C_t \subseteq \mathcal{X}_t$.

We recall that we denoted by $\mathcal{H}^{\mathcal{X}|\mathbf{D}}$ the relative Hilbert scheme of the family $\mathcal{X} \to \mathbf{D}$, by $\mathcal{H}$ the union of the irreducible components of $\mathcal{H}^{\mathcal{X}|\mathbf{D}}$ containing $[C]$, and by $\mathcal{H}_t$ the fibre of $\mathcal{H}$ over $t \in \mathbf{D}$. Let now $\mathbf{D}^o = \mathbf{D} \setminus \underline{0}$, $\mathcal{X}^o = \mathcal{X} \setminus \mathcal{X}_0$, and $\mathcal{H}^o = \mathcal{H} \setminus \mathcal{H}_0$. For all $\sigma \in \mathbf{N}$, we denote by $\mathcal{U}_\sigma^{\mathcal{X}^o|\mathbf{D}^o} \subseteq \mathcal{H}^o$ the relative Severi variety of $\sigma$-nodal curves, the fibre of which over $t \in \mathbf{D}^o$ is the locally closed (possibly empty) subscheme $\mathcal{U}_\sigma^{\mathcal{X}_t} \subseteq \mathcal{H}_t$ (endowed with its reduced structure) parametrizing reduced $\sigma$-nodal Cartier divisors on $\mathcal{X}_t$. We moreover denote by $\mathcal{V}_\sigma^{\mathcal{X}|\mathbf{D}} \subseteq \mathcal{H}$ the Zariski closure in $\mathcal{H}$ of $\mathcal{U}_\sigma^{\mathcal{X}^o|\mathbf{D}^o}$. Observe that, for every $t \neq 0$, the fiber $\mathcal{V}_\sigma^{\mathcal{X}_t}$ of $\mathcal{V}_\sigma^{\mathcal{X}|\mathbf{D}}$ over $t$ is the Zariski closure of $\mathcal{U}_\sigma^{\mathcal{X}_t}$, while $\mathcal{V}_\sigma^{\mathcal{X}|\mathbf{D}} \cap \mathcal{H}_0$ will have several irreducible components, whose general point corresponds to a curve on $\mathcal{X}_0$ with singularities possibly different than nodes.

**(3.1) Theorem.** *In the above setting, let*

$$W_{[C],es}^{|C|} \subseteq |\mathcal{O}_{\mathcal{X}_0}(C)|$$

*be the linear system of Cartier divisors in $|\mathcal{O}_{\mathcal{X}_0}(C)|$ passing through every node of $C$ on $\mathcal{X}_0 \setminus E$, and intersecting $E$ with multiplicity $m-1$ at every $m$-tacnode of $C$. Assume that:*

1. $h^0(C, \mathcal{N}_{C|\mathcal{X}_0}) = h^0(C, \mathcal{N}_{C|\mathcal{X}}) - 1$;

2. *the relative Hilbert scheme $\mathcal{H}^{\mathcal{X}|\mathbf{D}}$ is smooth at $[C]$ of dimension $h^0(C, \mathcal{N}_{C|\mathcal{X}})$. In particular, there is only one irreducible component $\mathcal{H} \subseteq \mathcal{H}^{\mathcal{X}|\mathbf{D}}$ containing $[C]$, which is smooth at $[C]$;*

3. $\dim(W_{[C],es}^{|C|}) = \dim(|\mathcal{O}_{\mathcal{X}_0}(C)|) - \delta - \sum_m \tau_m(m-1)$.

*Then $C \subseteq \mathcal{X}_0$ may be deformed to $(\delta + \sum_m \tau_m(m-1))$-nodal curves $C_t \subseteq \mathcal{X}_t$, in such a way every node of $C$ on $\mathcal{X}_0 \setminus E$ is preserved, and every $m$-tacnode of $C$ on $E$ is smoothed to $m-1$ nodes. All nodal curves $C_t$ obtained in this way are such that*

$$(3.1.1) \qquad h^0(C_t, \mathcal{N}'_{C_t|\mathcal{X}_t}) = h^0(C_t, \mathcal{N}_{C_t|\mathcal{X}_t}) - \delta - \sum_m \tau_m(m-1),$$

*thus, the corresponding point $[C_t] \in \mathcal{H}_t$ is a smooth point of the Severi variety $\mathcal{V}_{\delta+\sum_m \tau_m(m-1)}^{\mathcal{X}_t}$. In particular the nodes of $C_t$ may be smoothed independently, i.e. the versal map (2.4.4) for $C_t$ is surjective. This implies that*

$$[C] \in \mathcal{V}_\sigma^{\mathcal{X}|\mathbf{D}}, \text{ for every } \sigma \leqslant \delta + \sum_m \tau_m(m-1).$$

*Moreover, the equisingular deformation locus $\mathcal{ES}(C) \subseteq \mathcal{X}_0$ of $[C]$ in $\mathcal{X}_0$ (cf. Definition (2.5)) is an analytic open set of the smooth locus of an irreducible component $V_0$ of the central fibre*

$$\mathcal{V}_0 = \mathcal{V}_{\delta+\sum_m \tau_m(m-1)}^{\mathcal{X}|\mathbf{D}} \cap \mathcal{H}_0 = \mu V_0 + \cdots$$

*of $\mathcal{V}_{\delta+\sum_m \tau_m(m-1)}^{\mathcal{X}|\mathbf{D}} \subseteq \mathcal{H}$, of scheme theoretic multiplicity*

$$\mu = \prod_{p \text{ tacnode of } C} m_p.$$

*More precisely, at an analytic neighborhood of $[C]$, the universal Severi variety $\mathcal{V}_{\delta+\sum_m \tau_m(m-1)}^{\mathcal{X}|\mathbf{D}}$ will be the union of $k$ branches, each intersecting $\mathcal{H}_0$ at $[C]$ with multiplicity $\lambda$, where $k$ and $\lambda$ are defined by (3.0.2). These $k$ local analytic branches of $\mathcal{V}_{\delta+\sum_m \tau_m(m-1)}^{\mathcal{X}|\mathbf{D}}$ will all be smooth at $[C]$ if $\lambda = m_p$, for some tacnode $p$ of $C$; otherwise they will all be singular at $[C]$, with multiplicity $\frac{\lambda}{\max\{m_p | p \text{ tacnode of } C\}}$.*

**(3.2) Remark.** Hypothesis 3 in Theorem (3.1) can be weakened by asking

$$(3.2.1) \qquad h^0(C, \mathcal{N}'_{C|\mathcal{X}}) = h^0(C, \mathcal{N}'_{C|\mathcal{X}_0}) = h^0(C, \mathcal{N}_{C|\mathcal{X}_0}) - \delta - \sum_m \tau_m(m-1).$$

As we will observe in the proof of the theorem, under hypotheses 1 and 2, we have that hypothesis 3 is a sufficient condition for the equality (3.2.1). The reason we preferred this statement of Theorem (3.1) is that hypothesis 3 in the theorem is usually easier to be verified than (3.2.1).

We finally observe that, if

$$(3.2.2) \qquad H^0(C, \mathcal{N}_{C|\mathcal{X}_0}) = H^0(C, \mathcal{O}_C(C)),$$

which happens in particular if the family $\mathcal{X} \to \mathbf{D}$ is a family of regular surfaces (i.e., $h^1(\mathcal{X}_t, \mathcal{O}_{\mathcal{X}_t}) = 0$ for all $t$), then, under hypotheses 1, 2, and 3 of Theorem (3.1), we have that $H^0(C, \mathcal{N}'_{C|\mathcal{X}}) = W^{|C|}_{[C],es}$. Indeed, if (3.2.2) holds, then $H^0(C, \mathcal{N}'_{C|\mathcal{X}}) \subseteq W^{|C|}_{[C],es}$ by Lemma (2.8), and this inclusion is in fact an equality for dimensional reason.

**(3.3) Remark.** In Theorem (3.1), the nodes of $C$ on $E$ cannot be preserved by deforming $C$ to a curve in $\mathcal{X}_t$, with $t \neq 0$. This is a very general result that holds for any reduced Cartier divisor on $X_0$ with a node on $E$ under our hypotheses on $\mathcal{X} \to \mathbf{D}$, i.e., just assuming $\mathcal{X} \to \mathbf{D}$ a flat projective family of surfaces with smooth total space, smooth general fibre and central fibre $\mathcal{X}_0$ with normal crossing singularities (or more generally, with normal crossing singularity at a neighborhood of the node). A proof of this can be found in [9, Section 2]. This also follows from the hypotheses 1, 2 of Theorem (3.1), by Lemma (2.8). In order to see this, we fix a node $p$ of $C$ on $E$, and we consider the diagram (2.4.3). What we want to prove is that

$$(3.3.1) \qquad \phi_p^{-1}(\underline{0}) = U_p \cap \mathcal{H}_0.$$

By Lemma (2.8), we know that

$$T_{\underline{0}}(\phi_p^{-1}(\underline{0})) = \ker(d\phi_{p,[C]}) \subseteq H^0(C, \mathcal{N}_{C|\mathcal{X}_0}),$$

and, by the hypothesis 1 of Theorem (3.1), we get that $\dim(T_{\underline{0}}(\phi_p^{-1}(\underline{0}))) = h^0(C, \mathcal{N}_{C|\mathcal{X}_0}) = h^0(C, \mathcal{N}_{C|\mathcal{X}}) - 1$. Now observe that $U_p \cap \mathcal{H}_0 \subseteq \phi_p^{-1}(\underline{0})$ because every deformation of $C$ in $\mathcal{X}_0$ preserves the node of $C$ at $p$. Finally, by the hypotheses 1 and 2 in Theorem (3.1), we have that $\mathcal{H}_0$ is smooth at $[C]$ of dimension $h^0(C, \mathcal{N}_{C|\mathcal{X}_0}) = h^0(C, \mathcal{N}_{C|\mathcal{X}}) - 1$. This proves (3.3.1), which in particular says that Theorem (3.1) completely answers Problem (2.2) for a curve $C$ verifying the hypothesis (3.2.1).

*Proof of Theorem (3.1).* Let $\mathrm{Sing}(C)$ be the singular locus of $C$. For every singular point $p \in \mathrm{Sing}(C)$, we consider the mini-versal deformation space $\Delta_{(C,p)}$ described in Section 2.1, the versal map (2.4.4)

$$\phi = \prod_{p \in \mathrm{Sing}(C)} \phi_p : \bigcap_{p \in \mathrm{Sing}(C)} U_p \to \prod_{p \in \mathrm{Sing}(C)} \Delta_{(C,p)},$$

and its differential

$$d\phi_{[C]} = \bigoplus_{p \in \mathrm{Sing}(C)} d\phi_{p,[C]} : H^0(C, \mathcal{N}_{C|\mathcal{X}}) \to H^0(C, T^1_C) = \bigoplus_{p \in \mathrm{Sing}(C)} T^1_{C,p}$$

at $[C] \in \mathcal{H}$, described in Sections 2.2 and 2.3.

Arguing exactly as in the proof of Lemma (2.12), we see that

$$W^{|C|}_{[C],es} \subseteq |\mathcal{O}_{\mathcal{X}_0}(C)| \subseteq H^0(C, \mathcal{N}_{C|\mathcal{X}_0}),$$

and the hypothesis 3 of the theorem implies that the image of the differential $d\phi_{[C]}$ has dimension $\delta + \sum_m \tau_m(m-1)$, hence

$$(3.3.2) \quad \dim\bigl(\ker(d\phi_{[C]})\bigr) = h^0(C, \mathcal{N}'_{C|\mathcal{X}}) = h^0(C, \mathcal{N}'_{C|\mathcal{X}_0}) = h^0(C, \mathcal{N}_{C|\mathcal{X}_0}) - \delta - \sum_m \tau_m(m-1).$$

This in turn implies that

$$\dim\bigl(\ker(d\phi_{p,[C]})\bigr) = h^0(C, \mathcal{N}_{C|\mathcal{X}_0}) - m + 1, \quad \text{for every } m\text{-tacnode } p \in E \cap C \text{ of } C$$

and

$$\dim\bigl(\ker(d\phi_{p,[C]})\bigr) = h^0(C, \mathcal{N}_{C|\mathcal{X}_0}) - 1, \quad \text{for every node of } C \text{ in } \mathcal{X}_0 \setminus E.$$

In particular, Corollary (2.10) applies at every $m$-tacnode $p$ of $C$ on $E$.

Now, by Remark (3.3), we know that every 1-tacnode (i.e., node) of $C$ on $E$ is necessarily smoothed as we deform $C \subseteq \mathcal{X}_0$ outside the central fibre of $\mathcal{X} \to \mathbf{D}$. Thus, denoting by

$$\text{Sing}^\circ(C) = \text{Sing}(C) \setminus \{\text{nodes on } E\},$$

we may restrict our attention to the versal map

$$(3.3.3) \quad \varphi = \prod_{p \in \text{Sing}^\circ(C)} \phi_p : \bigcap_{p \in \text{Sing}^\circ(C)} U_p \to \left( \prod_{\substack{m \geqslant 2 \\ m\text{-tacnode}}} \prod_{p \in E \cap C} \Delta_{(C,p)} \right) \times \prod_{\substack{p \in C \setminus E \\ \text{node}}} \Delta_{(C,p)}$$

and its differential at $[C]$,

$$d\varphi_{[C]} : H^0(C, \mathcal{N}_{C|\mathcal{X}}) \to \left( \bigoplus_{\substack{m \geqslant 2 \\ m\text{-tacnode}}} \bigoplus_{p \in E \cap C} T^1_{C,p} \right) \oplus \bigoplus_{\substack{p \in C \setminus E \\ \text{node}}} T^1_{C,p},$$

obtained by composing $d\phi_{[C]}$ with the natural projection map. By Lemma (2.9) and Corollary (2.10), and using the notation therein, we know that

$$d\varphi_{[C]}\bigl(H^0(C, \mathcal{N}_{C|\mathcal{X}})\bigr) \subseteq \left( \bigoplus_{\substack{m \geqslant 2 \\ m\text{-tacnode}}} \bigoplus_{p \in E \cap C} H_p \right) \oplus \bigoplus_{\substack{p \in C \setminus E \\ \text{node}}} T^1_{C,p}$$

and

$$d\varphi_{[C]}(H^0(C, \mathcal{N}_{C|\mathcal{X}_0})) \subseteq \left( \bigoplus_{\substack{m \geqslant 2 \\ m\text{-tacnode}}} \bigoplus_{p \in E \cap C} \Gamma_p \right) \oplus \bigoplus_{\substack{p \in C \setminus E \\ \text{node}}} T^1_{C,p}.$$

More precisely, observing that

$$\ker(d\varphi_{[C]}) = H^0(C, \mathcal{N}'_{C|\mathcal{X}_0}) \subseteq H^0(C, \mathcal{N}_{C|\mathcal{X}_0})$$

because every node on $E$ is trivially preserved if we deform $C$ in $\mathcal{X}_0$, by (3.3.2) we get that $d\varphi_{[C]}(H^0(C, \mathcal{N}_{C|\mathcal{X}_0}))$ has dimension

$$\sum_m \tau_m(m-1) + \delta = \dim\left[ \left( \bigoplus_{\substack{m \geqslant 2 \\ m\text{-tacnode}}} \bigoplus_{p \in E \cap C} \Gamma_p \right) \oplus \bigoplus_{\substack{p \in C \setminus E \\ \text{node}}} T^1_{C,p} \right].$$

Thus we obtain that

(3.3.4) $$d\varphi_{[C]}\big(H^0(C,\mathcal{N}_{C|\mathcal{X}_0})\big) = \Big(\bigoplus_{m\geqslant 2}\underset{m\text{-tacnode}}{\bigoplus_{p\in E\cap C}} \Gamma_p\Big) \oplus \bigoplus_{\underset{\text{node}}{p\in C\setminus E}} T^1_{C,p},$$

and

(3.3.5) $$d\varphi_{[C]}(H^0(C,\mathcal{N}_{C|\mathcal{X}})) = \Omega \oplus \bigoplus_{\underset{\text{node}}{p\in C\setminus E}} T^1_{C,p},$$

for some linear subspace $\Omega \subseteq \bigoplus_{m\geqslant 2}\bigoplus_p H_p$ containing $\bigoplus_{m\geqslant 2}\bigoplus_p \Gamma_p$ in codimension 1. Arguing as in the proof of Corollary (2.10), we get that:

- the versal map $\varphi$ defined in (3.3.3) maps $\mathcal{H}_0$ surjectively onto the product

$$\Big(\prod_{m\geqslant 2}\underset{m\text{-tacnode}}{\prod_{p\in E\cap C}}\Gamma_p\Big)\times \underline{0}\times ... \times \underline{0},$$

- the image of $\varphi$ is an analytic variety, smooth of dimension $\sum_m \tau_m(m-1)+\delta+1$ at $\underline{0}$, which is a product

$$\mathcal{W}\times \prod_{p\in C\setminus E\ \text{node}} \Delta_{(C,p)},$$

where $\mathcal{W}\subseteq \prod_{m\geqslant 2}\big(\prod_{p\in E\cap C\ m\text{-tacnode}}\Delta_{(C,p)}\big)$ may be identified with the image of the versal map

(3.3.6) $$\psi = \prod_{p\in \mathrm{Sing}^\circ(C)\cap E} \phi_p : \bigcap_{p\in \mathrm{Sing}^\circ(C)\cap E} U_p \to \prod_{m\geqslant 2}\underset{m\text{-tacnode}}{\prod_{p\in E\cap C}} \Delta_{(C,p)}.$$

In particular, $\mathcal{W}$ is a smooth analytic variety of dimension $\sum_m \tau_m(m-1)+1$, with tangent space $T_0\mathcal{W}=\Omega$ at $\underline{0}$, where $\Omega$ is as in (3.3.5). Keeping in mind the affine equations of the subspaces $\Gamma_p \subseteq H_p \subseteq \Delta_{(C,p)}$ in Lemma (2.9), the theorem follows by [1, Lemma 4.4] (see [V, Lemma (1.3)]), as we shall now explain.

To help the reader comparing our notation with that in [1], we observe that, for every $m$-tacnode $(C,p)$, our subspace $\Gamma_p \subseteq \Delta_{(C,p)}$ is denoted in [1, Lemma 4.4] by $\Delta_{i,m_i}$. Moreover, the product of the subspaces $\Gamma_p$ is denoted in [1, Lemma 4.4] by $\Delta_m$. Finally, the curves denoted by $\Gamma_1,...,\Gamma_k$ in the statement of [1, Lemma 4.4], will be denoted here by $\gamma_1,...,\gamma_k$.

By [1, Lemma 4.4], we have that the intersection of $\mathcal{W}$ with the locus of $\sum_m \tau_m(m-1)$-nodal curves in $\prod_{m\geqslant 2}\prod_{p\in E\cap C\ m\text{-tacnode}}\Delta_{(C,p)}$ is given by

$$\Big(\prod_{m\geqslant 2}\prod_{p\in E\cap C\ m\text{-tacnode}}\Gamma_p\Big)\cup \gamma_1 \cup \cdots \cup \gamma_k,$$

where $\gamma_1,\ldots,\gamma_k \subseteq \mathcal{W}$ are distinct, reduced, unibranched curves having intersection multiplicity exactly $\lambda$ with $\prod_{m\geqslant 2}\prod_{p\in E\cap C\ m\text{-tacnode}}\Gamma_p$ at the origin, with $\lambda$ and $k$ as defined in (3.0.1) and (3.0.2). The curves $\gamma_1,...,\gamma_k$ are all smooth at $[C]$ if $\lambda = m_p$, for some tacnode $p$ of $C$; otherwise they are all singular at $[C]$, with multiplicity $\lambda/\max\{m_p|\,p\ \text{tacnode of}\ C\}$.

The proof of the theorem is now complete if we observe that the versal map $\varphi$ defined by (3.3.3) has differential of maximal rank at $[C]$, thus

$$\varphi^{-1}(\gamma_1 \cup ... \cup \gamma_k \times \underline{0}\times ... \times \underline{0})$$

is locally a fibration over $\gamma_1 \cup ... \cup \gamma_k \times \underline{0} \times ... \times \underline{0}$ with smooth fibre of dimension $h^0(C, \mathcal{N}_{C|\mathcal{X}}) - \delta - \sum_m \tau_m(m-1)$, and it actually consists of $k$ analytic branches of $\mathcal{V}^{\mathcal{X}|\mathbf{D}}_{\delta + \sum_m \tau_m(m-1)}$ at $[C]$. The fact that, if $[C_t]$ is general in one of these analytic branches of $\mathcal{V}^{\mathcal{X}|\mathbf{D}}_{\delta + \sum_m \tau_m(m-1)}$, then (3.1.1) holds, follows by semicontinuity. □

**(3.4) Remark.** Theorem (3.1) proves in particular [I, Proposition (5.5)]; besides, [I] contains many applications of Theorem (3.1). In order to help the reader compare Theorem (3.1) with [I, Proposition (5.5)], we observe that [I, Proposition (5.5)] is stated under the assumption that Hilb($\mathcal{L}$) (notation from loc. cit.) is an irreducible component of the relative Hilbert scheme of curves in the family of surfaces $S \to \mathbf{D}$, see [I, Section 4].

**(3.5) Remark.** The failure of hypothesis 3 in Theorem (3.1) is not always an obstruction to study deformations of $C$ outside $\mathcal{X}_0$. It may happen for example that the singularities of $C$ do not impose the expected number of conditions because several nodes of $\mathcal{X}_0 \setminus E$ are "disconnecting nodes", i.e., normalizing $C$ at these nodes produces a reducible curve, and for some geometric reason it is not possible to deform $C$ on $\mathcal{X}_t$ to a reducible curve. In this case the disconnecting nodes cannot be preserved by deforming $C$ outside the central fibre, and the singularities of $C$ do not impose the expected number of conditions. This happens for example for curves constructed in [10, Section 4].

When hypothesis 3 in (3.1) Theorem fails, one has to restrict the attention on the subset of singularities imposing the right number of conditions (namely, 1 condition for a node on $\mathcal{X}_0 \setminus E$, and $m - 1$ conditions for an $m$-tacnode on $E$) to the linear system $|\mathcal{O}_{\mathcal{X}_0}(C)|$, and study the versal map only for this subset of singularities.

**(3.6) Remark.** Under the hypotheses of Theorem (3.1), by [10, Proposition 3.7 and Corollary 3.12], one has more generally that there exist deformations $C_t \subseteq \mathcal{X}_t$ of $C \subseteq \mathcal{X}_0$ so that every $m$-tacnode $p$ of $C$ is deformed to $d_{i,p}$ singularity of type $A_i$, for every multi-index $(d_{1,p},...,d_{r,p})$ such that $\sum_i i d_{i,p} = m - 1$, and every node of $C \subseteq \mathcal{X}_0$ is preserved. For these curves $C_t$ the equality (3.1.1) still holds, thus the corresponding generalized universal Severi variety is smooth at $[C_t]$, but describing its local geometry at $[C_0]$ is much more complicated.

# References

**Chapters in this Volume**

**External References**

Dipartimento di Matematica e Informatica, Università degli Studi della Calabria, Via Pietro Bucci, cubo 31B, 87036 Arcavacata di Rende (CS), Italy
galati@mat.unical.it

# Appendix B
# Flag semi-stable reduction of tacnodal curves at the limit of nodal curves

## by Ciro Ciliberto, Thomas Dedieu, and Concettina Galati

## 0 – Introduction

The objective of this note is to provide down-to-earth explanations of the following phenomenon, which is a central theme of the present collection of articles. It is studied in detail in [IV], following Caporaso and Harris [2, 3], and in [VIII] in the same context that we are considering here; it also appears in various guises in many other places.

Let $\mathcal{S} \to \mathbf{A}^1$ be a family of surfaces such that the general member $S_t$, $t \neq 0$, is a smooth projective surface, and the central member $S_0$ is the union of two smooth projective surfaces intersecting transversely along a smooth curve $B$. Let $\mathcal{L}$ be a line bundle on $\mathcal{S}$, and consider the family of Severi varieties

$$\mathring{\mathcal{V}} = \coprod_{t \neq 0} V^1(S_t, L_t) \longrightarrow \mathring{\mathbf{A}}^1 = \mathbf{A}^1 - \{0\},$$

where $V^1(S_t, L_t)$ is the family of 1-nodal curves in the linear system $|L_t|$ on the surface $S_t$, with $L_t = \mathcal{L}|_{S_t}$. Let $V_0$ be the family of curves in the linear system $|L_0|$ on $S_0$ having a *tacnode* along $B$, i.e., $V_0$ is the family of curves $C_0 \in |L_0|$, such that there exists a point $p \in B \cap C_0$, at which $C_0$ locally consists of two smooth branches $C_0' \subseteq S_0'$ and $C_0'' \subseteq S_0''$, both of which are tangent to $B$ at $p$. Then, at least in principle, $V_0$ appears in the central fibre of the closure $\overline{\mathcal{V}}$ of $\mathring{\mathcal{V}}$, with multiplicity 2.

We carry out a detailed local study of the nodes degenerating to a tacnode in the above situation, with stable reduction as a central concept. We investigate in particular the following questions, with various points of view: (i) in which way is such a tacnodal curve of the appropriate genus? (ii) why does it count with multiplicity 2? (iii) how should one modify the family of surfaces in order to see a tacnodal limit curve in $S_0$ as a nodal curve in an alternative semi-stable limit of the $S_t$'s?

The multiplicity 2 is fairly straightforward to understand, cf. (2.5). Seeing the tacnodal curve as a curve of the appropriate genus may be done in various ways, according to the chosen point of view, and is not always straightforward. We do provide a solution to problem (iii), but the conclusion is that although it is possible to see one tacnodal limit curve as a nodal curve in some semi-stable limit of the $S_t$'s, it is impossible to do it simultaneously for all curves of the family $V_0$. Therefore, tacnodal curves in the degeneration cannot be avoided, and have to be considered as definitely acceptable limits of nodal curves in their own right: in [I, Section 5], a degeneration of surfaces is said to be well-behaved if the limits of Severi varieties can all be accounted for, up to some reduction, by families of nodal and tacnodal curves: the upshot of our discussion here is that it is hopeless to ask that the limits of Severi varieties be accounted for only by families of nodal curves.

The text is organized as follows. Section 1 is meant as an overview of this note, and hopefully clarifies in what sense we study the questions (i-ii-iii) mentioned above: in this section, we outline various ways of viewing a tacnodal curve as an appropriate limit of 1-nodal curves, in order to illustrate all the local computations performed in the subsequent sections. In Section 2, we give the equations of the local situation we want to study, namely 1-nodal curves in smooth surfaces degenerating to a tacnodal curve in the transverse union of two smooth surfaces. In Section 3, we perform the stable reduction of a family of 1-nodal curves degenerating to a tacnodal curve, forgetting about the ambient degeneration of surfaces. The stable reduction process in this section is rather simple (we first normalize the total space, which normalizes the general fibres, then do some modifications, and finally return to 1-nodal general fibres), but it is unfortunately not suitable to answer question (iii) above. Therefore, in Section 4, we compute the same stable reduction without normalizing the total space. Finally, we are able to answer question (iii) in Section 5, providing what we call a *flag semi-stable reduction* of the nodes degenerating to a tacnode inside surface sheets degenerating to two smooth transverse sheets. The name refers to the fact that we perform semi-stable reduction simultaneously for the family of curves and for the ambient family of surfaces.

# 1 – Synthetic descriptions of nodal curves degenerating to tacnode

In this section we outline various ways to see a tacnodal curve as the (equigeneric) limit of nodal curves, in the context described in the introduction above. We do not enter in too many details for the moment; a fully rigorous treatment will be given in the next sections.

**(1.1) An example of ambient deformation space.** We will follow throughout this section the guiding example of smooth quartics degenerating to two quadrics, and their hyperplane sections.

Let $\mathcal{S} \to \mathbf{A}^1$ be a pencil of quartic hypersurfaces in $\mathbf{P}^3$, generated by a general quartic $S_\infty$, and the sum $S_0$ of two smooth quadrics $S_0'$ and $S_0''$ intersecting transversely along a quartic elliptic normal curve $B$. The total space $\mathcal{S}$ has ordinary double points at the finitely many points of $S_\infty \cap B$, which one usually prefers to resolve in one way or another; in our context we can safely ignore this, as our study is local, away from the incriminated points. This example is described in more detail in [I, 7.1].

**(1.2) Degeneration of smooth curves.** A general hyperplane section $C_t$ of $S_t$, $t \neq 0$, is a smooth plane quartic, and thus has genus 3. On the other hand, a general hyperplane section $C_0$ of $S_0$ is the union of two non-degenerate conics meeting transversely at 4 points: it has arithmetic genus 3, and it is a stable curve of genus 3.

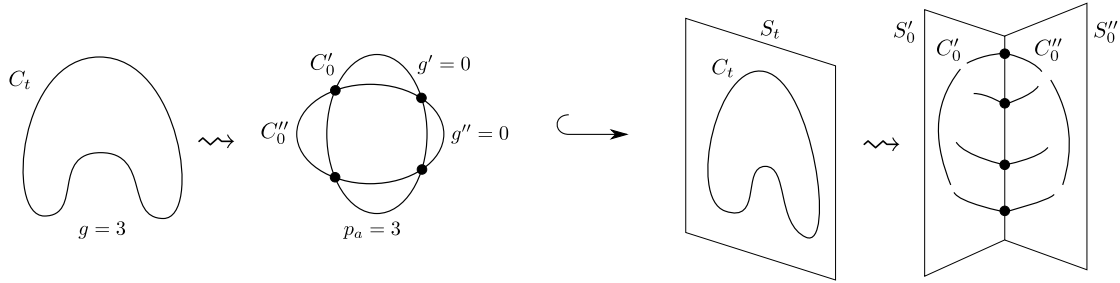

Figure 1: Degeneration of smooth curve sections

On Figure 1 above, the curves $C'_0$ and $C''_0$ are pictured in $S_0 = S'_0 \cup S''_0$ following a tropically inspired design, which is well suited for our drawings, see [**??**].

**(1.3) Degeneration of 1-nodal curves.** Continuing with our example of a pencil of quartic surfaces, we now consider for the rest of this section a family of 1-nodal hyperplane sections $C_t$ of $S_t$, $t \neq 0$, degenerating to a tacnodal hyperplane section $C'_0 \cup C''_0$ of $S_0$, i.e., $C'_0$ and $C''_0$ are non-degenerate conics on $S'_0$ and $S''_0$ respectively, tangent to $B$ at some point $q$, and meeting transversely at two other points $p_1, p_2 \in B$ (recall that $B$ is the curve $S'_0 \cap S''_0$).

*(1.3.1) An abuse of presentation.* In fact, as we shall see in Section 2, such a family must be 2-valued, i.e., the curve $C_t$ must consist of two 1-nodal hyperplane sections for $t \neq 0$, and the central fibre $C_0$ must be the double of $C'_0 \cup C''_0$ on $S_0$. This is an important point but we shall ignore it for the moment, while we are outlining the situation, as it is not the aspect we want to focus on in this note. Thus, for the sake of keeping this outline simple, we shall pretend for the rest of this section that the curves $C_t$ are irreducible 1-nodal curves, and $C_0$ is a reduced tacnodal curve, without any further mention. In particular, on Figure 2 and the other figures in this section, we only picture one 1-nodal curve tending to $C'_0 \cup C''_0$, whereas in fact there are two. From Section 2 on, our treatment will be fully rigorous without any abuse of presentation.

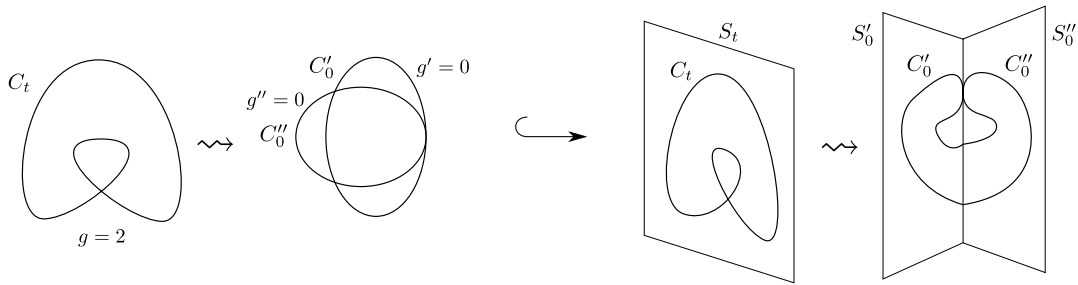

Figure 2: Degeneration of 1-nodal curves to a tacnodal curve

The curves $C_t$, $t \neq 0$, have geometric genus 2. In the next paragraphs we give various ways of appropriately seeing $C_0$ as a genus 2 limit of the $C_t$'s.

**(1.4) Degeneration of the normalizations.** We may incarnate the 1-nodal curves $C_t$, $t \neq 0$, as the stable maps $\bar{C}_t \to S_t$ given by the normalizations $\bar{C}_t \to C_t$; these are maps from smooth, genus 2, sources. From this perspective, the relevant incarnation of the curve $C_0$ is as the image of a stable map as pictured on Figure 3 below, with source a stable, genus 2, curve consisting of two smooth rational components meeting transversely at three points (this source is a partial normalization of the tacnode of $C_0$).

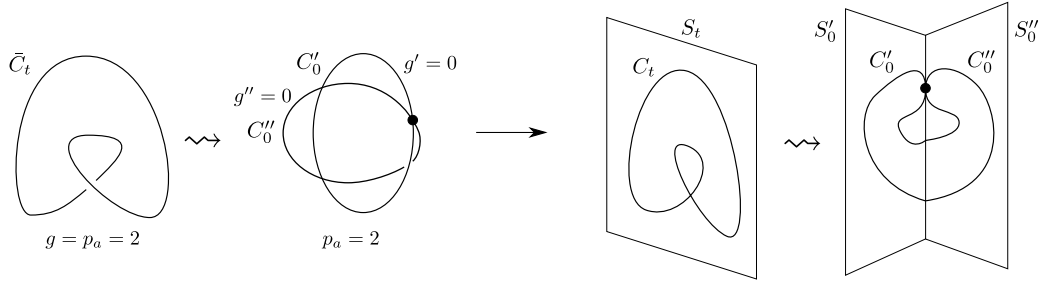

Figure 3: Degeneration of the normalizations of 1-nodal curve sections

To see this, one considers the family $\mathcal{C}$ of the curves $C_t$, $t \in \mathbf{A}^1$, and normalize it. It follows from Teissier's simultaneous resolution theorem [6] that the fiber of the normalization $\bar{\mathcal{C}}$ over $t \neq 0$ is the normalization $\bar{C}_t$. The fibre over 0 is computed in (3.2) below; it is singular over the tacnode, because there is a jump in $\delta$-invariant when the node of $C_t$ degenerates to the tacnode of $C_0$. This is compensated for by the jump in the number of irreducible components, which allows for the genus to be constant.

The latter point of view is arguably the more transparent. It is yet also desirable to understand the situation by seeing the curves $C_t$ as singular curves in the surfaces $S_t$.

**(1.5) Degeneration of the stable maps.** The 1-nodal curves $C_t$ are stable curves of arithmetic genus 3. As such, their limit $C_0$ can be modified into a stable curve of genus 3. In other words, the closed immersions $C_t \hookrightarrow S_t$ are stable maps of genus 3, and have a limit as such. These limits are pictured on Figure 4 below.

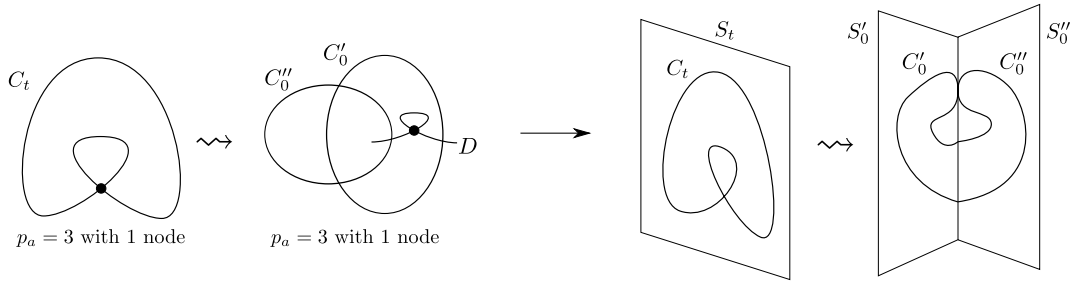

Figure 4: Degeneration as stable maps of the immersions of 1-nodal curve sections

The component $D$ on Figure 4 is a 1-nodal rational curve ($p_a = 1$ and $p_g = 0$); it is contracted to the tacnode of $C_0$ by the stable map with target $S_0$ which is the limit of the maps $C_t \hookrightarrow S_t$. The computations are carried out in Section 3.

One nice feature of this point of view is that we see our limit of curves of arithmetic genus 3 with 1 node as a curve of arithmetic genus 3 with 1 node as well; the node in question is that of the component $D$ (of course, the source of the limit stable map has other nodes; what we are saying is that the node of $D$ is the limit of the node of $C_t$ as $t$ tends to 0). The price one pays is that this stable limit of maps has a contracted component.

*(1.5.1) Observation.* The above model may be used to recover that described in (1.4). To do so, one simultaneously normalizes the nodes of all the curves $C_t$: on the central fibre, this normalizes the node of $D$, which then becomes a destabilizing smooth rational component; the stable model of the central fibre is obtained by contracting this component, which gives the central fibre of the family of curves described in (1.4).

**(1.6) Degeneration of the stable maps, II.** It frequently happens in practice that one has to modify the family of surfaces $\mathcal{S} \to \mathbf{A}^1$ by adding one or more ruled surfaces over $B$ between $S'_0$ and $S_0$. Typically, one performs a base change $t \in \mathbf{A}^1 \mapsto t^m \in \mathbf{A}^1$, and then resolves the singularities of the obtained family. In particular, the description below is very useful to understand [**??**].

Let us assume here that there is one ruled surface $\Sigma_B$ over $B$ between our two quadrics $S'_0$ and $S''_0$, and that our family of surfaces has reduced central fibre. In this model of $\mathcal{S}$, our tacnodal curve $C_0$ is modified as indicated on Figure 5, by the adding of (i) reduced fibres of $\Sigma_B$ to connect the points on $C'_0$ and $C''_0$ that previously were transverse intersection points of $C'_0$ and $C''_0$, and (ii) one double fibre to connect the two points that previously were the tacnode.

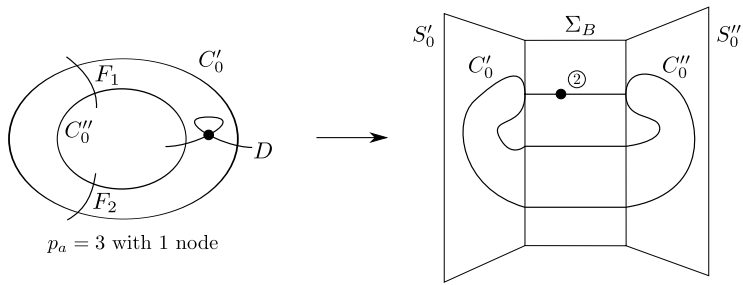

Figure 5: Degeneration as stable maps of the immersions of 1-nodal curve sections, II

The stable map limit of the immersions $C_t \hookrightarrow S_t$ in this context is easily deduced from the limit in paragraph (1.5). The two intersection points between $C'_0$ and $C''_0$ are replaced by smooth rational curves $F_1$ and $F_2$, mapped isomorphically to the corresponding fibres of $\Sigma_B$. The nodal curve $D$ on the other hand is mapped $2:1$ to the double fibre "at the tacnode", with ramification divisor the sum of the two points $D \cap C'_0$ and $D \cap C''_0$; the node of $D$ is mapped to the limit of the nodes of the curves $C_t$ as $t$ approaches 0, which lies on the double fibre of $\Sigma_D$; the pre-image of this point in $D$ consists only of the node of $D$.

Alternatively, this limit may be computed using the reduction which we discuss next, and which is constructed in Section 5.

**(1.7) Flag semi-stable reduction.** Ultimately, one may desire to see the stable limit of the 1-nodal curves $C_t$, depicted in Figure 4, not mapped to $S_0$ with some contracted components, but rather immersed in some semi-stable limit of the surfaces $S_t$. It is of course asserted by the general stable reduction theorems that this is possible, provided one allows the replacement of the stable limit of the $C_t$'s by another semi-stable model. We carry this out explicitly in Section 5, and the result is depicted in Figure 13 there. We dub this a *flag semi-stable reduction* of the pairs $(C_t, S_t)$, $t \neq 0$.

## 2 – The local situation

In this section, we set up an explicit incarnation of the situation that we study throughout this text, of a family of 1-nodal curves degenerating to a tacnodal curve inside a degeneration of smooth surfaces to the transverse union of two smooth surfaces. This is the same setup as in [IV].

**(2.1) Notation.** The notation $\mathbf{A}^4_{xyzt}$ denotes the affine 4-space equipped with affine coordinates $(x, y, z, t)$; we will also use obvious variations, e.g., $\mathbf{A}^1_t$ denotes the affine line with affine coordinate $t$. Similarly, $\mathbf{P}^1_{x:y}$ denotes the projective line with homogeneous coordinates $(x:y)$, and obvious variations will be used.

The arrow $\mathbf{A}^4_{xyzt} \to \mathbf{A}^1_t$ without further indication denotes the affine projection map $(x,y,z,t) \in \mathbf{A}^4 \to t \in \mathbf{A}^1$.

**(2.2)** We consider the degeneration of surfaces $\mathcal{S} \to \mathbf{A}^1_t$, where the total space $\mathcal{S}$ is the hypersurface of $\mathbf{A}^4_{xyzt}$ defined by the equation

$$\mathcal{S} : xy = t,$$

and the map $\mathcal{S} \to \mathbf{A}^1_t$ is the restriction of the affine projection $\mathbf{A}^4_{xyzt} \to \mathbf{A}^1_t$.

For all $t \neq 0$, the fibre $S_t$ of $\mathcal{S}$ over $t$ is smooth and irreducible, whereas the central fibre $S_0$ is the union of two planes $S'_0$ and $S''_0$, defined by the equations $t = y = 0$ and $t = x = 0$ respectively, which intersect transversely along the line defined by the equations $x = y = t = 0$.

**(2.3)** We want to study degenerations of nodal curves to a tacnodal curve inside this degeneration of surfaces, i.e., families of curves $\mathcal{C} \to \mathbf{A}^1$ fitting in a commutative diagram

$$\begin{array}{c} \mathcal{C} \subseteq \mathcal{S} \\ \searrow \downarrow \\ \mathbf{A}^1 \end{array}$$

such that the fibres $C_t$ with $t \neq 0$ are 1-nodal curves, and the central fibre $C_0$ is a curve with a tacnode along the double curve of $S_0$, by which we mean that $C_0$ consists of two smooth branches $C'_0$ and $C''_0$ contained in $S'_0$ and $S''_0$ respectively, which intersect in only one point where they are both tangent to the curve $S'_0 \cap S''_0$.

**(2.4)** As a matter of fact, families of curves as in (2.3) cannot exist, because if they would, then the singularities of the curves would form a section of $\mathcal{S} \to \mathbf{A}^1$ meeting the central fibre along its double curve, which is impossible since the total space $\mathcal{S}$ is smooth.

As we shall now see, there exist instead families of curves with central fibre the double of the curve $C_0$ as in (2.3) above, and fibre over $t \neq 0$ the sum of two 1-nodal curves. We let, indeed, $^\sharp\mathcal{C} \to \mathbf{A}^1_t$ be defined in $\mathcal{S}$ by the equation

$$^\sharp\mathcal{C} : (x + y - z^2)^2 = 4t.$$

Its central fibre is indeed the double of a tacnodal curve $C_0$ as described in (2.3): it is defined by $t = 0$, and thus it is the double of the curve with the following equations in $\mathbf{A}^3_{xyz}$ (the latter being the central fibre of $\mathbf{A}^4_{xyzt} \to \mathbf{A}^1_t$):

$$\begin{cases} xy = 0 \\ x + y - z^2 = 0 \end{cases} \iff \begin{cases} y(y-z)^2 = 0 \\ -x = y - z^2. \end{cases}$$

On the other hand, for all $t_0 \neq 0$, the fibre $^\sharp C_{t_0}$ consists of two irreducible 1-nodal curves $C^\pm_{t_0}$, with the following equations in $\mathbf{A}^3_{xyz}$, the fibre of $\mathbf{A}^4_{xyzt} \to \mathbf{A}^1_t$ over $t_0$:

(2.4.1) $$\begin{cases} -xy + t_0 = 0 \\ x + y - z^2 \mp 2\sqrt{t_0} = 0 \end{cases} \iff \begin{cases} (y \mp \sqrt{t_0})^2 - yz^2 = 0 \\ -x = y - z^2 \mp 2\sqrt{t_0}. \end{cases}$$

**(2.5)** The family $^\sharp\mathcal{C}$ inside $\mathcal{S}$ over $\mathbf{A}^1$ is the universal local model for 1-nodal curves in smooth surfaces degenerating to a tacnodal curve along the transverse intersection of two smooth surface sheets.

The necessity of having two 1-nodal curves in each fibre $S_t$, $t \neq 0$, tending to the tacnodal curve $C_0$ explains why a tacnodal curve in the limit amounts for (at least) two 1-nodal curves.

**(2.6)** Next, we shall perform a degree 2 base change in order to separate $^\sharp\mathcal{C}$ in two distinct families $\mathcal{C}^+$ and $\mathcal{C}^-$ with the same central fibre $C_0$.

Thus, let $\tilde{\mathcal{S}} \to \mathbf{A}_t^1$ be the family of surfaces with total space defined by

$$\tilde{\mathcal{S}} : xy = t^2$$

in $\mathbf{A}^4_{xyzt}$, and with the map $\tilde{\mathcal{S}} \to \mathbf{A}_t^1$ given by the affine projection $\mathbf{A}^4_{xyzt} \to \mathbf{A}_t^1$. Note that $\tilde{\mathcal{S}}$ is singular along the double curve of its central fibre $S_0$. There are two families of 1-nodal curves tending to a curve with a tacnode along the double curve of $S_0$, namely $\mathcal{C}^\pm$ defined by

$$\mathcal{C}^\pm : x + y - z^2 = \pm 2t$$

inside $\tilde{\mathcal{S}}$. We will usually only look at one of these two families, namely $\mathcal{C}^-$, which we will simply denote by $\mathcal{C}$, in spite of the conflict of notation with (2.3).

## 3 – Abstract reduction of nodes degenerating to a tacnode, via normalization of the nodes

In this section, we compute the stable reduction of the family $\mathcal{C} \to \mathbf{A}^1$ of 1-nodal curves degenerating to a tacnodal curve described in the previous section. We start by normalizing the total space $\mathcal{C}$, which has the effect of normalizing the fibres $C_t$ with $t \neq 0$. Then, we make some modifications in order for the pre-images of the nodes to form two disjoint sections. The last operation is to glue these two sections, so that we come back to a family with general member a 1-nodal curve.

**(3.1) Setup.** We consider the family $\mathcal{C}$ defined as $\mathcal{C}^-$ in (2.6) (recall also Notation (2.1)). It may be identified with its affine projection in $\mathbf{A}^3_{yzt}$, which is defined by the single equation

(3.1.1) $$\mathcal{C} : \quad y(y - z^2 + 2t) + t^2 = 0 \iff (y+t)^2 - yz^2 = 0,$$

obtained by eliminating $x$ from the two equations $xy = t^2$ and $x + y - z^2 + 2t = 0$ defining $\mathcal{C}^-$.

The total space $\mathcal{C}$ is a surface with a double curve along

$$R : \ y + t = z = 0,$$

and a pinch point at the origin (which is the intersection point of $R$ with the central fibre), and is otherwise smooth. For all $t \neq 0$, $C_t$ is a 1-nodal curve with its node at $C_t \cap R$, and $C_0$ is a 1-tacnodal curve with its tacnode at the pinch point $C_0 \cap R$. We let $C_0', C_0''$ be the two irreducible components of $C_0$, defined respectively by $y = t = 0$ and $y - z^2 = t = 0$.

**(3.2) Normalization of the total space.** The first move in the present approach is to normalize the total space $\mathcal{C}$. We do this by blowing-up the double curve $R$ in $\mathcal{C}$. Let $\varepsilon_1 : \mathcal{C}_1 \to \mathcal{C}$ be the corresponding morphism.

The surface $\mathcal{C}_1$ lives in the subscheme of the product $\mathbf{A}^3_{yzt} \times \mathbf{P}^1_{t_1:z_1}$ defined by the equation

(3.2.1) $$z_1(y+t) = t_1 z,$$

in which it is defined by the equation

(3.2.2) $$t_1^2 - yz_1^2 = 0.$$

First note that the surface $\mathcal{C}_1$ is smooth. The exceptional locus of $\varepsilon_1$ (i.e., the total transform via $\varepsilon_1$ of $R$) is the curve $R_1$ defined by

$$y + t = z = 0$$

in $\mathcal{C}_1$; it is a bisection of $\mathcal{C}_1 \to \mathbf{A}_t^1$ (which is the composed map $\mathcal{C}_1 \to \mathcal{C} \to \mathbf{A}_t^1$), mapped $2:1$ onto $R$ by $\varepsilon_1$ with ramification over the tacnode. The central fibre of $\mathcal{C}_1$ is defined by $t = 0$, which gives $z_1 y = t_1 z$ by (3.2.1), hence (3.2.2) becomes

$$t_1(t_1 - zz_1) = 0;$$

thus the central fibre of $\mathcal{C}_1$ consists of two branches intersecting transversely, which we shall still denote by $C_0'$ and $C_0''$. The situation is summed up in Figure 6.

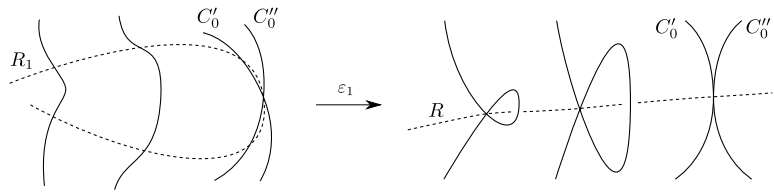

Figure 6: Normalization of the total space, $\varepsilon_1 : \mathcal{C}_1 \to \mathcal{C}$

**(3.3) Base change.** The next step is to perform the base change $t \in \mathbf{A}^1 \mapsto t^2 \in \mathbf{A}^1$ in order to separate the pre-images of the nodes of the $C_t$'s in two sections. Let $\tilde{\mathcal{C}}_1 \to \mathcal{C}_1$ be the corresponding double cover. In practice, this amounts to replacing $t$ by $t^2$ everywhere in the equations defining $\mathcal{C}_1$. The notable fact is that this yields an ordinary double point on the surface $\tilde{\mathcal{C}}_1$ at the point where $C_0'$ and $C_0''$ intersect. We call $R_1'$ and $R_1''$ the two irreducible components of the pre-image of $R_1$.

**(3.4) Resolution.** The next thing we do is to resolve the ordinary double point of $\tilde{\mathcal{C}}_1$. We call $\bar{\mathcal{C}}_1 \to \tilde{\mathcal{C}}_1$ the corresponding map. It is a birational morphism, with exceptional locus a $(-2)$-curve which we call $E$. The resulting family of curves is depicted on the right-hand-side of Figure 7.

**(3.5) Reconstruction of the nodes.** Finally, since what we are looking for is the stable reduction of the 1-nodal curves $C_t$, $t \neq 0$, we reconstruct the nodes that have been normalized in the first step by identifying the two disjoint sections $R_1'$ and $R_1''$ of the family $\bar{\mathcal{C}}_1 \to \mathbf{A}^1$. We call $\mathcal{C}_2 \to \mathbf{A}^1$ the obtained family of curves, see the left-hand-side of Figure 7. Its central fibre is reduced, with a 1-nodal rational curve $D$ (image of $E$) connecting $C_0'$ and $C_0''$; the node of $D$ is the limit position of the nodes of the curves $C_t$, $t \neq 0$.

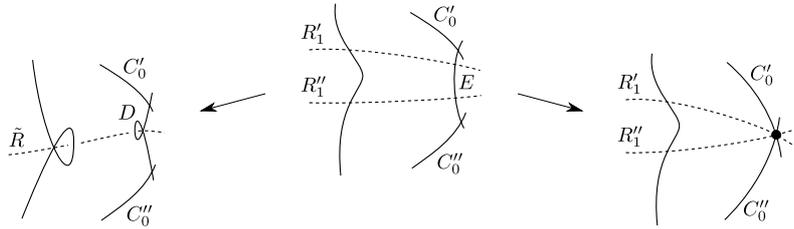

Figure 7: After the base change, $\mathcal{C}_2 \leftarrow \bar{\mathcal{C}}_1 \to \tilde{\mathcal{C}}_1$

**(3.6) Why we're not satisfied with this.** We want to perform the stable reduction of the nodal curves $C_t$ (or possibly a semi-stable variant of it), while simultaneously performing the semi-stable reduction of the degeneration of surfaces $\tilde{\mathcal{S}} \to \mathbf{A}^1$, with notation as in Section 2. The main problem with the above described stable reduction process is that the normalization of the family of curves and later the reconstruction of the nodes necessitate unwanted modifications of the surfaces $S_t$, $t \neq 0$.

Indeed, the normalization of $\mathcal{C}$ inside $\tilde{\mathcal{S}}$ requires to blow-up every surface $S_t$, $t \neq 0$ at the point where the curve $C_t$ has a node, which introduces a $(-1)$-curve. Finally, once the remaining steps of the reduction of the $C_t$'s are performed, we have to recontract all these $(-1)$-curves, and this messes everything up, in particular this destabilizes the central fibre. We leave it to the reader to check for himself what indeed happens.

Our solution is to follow another path to the stable reduction of $\mathcal{C}$, which doesn't involve the normalization of the nodal $C_t$'s. This is what we do in the next section.

## 4 – Abstract reduction of nodes degenerating to a tacnode, maintaining the nodes

### 4.1 – The construction

In this section we present our construction of the stable reduction of the family of 1-nodal curves $C_t$ without normalizing the nodes. The setup is the same as that formulated in (3.1), and we don't repeat it here. We use Notation (2.1), as always.

**(4.1) Blow-up of the tacnode.** Let $\varepsilon_1 : \mathcal{C}_1 \to \mathcal{C}$ be the blow-up of $\mathcal{C}$ at the origin of $\mathbf{A}^3_{yzt}$ (here we use the same notation as in paragraph (3.2) for a different blow-up, but we believe this will not cause any confusion). We will see that the reduced exceptional locus of $\varepsilon_1$ is a smooth $\mathbf{P}^1$, which we shall call $E_1$, and the central fibre of $\mathcal{C}_1$ consists of the proper transform of $C_0$, which is the transverse union of two smooth curves, namely the proper transforms of $C'_0$ and $C''_0$, which we continue to denote by the same symbols, plus $E_1$ with multiplicity 2; the total space $\mathcal{C}_1$ has a double curve along the proper transform of $R$, still with a pinch point at its intersection with the central fibre; besides this double curve, $\mathcal{C}_1$ has an ordinary double point at $C'_0 \cap C''_0$, which lies on $E_1$, and no other singularities. See Figure 8.

We realize $\mathcal{C}_1$ as the proper transform of $\mathcal{C}$ in the subvariety of $\mathbf{A}^3_{yzt} \times \mathbf{P}^2_{y_1:z_1:t_1}$ defined by the equations

$$\mathrm{rk} \begin{pmatrix} y & z & t \\ y_1 & z_1 & t_1 \end{pmatrix} < 2.$$

In the affine chart $t_1 = 1$, $\mathcal{C}_1$ has the equation

(4.1.1) $\qquad \mathcal{C}_1 \text{ in } [t_1 = 1]: \quad (y_1 + 1)^2 - y_1 z_1^2 t = 0.$

We thus see that $\mathcal{C}_1$ has a double curve along the proper transform $R_1$ of $R$, which has equations $y_1 + 1 = z_1 = 0$, with a pinch point at the point $(-1 : 0 : 1)$ on the exceptional divisor (the latter is the fibre of $\mathcal{C}_1$ over the origin in the projection $\mathbf{A}^3_{yzt} \times \mathbf{P}^2_{y_1:z_1:t_1} \to \mathbf{A}^3_{yzt}$). In this chart, the proper transform of $C_0$ is at infinity, and the central fibre $t = 0$ is the exceptional divisor, namely the double line $(y_1 + t_1)^2$ in $\mathbf{P}^2$.

In the affine chart $z_1 = 1$, we get the equation

(4.1.2) $\qquad\qquad \mathcal{C}_1 \text{ in } [z_1 = 1]: \quad (y_1 + t_1)^2 - y_1 z = 0.$

We thus see that $\mathcal{C}_1$ has an ordinary double point at the point $(0 : 1 : 0)$ on the exceptional divisor. The central fibre is defined by the equation $t = 0$, which gives $t_1 z = 0$: $t_1 = 0$ is the proper transform of $C_0$, and has equation $y_1(y_1 - z) = 0$, while $z = 0$ defines the exceptional divisor: it is the divisor defined by $(y_1 + t_1)^2 = 0$, i.e., $E_1$ with multiplicity 2.

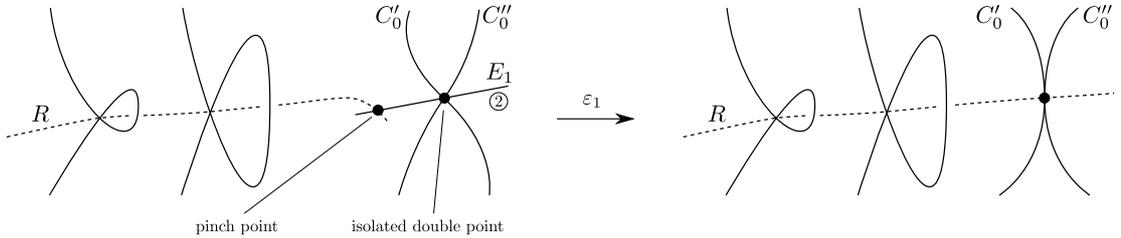

Figure 8: Blow-up of the tacnode, $\varepsilon_1 : \mathcal{C}_1 \to \mathcal{C}$

**(4.2) Blow-up of the ordinary double point.** Let $\varepsilon_2 : \mathcal{C}_2 \to \mathcal{C}_1$ be the blow-up at the point $C_0' \cap C_0'' \cap E_1$. The three curves $C_0', C_0'', E_1$ in $\mathcal{C}_1$ intersect transversely at this point, hence they are separated by this blow-up. The reduced exceptional locus is a smooth $\mathbf{P}^1$, which we shall call $E_2$; it appears with multiplicity 2 in the central fibre. See the right-hand-side of Figure 9 for a picture of the situation.

Note that the central fibre of $\mathcal{C}_1$ is $C_0' + C_0'' + 2E_1$, hence it has multiplicity 4 at the centre of the blow-up $\varepsilon_2$. However, as we shall now see, $E_2$ indeed appears only with multiplicity 2 in the central fibre of $\mathcal{C}_1$, because the centre of $\varepsilon_2$ is a singular point (of multiplicity 2) of $\mathcal{C}_1$.

This is fairly elementary, so we'll only sketch the computations. It is enough to consider the part of $\mathcal{C}_1$ in the affine chart $z_1 = 1$, which has equation (4.1.2). We add the new set of homogeneous variables $(y_2 : z_2 : t_2)$ with the new equations

$$\operatorname{rk} \begin{pmatrix} y_1 & z & y_1 + t_1 \\ y_2 & z_2 & t_2 \end{pmatrix} < 2.$$

Then, $\mathcal{C}_2$ has equation $t_2^2 - y_2 z_2 = 0$. Its central fibre has equation $t = 0$, which, in the affine chart $y_2 = 1$, since $t = t_1 z$, $t_1 = y_1(t_2 - 1)$, and $z = y_1 z_2$, reads

$$(t_2 - 1) z_2 y_1^2 = 0.$$

Now, $t_2 - 1 = 0$ defines one of the two branches $C_0'$ and $C_0''$ (the other is invisible in this affine chart), $z_2 = 0$ is the equation of the proper transform of $E_1$, and $y_1^2 = 0$ that of the exceptional curve $E_2$ with multiplicity two.

**(4.3) Base change and normalization.** Let $\tilde{\mathcal{C}}_2 \to \mathcal{C}_2$ be the double covering of $\mathcal{C}_2$ ramified over the central fibre, and $\bar{\mathcal{C}}_2 \to \tilde{\mathcal{C}}_2$ the partial normalization obtained by leaving only the singularity along the proper transform of $R$. As we shall see, the family $\tilde{\mathcal{C}}_2 \to \mathbf{A}^1$ is a semi-stable reduction of $\mathcal{C} \to \mathbf{A}^1$.

Algebraically, $\tilde{\mathcal{C}}_2$ is defined by adding a square root of $t$, i.e., a new variable $s$ with the new equation $s^2 = t$, to the equations of $\mathcal{C}_2$. The normalization process may be chopped in several local computation units, which we gather in Section 4.2 below. The upshot is as follows. Denote by abuse of notation $E_1, E_2 \subseteq \tilde{\mathcal{C}}_2$ the reduced pre-images of $E_1, E_2 \subseteq \mathcal{C}_2$. The normalization may be obtained by successively blowing-up the Weil divisors $E_1$ and $E_2$.

First, the blow-up of $E_1$ produces two smooth rational curves $E_1'$ and $E_1''$, disposed as pictured in Figure 9, and both mapped isomorphically onto $E_1$ (to see this, use (4.7) at the pinch point, and (4.6.1) at $E_1 \cap E_2$). Then, the blow-up of $E_2$ produces a double covering $\bar{E}_2 \to E_2$, with $\bar{E}_2$ a smooth rational curve disposed as pictured in Figure 9; this double covering branches at the two points $C_0' \cap E_2$ and $C_0'' \cap \bar{E}_2$, and the associated involution exchanges the two points $E_1' \cap \bar{E}_2$ and $E_1'' \cap \bar{E}_2$ (to see this use (4.5) at the two points $C_0' \cap E_2$ and $C_0'' \cap \bar{E}_2$, and (4.6.2) at the point $(E_1' + E_1'') \cap E_2$).

*(4.3.1) A useful variant.* Although arguably less natural, blowing-up first $E_2$ and then $E_1$ produces the same result. The blow-up of $E_2$ produces a double cover $\tilde{E}_2 \to E_2$, ramified at the two points $C_0' \cap \tilde{E}_2$ and $C_0'' \cap \tilde{E}_2$, and with a node over the point $E_1 \cap E_2$, as depicted in Figure 9. Then, the blow-up along $E_1$ produces two curves $E_1'$ and $E_1''$ over $E_1$, and resolves the node of $\tilde{E}_2$.

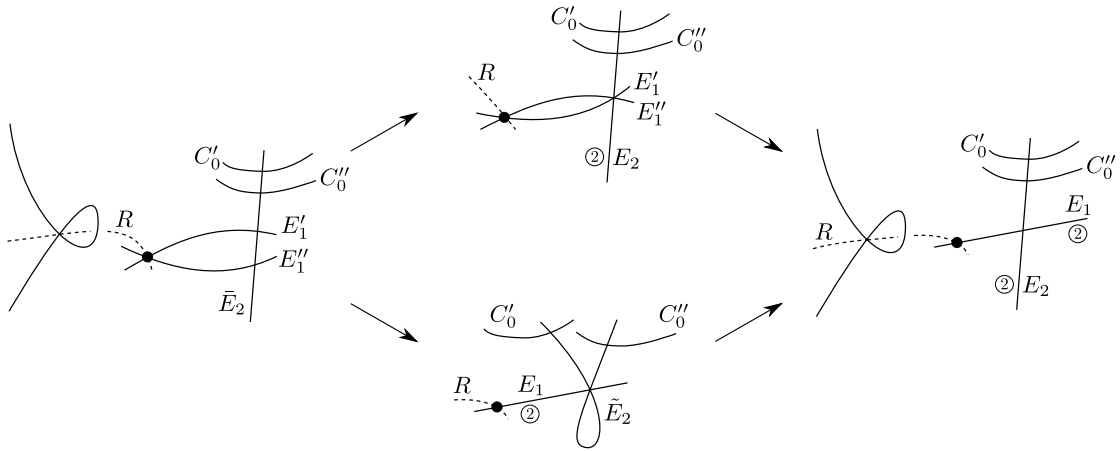

Figure 9: Removing the multiplicities, $\bar{\mathcal{C}}_2 \to \tilde{\mathcal{C}}_2$

**(4.4) Conclusion.** It follows from the previous considerations that the family $\bar{\mathcal{C}}_2 \to \mathbf{A}^1$ is indeed semi-stable. One may check that the final family has an ordinary double point all along the proper transform of $R$, i.e., the pinch point has disappeared. To build a stable model out of this semi-stable model, one has to contract the two rational curves $E_1'$ and $E_1''$, which turns the curve $\bar{E}_2$ into a 1-nodal rational curve, with its node at the limit of the nodes of the curves $C_t$, $t \neq 0$.

## 4.2 – Ready for use toroidal computations

In the following paragraphs we gather some local computation, useful to normalize families of curves gotten by base change from families with non-reduced but normal crossing central fibre; the normalization will be expressed as a succession of blow-ups. We only include the computations needed for (4.3) above, but the reader will undoubtedly figure out how to generalize them to similar situations. [5, p. 125] is a useful reference on this subject.

**(4.5) Local situation $x^2 y = t^2$.** We consider the surface $\mathcal{C}$ defined by this equation in $\mathbf{A}^3_{xyt}$, together with the fibration over $\mathbf{A}^1_t$ induced by the projection $\mathbf{A}^3_{xyt} \to \mathbf{A}^1_t$. This is the family gotten by a degree two base change totally ramified at the origin from a central fibre with one double component $2A$ and one reduced component $B$ meeting transversely, defined by the equations $x^2 = 0$ and $y = 0$ respectively.

The normalization of $\mathcal{C}$ amounts to the blow-up $\bar{\mathcal{C}} \to \mathcal{C}$ of the Weil divisor $A$, defined by the two equations $x = t = 0$ (note that $\mathcal{C}$ has an ordinary double point at the generic point of $A$, and a pinch point at the origin). To write down this blow-up, we add the set of homogeneous coordinates $(x_1 : t_1)$ with the new equation $x_1 t = t_1 x$.

The surface $\bar{\mathcal{C}}$ is defined by the homogeneous equation

(4.5.1) $$x_1^2 y = t_1^2.$$

One thus sees that $\bar{\mathcal{C}}$ is smooth. We claim that its central fibre has the form $\bar{A} + \bar{B}$, with two maps $\bar{A} \to A$ and $\bar{B} \to B$ induced by $\bar{\mathcal{C}} \to \mathcal{C}$, respectively a double covering of $A$ ramified at the point $\bar{A} \cap \bar{B}$, and an isomorphism.

Indeed, the central fibre is defined by $t = 0$, which gives $t_1 x = 0$. The equations $x = t = 0$ define the exceptional divisor in $\bar{\mathcal{C}}$; using (4.5.1), one sees that the exceptional divisor is a double cover of the line $x = t = 0$ in $\mathbf{A}^3_{xyt}$, which is the divisor $A$. The equation $t_1 = 0$, on the other hand, implies $y = 0$ by (4.5.1), and one sees that it defines a divisor projecting isomorphically to $B$.

**(4.6) Local situation $x^2 y^2 = t^2$.** We consider the surface $\mathcal{C}$ defined by this equation in $\mathbf{A}^3_{xyt}$, together with the fibration $\mathcal{C} \to \mathbf{A}^1_t$ defined by affine projection. The central fibre of $\mathcal{C}$ consists of the two Weil divisors $A$ ($y = t = 0$) and $B$ ($x = t = 0$), each with multiplicity 2.

*(4.6.1) Blow-up of the Weil divisor $A$.* We add the homogenous coordinates $(y_1 : t_1)$, with the equation $y_1 t = t_1 y$; the blow-up $\mathcal{C}_1$ of $\mathcal{C}$ is defined by the equation $x^2 y_1^2 = t_1^2$. Let us restrict to the chart $y_1 = 1$, to fix ideas; then $\mathcal{C}_1$ is defined by the equation $x^2 = t_1^2$. The central fibre is defined by $t = 0$, which implies $t_1 y = 0$, hence it consists of the proper transform $B_1$ of $B$, still with multiplicity 2 (defined by $t = t_1 = 0$, hence $x^2 = 0$), plus the exceptional divisor, defined by $t = y = 0$ and $(x - t_1)(x + t_1) = 0$, which is thus a sum $A' + A''$ of two divisors mapping isomorphically to $A$.

*(4.6.2) Successive blow-up of the second Weil divisor $B$.* Now, we consider the blow-up $\mathcal{C}_2$ of $\mathcal{C}_1$ along the Weil divisor $B_1$, which is the proper transform of $B$ in $\mathcal{C}_1$. Thus we blow-up along the subscheme defined by $t_1 = x = 0$, and take the proper transform $\mathcal{C}_2$ of $\mathcal{C}_1$: we add the homogenous coordinates $(x_2 : t_2)$, with the equation $x_2 t_1 = t_2 x$; the blow-up $\mathcal{C}_2$ of $\mathcal{C}_1$ is defined by the equation $x_2^2 = t_2^2$ which gives $(x_2 : t_2) = (1 : \pm 1)$. The central fibre is defined by $t = 0$, which implies $t_1 y = 0$. If $t_1 = 0$, we get $x = 0$ and we obtain two reduced curves since $(x_2 : t_2) = (1 : \pm 1)$. These are the curves $B' + B''$ because $x = 0$. If $y = 0$, we get the equations $x = \pm t_1$, which give the two proper transforms of $A'$ and $A''$ respectively.

Of course the situation in the present paragraph is more or less trivial, as the surface $\mathcal{C}$ has two smooth irreducible components, defined by $xy = t$ and $xy = -t$ respectively, hence its normalization is the disjoint union of these two components. However, we find it useful to describe the normalization as an explicit composition of blow-ups.

**(4.7) Local situation $x^2 = y^2 t^2$.** We consider the surface $\mathcal{C}$ defined by this equation in $\mathbf{A}^3_{xyt}$, together with the fibration over $\mathbf{A}^1_t$ defined by the projection $\mathbf{A}^3_{xyt} \to \mathbf{A}^1_t$. The family $\mathcal{C} \to \mathbf{A}^1$ is gotten by a degree 2 base change from the family defined by $x^2 = y^2 t$, which has a double curve dominating $\mathbf{A}^1$, and a pinch point on the central fibre; the central fibre is a smooth curve with multiplicity 2.

We blow-up the Weil divisor $x = t = 0$: this introduces homogeneous coordinates $(x_1 : t_1) \in \mathbf{P}^1$, with the equation $x_1 t = t_1 x$, and the blown-up surface is defined by the homogeneous equation $x_1^2 = y^2 t_1^2$. The exceptional divisor is given by the latter equation together with $x = t = 0$, hence it is the union of two lines meeting at the point $(0, 0, 0, (0 : 1))$. The central fibre is defined by $t = 0$, hence also $t_1 x = 0$, and one sees that it equals the exceptional divisor.

The total space is singular along the curve $y = x_1 = 0$, where it uniformly has an ordinary double point. In particular, the pinch point has been resolved by the blow-up.

## 5 – The flag reduction

We now use the construction of Section 4 in order to perform the semi-stable reduction of both families $\mathcal{C} \to \mathbf{A}^1$ and $\tilde{\mathcal{S}} \to \mathbf{A}^1$ at the same time (in the notation of (2.6)).

**(5.1) Setup.** We briefly recapitulate the situation. We consider the threefold $\tilde{\mathcal{S}}$ defined by the equation

$$(5.1.1) \qquad \tilde{\mathcal{S}}: \quad xy = t^2$$

in $\mathbf{A}^4_{xyzt}$, and the map $\tilde{\mathcal{S}} \to \mathbf{A}^1$ induced by the projection $\mathbf{A}^4_{xyzt} \to \mathbf{A}^1_t$. Thus $\tilde{\mathcal{S}} \to \mathbf{A}^1$ is a family of surfaces. We call $S'_0$ and $S''_0$ the two components of its central fibre, defined by $y = t = 0$ and $x = t = 0$ respectively. The family of curves $\mathcal{C}$ studied in Sections 3 and 4 is a Cartier divisor of $\tilde{\mathcal{S}}$, defined in $\mathbf{A}^4_{xyzt}$ by the equations

$$(5.1.2) \qquad \mathcal{C}: \begin{cases} y + x - z^2 + 2t = 0 \\ xy = t^2. \end{cases}$$

**(5.2) Blow-up of the tacnode.** Let $\varepsilon_1 : \mathcal{S}_1 \to \tilde{\mathcal{S}}$ be the blow-up of $\tilde{\mathcal{S}}$ at the point $\mathbf{z} = (0, 0, 0, 0) \in \mathbf{A}^4$, and let $\mathcal{C}_1$ be the proper transform of $\mathcal{C}$ in $\mathcal{S}_1$; the restriction $\varepsilon_1|_{\mathcal{C}_1}$ is the same map $\mathcal{C}_1 \to \mathcal{C}$ as in (4.1). We call $\Sigma_1$ the exceptional divisor of $\varepsilon_1$.

The blow-up amounts to adding the new homogeneous variables $(x_1 : y_1 : z_1 : t_1) \in \mathbf{P}^3$, with the equations

$$\mathrm{rk} \begin{pmatrix} x & y & z & t \\ x_1 & y_1 & z_1 & t_1 \end{pmatrix} < 2.$$

The total space $\mathcal{S}_1$ is defined by these equations plus the homogeneous equation

$$(5.2.1) \qquad x_1 y_1 = t_1^2$$

in $\mathbf{A}^4_{xyzt} \times \mathbf{P}^3_{x_1:y_1:z_1:t_1}$. The exceptional divisor $\Sigma_1$ is the corank 1 quadric defined by the equation (5.2.1) in $\mathbf{P}^3_{x_1:y_1:z_1:t_1}$ (and it sits over the origin of $\mathbf{A}^4_{xyzt}$). The central fibre of $\mathcal{S}_1$ is defined by the equation $t = 0$, which implies $t_1 x = t_1 y = t_1 z = 0$, hence it is the reduced sum of divisors $S'_0 + S''_0 + \Sigma_1$, where we denote by abuse of notation $S'_0, S''_0 \subseteq \mathcal{S}_1$ the respective proper transforms of $S'_0, S''_0 \subseteq \tilde{\mathcal{S}}$ (the former are the respective blow-ups of the latter at their smooth point $\mathbf{z}$).

The surface $\mathcal{C}_1$ on the other hand is defined inside $\mathcal{S}_1$ by the same equations as in (4.1). We leave it to the reader to check that the central fibres of $\mathcal{S}_1$ and $\mathcal{C}_1$ are arranged as indicated on Figure 10 (compare also with Figure 8, right-hand-side).

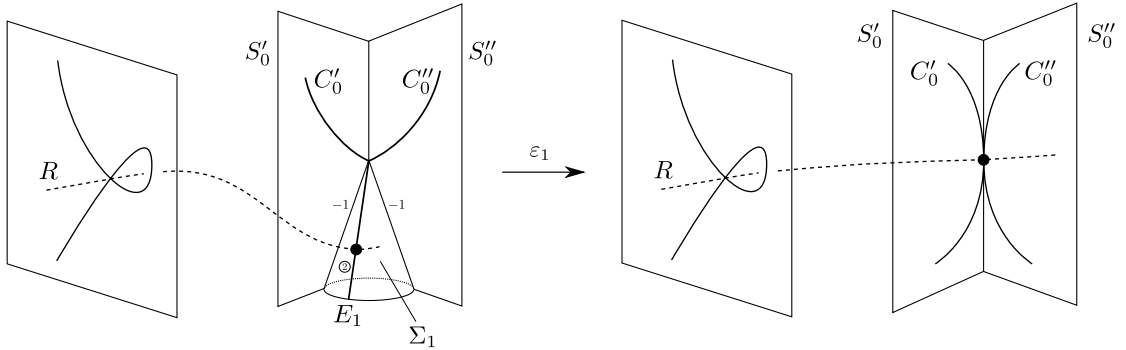

Figure 10: Blow-up of the tacnode, $\varepsilon_1 : \mathcal{S}_1 \to \mathcal{S}$

**(5.3) Blow-up of the vertex of $\Sigma_1$.** Let $\varepsilon_2 : \mathcal{S}_2 \to \mathcal{S}_1$ be the blow-up at the point $\mathbf{z}_1 = (0, 0, 0, 0, (0 : 0 : 1 : 0)) \in \Sigma_1$, $\mathcal{C}_2$ be the proper transform of $\mathcal{C}_1$, and $\Sigma_2$ be the exceptional divisor. As usual, we call $\Sigma_1$ its own proper transform in $\mathcal{S}_2$.

To fix ideas, we consider the affine chart $z_1 = 1$ with respect to the homogeneous coordinates introduced in (5.2) for the previous blow-up; there, $\mathcal{S}_1$ is defined by the affine equation $x_1 y_1 = t_1^2$ (compare (5.2.1)) in $\mathbf{A}^4$ with coordinates $(x_1, y_1, z, t_1)$. We add the new homogeneous variables $(x_2 : y_2 : z_2 : t_2) \in \mathbf{P}^3$, with the new equations

$$\mathrm{rk} \begin{pmatrix} x_1 & y_1 & z & t_1 \\ x_2 & y_2 & z_2 & t_2 \end{pmatrix} < 2.$$

Then $\mathcal{S}_2$ is defined by the homogeneous equation

(5.3.1) $$x_2 y_2 = t_2^2.$$

Exactly as in (5.2), the exceptional divisor $\Sigma_2$ is a corank 1 quadric surface. The central fibre of $\mathcal{S}_2$ is defined by $t = 0$, or equivalently (in the affine chart $z_1 = 1$) by $t_1 z = 0$. Let us work in the chart $x_2 = 1$. Then we have $z = x_1 z_2, t_1 = x_1 t_2$ and therefore $t_1 z = x_1^2 z_2 t_2$. Hence the central fibre of $\mathcal{S}_2$ is the sum of the three divisors respectively defined in $\mathcal{S}_2$ by the three equations $x_1^2 = 0, z_2 = 0, t_2 = 0$. The first equation gives a double component; the second a simple component; the third gives $x_2 y_2 = 0$, hence two reduced components. The double component is necessarily the exceptional divisor over the vertex of the cone. The upshot is that the central fibre is the non-reduced sum $S'_0 + S''_0 + \Sigma_1 + 2\Sigma_2$, with our usual abuse of notation: $S'_0, S''_0 \subseteq \mathcal{S}_2$ are the blow-ups of $S'_0, S''_0 \subseteq \mathcal{S}_1$ at their smooth point $\mathbf{z}_1$, and $\Sigma_1 \subseteq \mathcal{S}_2$ is the blow-up of $\Sigma_1 \subseteq \mathcal{S}_1$ at its vertex (in particular, it is isomorphic to an $\mathbf{F}_2$ minimal rational ruled surface).

Again, the surface $\mathcal{C}_2$ is defined inside $\mathcal{S}_2$ by the same equations as in (4.2), and we leave it to the reader to check that the central fibres of $\mathcal{S}_2$ and $\mathcal{C}_2$ are arranged as indicated on Figure 11, right-hand-side (compare also with Figure 9, right-hand-side). Note in particular that the curve $E_2$ is contained in $\Sigma_2$ and tangent to the curve $\Sigma_2 \cap \Sigma_1$.

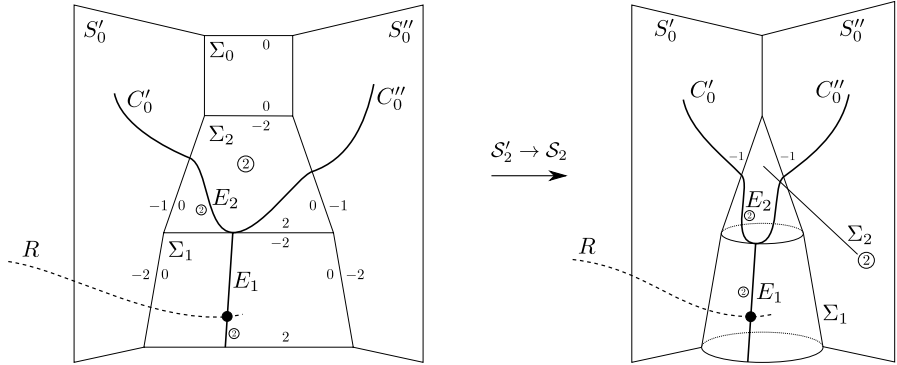

Figure 11: Blow-up of the vertex of $\Sigma_1$, $\varepsilon_2 : \mathcal{S}_2 \to \mathcal{S}_1$, and its resolution $\mathcal{S}'_2 \to \mathcal{S}_2$

**(5.4) Resolution of singularities of $\mathcal{S}_2$.** The threefold $\mathcal{S}_2$ is singular along the curve $S'_0 \cap S''_0$, at the generic point of which it has an $A_1$ singularity. Let $\mathcal{S}'_2 \to \mathcal{S}_2$ be the blow-up along this curve, with exceptional divisor $\Sigma_0$ a $\mathbf{P}^1$-bundle over $S'_0 \cap S''_0$. Then, $\mathcal{S}'_2$ is non-singular, with central fibre

(5.4.1) $$S'_0 + S''_0 + \Sigma_1 + 2\Sigma_2 + \Sigma_0;$$

all summands are isomorphic to their counterparts in $\mathcal{S}_2$ except $\Sigma_2 \subseteq \mathcal{S}'_2$, which is the minimal resolution of the quadric cone $\Sigma_2 \subseteq \mathcal{S}_2$. On Figure 11, left-hand-side, we picture the central fibre of $\mathcal{S}'_2$ with the self-intersections of the intersection curves of any two of its irreducible components. One may check that they concord with the triple-point-formula, which we recall as Lemma (5.5) below.

One may prefer to skip step (5.4); this induces only minor changes.

**(5.5) Lemma** (Triple Point Formula, see [1, 4]). *Let $f : \mathcal{X} \to \mathbf{A}^1$ be a family of surfaces with smooth total space, such that the central fibre $X_0$ has simple normal crossing support but may be non-reduced. Let $Q, Q'$ be irreducible components of $X_0$, with respective multiplicities $m$ and $m'$ in $X_0$, intersecting along the* double curve $B$. *Then*

$$m' B_Q^2 + m B_{Q'}^2 + \sum_{Q'' \neq Q, Q'} \mathrm{mult}_{X_0}(Q'') \cdot \mathrm{Card}(B \cap Q'') = 0,$$

*where $Q''$ ranges through all irreducible components of $X_0$ besides $Q'$ and $Q''$.*

In the above statement, $\mathrm{mult}_{X_0}(Q'')$ denotes the multiplicity of $Q''$ in the central fibre $X_0$, and $B_Q^2 = \deg(N_{B/Q})$ (resp. $B_{Q'}^2 = \deg(N_{B/Q'})$) is the self-intersection of $B$ as a curve in the surface $Q$ (resp. $Q'$).

**(5.6) Base change and resolution.** Let $\tilde{\mathcal{S}}_2 \to \mathcal{S}'_2$ be the double covering ramified over the central fibre. It is defined by adding a square root of $t$, i.e., we introduce a new variable $s$ with the new equation $s^2 = t$, and consider the family of surfaces $\tilde{\mathcal{S}}_2 \to \mathbf{A}^1_s$ defined by the projection map. The proper transform of $\mathcal{C}_2$ in $\tilde{\mathcal{S}}_2$ is isomorphic to the family $\tilde{\mathcal{C}}_2$ of (4.3).

The central fibre of $\tilde{\mathcal{S}}_2$ is, of course, the same as that of $\mathcal{S}'_2$. We let $\tilde{\mathcal{S}}'_2 \to \tilde{\mathcal{S}}_2$ be the blow-up of the Weil divisor $\Sigma_2$. This has the effect of replacing $\Sigma_2$ by its double cover $\tilde{\Sigma}_2$ branched over the anticanonical divisor $(S'_0 + S''_0 + \Sigma_1 + \Sigma_0) \cap \Sigma_2$, and doesn't change anything else.

Let $\tilde{\mathcal{C}}'_2$ be the proper transform of $\tilde{\mathcal{C}}_2$ in $\tilde{\mathcal{S}}'_2$. It is obtained by replacing the curve $E_2$ with its double cover $\tilde{E}_2$ branched along the divisor $(C'_0 + C''_0 + 2E_1) \cap E_2$: this is a rational curve with an ordinary double point at the intersection with $E_1$ (see (4.3.1)). The central fibre is depicted on Figure 12, right-hand-side.

One checks with local computations in the style of those in subsection $4.2^1$ that $\tilde{\mathcal{S}}'_2$ is singular only along the four curves at the intersection of one of $S'_0, S''_0$ with one of $\Sigma_0, \Sigma_1$, where it has lines of $A_1$ double points. Correspondingly, the surface $\tilde{\Sigma}_2$ has four ordinary double points over the double points of the anticanonical divisor of $\Sigma_2$.

It follows that the blow-up of $\tilde{\mathcal{S}}'_2$ along the four aforementioned curves is non-singular. We let $\bar{\mathcal{S}}_2 \to \tilde{\mathcal{S}}'_2$ be the corresponding blow-up morphism. The central fibre of $\bar{\mathcal{S}}_2 \to \mathbf{A}^1_s$ is pictured on the left-hand-side of Figure 12 below (see also Figure 9 central-bottom). The self-intersections of the double curves of the central fibre are indicated as well.

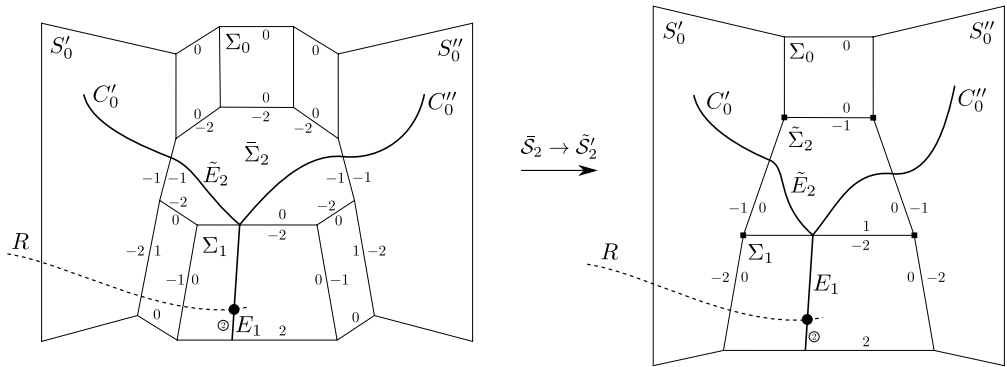

Figure 12: Base change and resolution

**(5.7) Conclusion.** Eventually, let $\mathcal{S}_3 \to \bar{\mathcal{S}}_2$ be the blow-up along the curve $E_1$, and $\Sigma_1^1$ its exceptional divisor. The curve $E_1$ is contained in $\tilde{\mathcal{C}}'_2$, hence the proper transform of the latter family of curves in $\mathcal{S}_3$ is isomorphic to $\bar{\mathcal{C}}_2$ in (4.3), see (4.3.1). The result is indicated on Figure 13 below (see also Figure 9, left-hand-side).

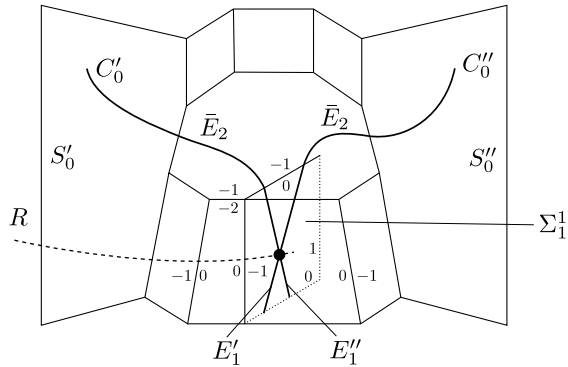

Figure 13: Flag semi-stable reduction, $\mathcal{S}_3$

This is the flag semi-stable reduction we were looking for: the family of surfaces $\mathcal{S}_3$ is semi-stable, and so is the family of curves $\bar{\mathcal{C}}_2$ it contains; thus, we have realized a tacnodal limit curve of 1-nodal curves as a curve with one distinguished node, in a semi-stable limit of the ambient surfaces.

---

[1]for instance, at the generic point of $S'_0 \cap \Sigma_2$, $\tilde{\mathcal{S}}_2$ has an equation equivalent to $a^2 b = u^2$, and $\tilde{\mathcal{S}}'_2 \to \tilde{\mathcal{S}}_2$ is locally the blow-up of $a = u = 0$.

One may in turn modify this flag semi-stable reduction in order to obtain whatever variant we fancy. For instance, it is possible to understand in this way the model described in (1.6).

## References

**Chapters in this Volume**

**External References**